\documentclass[a4paper,10pt]{article}
\usepackage{amsmath, amsthm, amssymb}
\usepackage{url}

\usepackage{physics}
\usepackage{graphicx,lipsum}
\usepackage{adjustbox}
\graphicspath{ {./downloads/} }

\usepackage[colorlinks,citecolor=red,urlcolor=blue,bookmarks=true,hypertexnames=true]{hyperref}
\usepackage{tikz}
\usetikzlibrary{calc}
\usetikzlibrary{shapes}
\usepackage[autostyle]{csquotes}
\makeatletter

\usepackage[colorlinks]{hyperref}
\usepackage[nameinlink,capitalize]{cleveref}
\newtheorem{theorem}{Theorem}[section]
\newtheorem{corollary}{Corollary}[theorem]
\newtheorem{lemma}[theorem]{Lemma}
\newtheorem{proposition}{Proposition}[section]
\newtheorem{conjecture}{Conjecture}[section]
\newtheoremstyle{named}{}{}{\itshape}{}{\bfseries}{.}{.5em}{\thmnote{#3's }#1}
\theoremstyle{named}

\newtheorem{remark}{Remark}
\theoremstyle{definition}
\newtheorem{definition}{Definition}[section]
\newtheorem{example}{Example}

\theoremstyle{remark}

\theoremstyle{conclusion}

\theoremstyle{observation}

\usepackage{xspace}
\usepackage[margin=1.0in]{geometry}

\newcommand\isomto{\stackrel{\textstyle\sim}{\smash{\longrightarrow}\rule{0pt}{0.4ex}}}

\usepackage{titlesec}

\usepackage{mathtools}

\DeclarePairedDelimiterX{\inp}[2]{\langle}{\rangle}{#1, #2}
\titleformat{\chapter}
  {\Large\bfseries}
  {}
  {0pt}
  {\huge}

\usepackage{makeidx}  
\usepackage{dsfont}
\usepackage{amsmath,amssymb,amsfonts,amscd}
\usepackage{hyperref}
\usepackage{booktabs} 
\usepackage{appendix}
\usepackage[utf8]{inputenc}

\usepackage{tikz,tikz-cd}
\tikzset{
    invisible/.style={opacity=0},
    visible on/.style={alt={#1{}{invisible}}},
    alt/.code args={<#1>#2#3}{%
      \alt<#1>{\pgfkeysalso{#2}}{\pgfkeysalso{#3}}%
  }
}
\tikzset{
  symbol/.style={
    draw=none,
    every to/.append style={
      edge node={node [sloped, allow upside down, auto=false]{$#1$}}}
  }
}

\usepackage{xcolor}

\begin{document}

\begin{center}
\fontsize{13pt}{10pt}\selectfont
    \textsc{\textbf{NONCOMMUTATIVE HOURGLASSES I: ON CLASSIFICATION OF THE $\mathds{Q}$-FANO 3-FOLDS GORENSTEIN INDEX $2$ VIA DERIVED CATEGORY}}
    \end{center}
\vspace{0.1cm}
\begin{center}
   \fontsize{10pt}{8pt}\selectfont
    \textsc{HAO XINGBANG}
\end{center}
\vspace{0.2cm}
\begin{center}
   \fontsize{12pt}{10pt}\selectfont
    \textsc{Abstract}
\end{center}

In previous work, Takagi used the methods of solving the Sarkisov links by calculating the corresponding Diophantine equations and the construction of key varieties to give all possible classifications and some implementations of a class $\mathds{Q}$-Fano 3-fold with Fano index 1/2 and at worst $(1, 1, 1)/2$ or QODP singularities. Firstly, we use a method different from Kawamata's work to give the derived category formulas for general weighted blow-up and Kawamata weighted blow-up. On this basis, we study the changing behavior of the derived category of Takagi's varieties under Sarkisov links. Finally, by studying non-commutative projections, we give exceptional collections on the derived category of Takagi's varieties and their corresponding geometric meanings.
\tableofcontents

%

\section{Introduction}
\label{sec:introduction}

The study of derived categories of algebraic varieties traces its origins to Beilinson's seminal work on diagonal resolutions for projective spaces \cite{Bei}. Kapranov's subsequent generalization to homogeneous spaces and quadrics \cite{Kap} established the fundamental role of exceptional collections and their duals in this theory. These developments laid the foundation for categorical approaches to studying algebraic varieties, particularly in low dimensions.

Mukai's work \cite{MuF} revealed that three-dimensional smooth varieties often admit highly symmetric ambient spaces - typically homogeneous varieties - whose derived categories are generated by complete exceptional collections. Geometrically realized as complete intersections, these spaces exhibit a remarkable categorical correspondence: the process of forming complete intersections manifests categorically as partitions of exceptional collections, naturally yielding matrix factorization categories as byproducts \cite{Orlh}. Mukai further developed this program through projections along rational loci, showing how birational transformations preserving singularity types encode essential information about the varieties.

The categorical counterpart of these geometric constructions was pioneered by Kuznetsov, who established the existence of full exceptional collections on $A_{22}$ varieties \cite{KuV22}. Subsequent work \cite{KuV14,KuGM} extended these results to broader classes of Sarkisov links, requiring modifications to classical constructions through techniques including blow-ups, projective bundle formations, standard flips and flops, and contractions of linear fibers. Bondal and Orlov's foundational work \cite{BO,Orltp} on the derived categories of these transformations provides crucial tools for our current investigation.

This article presents a classification of derived categories for singular three-dimensional algebraic varieties, building on modern adaptations of Sarkisov links for singular cases. Our approach synthesizes Takagi's extension \cite{Tak,Tak1} of Takeuchi's methods to Fano index $\frac{1}{2}$ $\mathbb{Q}$-Fano 3-folds with mild singularities, and the projection/unprojection techniques for orbifold singularities developed by Brown-Kerber-Reid \cite{BKR}. The key geometric process involves weighted blow-ups at orbifold singular points followed by flop transformations and sequences of local toric flips.

Adopting Bondal-Orlov's perspective \cite{BO}, we focus on noncommutative aspects of these constructions. To establish rigorous foundations, we first provide detailed proofs for two fundamental results that, while known to experts, lack complete documentation in the literature.

\begin{theorem}[\cite{Kalog}\cite{KaTor}, Thm. \ref{thm:smoohtweightedblowup}]
If $X$ is a smooth variety, after $(a_{1},..,a_{c})$ \footnote{\, A sequence of positive integers has no non-trivial common divisor \ref{no:ai}, to work with the canonical stack we assume it is well-formed.\label{foot:ai:intro}}-weighted blowing up at a codimension $c$ smooth center $S$, we get another variety (algebraic space) $Y$ with only cyclic quotient singularities and we denote its canonical stack by $\widetilde{Y}$.
$$\begin{tikzcd}[ampersand replacement=\&]
E_{\widetilde{Y}}\arrow[r,"\iota"]\arrow[d,"\pi_{E_{\widetilde{Y}}}"] \&\widetilde{Y}\arrow[d,"\pi"]  \\
S \arrow[r,"\kappa"] \&X
\end{tikzcd}$$
Then we have a semi-orthogonal decomposition of $\widetilde{Y}$ with respect to $X$ as following:
$$\mathrm{D}^{b}(\widetilde{Y})\cong\langle \mathrm{Im}\,\mathrm{\Psi}_{-\sum_{i}{a_{i}}+1},...,\mathrm{Im}\,\mathrm{\Psi}_{-1},\mathrm{Im}\, \mathrm{\Psi}\rangle$$
where $\mathrm{\Psi}$ and $\mathrm{\Psi}_{i}$ are defined by $\mathrm{\Psi}(-):=\pi^{*}(-)$ and $\mathrm{\Psi}_{i}(-):=\iota_{*}\bigl(\pi_{E_{\widetilde{Y}}}^{*}(-)\otimes\mathcal{O}_{E_{\widetilde{Y}}}(i)\bigr)$ for any integer $i$.

\end{theorem}

\begin{theorem}[Thm. \ref{thm:singularweightedblowup}]
If $X$ is a variety with some codimension $c$ and $(a_{1},..,a_{c})/r$  type \footnote{\, Same as above \ref{foot:ai:intro}, additionally condition $a_{i}<r$ is unnecessary.} cyclic quotient center $S$, we denotes $\widetilde{X}$ as its canonical stack, with the global split character and weighted normal data specified in Section~4.3, after weighted blowing up with the compatible weight at $S$, we get another variety (algebraic space) $Y$ with some cyclic quotient center, and denotes its canonical stack by $\widetilde{Y}$.\\
\begin{enumerate}
\item If $r<\sum_{i=1}^{c} a_{i}$:\\

There is a semi-orthogonal decomposition of $\widetilde{Y}$ with respect to $\widetilde{X}$ as following:
$$\mathrm{D}^{b}(\widetilde{Y})\cong\langle \mathrm{Im}\,\mathrm{\Theta}_{-\sum_{i}{a_{i}}+r},..., \mathrm{Im}\,\mathrm{\Theta}_{-1},\mathrm{Im}\, \mathrm{\Theta}\rangle$$
where $\mathrm{\Theta}$ and $\mathrm{\Theta}_{i}$ are defined by considering the following diagrams:
$$\begin{tikzcd}[ampersand replacement=\&]
E_{\widetilde{Y}} \arrow[r,"\iota"]\arrow[d,"f_{E_{\widetilde{Y}}}"] \&\widetilde{Y}\arrow[d,"f"]  \&\sqrt[ r]{E_{\widetilde{Y}}/\widetilde{Y}} \arrow[r,"\varpi"]\arrow[d,"\pi"] \&\widetilde{Y}\arrow[d,"f"] \\
S \arrow[r,"\kappa"] \&X  \&\widetilde{X} \arrow[r,"\theta"] \& X
\end{tikzcd}$$
we have $\mathrm{\Theta}(-):=\varpi_{*}\pi^{*}(-)$ and $\mathrm{\Theta}_{i}(-):=\iota_{*}\bigl(f_{E_{\widetilde{Y}}}^{*}(-)\otimes\mathcal{O}_{E_{\widetilde{Y}}}(i)\bigr)$
for any integer $i$.\\

\item If $r=\sum_{i=1}^{c} a_{i}$:\\

There is an equivalence induced by either $\mathrm{\Theta}$ or $\mathrm{\Phi}$:

$$\mathrm{D}^{b}(\widetilde{Y})\simeq \mathrm{D}^{b}(\widetilde{X})$$

\item If $r>\sum_{i=1}^{c} a_{i}$:\\

There is a semi-orthogonal decomposition of $\widetilde{X}$ with respect to $\widetilde{Y}$ as following:
$$\mathrm{D}^{b}(\widetilde{X})=\langle\mathrm{Im}\,\Phi_{r-\sum_{i=1}^{c}a_i},\ldots,\mathrm{Im}\,\Phi_{1},\mathrm{Im}\, \mathrm{\Phi}\rangle$$
where $\mathrm{\Phi}$ and $\mathrm{\Phi}_{i}$ are defined by considering the following diagrams:
$$\begin{tikzcd}[ampersand replacement=\&]
[S/\mu_{r}] \arrow[r,"\varsigma"]\arrow[d,"\tau"] \&\widetilde{X}  \&\sqrt[r]{E_{\widetilde{Y}}/\widetilde{Y}} \arrow[r,"\varpi"]\arrow[d,"\pi"] \&\widetilde{Y}\arrow[d,"f"] \\
S  \& \&\widetilde{X} \arrow[r,"\theta"] \& X
\end{tikzcd}$$
we have $\mathrm{\Phi}(-):=\pi_{*}\varpi^{*}(-)$ and $\mathrm{\Phi}_{i}(-):=\varsigma_{*}(\tau^{*}(-)\otimes \nu^{i})$ for any integer $i$.
\end{enumerate}
\end{theorem}

We compare our approach with Kawamata's method for analyzing toric divisorial contractions. In \cite[Theorem 4.2]{Kalog}, Kawamata established that for quasi-smooth toroidal divisorial contractions $Y \to X$, the canonical stack lift diagram
\[
\begin{tikzcd}[column sep=small]
& (\widetilde{Y}\times_X\widetilde{X})^{\mathrm{nor}} \arrow[rr, "\nu"] \arrow[d, "\mu"'] & & \widetilde{Y} \\
& \widetilde{X} & &
\end{tikzcd}
\]
induces a fully faithful functor $\mu_*\nu^*$ or $\nu_*\mu^*$, depending on their $K$-relationship. Subsequently, in \cite[Theorem 5.2]{KaTor}, Kawamata showed that for certain quasi-smooth toric divisorial contractions, there exists an orthogonal complement to $\mathrm{Im}\,\mu_*\nu^*$ (or $\mathrm{Im}\,\nu_*\mu^*$) generated by fully faithful functors associated to ray divisors on $\widetilde{X}$ (or $\widetilde{Y}$). We observe that these results could alternatively be derived from the fundamental property that skyscraper sheaves form a spanning class in the toric case.

The primary objective of our work is to refine and elaborate on \cite[Theorem 5.2]{KaTor}, providing a more intrinsic and concrete explanation of the underlying method. Rather than examining the birational transformations of ray divisors in analytically toric coordinates (which naturally constitute a spanning class), we concentrate directly on the behavior of skyscraper sheaves under triangulated functors. This approach leads to our key proposition (see also Proposition~\ref{prop:pullbackskycaperalongsmoothcenter} for the weighted blow-up version).

\begin{proposition}[Prop. \ref{prop:pullbackcenterlocal}]
If we Kawamata blow up $X$ along the cyclic quotient center $S$ and consider its stacky lift picture:
$$\begin{tikzcd}[ampersand replacement=\&]
 \&\sqrt[r]{E_{\widetilde{Y}}/\widetilde{Y}}\arrow[d,"\pi"] \arrow[r,"\varpi"]\&\widetilde{Y}:=\widetilde{w_{1/r}Bl_{S}X}\\
 \&  {\widetilde{X}}
\end{tikzcd}$$
where $E_{\widetilde{Y}}$ is the stacky lift of the exceptional divisor of the Kawamata blow-up. Then we have for any generalized geometrical point $\mathrm{k}_{x}\otimes\nu^{i}$ on $[S/\mu_{r}]$, $\overline{\mathcal{D}}_{x,i}:=\varpi_{*}\pi^{*}(\mathrm{k}_{x}\otimes\nu^{i})$ has the following finite Postnikov filtration:
\[\begin{adjustbox}{max width=0.98\linewidth}\(\begin{tikzcd}[column sep=0.5em]
 & 0 \arrow{rr}&& E_{-i-kr}\arrow{dl} \\
&& \iota_{x*}\mathcal{D}_{i+kr}\arrow[ul,dashed,"\Delta"]
\end{tikzcd}......\begin{tikzcd}[column sep=0.5em]
 &   E_{-i-r}\arrow{rr}&&  E_{-i}\arrow{rr}\arrow{dl}&& \overline{\mathcal{D}}_{x,i}\arrow{dl} \\
 && \iota_{x*}\mathcal{D}_{i+r}\arrow[ul,dashed,"\Delta"]  && \iota_{x*}\mathcal{D}_{i}\arrow[ul,dashed,"\Delta"]
 \end{tikzcd}\)\end{adjustbox}\]
and $\mathcal{D}_{x,i}$ is the $i$-th dual exceptional object on fiber over $x$. For $\sum_{j=1}^c a_j\le r$ and $0\le i<r$, the image is $\iota_{x*}\mathcal D_i$ if $i<\sum_j a_j$, and is zero otherwise.\\
\end{proposition}
This illustrates a correspondence between local finite group representations and the dual exceptional collections on the exceptional linear fibers produced by Kawamata blow-up (or weighted blow-up). This method is analogous to and generalizes the classical Bondal-Orlov derived blow-up formula \cite{Orltp} and Krug, Ploog, and Sosna's proof for cyclic quotient singularities \cite[Theorem 4.1]{KPS}.

Our approach's advantage lies in naturally revealing all exceptional collections and their geometric implications in the semi-orthogonal decomposition for local generalized McKay correspondences. Extensions to non-abelian groups remain to be studied.

Direct applications of our non-commutative weighted blow-up formula yield categorical behaviors for low-dimensional toric flips (Theorem~\ref{thm:fraciaflip}), generalizing classical standard flips. We also reproduce a special case of this result alternatively for later use.
\begin{proposition}[Prop. \ref{prop:secondshoukufr}]
For $(1,2:1,1)$ type Francia's flip, we have another categorical description under Shokurov's link:
\begin{equation}\mathrm{D}^{b}(\widetilde{X})\simeq\langle\mathcal{O}_{C}(-1), \mathrm{Im}\,\Gamma\big(D^{b}(X^{+})\big) \rangle\end{equation}
where $\Gamma(-):=\pi_{*}\varpi^{!}\overline{\phi}_{*}\overline{\varphi}^{!}f^{+*}(-)$  is a fully-faithful functor, and we view the flipping curve $C$ as $\mathcal{P}_{2}^{1}$.
\end{proposition}

Then, we apply the derived category behavior in birational transformations obtained above to study the categorical classification of $\mathbb{Q}$-Fano $3$-folds given by Takagi \cite{Tak,Tak1,Tak2}.
\begin{definition}[\cite{Tak} Main Assumption 0.1] Takagi's varieties $\mathrm{X}$ are $\mathds{Q}$-factorial $\mathds{Q}$-Fano 3-fold $\mathrm{X}$ satisfy the following conditions:
\begin{enumerate}
  \item picard number $\rho(\mathrm{X})=1$,
  \item Gorenstein index $I(\mathrm{X})=2$,
  \item Cartier Fano index $F(\mathrm{X}) = 1/2$,
  \item $h_{0}(-\mathrm{K}_{\mathrm{X}})\geq 4$,
  \item there exists an index $2$ point $p$ such that it is analytically isomorphic to
$$(xy+z^{2}+u^{a}=0)\Big/\mu_{2}^{(1,1,1,0)}$$
for some $a\geq 0$.
\end{enumerate}
\end{definition}

Takagi performed a weighted blow-up at the singular point $\mathrm{p}$ and ran a minimal model program. By utilizing Mori-type extremal contraction classification and calculating divisor intersection numbers, the complete picture is described by:
$$\begin{tikzcd}[ampersand replacement=\&]
\mathrm{Y}=w_{\frac{1}{2}}Bl^{(1,1,1,0)}_{\mathrm{p}}\mathrm{X} \arrow[rr,dashrightarrow]\arrow[rd,"\mathrm{g}"] \arrow[d,"\mathrm{f}"]\&\& \mathrm{Y}_{1}\arrow[rr,dashrightarrow]\arrow[ld,"\mathrm{g}'"]\arrow[rd,"\mathrm{g}_{1}"] \& \&\mathrm{Y_{2}}=\mathrm{Y'}\arrow[d,"\mathrm{f'}"]\arrow[ld,"\mathrm{g}'_{1}"] \\
\mathrm{X}\&\mathrm{Z}\&\&\mathrm{Z}_{1}\& \mathrm{X'}
\end{tikzcd}$$
where the small $\mathbb{Q}$-factorial modifications ($\mathrm{SQM}$s) consist of first a flop followed by a flip. Takagi provides complete classifications of $\mathrm{f}'$ and $\mathrm{X}'$ in \cite[Main Theorem 0.3]{Tak}, Tables 1 to 5.

However, not all classifications can be realized concretely via inverse Sarkisov links. Takagi establishes implementation results for cases where X has only orbifold or QODP singularities in \cite[Theorem 0.10 to Theorem 0.20]{Tak1}, proving existence of such varieties.

We now consider a variety with only orbifold singularities under Takagi's classification in Table 1. In these cases, $\mathrm{f}'$ is simply a blow-up along a smooth curve $\mathrm{C}_{1,i}$ on $\mathrm{X}'_{1,i}$, and all varieties $\mathrm{X}_{1,i}$ are realizable. In Section~\ref{Geometryofhourglasses}, we show these varieties are geometrically similar to ``hourglasses'' of varying shapes, and we adopt this term for brevity. First, we prove the residue categories of $\mathrm{X}_{1,i}'$ remain invariant under Sarkisov links.

\begin{theorem}[Thm. \ref{thm:SODh}]
For any $i$, the canonical stack $\widetilde{\mathrm{X}}_{1,i}$ of hourglass $\mathrm{X}_{1,i}$  admits semi-orthogonal decomposition as following:
        $$\mathrm{D}^{b}(\widetilde{\mathrm{X}}_{1,i})\cong \langle\,  {\mathcal{R}}es_{1,i},\mathrm{Exc}_{1,i} \,\rangle$$
        where $\mathrm{Exc}_{1,i}$ is an exceptional collection,
        and the residue part ${\mathcal{R}}es_{1,i} $ is equivalent to a span of $\mathrm{D}^{b}(\mathrm{C_{1,i}} )$ and the matrix factorization category $\mathcal{MF}(\mathrm{\widetilde{X}}'_{1,i})$ as in Corollary \ref{cor:mfx'}, that is
        $${\mathcal{R}}es_{1,i}\simeq\langle\mathrm{D}^{b}(\mathrm{C_{1,i}} ),\mathcal{MF}(\mathrm{\widetilde{X}}'_{1,i})\rangle $$
Specifically, if $g(\mathrm{C_{1,i}} )=0$ and the matrix factorization category is degenerated, the noncommutative hourglass admits a full exceptional collection.
\end{theorem}

We investigate whether the exceptional collection in Theorem~\ref{thm:SODh} admits geometric interpretations similar to varieties with projective spaces or Grassmannians as ambient spaces, as seen in smooth Fano 3-folds \cite{KuFano}. For this purpose, we examine the specific case of $\mathrm{X}_{1.9}$ with $\mathrm{X}'_{1.9} = B_3$, which exhibits cubic geometry. The geometric description is as follows:

\begin{proposition}[Prop. \ref{prop:X1.9}]
We have a semi-orthogonal decomposition of $\widetilde{\mathrm{X}}_{1.9}$:
$$\mathrm{D}^{b}(\widetilde{\mathrm{X}}_{1.9})\simeq \langle \mathcal{O}_{\widetilde{\mathrm{X}}_{1.9}}(-2\mathrm{p})\rightarrow \bigoplus_{i=1}^{6}\mathcal{O}_{L_{i}}(-1),\mathcal{O}_{\widetilde{\mathrm{X}}_{1.9}},\mathrm{U}^{\vee}_{\widetilde{\mathrm{X}}_{1.9}},\mathcal{MF}(\widetilde{\mathrm{X}}_{1.9})\rangle$$
and there is an equivalence between the residue component $\mathcal{MF}(\widetilde{\mathrm{X}}_{1.9})$ and $\mathcal{MF}(B_{3})$ via a functor
$$\Gamma_{5}:\mathcal{MF}(B_{3})\longrightarrow\mathcal{MF}(\widetilde{\mathrm{X}}_{1.9})$$
$$\Gamma_{5}(-):=\pi_{*}\varpi^{*}\big[\phi_{*}\varphi^{*}\mathbb{R}_{\langle\mathcal{O}_{\mathrm{S}}(\mathrm{S}),\mathcal{O}_{\mathrm{Y}'},g'^{*}\mathrm{U}_{\mathrm{Z}}^{\vee}\rangle}\big(\mathrm{f}'^{*}(-)(\mathrm{H})\big)\big]$$
\end{proposition}

This conclusion suggests that the stacky line on $\mathrm{X}_{1,i}$ and the exceptional objects in $\mathrm{D}^b(\mathrm{Z})$ naturally induce exceptional objects on $\mathrm{X}_{1,i}$. By analyzing the projection $\mathrm{X}_{1,i} \to \mathrm{Z}_{1,i}$ through derived categories, we provide the following categorical explanation:

\begin{proposition}[Non-commutative projection, Prop. \ref{prop:projection}-\ref{prop:sgequi}]

We consider the first midpoint $\mathrm{Z}$ of our Sarkisov link, it admits an equivalence up to quotient via projection:
$$g_{*}\mathrm{\Theta}:\mathrm{D}^{b}(\widetilde{\mathrm{X}})/\mathcal{L}\simeq\mathrm{D}^{b}(\mathrm{\widetilde{Z}})/\langle\mathcal{O}_{\mathrm{A}}(-1)\rangle $$
while it also induces an equivalence:
$$g_{*}\mathrm{\Theta}:\mathcal{B}\isomto [^{\perp}\langle\mathcal{O}_{\mathrm{A}}(-1)\rangle]^{perf}$$
where $\mathrm{\widetilde{Z}}$ is stacky smoothing of $\mathrm{Z}$ at orbifold points, $\mathcal{L}$ is the subcategory generated by $\mathcal{O}_{L_{i}}(-1/2)$ for all stacky lines  $L_{i}$ passing the projection point $\mathrm{p}$ on $\widetilde{\mathrm{X}}_{1,i}$ and $\mathcal{B}$ is its left orthogonal component.
\end{proposition}

So we have a general description of the exceptional collections on a non-commutative hourglass.

\begin{corollary}
  Assume the common-model and local-incidence hypotheses of
Corollary~\ref{cor:multipleunprojection}. If we have an exceptional collection consisting of perfect objects on the \textit{extended global midpoint} $\mathrm{\overline{Z}}_{1,i}$  of (\ref{exsark3}), e.g.
$$\langle e_{1},...,e_{l}\rangle$$
and they are left orthogonal to the component generated by all unprojection planes $\mathrm{A}_{1}$,...,$\mathrm{A}_{N}$ $$\langle \mathcal{O}_{\mathrm{A}_{1}}(-1),...,\mathcal{O}_{\mathrm{A}_{N}}(-1)\rangle$$
then their unprojection
   $$\langle \mathcal{O}_{L_{i_{1}},\mathrm{p}_{a(1)}}(-1/2),..., \mathcal{O}_{L_{i_{m}},\mathrm{p}_{a(m)}}(-1/2), \mathrm{\Theta}^{!}g^{*}e_{1},...,\mathrm{\Theta}^{!}g^{*}e_{l} \rangle$$
  form an exceptional collection on $\widetilde{\mathrm{X}}_{1,i}$ provided the chosen curve objects have the flop-image local form considered above and are mutually semi-orthogonal (for instance, their supports are pairwise disjoint), where $N=N(\mathrm{X}_{1,i})$, $m$ is the number of selected curve objects and $L_{i_j,\mathrm{p}_{a(j)}}$ is a chosen stacky line through its marked centre $\mathrm{p}_{a(j)}$, with $a(j)\in\{1,\ldots,N\}$ on $\widetilde{\mathrm{X}}_{1,i}$.
\end{corollary}

Under these geometric hypotheses, this reduces the search for such exceptional collections to perfect exceptional collections on an extended global midpoint $\overline{\mathrm{Z}}_{1,i}^{1}$, which has only non-factorial ordinary double points as singularities. We are thus motivated to identify both the key varieties of these global extended midpoints and their exceptional collections, as their pullback should yield perfect exceptional collections on $\overline{\mathrm{Z}}_{1,i}^{1}$.

Based on specific examples (see Conjecture~\ref{conj:noncommuativehourglasses}), we propose that proving this may require considering additional Sarkisov links that clearly distinguish symmetries between stacky lines on $\widetilde{\mathrm{X}}_{1,i}$.
\section{Notations}

\begin{itemize}
  \item $\mathrm{k}:=\mathbb{C}$ - base field for all triangulated categories and varieties
  \item $V_{i}:=\mathrm{k}^{\oplus i}$ - $i$-dimensional vector space
  \item $\widetilde{X}$ - canonical stack of variety $X$ with quotient singularities
  \item $\widehat{Z}$ - strict transform of cycle $Z$ under birational morphism
  \item $\mathrm{k}_{x}$ - skyscraper sheaf at point $x$
  \item $\mu_{r}$ - cyclic group $\mathbb{Z}/r\mathbb{Z}$ of order $r$
  \item $\zeta_{r}$ - primitive $r$-th root of unity
  \item $i,j,k,\lambda,a,b$ - summation indices (ranges sometimes omitted when clear)
  \item $\chi$ - group character element
  \item $\nu_{i},\nu$ - cyclic group character generator
  \item $\nu_{x}^{i}:=\mathrm{k}_{x}\otimes\nu^{i}$ - character-twisted skyscraper sheaf
  \item $\mathrm{\Upsilon}^{*},\mathrm{\Upsilon}^{!}$ - left/right adjoints of admissible functor $\mathrm{\Upsilon}$
  \item $\mathcal{H}^{r}(-)$ - $r$-th cohomology under standard filtration
  \item $\mathcal{H}^{\mathrm{\Upsilon}}(-)$ - projection to $\mathrm{Im}\,\mathrm{\Upsilon}$ subcategory
  \item $(a_{0},..,a_{n})$, $(a_{1},..,a_{c})$, $(1,a+b;a,b)$ - integer sequences with $\gcd=1$ \label{no:ai}
  \item $\mathbb{P}(a_{0},..,a_{n})$ - weighted projective space
  \item $\mathcal{P}(a_{0},..,a_{n})$ - weighted projective stack (or $\mathcal{P}$ when unambiguous)
  \item $\mathcal{P}^{1}_{a}$ - stacky $\mathbb{P}^{1}$ with degree $a$ point; $\mathcal{O}_{\mathcal{P}^{1}_{a}}(1/a)$ associated line bundle
  \item $\mathbb{P}_{X}(E)$ - projectivization of sheaf $E$
  \item $wBl_{S}^{(a_{1},..,a_{c})}X$ - $(a_{1},..,a_{c})$-weighted blow-up of smooth center $S$
  \item $w_{1/r}Bl_{S}^{(a_{1},..,a_{c})}X$ - degree $r$ cyclic quotient weighted blow-up
  \item $f^{*},f_{*}, f^{!}, f_{!}$ - derived functors (pullback, pushforward, etc.)
  \item $\mathbf{L}_{0}f^{*},\mathbf{R}^{0}f_{*}$ - ordinary sheaf pullback/pushforward
  \item $\mathrm{H}^{*}$, $\mathrm{Ext}^{*}$ - cohomology and extension groups
  \item $\mathbb{R}_{E}X$, $\mathbb{L}_{E}X$ - right/left mutation with respect to $E$
  \item $\langle e_{1},...,e_{n} \rangle$ - triangulated subcategory generated by $e_{1},..,e_{n}$
  \item $(a_{1},..,a_{c})/r$ cyclic quotient singularity - $\mathbb{C}^{n}$ quotient by $x_{i}\mapsto \zeta_{r}^{a_{i}}x_{i}$ \label{cqs}
  \item $(1^{3}/2)$-singularity - $\mathbb{C}^{3}$ quotient by $x_{i}\mapsto -x_{i}$
  \item $\sqrt[n]{E/X}$ - root stack of $(\mathcal{O}(-E),s_{E})/X$
  \item $\mathcal{M}_{E}=\mathcal{O}(M_E)$ - positive tautological line bundle on $\sqrt[n]{E/X}$, dual to the universal ideal root
  \item $M_{E}$ - tautological divisor on $\sqrt[n]{E/X}$
  \item $\mathcal{N}$, $\mathcal{I}$ - normal and ideal sheaves
  \item $A_{2g-2}$ ($1\leq g \leq12$, $g \neq 11$) - Terminal Fano 3-fold index 1, $\rho=1$, genus $g$
  \item $B_{i}$ ($1 \leq i \leq 5$) - Terminal Fano 3-fold index 2, $\rho=1$, degree $8i$
  \item $\mathrm{X}_{1.i}$ ($1\leq i \leq14$) - Takagi's $\mathbb{Q}$-Fano 3-folds index 2 (Table 1)
  \item Hourglass - $\mathbb{Q}$-Fano 3-fold with orbifold singularities in Takagi's classification
\end{itemize}

\section{Preliminaries}

\subsection{Stack elementary}

\begin{definition}
A stack $\widetilde{X} /S$ is a Delinge-Mumford stack if the following hold:
\begin{enumerate}
  \item The diagonal
$$\Delta:\widetilde{X}\longrightarrow\widetilde{X}\times_{S}\widetilde{X}$$
is representable.
  \item There exists a $\acute{e}$tale surjective morphism $\rho: U\longrightarrow \widetilde{X} $ with $U$ a scheme.
\end{enumerate}
\end{definition}

We assume $\widetilde{X}/S$ is a Deligne-Mumford stack locally of finite type and with finite diagonal. By Keel-Mori theorem always admits a coarse Moduli space
$X$. For any geometric point $x$, the stabilizer group $G_{x}$ of $x$ is always a finite group. There exists an $\acute{e}$tale neighborhood $V\longrightarrow X$ of $x$ and a
finite $V$-scheme U with action of $G_{x}$ such that
$$\widetilde{X}\times_{X}V\simeq[U/G_{x}]$$
We usually say that a smooth stack has orbifold singularity locally.\\

On the contrary, if we have a variety (or algebraic space) $X$ with orbifold singularity, we expect that we have an intrinsic Deligne-Mumford stack $\widetilde{X}$ whose coarse moduli space $X$ is this scheme itself. Let $\widetilde{X}$ be an irreducible smooth Deligne-Mumford stack, and $\rho: \widetilde{X} \longrightarrow X$ be the structure morphism to the coarse moduli space. The Deligne-Mumford stack $\widetilde{X}$ will be called \textit{canonical} if $\rho$ is an isomorphism outside a small locus  \cite[Definition 4.4]{FMN}. The canonical property is universal among all dominant codimension preserving morphism from another smooth Deligne-Mumford stack $\rho': \widetilde{Y} \longrightarrow X$. For example, $X:=U/G$ is a variety with finite quotient singularities, if the action of $G$ on $U$ admits no pseudo-reflective divisors, we can see the quotient map
$$\rho:[U/G]\longrightarrow U/G$$
defines the coarse Moduli morphism from the canonical Deligne-Mumford stack $[U/G]$. Similarly, Vistoli \cite[Proposition 2.8]{Vis} gives the construction of canonical Deligne-Mumford stack for general variety with quotient singularities.\\

Using the stack method, we can smooth the general classical ramified morphism, this method is commonly used in mathematical physics when dealing with counting curves passing through a fixed divisor. First, we only consider a \textit{generalized effective Cartier divisor} $(L,\rho)$ on $X$, where $X$ is a variety, $L$ is a line bundle on $X$, while $\rho$ is a mapping from $L$ to $\mathcal{O}_{X}$.
\begin{definition}[Root stack]
Fix a generalized effective Cartier divisor $(L,\rho)$ and an integer $n\geq 1$.
Let $\sqrt[n]{(L,\rho)/X}$ be the fibered category over the category of schemes, whose objects defined over scheme $T$ are
triples
$$\big(f,(M,\rho_{M}),\sigma\big)$$
where  \begin{enumerate}
  \item $f$ is a $T$-point on $X$,
  \item $(M,\rho_{M})$ is a generalized effective divisor on $T$,
  \item $\sigma:M^{\otimes n}\xrightarrow{\sim}f^{*}L$ satisfies $f^{*}\rho\circ\sigma=\rho_M^{\otimes n}$.
\end{enumerate}
whose morphisms defined between triples $\big(f,(M,\rho_{M}),\sigma\big)$ and $\big(f',(M',\rho_{M'}),\sigma'\big)$ are pairs
$$(h,H)$$
where \begin{enumerate}
  \item $h:T'\to T$ satisfies $f'=f\circ h$,
  \item $H:M'\xrightarrow{\sim}h^{*}M$ satisfies $h^{*}\rho_M\circ H=\rho_{M'}$ and $h^{*}\sigma\circ H^{\otimes n}=\sigma'$.
\end{enumerate}
For an effective prime Cartier divisor $D$ on $X$, we write
\[
\sqrt[n]{D/X}:=\sqrt[n]{(\mathcal{O}_{X}(-D),s_D)/X},\qquad
s_D:\mathcal{O}_X(-D)\longrightarrow\mathcal{O}_X.
\]
Here $s_D$ is the inclusion defining $D$.
\end{definition}

We can extend this definition to the general case when $X$ is a  Deligne-Mumford stack, for more constructions and properties of the root stack we refer to \cite[Section 1.3.b]{FMN}.

\subsection{Weighted blowing up}

Let $\widetilde{X}$ be a smooth Delinge-Mumford stack, and $\widetilde{S}$ is a smooth center of codimension $c$. A weighted blow-up of weight $(a_{1},...,a_{c})$, where $a_{i}$ are positive integers has no common divisor, with the center $\widetilde{S}$ on $\widetilde{X}$ is the stacky projectivization

$$wBl_{\widetilde{S}}\widetilde{X}:=\mathcal{P}roj_{\widetilde{X}}(\bigoplus_{k\geq 0}\mathcal{\widetilde{I}}_{k})\longrightarrow \widetilde{X}$$
where $\mathcal{\widetilde{I}}_{k}$ are sheaves of ideals on $\widetilde{X}$, such that for any geometrical point $x$ on $\widetilde{S}$ there are $\acute{e}$tale lcoal smooth coordinate functions $x_{1}$,...,$x_{c}$ defining the center $\widetilde{S}$, and $\mathcal{\widetilde{I}}_{k}$ is generated by the monomials
$$<x_{1}^{b_{1}}\cdot,...,x_{c}^{b_{c}}|\sum_{i}a_{i}b_{i}\geq k>$$

We choose the local coordinates so that the weighted ideals agree on overlaps of the $\acute{e}$tale charts. The local weighted blow-ups then glue, and we denote their coarse moduli space by $wBl_{S}X$.

If we further assume $\widetilde{X}$ locally of finite type and with finite diagonal,  under the local $\acute{e}$tale neighborhood $V$ of the moduli $X$ such that $\widetilde{X} \times_{ X} V = [U/G]$
for some finite group acting on a smooth affine variety  $U$. If $I_{k}$ are ideal sheaves correspondences to $\mathcal{I}_{k}|_{[U/G]}$ on
$U$, we have a local discerption for $wBl_{\widetilde{S}}\widetilde{X}\times_{ X} V$ as quotient of the projectivization

$$[\mathcal{P}roj_{U}(\bigoplus_{k\geq 0}I_{k})/G]\longrightarrow [U/G]$$
while taking coarse moduli gives the local discerption for $wBl_{S}X$.\\

\subsubsection{Smooth type}
A trivial case ($G$ is trivial as above) is a weighted blow-up with weight $(a_{1},...,a_{c})$ at a smooth center $S$ on a smooth variety $X$ (or admitting mild singularities outside the center), the description is consistent with above
$$wBl_{S}X:=\mathcal{P}roj_{X}(\bigoplus_{k\geq 0}\mathcal{I}_{k})\longrightarrow X$$
$$<x_{1}^{b_{1}}\cdot,...,x_{c}^{b_{c}}|\sum_{i}a_{i}b_{i}\geq k>$$

We consider $\acute{e}$tale locally that there exists a smooth variety $U$  and a iterated cyclic cover $\rho:U\longrightarrow V$ branched at smooth divisor $D_{1},...,D_{c}$ with degree $a_{1},..,a_{c}$, such that $D_{1},...,D_{c}$ define the center $S$ and $U/\mu_{a_{1}}\times,..,\times\mu_{a_{c}}=V$. Then we have a local discerption for $wBl_{S}X\times_{ X} V$ as the projectivization

$$\mathcal{P}roj_{U}(\bigoplus_{k\geq 0}I_{k})/\mu_{a_{1}}\times,..,\times\mu_{a_{c}}\longrightarrow U/\mu_{a_{1}}\times,..,\times\mu_{a_{c}}$$
while $I_{i}$ is generated by the monomials with all variables contributing only a weight of $1$, $\mathcal{P}roj_{U}(\bigoplus_{k\geq 0}I_{k})$ is the ordinary blow-up. \\

We denote such a smooth $(a_{1},..,a_{c})$-weighted blow-up as $wBl^{(a_{1},..,a_{c})}_{S}X$ which is in general an algebraic space and its canonical stack $\widetilde{wBl^{(a_{1},..,a_{c})}_{S}X}$.

\subsubsection{Kawamata type}

We are especially interested in the case $S$ is a smooth variety lies on the $\mu_{r}$ cyclic quotient singular locus of weight $(a_{1},...,a_{c})$ on $X$ see (\ref{cqs}), the classic definition requires condition $a_{i}<r$ for any $i$, but we ignore this condition by considering a more general definition. That is to say, under the local $\acute{e}$tale neighborhood $V$ of $X$ such that $\widetilde{X} \times_{ X} V = [U/\mu_{r}]$, we can choose a set coordinates $x_{1}$, ... , $x_{c}$ of $U$ (they correspond to part of the affine coordinate functions of $\mathds{A}^{n}$ through an $\acute{e}$tale map) defining the center $S$, such that the action follows:

$$\mu_{r}\times U\longrightarrow U $$
$$\zeta_{r}\times x_{i}\mapsto \zeta_{r}^{a_{i}} x_{i} $$
We use such equivariant coordinates subject to the compatibility condition above.\\

From the local common $\acute{e}$tale diagram we can see the divisors associating with  $x_{i}$  are all smooth divisors and completely intersected (SNC) if we shrink our $\acute{e}$tale atlas small enough.  By taking iterated cyclic covering branched along $D(x_{i})$ of degree $a_{i}$, we get a smooth variety $U'$ with compatible $\mu_{r}$ action

$$\mu_{r}\times U'\longrightarrow U' $$
$$\zeta_{r}\times y_{i}\mapsto \zeta_{r} y_{i} $$
where $y_{i}$ is a local coordinate of $U'$ defining the center $S$, such that $y_{i}^{a_{i}}=x_{i}$.\\

Now if we weighted blow up $X$ at the center
$S$ with weight ($a_{1},..,a_{c}$), we use the following auxiliary root-cover presentation mapping to $wBl_{\widetilde{S}}\widetilde{X}\times_X V$

$$[\mathcal{P}roj_{U'}(\bigoplus_{k\geq 0}I'_{k})/\mu_{a_{1}}\times,..,\times\mu_{a_{c}}\times\mu_{r}]\longrightarrow [U'/\mu_{a_{1}}\times,..,\times\mu_{a_{c}}\times\mu_{r}]$$
and $I'_{i}$ is generated by the monomials
$$<y_{1}^{b_{1}}\cdot,...,y_{c}^{b_{c}}|\sum_{i}b_{i}\geq k>$$
It's not difficult to see $\mathcal{P}roj_{U'}(\bigoplus_{k\geq 0}I'_{k})$ is just the ordinary blow-up $U'$ at center $S$, so $[Bl_{S}U'/\mu_{a_{1}}\times\cdots\times\mu_{a_{c}}\times\mu_{r}]$ is an auxiliary root stack mapping to the weighted blow-up. It retains root stabilizers along the strict transforms of the coordinate divisors and is not, in general, the canonical stack of that blow-up.\\

 On the other hand side, we consider $\widetilde{wBl_{S}X}$ (the canonical stack of $wBl_{S}X$)  where they all have the coarse moduli space $wBl_{S}X$, but it is coincident with $wBl_{\widetilde{S}}\widetilde{X}$ up to a ramification on exceptional divisor. We define Kawamata blow-up along a $(a_{1},...,a_{c})$ type $\mu_{r}$ cyclic quotient singularity center $S$ on $X$ of weight $(a_{1},...,a_{c})$  as $w_{1/r}Bl_{S}^{a_{1},...,a_{c}}X$ which is in general an algebraic space, and its canonical stack denotes by $\widetilde{w_{1/r}Bl_{S}^{a_{1},...,a_{c}}X}$.

\subsection{Semi-orthogonal decomposition}

 Let $\mathcal{D}$ be a triangulated category, we refer to some basic definitions and properties of triangulated categories as \cite[Chapter 1]{H}.

\begin{definition}[Exceptional collection]
  We say an ordered collection of $$(E_{1},E_{2},\cdots,E_{n})$$ of objects of $\mathcal{D}$ is an exceptional collection if
  $$\mathrm{Ext}^{*}(E_{j},E_{i})=0\quad \forall j>i \quad\text{and} \quad \mathrm{Ext}^{*}(E_{i},E_{i})\simeq \mathrm{k}\quad \forall i$$
 an exceptional collection is full if it generates category $\mathcal{D}$.
\end{definition}

\begin{definition}[Semi-orthogonal decomposition]\label{def:SOD}
\begin{enumerate}
  \item A sequence $\mathcal{A}_{0},...,\mathcal{A}_{n}$ of full subcategories of $\mathcal{D}$ is semi-orthogonal if
$$\mathrm{Ext}^{*}(a_{j},a_{i})=0$$
for any $a_{i},a_{j}$ contained in $\mathcal{A}_{i},\mathcal{A}_{j}$ and $j>i$.
  \item The sequence is a semi-orthogonal decomposition of category $\mathcal{D}$, if
for any object $a$ in  $\mathcal{D}$ we have a filtration:
  $$\begin{tikzcd}[column sep=0.5em]
 & 0 \arrow{rr}&& E_{n-1}\arrow{dl}\arrow{rr}&& E_{n-2}  \arrow{dl} \\
&& a_{n}\arrow[ul,dashed,"\Delta"]&& a_{n-1}\arrow[ul,dashed,"\Delta"]
\end{tikzcd}
...\begin{tikzcd}[column sep=0.5em]
 &     E_{0}\arrow{rr}&& a\arrow{dl} \\
&&a_{0}\arrow[ul,dashed,"\Delta"]
\end{tikzcd}$$
  where $a_{i}$ is contained in $\mathcal{A}_{i}$, for any $i$. We write the decomposition as:
  $$\mathcal{D}\simeq\langle\mathcal{A}_{0},..,\mathcal{A}_{n}\rangle$$
\end{enumerate}
\end{definition}

For a full subcategory $\mathcal{A}$ in $\mathcal{D}$, if there is no special statement, we generally refer to the embedding functor of $\mathcal{A}$ to $\mathcal{D}$ by $i_{\mathcal{A}}$.  We say $\mathcal{A}$ is left (resp. right) admissible if $i_{\mathcal{A}}$ admits a left adjoint $i^{*}_{\mathcal{A}}$ (resp. a right adjoint $i^{!}_{\mathcal{A}}$).

\begin{definition}[Mutation]\label{def:mutataion}
Let's assume that $\mathcal{A}$ is a left or right admissible full subcategory of $\mathcal{D}$, then we can define the right or left mutation functor for its embedding $i_{\mathcal{A}}$ by:
$$\mathbb{R}_{\mathcal{A}}(-):=\mathrm{Cone} \big((-)\longrightarrow i_{\mathcal{A}}i^{*}_{\mathcal{A}}(-)\big)[-1]$$
or
$$\mathbb{L}_{\mathcal{A}}(-):=\mathrm{Cone} \big( i_{\mathcal{A}}i^{!}_{\mathcal{A}}(-)\longrightarrow(-)\big)$$
\end{definition}

Especially, if $\mathcal{A}$ is generated by a single exceptional object $a$, the embedding functor is defined by
$$i_{\mathcal{A}}:\mathrm{D}^{b}(\mathrm{k})\longrightarrow\mathcal{D}\quad (-)\mapsto a\otimes (-)$$
and its left or right adjoint:
$$i^{*}_{\mathcal{A}}:\mathcal{D}\longrightarrow\mathrm{D}^{b}(\mathrm{k})\quad (-)\mapsto \mathrm{Ext}^{*}((-),a)^{\vee}$$
$$i^{!}_{\mathcal{A}}:\mathcal{D}\longrightarrow\mathrm{D}^{b}(\mathrm{k})\quad (-)\mapsto \mathrm{Ext}^{*}(a,(-))$$
hence we have the mutation distinguished triangle concerning a single exceptional object $a$,
$$\mathbb{R}_{a}(-)\longrightarrow(-)\longrightarrow  \mathrm{Ext}^{*}((-),a)^{\vee}\otimes a\longrightarrow\mathbb{R}_{a}(-)[1]$$
$$\mathbb{L}_{a}(-)[-1]\longrightarrow\mathrm{Ext}^{*}(a,(-))\otimes a\longrightarrow (-) \longrightarrow\mathbb{L}_{a}(-)$$
\begin{definition}[Left and Right dual]\label{def:dual}
If $\mathbf{E}:=(E_{0},\cdots,E_{n})$ is a full exceptional collection in a triangulated category $\mathcal{D}$, we have its left or right dual collection by giving
$$\mathbf{E}^{\vee}=(E_{n}^{\vee},\cdots,E_{0}^{\vee})\quad\quad^{\vee}\mathbf{E}=(^{\vee}E_{n},\cdots,^{\vee}E_{0})$$
where
$$E_{i}^{\vee}:=\mathbb{L}_{E_{0}}\cdots\mathbb{L}_{E_{i-1}}E_{i}\quad\quad^{\vee}E_{i}:=\mathbb{R}_{E_{n}}\cdots\mathbb{R}_{E_{i+1}}E_{i}$$
\end{definition}

By our definition, the left of the right dual collection is again a full exceptional collection of $\mathcal{D}$, because the transformations between them are determined by a set of braid groups, see \cite[Section 2.4]{GK}.

\begin{proposition}\label{Prop:dual.collection}
We have
$$\mathrm{Ext}^{*}(E_{i},E_{j}^{\vee})=\mathrm{k}\cdot\delta_{i,j}\text{\quad and\quad}
\mathrm{Ext}^{*}(^{\vee}E_{i},E_{j})=\mathrm{k}\cdot\delta_{i,j}$$
for any $i,j$ in the range.
\end{proposition}
\begin{proof}
If $i>j$, noticing by definition $E_{j}^{\vee}$ is contained in $\langle E_{0}...,E_{j}\rangle$, and from the semi-orthogonality of  $\mathbf{E}$ we have $\mathrm{Ext}^{*}(E_{i},E_{j}^{\vee})=0$.\\
If $i=j$,  noticing
$$\mathbb{L}_{\langle E_{0},..,E_{i-1}\rangle}E_{i}=\mathrm{Cone} \big( i_{\langle E_{0},..,E_{i-1}\rangle}i^{!}_{\langle E_{0},..,E_{i-1}\rangle}E_{i}\longrightarrow E_{i}\big)$$
so there is $\mathrm{Ext}^{*}(E_{i},E_{i}^{\vee})=\mathrm{Ext}^{*}(E_{i},\mathbb{L}_{E_{0}}\cdots\mathbb{L}_{E_{i-1}}E_{i})=\mathrm{Ext}^{*}(E_{i},E_{i})=\mathrm{k}$.\\
If $i<j$,  noticing
$$\mathbb{L}_{\langle E_{0},..,E_{i-1}\rangle}F=\mathrm{Cone} \big( i_{\langle E_{0},..,E_{i-1}\rangle}i^{!}_{\langle E_{0},..,E_{i-1}\rangle}F\longrightarrow F\big)$$
where $F:=\mathbb{L}_{\langle E_{i},..,E_{j-1}\rangle}E_{j}$, so we have $\mathrm{Ext}^{*}(E_{i},E_{j}^{\vee})= \mathrm{Ext}^{*}(E_{i},F)$. On the other hand side
$$\mathrm{Ext}^{*}(E_{i},F)= \mathrm{Ext}^{*}(\mathbb{L}_{\langle E_{0},..,E_{i-1}\rangle}E_{i},\mathbb{L}_{\langle E_{0},..,E_{i-1}\rangle}F)=\mathrm{Ext}^{*}(E_{i}^{\vee},E_{j}^{\vee})=0$$
the first equation is due to the fact  $\mathrm{Ext}^{*}(E,\mathbb{L}_{E}(-))=0$ for any exceptional object $E$, and by a simple computation expanding all terms using mutation triangles we can see
$$\mathrm{Ext}^{*}(E_{i},F)=\mathrm{Ext}^{*}(\mathbb{L}_{E_{k}}E_{i},\mathbb{L}_{E_{k}}F)$$
for any $k<i$.
\end{proof}
\begin{proposition}
If $\mathbf{E}=(E_{0},\cdots,E_{n})$ is a full exceptional collection in category $\mathcal{D}$, then for any $i=0,\cdots,n$ we have
$$\mathrm{Ext}^{*}(^{\vee}E_{i},X)=\mathrm{Ext}^{*}(E_{i},\mathbb{L}_{E_{i+1}}\cdots\mathbb{L}_{E_{n}}X)$$
$$\mathrm{Ext}^{*}(X,E_{i}^{\vee})=\mathrm{Ext}^{*}(E_{i},\mathbb{L}_{E_{i+1}}\cdots\mathbb{L}_{E_{n}}X)^{\vee}$$
for any $X$ contained in $\mathcal{D}$.
\end{proposition}
\begin{proof}
See {\cite[Proposition 2.6.1]{GK}}. Because $\mathbf{E}=(E_{0},\cdots,E_{n})$ is full, $\mathbb{L}_{E_{0}}\cdots\mathbb{L}_{E_{n}}X=\mathbb{L}_{\mathcal{D}}X=0$.\\
Then we use induction for triangles traversing $j=0,\cdots,i$, by mutation triangle
$$\mathrm{Ext}^{*}(E_{j},\mathbb{L}_{E_{j+1}}\cdots\mathbb{L}_{E_{n}}X)\otimes E_{j}\rightarrow\mathbb{L}_{E_{j+1}}\cdots\mathbb{L}_{E_{n}}X\rightarrow\mathbb{L}_{E_{j}}\cdots\mathbb{L}_{E_{n}}X$$
so for $j=0,\cdots,i-1$, we have $\mathrm{Ext}^{*}(^{\vee}E_{i},\mathbb{L}_{E_{j+1}}\cdots\mathbb{L}_{E_{n}}X)=0$ by  Proposition \ref{Prop:dual.collection}.\\
If we apply $\mathrm{Ext}^{*}(^{\vee}E_{i},-)$ to the $j=i$ triangle, we have
$$\mathrm{Ext}^{*}(^{\vee}E_{i},X)=\mathrm{Ext}^{*}(E_{i},\mathbb{L}_{E_{i+1}}\cdots\mathbb{L}_{E_{n}}X)$$
since $\mathrm{Ext}^{*}(^{\vee}E_{i},E_{i})=\mathrm{k}$ by  Proposition \ref{Prop:dual.collection}.
\end{proof}

Let $\mathcal{B}$ be a full triangulated subcategory of $\mathcal{D}$, consider morphisms in $\mathcal{D}$ fitting into a distinguished triangle:
$$a\longrightarrow b\longrightarrow c\longrightarrow a[1]$$
with $c$ contained in $\mathcal{B}$, these morphisms actually form a multiplicative system $\Sigma (\mathcal{B})$, e.g.\cite{Orlh}.

\begin{definition}[Verdier quotient]
We define a quotient category $\mathcal{D}/\mathcal{B}$ as the localization of $\mathcal{D}$ for the multiplicative system $\Sigma (\mathcal{B})$.
\end{definition}

\subsection{Derived category}

For a variety (algebraic space) or a Deligne-Mumford stack $X$, we denote $\mathrm{D}^{b}(X)$ by its bounded derived category of coherent sheaves, see \cite[Chapter 2]{H}. We denote by $\mathrm{D}^{perf}(X)$ the full subcategory of $\mathrm{D}^{b}(X)$ consisting
of objects locally quasi-isomorphic to bounded complexes of finite-rank locally free sheaves, and $\mathrm{D}_{sg}(X)$ is the Verdier quotient of $\mathrm{D}^{b}(X)$ over $\mathrm{D}^{perf}(X)$.

\section{Derived category of blowing up}

\subsection{Dual collection on weighted projective space}

\begin{definition} For positive integers $a_{0},a_{1},...,a_{n}$ such that $\mathrm{gcd}\,(a_{0},...,a_{n})=1$, we consider weighted projective space $\mathds{P}(a_{0},a_{1},...,a_{n})$  defined as
$$\mathcal{P}roj (\mathrm{k}[x_{0},..,x_{n}])$$
for graded ring $\mathrm{k}[x_{0},..,x_{n}]$ with $\deg(x_{i}):=a_{i}$ for $0\leq i\leq n$, and its quotient stack lift
$$[\mathds{A}^{n+1}-\{0\}/\mathbb{G}^{a_{0},a_{1},...,a_{n}}_{m}]$$
with $\mathbb{G}_{m}$ acting on coordinates of $\mathds{A}^{n+1}$ with weight $a_{0},a_{1},...,a_{n}$, and denote it by $\mathcal{P}(a_{0},a_{1},...,a_{n})$.
\end{definition}

For weighted projective stack $\mathcal{P}(a_{0},a_{1},...,a_{n})$, for certain $i$ with $\mathrm{gcd}\,(a_{0},...,\hat{a_{i}},...,a_{n})=d\neq 1 $, the coarse moduli map has ramified structure over $\mathcal{P}(a_{0}/d,...,\hat{a_{i}},...,a_{n}/d)$ on the canonical stack of $\mathds{P}(a_{0},a_{1},...,a_{n})$. If $ \mathrm{gcd}\,(a_{0},...,\hat{a_{i}},...,a_{n})=1$ for any $i$, we call this sequence is \textit{well-formed}, then $\mathds{P}(a_{0},a_{1},...,a_{n})$ has only codimension 2 orbifold singularities and admits unique canonical stack  \cite[Example 7.27]{FMN}, so $\widetilde{\mathds{P}}(a_{0},a_{1},...,a_{n})=\mathcal{P}(a_{0},a_{1},...,a_{n})$.

\begin{example}
$\mathcal{P}(1,1,2)$ is defined to be the canonical stack of $\mathds{P}(1,1,2)$, at two atlas with trivial coarse moduli, while one atlas with cyclic $(1,1)$ weight quotient morphism $[\mathds{A}^{2}/\mu_{2}]\rightarrow \mathds{A}^{2}/ \mu_{2}$.
\end{example}

We should mention the following basic properties which are direct generalizations of projective spaces' cases.
\begin{lemma}\label{lemma:wpscoho}
\begin{enumerate}
  \item  $\mathrm{K}_{\mathcal{P}(a_{0},a_{1},...,a_{n})}=\mathcal{O}_{\mathcal{P}(a_{0},a_{1},...,a_{n})}(-\sum_{i=0}^{n} a_{i})$\\
\item \begin{equation}
 \mathrm{H}^{j}(\mathcal{P}(a_{0},a_{1},...,a_{n}),\mathcal{O}_{\mathcal{P}(a_{0},a_{1},...,a_{n})}(i)) =
    \begin{cases}
      \mathrm{dim}_{\mathrm{k}}\,\mathrm{k}[x_{0},..,x_{n}]^{(a_{0},..,a_{n})}_{(i)}  & \text{$j=0$}\\
     0 & \text{$0<j<n$}\\
       \mathrm{dim}_{\mathrm{k}}\,\mathrm{k}[x_{0},..,x_{n}]^{(a_{0},..,a_{n})}_{(-\sum_{i}a_{i})-i}  & \text{$j=n$}\\
    \end{cases}
\end{equation}
where $\mathrm{k}[x_{0},..,x_{n}]^{(a_{0},..,a_{n})}_{(i)}$ denotes the degree $i$ part of homogeneous ring $\mathrm{k}[x_{0},..,x_{n}]$ with homogeneous degree of $x_{i}$ as $a_{i}$.\end{enumerate}
\end{lemma}

This part explains two equations in \cite[Section 5, p.15]{KaEq}. Write
$G\simeq\mu_{a_0}\times\cdots\times\mu_{a_n}$. We represent a character by
$\chi=(\chi_0,\ldots,\chi_n)$, where $0\leq\chi_i<a_i$ is the integer
exponent of the chosen character generator of $\mu_{a_i}$, and set
$|\chi|:=\sum_{i=0}^n\chi_i$. For $-\chi$, we again take the least
nonnegative representatives before computing $|-\chi|$.
We write $\mathrm D_\chi$ for the union of coordinate divisors $(x_i=0)$
for which $\chi_i\ne0$.
\begin{lemma}\label{lemma:Kawambi}
If $a_{0},...,a_{n}$ is a well-formed positive integer sequence, consider the natural map
$$\overline{\pi}:[\mathds{P}^{n}/G]\longrightarrow\mathcal{P}(a_{0},...,a_{n})$$
with the action of $G$ on $\mathds{P}^{n}$ is given by
\begin{equation}\label{equ:actclassical}
  (\zeta_{a_{0}},..,\zeta_{a_{n}})\times[x_{0},..,x_{n}]\longmapsto[\zeta_{a_{0}}x_{0},..,\zeta_{a_{n}}x_{n}]
\end{equation}
and $\overline{\pi}$ is defined via $[(\mathbb{A}^{n+1}\setminus\{0\})/\mathbb{G}^{1^{n+1}}_{m}\times G]\longrightarrow[(\mathbb{A}^{n+1}\setminus\{0\})/\mathbb{G}^{a_{0},...,a_{n}}_{m}]$, we have:
$$\overline{\pi}_{*}(\mathcal{O}(i)_{[\mathds{P}^{n}/G]}\otimes\chi) =\mathcal{O}_{\mathcal{P}(a_{0},..,a_{n})}(i-|\chi|)\quad\text{and}\quad \overline{\pi}_{*}(\mathrm{\Omega}^{k}_{[\mathds{P}^{n}/G]}(i)\otimes\chi) =\mathrm{\Omega}^{k}_{\mathcal{P}(a_{0},..,a_{n})}(\mathrm{log}\,D_{\chi})(i-|\chi|)$$
for any $i$, $k$ and $\chi$.
\end{lemma}
\begin{proof}

\begin{enumerate}
\item\label{enu:wpscon1} The well-formed condition ensures that the coarse moduli map
$$\alpha:\mathcal{P}(a_{0},\ldots,a_{n})\longrightarrow\mathds{P}(a_{0},\ldots,a_{n})$$
is an isomorphism over the smooth locus of its coarse space, with complements of codimension at least two. Pushforward along $\alpha$ induces an equivalence between the categories of coherent reflexive sheaves on the stack and its coarse space. Its inverse sends $\mathcal F$ to $(\alpha^*\mathcal F)^{**}$, where the pullback is ordinary. Vector bundles on the stack correspond to those reflexive sheaves whose reflexive pullbacks are locally free. Reflexive sheaves on either side are uniquely determined by their restrictions to the common open subset.\\
\item In proof of \cite{Ca}, actually A.Canonaco considered a morphism
$$\pi:[\mathds{P}^{n}/G]\longrightarrow\mathds{P}(a_{0},...,a_{n})$$
where $G:=\mu_{a_{0}}\times...\times \mu_{a_{n}}$, and gave two important results in \cite[Corollary 1.4 and Lemma 2.6]{Ca} under pure (equivariant) module setting:
\begin{equation}\label{Kawkeyeq}
  \big(\Omega^{k}_{S}(i)\big)^{\chi}=\Omega^{k}_{S^{G}}(\mathrm{log}\,D_{\chi})(i-|\chi|)\footnote{\,Our $\Omega$ is $\bar{\Omega}$ under the notation of \cite[p.34]{Ca}.}
\end{equation}
for any $i$, $k$ and $\chi$, where $S$ is the ordinary polynomial ring with a group $G$ (co)action on it.\\
\item\label{enu:wpscon3} We have some \lq\lq non-commutative" description for these weighted projective spaces,
\begin{equation}\label{equ:serreequa}
 \approx_{s}:\mathrm{Gr}(S^{(a_{0},...,a_{n})})/\mathrm{Tor}\simeq\mathrm{Coh}(\mathcal{P}(a_{0},...,a_{n}))
\end{equation}
where $S^{(a_{0},...,a_{n})}$ is $k[x_{0},..,x_{n}]$ with $\deg\,x_{i}:=a_{i}$, e.g. \cite[Proposition 2.3]{AKO}. And Serre's functor $\lq\lq \sim_{s}"$ induces a morphism (which is in general not an equivalent) from it to $\mathrm{Coh}(\mathds{P}(a_{0},...,a_{n}))$, which is coincident with the pushforward along coarse moduli map $\alpha$.\\
\item We consider the following commutative diagram:
$$\begin{tikzcd}[column sep=0.5em]
 & {[\mathds{P}^{n}/G] }\arrow[rr,"\pi"]\arrow[dr,"\overline{\pi}"]&& \mathds{P}(a_{0},...,a_{n}) \\
&& \mathcal{P}(a_{0},...,a_{n})\arrow[ur,"\alpha"]&&
\end{tikzcd}$$
we assume $\overline{\pi}^{\chi}_{*}\big(\Omega^{k}_{S}(i)^{\sim_{s}}\big)=N^{\approx_{s}}$ for certain $G$ graded module $N$, it can be easily seen than $N^{\approx_{s}}$ is locally free, e.g.\cite[Proposition 1.1]{Ca}.  Then pushforward it along $\alpha$ and compare with \cite[Lemma 1.3]{Ca} and (\ref{Kawkeyeq}), we have
$$\alpha_{*}N^{\approx_{s}}=\pi^{\chi}_{*}\big(\Omega^{k}_{S}(i)^{\sim_{s}}\big)=\Omega^{k}_{S^{G}}(\mathrm{log}\,D_{\chi})(i-|\chi|)^{\sim_{s}}$$
On the common smooth locus, the preceding identity identifies $N^{\approx_s}$ with
$\mathrm{\Omega}^{k}_{\mathcal{P}(a_{0},\ldots,a_{n})}(\mathrm{log}\,D_{\chi})(i-|\chi|)$.
The latter is locally free because $D_{\chi}$ is a union of coordinate divisors with simple normal crossings on the smooth stack. Both sheaves are therefore reflexive, and uniqueness of reflexive extension in (\ref{enu:wpscon1}) gives
$$N^{\approx_s}\simeq\mathrm{\Omega}^{k}_{\mathcal{P}(a_{0},\ldots,a_{n})}(\mathrm{log}\,D_{\chi})(i-|\chi|).$$
This proves the asserted pushforward formula.
\end{enumerate}
\end{proof}

In fact for the Lemma \ref{lemma:Kawambi}, we can remove the well-formed condition.

\begin{theorem}[Beilinson type resolution {\cite[Section 5]{KaEq}}]\label{thm:beilinson}
For any positive integers $a_{0},...,a_{n}$ with $\mathrm{gcd}\,(a_{0},...,a_{n})=1$ and weighted projective stack $\mathcal{P}:=\mathcal{P}(a_{0},a_{1},...,a_{n})$, we have the following complex $\mathbf{e}$ to diagonal $\mathcal{O}_{\mathrm{\Delta}}$ on $\mathcal{P}\times\mathcal{P}$:\\
$$0\longrightarrow
\bigoplus_{\chi}
\mathcal{O}_{\mathcal{P}}(-n-|\chi|)\boxtimes \mathrm{\Omega}^{n}(\mathrm{log}\, D_{-\chi})(n-|-\chi|)...\longrightarrow $$
$$\bigoplus_{\chi}
\mathcal{O}_{\mathcal{P}}(-i-|\chi|)\boxtimes \mathrm{\Omega}^{i}(\mathrm{log}\, D_{-\chi})(i-|-\chi|)...\longrightarrow \bigoplus_{\chi}
\mathcal{O}_{\mathcal{P}}(-|\chi|)\boxtimes \mathcal{O}_{\mathcal{P}}(-|-\chi|)\longrightarrow 0$$
for any object $a$ in $\mathrm{D}^{b}(\mathcal{P})$, we have $p_{1*}(p_{2}^{*}a\otimes \mathbf{e})=p_{2*}(p_{1}^{*}a\otimes \mathbf{e})=a$, where $p_{i}$  $i=1,2$ is projection of $\mathcal{P}\times\mathcal{P}$ to any component.
\end{theorem}

\begin{proof}
If we assume $ \mathrm{gcd}\,(a_{0},...,\hat{a_{i}},...,a_{n})=1$ for any $i$ for positive integers $a_{0},..,a_{n}$ (\textit{well formed}), it follows the result of {\cite[Section 5]{KaEq}}.\\

If for certain $i$, $ \mathrm{gcd}\,(a_{0},...,\hat{a_{i}},...,a_{n})\neq1$, we enlarge the ambient space of quotient as:
$$\begin{tikzcd}[ampersand replacement=\&]
\&{[\mathds{P}^{n}/G]}\arrow[r,"\pi"]\arrow[d,"i"] \&\mathcal{P}(a_{0},..,a_{n})\arrow[d,"i'"] \\
 \&{[\mathds{P}^{n+1}/G]} \arrow[r,"\overline{\pi}"] \&  \mathcal{P}(a_{0},..,a_{n},1)
\end{tikzcd}$$
where $G:=\mu_{a_{0}}\times...\times \mu_{a_{n}}$. So the new sequence of $a_{0},..,a_{n},1$ is well-formed by our assumption, especially by the previous Lemma \ref{lemma:Kawambi} we have
$$\overline{\pi}_{*}(\mathcal{O}(i)_{[\mathds{P}^{n+1}/G]}\otimes\chi) =\mathcal{O}_{\mathcal{P}(a_{0},..,a_{n},1)}(i-|\chi|)\quad\text{and}\quad \overline{\pi}_{*}(\mathrm{\Omega}^{k}_{[\mathds{P}^{n+1}/G]}(i)\otimes\chi) =\mathrm{\Omega}^{k}_{\mathcal{P}(a_{0},..,a_{n},1)}(\mathrm{log}\,D_{\chi})(i-|\chi|)$$
for any $i$, $k$ and $\chi$, noticing short exact sequences:
$$0\longrightarrow\mathcal{O}_{[\mathds{P}^{n+1}/G]}(-1)\longrightarrow\mathcal{O}_{[\mathds{P}^{n+1}/G]}\longrightarrow\mathcal{O}_{[\mathds{P}^{n}/G]}\longrightarrow0$$
$$0\longrightarrow\mathcal{O}_{\mathcal{P}(a_{0},..,a_{n},1)}(-1)\longrightarrow\mathcal{O}_{\mathcal{P}(a_{0},..,a_{n},1)}\longrightarrow\mathcal{O}_{\mathcal{P}(a_{0},..,a_{n})}\longrightarrow0$$
we tensor the first exact sequence with $\chi$ then  pushforward along $\overline{\pi}$ and compare with the second one, it's not difficult to see $\pi_{*}(\mathcal{O}_{[\mathds{P}^{n}/G]}\otimes\chi) =\mathcal{O}_{\mathcal{P}(a_{0},..,a_{n})}(-|\chi|)$ by Lemma \ref{lemma:Kawambi}.\\

Then notice for any $\chi$, $D_{\chi}$ is transversal with $\mathcal{P}(a_{0},..,a_{n})$ in $\mathcal{P}(a_{0},..,a_{n},1)$, and relative K$\ddot{a}$hler differential exact sequences by a local equivariant descent argument:
$$0\longrightarrow\mathcal{O}_{[\mathds{P}^{n}/G]}(-1)\longrightarrow\mathrm{\Omega}_{[\mathds{P}^{n+1}/G]}|_{[\mathds{P}^{n}/G]}\longrightarrow\mathrm{\Omega}_{[\mathds{P}^{n}/G]}\longrightarrow0$$
$$0\longrightarrow\mathcal{O}_{\mathcal{P}(a_{0},..,a_{n})}(-1)\longrightarrow\mathrm{\Omega}_{\mathcal{P}(a_{0},..,a_{n},1)}(\mathrm{log}\,D_{\chi})|_{\mathcal{P}(a_{0},..,a_{n})}\longrightarrow\mathrm{\Omega}_{\mathcal{P}(a_{0},..,a_{n})}(\mathrm{log}\,D_{\chi}\cap\mathcal{P}(a_{0},..,a_{n}))\longrightarrow0$$
similarly we tensor the first exact sequence with $\chi$ then  pushforward along $\overline{\pi}$ and compare with the second one, we can see $\overline{\pi}_{*}(\mathrm{\Omega}_{[\mathds{P}^{n}/G]}\otimes\chi) =\mathrm{\Omega}_{\mathcal{P}(a_{0},..,a_{n})}(\mathrm{log}\,D_{\chi})(-|\chi|)$ by Lemma \ref{lemma:Kawambi}.\\

By taking the wedge product of the previous two exact sequences since each term is locally free, we have:
$$0\longrightarrow\mathcal{O}_{[\mathds{P}^{n}/G]}(-1)\otimes\mathrm{\Omega}^{k-1}_{[\mathds{P}^{n}/G]}\longrightarrow\mathrm{\Omega}^{k}_{[\mathds{P}^{n+1}/G]}|_{[\mathds{P}^{n}/G]}\longrightarrow\mathrm{\Omega}^{k}_{[\mathds{P}^{n}/G]}\longrightarrow0$$
$$\begin{aligned}0\longrightarrow\mathcal{O}_{\mathcal{P}(a_{0},..,a_{n})}(-1)\otimes\mathrm{\Omega}^{k-1}_{\mathcal{P}(a_{0},..,a_{n})}(\mathrm{log}\,D_{\chi})&\longrightarrow\mathrm{\Omega}^{k}_{\mathcal{P}(a_{0},..,a_{n},1)}(\mathrm{log}\,D_{\chi})|_{\mathcal{P}(a_{0},..,a_{n})}\\&\longrightarrow\mathrm{\Omega}^{k}_{\mathcal{P}(a_{0},..,a_{n})}(\mathrm{log}\,D_{\chi})\longrightarrow0\end{aligned}$$
by induction if $\overline{\pi}_{*}(\mathrm{\Omega}^{k-1}_{[\mathds{P}^{n}/G]}(i)\otimes\chi) =\mathrm{\Omega}^{k-1}_{\mathcal{P}(a_{0},..,a_{n})}(\mathrm{log}\,D_{\chi})(i-|\chi|)$  holds for any $i$ and $\chi$, then tensor the first exact sequence with $\chi$ then  pushforward along $\overline{\pi}$ and compare with the second one, we have $\overline{\pi}_{*}(\mathrm{\Omega}^{k}_{[\mathds{P}^{n}/G]}(i)\otimes\chi) =\mathrm{\Omega}^{k}_{\mathcal{P}(a_{0},..,a_{n})}(\mathrm{log}\,D_{\chi})(i-|\chi|)$, for any $i$ and $\chi$.\\

The rest of the proof is consistent with the proof in {\cite[Section 5]{KaEq}}, first, consider the equivariant version of Beilinson type resolution on $[\mathds{P}^{n}/G]$, then pushforward it along $\pi$ gets the result.

\end{proof}

\begin{corollary}
The weighted projective stack $\mathcal{P}(a_{0},...,a_{n})$ has the semi-orthogonal decomposition:
$$\mathrm{D}^{b}(\mathcal{P})=\langle \mathcal{O}_{\mathcal{P}}(-\sum_{i}a_{i}+1),...,\mathcal{O}_{\mathcal{P}}(-1),\mathcal{O}_{\mathcal{P}}\rangle$$
\end{corollary}

Like the construction for projective space, we can construct dual collection using mutations:
\begin{corollary}\label{cor:sodofwps}
The dual semi-orthogonal decomposition of $\mathcal{P}(a_{0},...,a_{n})$ is:
$$\mathrm{D}^{b}(\mathcal{P})=\langle \mathcal {D}_{\sum_{i}a_{i}-1},...,\mathcal {D}_{1},\mathcal D_{0}\rangle$$
where $\mathcal D_{i}:=\mathbb{L}_{\langle\mathcal{O}_{\mathcal{P}},...,\mathcal{O}_{\mathcal{P}}(i-1)\rangle}\mathcal{O}_{\mathcal{P}}(i)$ for $0\leq i\leq \sum_{i}a_{i}-1$\footnote{\, For our convenience, we assume that $\mathcal D_{i}$ is $0$ for $i$ exceeds this range.}.
\end{corollary}
\begin{proof}
See Definition \ref{def:dual} or \cite[Section 2.4]{GK}.
\end{proof}

The following property holds for general dual collection of a full exceptional collection, see Proposition \ref{Prop:dual.collection}.

\begin{corollary}\label{cor:dualvanishing}
For $0\leq i,j<\sum_{k=0}^{n} a_{k}$,
\begin{equation}
 \mathrm{Ext}^{*}(\mathcal{O}_{\mathcal{P}}(i),\mathcal{D}_{j}) =
    \begin{cases}
      \mathrm{k}  & \text{if $i=j$}\\
      0 & \text{else}\\
    \end{cases}
\end{equation}
where $\mathrm{k}$ is a complex concentrates on degree $0$.
\end{corollary}

\begin{corollary}\label{cor:dualdes}
For any $1\leq k\leq \sum_{i}a_{i}-1$, $\mathcal{D}_{k}$ is quasi-isomorphic to a complex $\mathbf{D}_{k}$:
\begin{equation}\label{equ:dualcoll}
\mathrm{\Omega}^{k}(k)\longrightarrow\bigoplus_{\chi\big||\chi|=1}\mathrm{\Omega}^{k-1}(\mathrm{log}\, D_{-\chi})(k-1-|-\chi|)...\longrightarrow
\end{equation}
$$\bigoplus_{\chi\big||\chi|=i}\mathrm{\Omega}^{k-i}(\mathrm{log}\, D_{-\chi})(k-i-|-\chi|)\longrightarrow\bigoplus_{\chi\big||\chi|=\sum_{i}a_{i}-n-1}\mathrm{\Omega}^{k-(\sum_{i}a_{i}-n-1)}(\mathrm{log}\, D_{-\chi})(k-(\sum_{i}a_{i}-n-1)-|-\chi|)$$
where term $ \mathrm{\Omega}^{k}(k)$ lies in degree $-k$, and each morphism in the complex is non-trivially and uniquely induced by the quotient of morphism on projective space.
\end{corollary}
\begin{proof}
We prove it by induction. If $k=1$, $\mathcal{D}_{1}=\mathbb{L}_{\mathcal{O}_{\mathcal{P}}}\mathcal{O}_{\mathcal{P}}(1)$, notice that $\mathcal{O}_{\mathcal{P}}(1)$ is isomorphic to the following complex if we test Fourier-Mukai transformation of $\mathcal{O}_{\mathcal{P}}(1)$ to complex in Theorem \ref{thm:beilinson} i.e. $p_{2*}(p_{1}^{*}\mathcal{O}_{\mathcal{P}}(1)\otimes \mathbf{e})$:
$$0\longrightarrow
\bigoplus_{\chi}
\mathrm{H}^{*}(\mathcal{O}_{\mathcal{P}}(-n-|\chi|+1))\otimes \mathrm{\Omega}^{n}(\mathrm{log}\, D_{-\chi})(n-|-\chi|)...\longrightarrow $$
$$\bigoplus_{\chi}
\mathrm{H}^{*}(\mathcal{O}_{\mathcal{P}}(-i-|\chi|+1))\otimes \mathrm{\Omega}^{i}(\mathrm{log}\, D_{-\chi})(i-|-\chi|)...\longrightarrow \bigoplus_{\chi}
\mathrm{H}^{*}(\mathcal{O}_{\mathcal{P}}(-|\chi|+1))\otimes \mathcal{O}_{\mathcal{P}}(-|-\chi|)\longrightarrow 0$$
and $|\chi|\leq \sum_{i=0}^{n}(a_{i}-1)$ we have $-\sum_{i}a_{i}<-n-|\chi|+1$, by cohomology computation using Lemma \ref{lemma:wpscoho}
$$\mathrm{H}^{*}(\mathcal{O}_{\mathcal{P}}(-i-|\chi|+1))=\mathrm{H}^{0}(\mathcal{O}_{\mathcal{P}}(-i-|\chi|+1))$$
so we only need to consider terms whence $i+|\chi|=0,1$, the complex is just
$$\mathrm{\Omega}_{\mathcal{P}}^{1}(1)\longrightarrow \bigoplus_{i|a_{i}>1}\mathcal{O}_{\mathcal{P}}(-a_{i}+1)\bigoplus \mathrm{H}^{*}(\mathcal{P},\mathcal{O}_{\mathcal{P}}(1))\otimes\mathcal{O}_{\mathcal{P}}$$
hence by the definition of $\mathcal{D}_{1}$ which fits in distinguished triangle $\mathrm{H}^{*}(\mathcal{P},\mathcal{O}_{\mathcal{P}}(1))\otimes\mathcal{O}_{\mathcal{P}}\longrightarrow\mathcal{O}_{\mathcal{P}}(1)\longrightarrow\mathcal{D}_{1}$, there is $\mathcal{D}_{1}\simeq \big(\mathrm{\Omega}_{\mathcal{P}}^{1}(1)\longrightarrow \bigoplus_{i|a_{i}>1}\mathcal{O}_{\mathcal{P}}(-a_{i}+1)\big)$.\\

Then let's consider the morphism in this complex, recall the construction of Beilinson type resolution for weighted projective stack in Theorem \ref{thm:beilinson}, which is constructed from the equivariant version of Beilinson type resolution for projective space. So this morphism comes from the pushforward of a non-trivial morphism in $\mathrm{Hom}(\mathrm{\Omega}_{[\mathds{P}^{n}/G]}^{1}(1),\mathcal{O}_{[\mathds{P}^{n}/G]}\otimes\nu^{-1}_{j})\simeq[(\bigoplus_{i}\mathrm{k}\otimes\nu_{i})\otimes\nu_{j}^{-1}]^{G}\simeq\mathrm{k}$ by \cite[Lemma 2]{Bei}
for any $j$. The assertion holds for $k=1$.\\

Since $\mathcal{O}_{\mathcal{P}}(k)$ (for $1\leq k\leq \sum_{i}a_{i}-1$) is quasi-isomorphic to the following complex by similarly Fourier-Mukai testing for $\mathcal{O}_{\mathcal{P}}(k)$ to the complex in Theorem \ref{thm:beilinson} i.e. $p_{2*}(p_{1}^{*}\mathcal{O}_{\mathcal{P}}(k)\otimes \mathbf{e})$:
\begin{equation}\label{eq:conek}
0\longrightarrow
\bigoplus_{\chi}
\mathrm{H}^{*}(\mathcal{O}_{\mathcal{P}}(-n-|\chi|+k))\otimes \mathrm{\Omega}^{n}(\mathrm{log}\, D_{-\chi})(n-|-\chi|)...\longrightarrow \end{equation}
$$\bigoplus_{\chi}
\mathrm{H}^{*}(\mathcal{O}_{\mathcal{P}}(-i-|\chi|+k))\otimes \mathrm{\Omega}^{i}(\mathrm{log}\, D_{-\chi})(i-|-\chi|)...\longrightarrow \bigoplus_{\chi}
\mathrm{H}^{*}(\mathcal{O}_{\mathcal{P}}(-|\chi|+k))\otimes \mathcal{O}_{\mathcal{P}}(-|-\chi|)\longrightarrow 0$$
noticing $-\sum_{i}a_{i}<-n-|\chi|+k $, by cohomology computation using Lemma \ref{lemma:wpscoho} we have
$$\mathrm{H}^{*}(\mathcal{O}_{\mathcal{P}}(-i-|\chi|+k))=\mathrm{H}^{0}(\mathcal{O}_{\mathcal{P}}(-i-|\chi|+k))$$
and we only need to consider terms whence $i+|\chi|\leq k$, we can check that the subcomplex of (\ref{eq:conek})  with terms satisfying $i+|\chi|=i'$ is exactly $\mathrm{H}^{0}(\mathcal{O}_{\mathcal{P}}(k-i'))\otimes \mathbf{D}_{i'}$ according to the expression of (\ref{equ:dualcoll}), so we can consider complex (\ref{eq:conek}) as the iterated extensions of such subcomplexes:

\begin{equation}\label{DD1complex}\begin{tikzcd}[column sep=0.5em]
 & 0 \arrow{rr}&& E_{0}\arrow{dl}\arrow{rr}&& E_{1}  \arrow{dl} \\
&& V_{k}^{0}\otimes\mathbf{D}_{0}\arrow[ul,dashed,"\Delta"]&& V_{k}^{1}\otimes\mathbf{D}_{1}\arrow[ul,dashed,"\Delta"]
\end{tikzcd}
...\begin{tikzcd}[column sep=0.5em]
 &     E_{k-1}\arrow{rr}&& \mathcal{O}_{\mathcal{P}}(k)\arrow{dl} \\
&& V_{k}^{k}\otimes\mathbf{D}_{k}\arrow[ul,dashed,"\Delta"]
\end{tikzcd}\end{equation}
where $V_{k}^{i'}:=\mathrm{H}^{0}(\mathcal{P},\mathcal{O}_{\mathcal{P}}(k-i'))$.\\

While also notice for any $k$, $\mathcal{D}_{k}:=\mathbb{L}_{\langle\mathcal{O}_{\mathcal{P}},...,\mathcal{O}_{\mathcal{P}}(k-1)\rangle}\mathcal{O}_{\mathcal{P}}(k)=\mathbb{L}_{\langle\mathcal{D}_{k-1},...,\mathcal{D}_{0}\rangle}\mathcal{O}_{\mathcal{P}}(k)$, and by similar argument in Proposition \ref{Prop:dual.collection} since $\mathcal{D}_{i}=\mathbb{L}_{\langle\mathcal{D}_{i-1},...,\mathcal{D}_{0}\rangle}\mathcal{O}_{\mathcal{P}}(i)$ and $\mathrm{Ext}^{*}(\mathcal{O}_{\mathcal{P}}(i),\mathcal{D}_{\lambda})=0$ for $0\leq \lambda<i$ by Corollary \ref{cor:dualvanishing}, we have for any $0\leq i\leq k$ $$\mathrm{Ext}^{*}(\mathcal{D}_{i},\mathbb{L}_{\langle\mathcal{D}_{i-1},...,\mathcal{D}_{0}\rangle}\mathcal{O}_{\mathcal{P}}(k))=\mathrm{Ext}^{*}(\mathcal{O}_{\mathcal{P}}(i),\mathcal{O}_{\mathcal{P}}(k))=V_{k}^{i}$$ and $\mathbb{L}_{\langle\mathcal{D}_{i},...,\mathcal{D}_{0}\rangle}\mathcal{O}_{\mathcal{P}}(k)=\mathbb{L}_{\mathcal{D}_{i}}\mathbb{L}_{\langle\mathcal{D}_{i-1},...,\mathcal{D}_{0}\rangle}\mathcal{O}_{\mathcal{P}}(k)$ fits in the mutation distinguished triangle:
$$V^{i}_{k}\otimes\mathcal{D}_{i}\longrightarrow\mathbb{L}_{\langle\mathcal{D}_{i-1},...,\mathcal{D}_{0}\rangle}\mathcal{O}_{\mathcal{P}}(k)\longrightarrow\mathbb{L}_{\langle\mathcal{D}_{i},...,\mathcal{D}_{0}\rangle}\mathcal{O}_{\mathcal{P}}(k)$$
if we traverse $i$ from $0$ to $k$, a repeated extension of the above distinguished triangles by Octahedron axiom in triangulated category induces:

\begin{equation}\label{DD'1complex}\begin{tikzcd}[column sep=0.5em]
 & 0 \arrow{rr}&& \mathcal{E}_{0}\arrow{dl}\arrow{rr}&& \mathcal{E}_{1}  \arrow{dl} \\
&& V_{k}^{0}\otimes\mathcal{D}_{0}\arrow[ul,dashed,"\Delta"]&& V_{k}^{1}\otimes\mathcal{D}_{1}\arrow[ul,dashed,"\Delta"]
\end{tikzcd}
...\begin{tikzcd}[column sep=0.5em]
 &     \mathcal{E}_{k-1}\arrow{rr}&& \mathcal{O}_{\mathcal{P}}(k)\arrow{dl} \\
&& V_{k}^{k}\otimes\mathcal{D}_{k}\arrow[ul,dashed,"\Delta"]
\end{tikzcd}\end{equation}
using our assumption of induction on $\mathcal{D}_{i}\simeq\mathbf{D}_{i}$ for any $i\leq k-1$ and comparing the above two filtrations (\ref{DD1complex}) and (\ref{DD'1complex}) for $\mathcal{O}_{\mathcal{P}}(k)$ by Corollary \ref{cor:sodofwps}, we can see $\mathcal{D}_{k}\simeq\mathbf{D}_{k}$  which gives the result. Similarly all nontrivial morphisms in expression (\ref{equ:dualcoll}) of $\mathbf{D}_{k}$ come from the pushforward of morphisms among $$\mathrm{Hom}(\mathrm{\Omega}_{[\mathds{P}^{n}/G]}^{i}(i),\mathrm{\Omega}_{[\mathds{P}^{n}/G]}^{i-1}(i-1)\otimes\nu^{-1}_{j})\simeq[(\bigoplus_{i}\mathrm{k}\otimes\nu_{i})\otimes\nu_{j}^{-1}]^{G}\simeq\mathrm{k}$$ by \cite[Lemma 2]{Bei}, they are non-trivial and unique.
\end{proof}

\begin{corollary}\label{cor:okfmtest}
For $k>0$, $\mathcal{O}_{\mathcal{P}}(-k)$ has a distinguished triangle filtration:
\[\begin{adjustbox}{max width=0.98\linewidth}\(\begin{tikzcd}[column sep=0.5em]
 & 0 \arrow{rr}&& E_{k-1}\arrow{dl} \\
&& \mathcal{D}_{\sum_{i=0}^{n}a_{i}-k}\arrow[ul,dashed,"\Delta"]
\end{tikzcd}...
\begin{tikzcd}[column sep=0.5em]
 & E_{\lambda} \arrow{rr}&& E_{\lambda-1}\arrow{dl} \\
&& V_{\lambda}^{k}\otimes\mathcal{D}_{\sum_{i=0}^{n}a_{i}-\lambda}\arrow[ul,dashed,"\Delta"]
\end{tikzcd}...
\begin{tikzcd}[column sep=0.5em]
 & E_{1} \arrow{rr}&& \mathcal{O}_{\mathcal{P}}(-k)[n]\arrow{dl} \\
&& V_{1}^{k}\otimes\mathcal{D}_{\sum_{i=0}^{n}a_{i}-1}\arrow[ul,dashed,"\Delta"]
\end{tikzcd}\)\end{adjustbox}\]
where $V_{\lambda}^{k}:=\mathrm{H}^{0}(\mathcal{P},\mathcal{O}_{\mathcal{P}}(\lambda+k-\sum_{i=0}^{n}a_{i}))$.
\end{corollary}

\begin{proof}
If we test Fourier-Mukai transformation of $\mathcal{O}_{\mathcal{P}}(-k)$ with respect to complex in Theorem \ref{thm:beilinson} i.e. $p_{2*}(p_{1}^{*}\mathcal{O}_{\mathcal{P}}(-k)\otimes \mathbf{e})$:
$$0\longrightarrow
\bigoplus_{\chi}
\mathrm{H}^{*}(\mathcal{O}_{\mathcal{P}}(-n-|\chi|-k))\otimes \mathrm{\Omega}^{n}(\mathrm{log}\, D_{-\chi})(n-|-\chi|)...\longrightarrow $$
$$\bigoplus_{\chi}
\mathrm{H}^{*}(\mathcal{O}_{\mathcal{P}}(-i-|\chi|-k))\otimes \mathrm{\Omega}^{i}(\mathrm{log}\, D_{-\chi})(i-|-\chi|)...\longrightarrow \bigoplus_{\chi}
\mathrm{H}^{*}(\mathcal{O}_{\mathcal{P}}(-|\chi|-k))\otimes \mathcal{O}_{\mathcal{P}}(-|-\chi|)\longrightarrow 0$$
noticing $-i-|\chi|-k<0$,
$$\mathrm{H}^{*}(\mathcal{P},\mathcal{O}_{\mathcal{P}}(-i-|\chi|-k))=\mathrm{H}^{n}(\mathcal{P},\mathcal{O}_{\mathcal{P}}(-i-|\chi|-k))[-n]=\mathrm{H}^{0}(\mathcal{P},\mathcal{O}_{\mathcal{P}}(i+|\chi|+k-\sum_{i=0}^{n}a_{i}))[-n]$$ by Lemma \ref{lemma:wpscoho}, so $\mathcal{O}_{\mathcal{P}}(-k)[n]$ is isomorphic to the following complex
$$0\longrightarrow
\bigoplus_{\chi}
\mathrm{H}^{0}(\mathcal{P},\mathcal{O}_{\mathcal{P}}(n+|\chi|+k-\sum_{i=0}^{n}a_{i}))\otimes \mathrm{\Omega}^{n}(\mathrm{log}\, D_{-\chi})(n-|-\chi|)...\longrightarrow $$
\begin{multline*}
\bigoplus_{\chi}
\mathrm{H}^{0}(\mathcal{P},\mathcal{O}_{\mathcal{P}}(i+|\chi|+k-\sum_{i=0}^{n}a_{i}))\otimes \mathrm{\Omega}^{i}(\mathrm{log}\, D_{-\chi})(i-|-\chi|)...\longrightarrow\\
\bigoplus_{\chi}
\mathrm{H}^{0}(\mathcal{P},\mathcal{O}_{\mathcal{P}}(|\chi|+k-\sum_{i=0}^{n}a_{i}))\otimes \mathcal{O}_{\mathcal{P}}(-|-\chi|)\longrightarrow 0
\end{multline*}
observing the items with $\lambda=i+|\chi|$ is a constant in the above complex, by Corollary \ref{cor:dualdes} it is just $$\mathrm{H}^{0}(\mathcal{P},\mathcal{O}_{\mathcal{P}}(\lambda+k-\sum_{i=0}^{n}a_{i}))\otimes \mathcal{D}_{\lambda}$$
up to a shift and $\mathcal{O}_{\mathcal{P}}(-k)$ has a distinguished triangle filtration:
\[\begin{adjustbox}{max width=0.98\linewidth}\(\begin{tikzcd}[column sep=0.5em]
 & 0 \arrow{rr}&& E_{k-1}\arrow{dl} \\
&& \mathcal{D}_{\sum_{i=0}^{n}a_{i}-k}\arrow[ul,dashed,"\Delta"]
\end{tikzcd}...
\begin{tikzcd}[column sep=0.5em]
 & E_{\lambda} \arrow{rr}&& E_{\lambda-1}\arrow{dl} \\
&& V_{\lambda}^{k}\otimes\mathcal{D}_{\sum_{i=0}^{n}a_{i}-\lambda}\arrow[ul,dashed,"\Delta"]
\end{tikzcd}...
\begin{tikzcd}[column sep=0.5em]
 & E_{1} \arrow{rr}&& \mathcal{O}_{\mathcal{P}}(-k)[n]\arrow{dl} \\
&& V_{1}^{k}\otimes\mathcal{D}_{\sum_{i=0}^{n}a_{i}-1}\arrow[ul,dashed,"\Delta"]
\end{tikzcd}\)\end{adjustbox}\]
where $V_{\lambda}^{k}:=\mathrm{H}^{0}(\mathcal{P},\mathcal{O}_{\mathcal{P}}(\lambda+k-\sum_{i=0}^{n}a_{i}))$.

\end{proof}

\subsection{Derived category of weighted blow-up}

We assume $U$ is a $n$ dimensional smooth quasi-projective variety, which admits a $\mu_{s}$ cyclic covering
$$\rho: V\longrightarrow U$$
branched along a single prime smooth divisor $D$. Alternatively, there is a $\mu_{s}$ group action on V with pseudo-reflexive divisor $D_{V}$, since $\rho$ is ramified, we consider its stacky smoothing instead and denote it also by $\rho$:
$$\rho: [V/\mu_{s}]\longrightarrow U$$
under a local coordinate:
\begin{equation}\label{equ:actclassicals}
  (\zeta_{s})\times(x_{1},..,x_{n})\longmapsto(\zeta_{s}x_{1},x_{2},...,x_{n})
\end{equation}

In particular, we use the group action and character system different from the definition in \cite[Section 4.1]{KAC}, and we can see that such a character choice induces a $\mu_{s}$ action on the projective normal bundle which coincident with the action defined in (\ref{equ:actclassical}), we always follow this system in the rest part of the paper without mention it. Now we are expected to show how a skyscraper sheaf changes under pullback along $\rho$.\\

We will use the following lemma, which can be considered as another description of semi-orthogonal decomposition, see Definition \ref{def:SOD}.

\begin{lemma}\label{lem:projection}
If a triangulated category $\mathcal{D}$ admits a semi-orthogonal decomposition to subcomponents $\mathcal{A}_{i}$, where $0\leq i\leq n$, then $n$ steps of right mutations coincidences with the projection to the last component. That is:
$$\mathbb{R}_{\langle \mathcal{A}_{0},...,\mathcal{A}_{n-1} \rangle}\simeq i_{\mathcal{A}_{n}} i_{\mathcal{A}_{n}}^{!}$$
where $i_{\mathcal{A}_{n}}$ is the embedding functor of component $\mathcal{A}_{n}$ and $i_{\mathcal{A}_{n}}^{!}$ is its right adjoint.
\end{lemma}

\begin{lemma}\label{lem:pullbackskycraper}
Write $\mathcal G=[D/\mu_s]$, let $l:\mathcal G\hookrightarrow[V/\mu_s]$ be the reduced ramification divisor, and let $\pi_D:\mathcal G\to D$ be its projection. Put $\mathcal M_D=\mathcal{O}(\mathcal G)$ and $\mathcal I=\mathcal M_D^{-1}$. We use the character convention $\mathcal M_D|_{\mathcal G_x}\simeq\nu$.
\begin{enumerate}
\item\label{num:pullbackskycraper1}
There is an identification $[V/\mu_s]\simeq\sqrt[s]{D/U}$ and a semi-orthogonal decomposition
\[
\mathrm{D}^{b}([V/\mu_s])=
\left\langle\mathrm\Upsilon_{s-1},\ldots,\mathrm\Upsilon_1,
\rho^*\mathrm{D}^{b}(U)\right\rangle,
\]
where
\[
i_{\mathrm\Upsilon_i}(F)=l_*\pi_D^*F\otimes\mathcal M_D^i,
\qquad
i_{\mathrm\Upsilon_i}^*(F)=\pi_{D*}l^*(F\otimes\mathcal M_D^{-i}).
\]
\item\label{num:pullbackskycraper2}
For a geometric point $x\in D$, write $\mathcal G_x=l_*\pi_D^*\mathrm k_x$. The sheaf $\rho^*\mathrm k_x$ has the finite filtration by powers of $\mathcal I$, whose associated factors are
\[
\frac{\mathcal I^a\rho^*\mathrm k_x}
{\mathcal I^{a+1}\rho^*\mathrm k_x}
\simeq \mathcal G_x\otimes\nu^{-a},
\qquad 0\le a<s.
\]
Equivalently, put
\[
Q_j=\rho^*\mathrm k_x/\mathcal I^{s-j}\rho^*\mathrm k_x
\quad(0\le j\le s).
\]
Then $Q_0=\rho^*\mathrm k_x$, $Q_{s-1}=\mathcal G_x$, $Q_s=0$, and there are short exact sequences
\[
0\longrightarrow\mathcal G_x\otimes\nu^{j+1}
\longrightarrow Q_j\longrightarrow Q_{j+1}\longrightarrow0,
\qquad 0\le j<s-1.
\]
Each of these extensions is non-split and unique up to isomorphism and rescaling of its end terms. Their rotations, followed by $\mathcal G_x\to0\to\mathcal G_x[1]$, give the left Postnikov system. Its outgoing factors are $\mathcal G_x\otimes\nu[1],\ldots,\mathcal G_x\otimes\nu^{s-1}[1],\mathcal G_x[1]$.
These are associated factors of the indicated filtration, not projections of $\rho^*\mathrm k_x$ to the components $\mathrm\Upsilon_i$.
\end{enumerate}
\end{lemma}
\begin{proof}
The root identification and the decomposition in \ref{num:pullbackskycraper1} are the usual root-stack formula; see \cite[Proposition~6.1]{IU}. Here is its local mechanism. The projection $\rho$ is flat and $\rho_*\mathcal{O}=\mathcal{O}_U$, so $\rho^*$ is fully faithful. The self-intersection triangle for $l$ has conormal line $\mathcal M_D^{-1}$. Taking $\mu_s$-invariants therefore proves full faithfulness of the functors $i_{\mathrm\Upsilon_i}$ and their semi-orthogonality in the displayed order. For generation, the cone of $\rho^*\rho_*F\to F$ is supported on $\mathcal G$. A coherent sheaf with this support is killed by a power of $\mathcal I$, and hence is built from sheaves on $\mathcal G$. Split their characters. The nontrivial characters are the components $\mathrm\Upsilon_i$; the trivial-character component is obtained from $\rho^*\kappa_*\mathrm{D}^{b}(D)$ and the same finite ideal-power filtration. This proves generation and the stated decomposition. The formula for the left adjoint follows from adjunction for $l$ and $\pi_D$.

For \ref{num:pullbackskycraper2}, choose local parameters $(u,z_2,\ldots,z_n)$ with $D=(u=0)$ and root coordinate $t$ satisfying $t^s=u$. On this chart the root cover is finite flat and
\[
\rho^*\mathrm k_x=\mathrm k[t]/(t^s),
\qquad \mathcal I=(t),
\]
where $t$ has character $\nu^{-1}$. Thus its successive quotients are generated by $1,t,\ldots,t^{s-1}$ with the stated characters. Also $Q_j=\mathrm k[t]/(t^{s-j})$ and the kernel of $Q_j\to Q_{j+1}$ is generated by $t^{s-j-1}$, of character $\nu^{j+1}$. This proves the exact sequences.

They do not split: $t^{s-j-1}$ acts nontrivially on $Q_j$ but annihilates $Q_{j+1}\oplus(\mathcal G_x\otimes\nu^{j+1})$. Their extension space is one-dimensional. Indeed, the Koszul resolution of $Q_{j+1}$ on the smooth root chart uses the regular sequence $(t^{s-j-1},z_2,\ldots,z_n)$. Among its first syzygies only $t^{s-j-1}$ has character $\nu^{j+1}$; the tangent parameters have trivial character. Taking invariants gives
\[
\mathrm{Ext}^1(Q_{j+1},\mathcal G_x\otimes\nu^{j+1})=\mathrm k.
\]
The filtration itself is intrinsic because it is defined by powers of the ramification ideal, so these local descriptions glue. Rotating the short exact sequences gives the stated shifts. Since $\rho^*\mathrm k_x$ belongs to the last component of the root decomposition, its actual projections to $\mathrm\Upsilon_i$ are zero for $1\le i<s$.
\end{proof}

If $S$ is a smooth codimension $c$ subvariety contained in $U$, which can be realized as complete intersection of smooth divisors $D_{1}$,...$D_{c-1}$, $D$ then we consider $(1,...,1,s)$-type non-commutative weighted blow-up along $S$.
\begin{proposition}\label{prop:singlecover}
Consider the following commutative diagram, we have:
$$\begin{tikzcd}[ampersand replacement=\&]
\&{[Bl_{S}V/\mu_{s}]}\arrow[r,"q"]\arrow[d,"p"] \&\widetilde{wBl_{S}U}\arrow[d,"\pi"] \\
 \&{[V/\mu_{s}]} \arrow[r,"{\rho}"] \&  U
\end{tikzcd}$$
for any geometric point $x$ in $S$, $\pi^{*} (\mathrm{k}_{x})$ has a unique Postnikov filtration:
$$\begin{tikzcd}[column sep=0.5em]
 & 0 \arrow{rr}&& E_{-r}\arrow{dl}\arrow{rr}&& E_{-r+1}  \arrow{dl} \\
&& \iota_{x*}\mathcal{D}_{r}\arrow[ul,dashed,"\Delta"]&& \iota_{x*}\mathcal{D}_{r-1}\arrow[ul,dashed,"\Delta"]
\end{tikzcd}
...\begin{tikzcd}[column sep=0.5em]
 &   E_{-2}\arrow{rr}&&  E_{-1}\arrow{rr}\arrow{dl}&& \pi^{*} (\mathrm{k}_{x})\arrow{dl} \\
 && \iota_{x*}\mathcal{D}_{1}\arrow[ul,dashed,"\Delta"]  && \iota_{x*}\mathcal{D}_{0}\arrow[ul,dashed,"\Delta"]
\end{tikzcd}$$
where $\iota_{x}$ is the embedding the fiber of exceptional divisor above $x$:
$$\begin{tikzcd}[ampersand replacement=\&]
\&{E}|_{x}\simeq \mathcal{P}(1^{c-1},s)\arrow[r,"\iota_{x}"]\arrow[d] \&\widetilde{wBl_{S}U}\arrow[d,"\pi"] \\
 \&{x} \arrow[r] \&  U
\end{tikzcd}$$
and $\mathcal{D}_{i}$ is the dual exceptional object for $\mathcal{P}(1^{c-1},s)$, $r:=c+s-2$.\\
\end{proposition}
\begin{proof}

Step 1: Pull back $\mathrm k_x$ by the flat root projection.
Lemma~\ref{lem:pullbackskycraper}(\ref{num:pullbackskycraper2})
identifies it with $\mathrm k[t]/(t^s)$ and supplies its quotient
filtration. Rotating the short exact sequences gives outgoing factors
$l_*\mathrm k_x\otimes\nu[1],\ldots,
l_*\mathrm k_x\otimes\nu^{s-1}[1],l_*\mathrm k_x[1]$.
In the following diagrams write $E^j=Q_j$. We retain these actual extension maps. These are filtration factors,
not projections to the root semi-orthogonal components.\\

Step 2: By our assumption $S$ can be realized as complete intersection of divisors $D$ and $D_{i}$ in $U$, so on cyclic covering of
$U$ the reduced preimage of $S$ can be realized as the complete intersection of divisors $D_{V}$ and $\rho^{*}D_{i}$, it is not difficult to see that induces an equivariant version of Koszul resolution defined the center as follows:

\begin{equation}\label{equ:kouzcenter}
\mathcal{E}^{\vee}:=\bigoplus_{i}\mathcal{O}(-\rho^{*}D_{i})\bigoplus\mathcal{M}^{-1}_{D}\longrightarrow\mathcal{O}_{[V/\mu_{s}]}\longrightarrow\mathcal{O}_{[S/\mu_{s}]}\longrightarrow 0\end{equation}
and we can repeat the trick they did in the classical blow-up, see \cite[Section 11.1, p.252-253]{H}. We consider the surjection given by the Koszul resolution above, which induces a closed immersion of $Bl_{S}V$ to $\mathds{P}_{V}(E)$, where $E:=\mathcal{O}_{V}(D_{V})\bigoplus_{i}\mathcal{O}_{V}(\rho^{*}D_{i})$, we have a compatible $\mu_{s}$ action under this immersion, under one atlas of local coordinate it is like:
$$(x_{1},t_{x_{2}},...,t_{x_{c}},x_{c+1},...,x_{n})\longmapsto (e_{x_{2}}=t_{x_{2}},...,e_{x_{c}}=t_{x_{c}},x_{1},t_{x_{2}}x_{1},...,t_{x_{c}}x_{1}, x_{c+1},...,x_{n})$$
with compatible action:
$$(x_{1},t_{x_{2}},...,\zeta_{s}t_{x_{c}},x_{c+1},...,x_{n})\longmapsto (e_{x_{2}}=t_{x_{2}},...,\zeta_{s}e_{x_{c}}=\zeta_{s}t_{x_{c}},x_{1},t_{x_{2}}x_{1},...,\zeta_{s}t_{x_{c}}x_{1}, x_{c+1},...,x_{n})$$
where $t_{x_{i}}$ and $e_{x_{i}}$ is the projective ray corresponding to the divisor $\rho^{*}D_{i}$ and $D_{V}$ in $Bl_{S}V$ and $\mathds{P}_{V}(E)$. We have a stacky version description as follows:
$$\begin{tikzcd}[column sep=0.5em]
 & {[Bl_{S}V/\mu_{s}]} \arrow[dr,"p"]\arrow[rr,"\beta"]&& {\mathds{P}_{[V/\mu_{s}]}(\mathcal{E})}\arrow[dl,"\alpha"]\\
&&{[V/\mu_{s}]}
\end{tikzcd}$$
and noticing $\mathds{P}_{[V/\mu_{s}]}(\mathcal{E})$ restrict at a fiber of center is ${[\mathds{P}(V_{c})/\mu_{s}]}$, where $\mu_{s}$ act on the correspondent coordinate of  $\mathds{P}(V_{c})$ as (\ref{equ:actclassical}).  Then following the construction in \cite[Section 11.1]{H} we can realize this embedding as a Koszul immersion concerning the relative Euler exact sequence
$$...\longrightarrow\mathrm{\Omega}_{\mathds{P}_{V}(E)/V}(1)\longrightarrow\mathcal{O}_{\mathds{P}_{V}(E)}\longrightarrow\mathcal{O}_{Bl_{S}V}\longrightarrow 0$$
In particular, the global coordinates described by \cite[Section 2]{Ca} tell us this Koszul resolution is not merely exact, but also $\mu_{s}$-equivariant, so we have a stacky description as follows:
$$...\longrightarrow\mathrm{\Omega}_{\mathds{P}_{[V/\mu_{s}]}(\mathcal{E})}(1)\longrightarrow\mathcal{O}_{\mathds{P}_{[V/\mu_{s}]}(\mathcal{E})}\longrightarrow\mathcal{O}_{[Bl_{S}V/\mu_{s}]}\longrightarrow 0$$
for any geometric point $x$ on $[S/\mu_{s}]$, we have
$$
\beta_{*}p^{*}\mathrm{k}_{x}= \beta_{*}\beta^{*}\alpha^{*}\mathrm{k}_{x} =\alpha^{*}\mathrm{k}_{x}\otimes \mathcal{O}_{[Bl_{S}V/\mu_{s}]}=\alpha^{*}\mathrm{k}_{x}\otimes \big[\bigwedge\mathrm{\Omega}_{\mathds{P}_{[V/\mu_{s}]}(\mathcal{E})}(1)\big]
$$
where $\big[\bigwedge\mathrm{\Omega}_{\mathds{P}_{[V/\mu_{s}]}(\mathcal{E})}(1)\big]$ denotes the Koszul complex, noticing when we restrict at a fiber of $\mathds{P}_{[V/\mu_{s}]}$ on the zero locus of the section of (\ref{equ:kouzcenter}), $\big[\bigwedge\mathrm{\Omega}_{\mathds{P}_{[V/\mu_{s}]}(\mathcal{E})}(1)\big]$ degenerates and their morphisms are zero maps, so we get
 $$\mathcal{H}^{-r}(\beta_{*}p^{*} \mathrm{k}_{x})=\beta_{*}\mathcal{H}^{-r}(p^{*} \mathrm{k}_{x})=\beta_{*}\iota_{x*}\mathrm{\Omega}_{[\mathds{P}(V_{c})/\mu_{s}]}^{r}\otimes\mathcal{O}_{[\mathds{P}(V_{c})/\mu_{s}]}(r)$$
and
$$\mathcal{H}^{-r} (p^{*} \mathrm{k}_{x})=\iota_{x*}\mathrm{\Omega}_{[\mathds{P}(V_{c})/\mu_{s}]}^{r}\otimes\mathcal{O}_{[\mathds{P}(V_{c})/\mu_{s}]}(r)$$
 where $0\leq r\leq c-1$.\\

Step 3: We compute the action of functor $p^{*}$ on the left Postnikov system in Lemma \ref{lem:pullbackskycraper} (\ref{num:pullbackskycraper2}), we have
\begin{equation}\label{filtration22}
\begin{tikzcd}[column sep=0.5em]
 & p^{*}\rho^{*}(\mathrm{k}_{x}) \arrow{rr}&& p^{*}E^{1}\arrow{dl}\arrow{rr}&& p^{*}E^{2}  \arrow{dl} \\
&&p^{*}l_{*}\mathrm{k}_{x}\otimes\nu[1]\arrow[ul,dashed,"\Delta"]&& p^{*}l_{*}\mathrm{k}_{x}\otimes\nu^{2}[1]\arrow[ul,dashed,"\Delta"]
\end{tikzcd}...\end{equation}
$$...\begin{tikzcd}[column sep=0.5em]
 &  p^{*}E^{s-2}\arrow{rr}&& p^{*}E^{s-1}\arrow{rr} \arrow{dl}&& 0\arrow{dl} \\
 && p^{*}l_{*}\mathrm{k}_{x}\otimes\nu^{s-1}[1]\arrow[ul,dashed,"\Delta"]  && p^{*} l_{*}\mathrm{k}_{x}[1]\arrow[ul,dashed,"\Delta"]
\end{tikzcd}$$
notice $p^{*} D_{V}\sim \widehat{D}_{V}+E$ and $\widehat{D}_{V}$ is the ramified locus in $Bl_{S}(V)$, hence $p^{*}\mathcal{M}_{D}\sim\widehat{\mathcal{M}}_{D}+E$, then for any $k$ there is:
\begin{align*}
\mathcal{H}^{-\lambda} (p^{*} \mathrm{k}_{x}\otimes\nu^{k})&= \mathcal{H}^{-\lambda} (p^{*} \mathrm{k}_{x}\otimes\mathcal{M}_{D}^{k}) \\
&=\mathcal{H}^{-\lambda} (p^{*} \mathrm{k}_{x})\otimes\mathcal{M}_{\widehat{D}}^{k}\otimes\mathcal{O}(E)^{k}\\
&=\iota_{x*}(\mathrm{\Omega}_{[\mathds{P}(V_{c})/\mu_{s}]}^{\lambda}\otimes\mathcal{O}_{[\mathds{P}(V_{c})/\mu_{s}]}(\lambda)\otimes\nu^{k}) \stepcounter{equation}\tag{\theequation}\label{equ:filofpkx}
\end{align*}
the last step is due to $\mathcal{M}_{\widehat{D}}|_{[\mathds{P}(V_{c})/\mu_{s}]}=\mathcal{O}_{[\mathds{P}(V_{c})/\mu_{s}]}([\mathds{P}(V_{c-1})/\mu_{s}])$, while $\mathcal{O}(E)|_{[\mathds{P}(V_{c})/\mu_{s}]}=\mathcal{O}_{[\mathds{P}(V_{c})/\mu_{s}]}(-1)$, recall the equivariant idea sheaf defined the center $\mathds{P}(V_{c-1})$ e.g. Lemma \ref{lemma:Kawambi}, we have:
$$0\longrightarrow\mathcal{O}_{[\mathds{P}(V_{c})/\mu_{s}]}(-1)\otimes\nu^{-1}\longrightarrow\mathcal{O}_{[\mathds{P}(V_{c})/\mu_{s}]}\longrightarrow\mathcal{O}_{[\mathds{P}(V_{c-1})/\mu_{s}]}\longrightarrow0$$
since $\mathrm{Hom}(\mathcal{O}_{[\mathds{P}(V_{c})/\mu_{s}]}(-1)\otimes\nu^{-1},\mathcal{O}_{[\mathds{P}(V_{c})/\mu_{s}]})=\mathrm{k}$ and $\mathrm{Hom}(\mathcal{O}_{[\mathds{P}(V_{c})/\mu_{s}]}(-1),\mathcal{O}_{[\mathds{P}(V_{c})/\mu_{s}]})=V_{c-1}$, the zero locus of the section corresponding to different divisor. So $$\mathcal{O}_{[\mathds{P}(V_{c})/\mu_{s}]}([\mathds{P}(V_{c-1})/\mu_{s}])\simeq \mathcal{O}_{[\mathds{P}(V_{c})/\mu_{s}]}(1)\otimes\nu$$
that means $\mathcal{M}_{\widehat{D}}\otimes\mathcal{O}(E)|_{\mathds{P}([V_{c}/\mu_{s}])}\simeq \mathcal{O}_{\mathds{P}([V_{c}/\mu_{s}])}\otimes\nu$. \\

In order to simplify notation, we will refer to $\mathrm{\Omega}^{\lambda}(\lambda)$ as $\mathrm{\Omega}_{[\mathds{P}(V_{c})/\mu_{s}]}^{\lambda}\otimes\mathcal{O}_{[\mathds{P}(V_{c})/\mu_{s}]}(\lambda)$ only in the following proof part. Then by (\ref{equ:filofpkx}) we have a heart filtration for each $p^{*}l_{*}\mathrm{k}_{x}\otimes\nu^{i}$ as following:
\begin{equation}\label{filtra1}
\begin{tikzcd}[column sep=0.5em]
 & 0 \arrow{rr}&& E_{-c+1}\arrow{dl}\arrow{rr}&& E_{-c+2}  \arrow{dl} \\
&& \iota_{x*}\mathrm{\Omega}^{c-1}(c-1)[c-1]\otimes\nu^{i}\arrow[ul,dashed,"\Delta"]&&  \iota_{x*}\mathrm{\Omega}^{c-2}(c-2)[c-2]\otimes\nu^{i}\arrow[ul,dashed,"\Delta"]
\end{tikzcd}...
\end{equation}
$$...\begin{tikzcd}[column sep=0.5em]
 &  E_{-2}\arrow{rr}&& E_{-1}\arrow{rr} \arrow{dl}&& p^{*}l_{*}\mathrm{k}_{x}\otimes\nu^{i}\arrow{dl} \\
 &&  \iota_{x*}(\mathrm{\Omega}^{1}(1)\otimes\nu^{i})[1]\arrow[ul,dashed,"\Delta"]  &&  \iota_{x*}(\mathcal{O}_{[\mathds{P}(V_{c})/\mu_{s}]}\otimes\nu^{i})\arrow[ul,dashed,"\Delta"]
\end{tikzcd}$$
Combine the above filtrations (\ref{filtration22}) and (\ref{filtra1}), we have a double filtration for object $p^{*}\rho^{*}(\mathrm{k}_{x})$, conversely if we record each extension in this filtration we can also recover $p^{*}\rho^{*}(\mathrm{k}_{x})$ from it.
\\

Step 4:
We figure out all possible factors that might contribute to the unique extension
\begin{equation}\label{equ:extension}
\mathrm{Ext}^{1}(p^{*}(\mathrm{k}_{x}\otimes\nu^{i}),p^{*}(\mathrm{k}_{x}\otimes\nu^{i-1}))=\mathrm{Ext}^{1}(\mathrm{k}_{x}\otimes\nu^{i},\mathrm{k}_{x}\otimes\nu^{i-1})= \mathrm{k}\end{equation}
via showing behaviors of the morphisms in spectral sequence,
$$E^{2}_{m,n}:=\mathrm{Ext}^{m}(\mathcal{H}^{-n}p^{*}(\mathrm{k}_{x}\otimes\nu^{i}),p^{*}(\mathrm{k}_{x}\otimes\nu^{i-1}))\Longrightarrow \mathrm{Ext}^{m+n}(p^{*}(\mathrm{k}_{x}\otimes\nu^{i}),p^{*}(\mathrm{k}_{x}\otimes\nu^{i-1}))$$
$$E'^{2}_{s,t}:=\mathrm{Ext}^{s}(\mathcal{H}^{-n}p^{*}(\mathrm{k}_{x}\otimes\nu^{i}),\mathcal{H}^{t}p^{*}(\mathrm{k}_{x}\otimes\nu^{i-1}))\Longrightarrow \mathrm{Ext}^{s+t}(\mathcal{H}^{-n}p^{*}(\mathrm{k}_{x}\otimes\nu^{i}),p^{*}(\mathrm{k}_{x}\otimes\nu^{i-1}))$$
where the fully-faithfulness of $p^{*}$ is due the fact $p_{*}\mathcal{O}_{[Bl_{S}(V)/\mu_{s}]}=\mathcal{O}_{[V/\mu_{s}]}$. It's not difficult to see all the possible contribution factors come from terms of form \begin{equation}\label{equ:conterms}\mathrm{Ext}^{k}(\iota_{x*}\mathrm{\Omega}^{a}(a)\otimes\nu^{i},\iota_{x*}\mathrm{\Omega}^{a-1+k}(a-1+k)\otimes\nu^{i-1})\end{equation}
for any integer $a$ and $k$ in their range (fixing $i$). \\

We next compute all these possible terms. Let's consider the following diagram:
$$\begin{tikzcd}[ampersand replacement=\&]
\&{E}|_{x}\simeq {[\mathds{P}(V_{c})/\mu_{s}]}\arrow[r,"\iota'_{x}"]\arrow[d,"p_{E_{x}}"]\&E\arrow[d,"p_{E}"]\arrow[r,"\iota"] \&{[Bl_{S}V/\mu_{s}]}\arrow[d,"p"] \\
 \&{[x/\mu_{s}]} \arrow[r,"\theta"]\&  {[S/\mu_{s}]} \arrow[r,"\vartheta"]\&  {[V/\mu_{s}]}
\end{tikzcd}$$
by excess distinguished triangle:\\
\begin{equation}\label{equa:exceedofomega1}
\iota'_{x*}\mathrm{\Omega}^{a}(a)\otimes\mathcal{O}_{E}(-E)[1]\longrightarrow \iota^{*}\iota_{*}\iota'_{x*}\mathrm{\Omega}^{a}(a)\longrightarrow\iota'_{x*}\mathrm{\Omega}^{a}(a)\longrightarrow\iota'_{x*}\mathrm{\Omega}^{a}(a)\otimes\mathcal{O}_{E}(-E)[2]\end{equation}
we can see
$$\mathrm{Ext}^{k}(\iota_{x*}\mathrm{\Omega}^{a}(a)\otimes\nu^{i},\iota_{x*}\mathrm{\Omega}^{a-1+k}(a-1+k)\otimes\nu^{i-1})=\mathrm{Ext}^{k}(\iota^{*}\iota_{*}\iota'_{x*}\mathrm{\Omega}^{a}(a)\otimes\nu^{i},\iota'_{x*}\mathrm{\Omega}^{a-1+k}(a-1+k)\otimes\nu^{i-1})$$
which fits in long exact sequence
$$\mathrm{Ext}^{k}(\iota'_{x*}\mathrm{\Omega}^{a}(a),\iota'_{x*}\mathrm{\Omega}^{a-1+k}(a-1+k)\otimes\nu^{-1})\longrightarrow\mathrm{Ext}^{k}(\iota_{x*}\mathrm{\Omega}^{a}(a)\otimes\nu^{i},\iota_{x*}\mathrm{\Omega}^{a-1+k}(a-1+k)\otimes\nu^{i-1})$$
\begin{equation}\label{equ:3term}
\longrightarrow\mathrm{Ext}^{k-1}(\iota'_{x*}\mathrm{\Omega}^{a}(a+1),\iota'_{x*}\mathrm{\Omega}^{a-1+k}(a-1+k)
\otimes\nu^{-1})\longrightarrow...
\end{equation}
for the first term of (\ref{equ:3term}) use the adjoint formula
$$\mathrm{Ext}^{k}(\iota'_{x*}\mathrm{\Omega}^{a}(a),\iota'_{x*}\mathrm{\Omega}^{a-1+k}(a-1+k)\otimes\nu^{-1})=\mathrm{Ext}^{k}(\iota'^{x*}\iota'_{x*}\mathrm{\Omega}^{a}(a),\mathrm{\Omega}^{a-1+k}(a-1+k)\otimes\nu^{-1})$$
and notice $\mathcal{H}^{-\lambda}(\iota'^{x*}\iota'_{x*}\mathrm{\Omega}^{a}(a))=\bigwedge^{\lambda}p_{E_{x}}^{*}\mathcal{N}^{\vee}_{x|S}\otimes\mathrm{\Omega}^{a}(a)=(\bigwedge^{\lambda}V_{n-c})\otimes\mathrm{\Omega}^{a}(a)$, since $[x/\mu_{s}]$ is a smooth point on $[S/\mu_{s}]$ and $[S/\mu_{s}]$ is a gerbe hence there is no representation terms in $\mathcal{N}^{\vee}_{x|S}$, then we have spectral sequence:
\begin{multline}\label{equ:speccoh1}
E^{p,q}_{2}:=\mathrm{Ext}^{p}((\bigwedge^{q}V_{n-c})\otimes\mathrm{\Omega}^{a}(a),\mathrm{\Omega}^{a-1+k}(a-1+k)\otimes\nu^{-1})\\
\Longrightarrow\mathrm{Ext}^{p+q}(\iota'_{x*}\mathrm{\Omega}^{a}(a),\iota'_{x*}\mathrm{\Omega}^{a-1+k}(a-1+k)\otimes\nu^{-1})
\end{multline}

We need the following lemma about cohomology computation:
\begin{lemma}
  $\mathrm{Ext}^{*}(\mathrm{\Omega}^{a}(a),\mathrm{\Omega}^{a-1+k}(a-1+k)\otimes\nu^{-1})=\begin{cases}
     \mathrm{ k}  & \text{if $k=0$}\\
     0 & \text{if $k\geq 1$}\\
     \end{cases}$
\end{lemma}

\begin{proof}
Noticing if $k<0$ we have $\mathrm{Ext}^{k}(\iota_{x*}\mathrm{\Omega}^{a}(a)\otimes\nu^{i},\iota_{x*}\mathrm{\Omega}^{a-1+k}(a-1+k)\otimes\nu^{i-1})=0$ automatically by Grothendieck vanishing, hence we only consider the cases whence $k\geq 0$.\\

If $k=0$, we have
\begin{align*}
\mathrm{Ext}^{*}(\mathrm{\Omega}^{a}(a),\mathrm{\Omega}^{a-1+k}(a-1+k)\otimes\nu^{-1})&= \mathrm{Ext}^{*}(\mathrm{\Omega}^{a}(a),\mathrm{\Omega}^{a-1}(a-1)\otimes\nu^{-1}) \\
&=\mathrm{Hom}_{\mathds{P}(V_{c})}(\mathrm{\Omega}_{\mathds{P}(V_{c})}^{a}(a),\mathrm{\Omega}_{\mathds{P}(V_{c})}^{a-1}(a-1)\otimes\nu^{-1})^{\mu_{s}}\\
&=[(V_{c-1}\bigoplus\mathrm{k}\otimes\nu)\otimes\nu^{-1}]^{\mu_{s}}=\mathrm{k}
\end{align*}
where the second and third steps are due to \cite[Lemma 2]{Bei}.

If $k=1$, we have $\mathrm{Ext}^{*}(\mathrm{\Omega}^{a}(a),\mathrm{\Omega}^{a-1+k}(a-1+k)\otimes\nu^{-1})=\mathrm{Ext}^{*}(\mathrm{\Omega}^{a}(a),\mathrm{\Omega}^{a}(a)\otimes\nu^{-1})$, while by the fact $\mathrm{Ext}_{\mathds{P}(V_{c})}^{*}(\mathrm{\Omega}_{\mathds{P}(V_{c})}^{a}(a),\mathrm{\Omega}_{\mathds{P}(V_{c})}^{a}(a))=\mathrm{k}$ since $\mathrm{\Omega}_{\mathds{P}(V_{c})}^{a}(a)[a]$ is a dual exceptional object in $\mathrm{D}^{b}(\mathds{P}(V_{c}))$, so
\begin{align*}
\mathrm{Ext}^{*}(\mathrm{\Omega}^{a}(a),\mathrm{\Omega}^{a}(a)\otimes\nu^{-1})&= \mathrm{Ext}_{\mathds{P}(V_{c})}^{*}(\mathrm{\Omega}_{\mathds{P}(V_{c})}^{a}(a),\mathrm{\Omega}_{\mathds{P}(V_{c})}^{a}(a)\otimes\nu^{-1})^{\mu_{s}} \\
&=\mathrm{Hom}_{\mathds{P}(V_{c})}(\mathrm{\Omega}_{\mathds{P}(V_{c})}^{a}(a),\mathrm{\Omega}_{\mathds{P}(V_{c})}^{a}(a)\otimes\nu^{-1})^{\mu_{s}}=0
\end{align*}

If $k>1$, we have $\mathrm{Ext}^{*}(\mathrm{\Omega}^{a}(a),\mathrm{\Omega}^{a-1+k}(a-1+k)\otimes\nu^{-1})=0$ since we have a semi-orthogonal decomposition
\begin{equation}\label{equ:dualforpojectiveVc}
\mathrm{D}^{b}(\mathds{P}(V_{c}))=\langle\mathrm{\Omega}_{\mathds{P}(V_{c})}^{c-1}(c-1)[c-1],...,\mathrm{\Omega}_{\mathds{P}(V_{c})}(1)[1], \mathcal{O}_{\mathds{P}(V_{c})}\rangle
\end{equation}
by Corollary \ref{cor:sodofwps} and \ref{cor:dualdes}, so
$\mathrm{Ext}_{\mathds{P}(V_{c})}^{*}(\mathrm{\Omega}_{\mathds{P}(V_{c})}^{a}(a),\mathrm{\Omega}_{\mathds{P}(V_{c})}^{a-1+k}(a-1+k))=0$.
\end{proof}

Using the above lemma and (\ref{equ:speccoh1}) we know the first term in our long exact sequence (\ref{equ:3term}) is nontrivial only if $k=0$. \\

Then let's treat the third term in a long exact sequence (\ref{equ:3term}) similarly there is
$$\mathrm{Ext}^{k-1}(\iota'_{x*}\mathrm{\Omega}^{a}(a+1),\iota'_{x*}\mathrm{\Omega}^{a-1+k}(a-1+k)\otimes\nu^{-1})=\mathrm{Ext}^{k-1}(\iota'^{x*}\iota'_{x*}\mathrm{\Omega}^{a}(a+1),\mathrm{\Omega}^{a-1+k}(a-1+k)\otimes\nu^{-1})$$
and a spectral sequence
\begin{multline}\label{equ:speccoh2}
E^{p,q}_{2}:=\mathrm{Ext}^{p}((\bigwedge^{q}V_{n-c})\otimes\mathrm{\Omega}^{a}(a+1),\mathrm{\Omega}^{a-1+k}(a-1+k)\otimes\nu^{-1})\\
\Longrightarrow\mathrm{Ext}^{p+q}(\iota'_{x*}\mathrm{\Omega}^{a}(a+1),\iota'_{x*}\mathrm{\Omega}^{a-1+k}(a-1+k)\otimes\nu^{-1})
\end{multline}
so if we claim the last term of (\ref{equ:3term}) vanishes for any $a$ and $k$, it's only enough for us to verify
$$\mathrm{Ext}^{\leq k-1}(\mathrm{\Omega}^{a}(a+1),\mathrm{\Omega}^{a-1+k}(a-1+k)\otimes\nu^{-1})=0$$
for $k>0$, furthermore which can also be reduced to verify
$$\mathrm{Ext}_{\mathds{P}(V_{c})}^{\leq k-1}(\mathrm{\Omega}_{\mathds{P}(V_{c})}^{a}(a+1),\mathrm{\Omega}_{\mathds{P}(V_{c})}^{a-1+k}(a-1+k)\otimes\nu^{-1})=0$$
for $k>0$. On the other hand, we know from (\ref{equ:dualforpojectiveVc})
$$\mathrm{Ext}_{\mathds{P}(V_{c})}^{\leq k-1}(\mathrm{\Omega}_{\mathds{P}(V_{c})}^{a}(a+1),\mathrm{\Omega}_{\mathds{P}(V_{c})}^{a-1+k}(a-1+k)\otimes\nu^{-1})= \mathrm{Hom}_{\mathds{P}(V_{c})}^{\leq 0}(\mathcal{D}^{c}_{a}(1),\mathcal{D}^{c}_{a-1+k}
\otimes\nu^{-1})$$
where  we refer $\mathcal{D}^{c}_{i}$ as the $i$-th dual exceptional object on projective space $\mathds{P}(V_{c})$.
\begin{lemma}\label{lemma:vanshingcoho1}
$\mathrm{Hom}_{\mathds{P}(V_{c})}^{\leq 0}(\mathcal{D}^{c}_{a}(1),\mathcal{D}^{c}_{a-1+k}
\otimes\nu^{-1})=0$ if $k>0$.
\end{lemma}
\begin{proof}

If $0\leq a<c-1$, noticing $\mathbb{L}_{\mathcal{O}_{{\mathds{P}(V_{c})}}}\mathcal{D}^{c}_{a}(1)=\mathcal{D}^{c}_{a+1}$ by the periodicity of the semi-orthogonal decomposition of projective space, so the mutation triangle
$$\bigoplus_{k}\mathrm{Ext}^{k}(\mathcal{O}_{{\mathds{P}(V_{c})}},\mathcal{D}^{c}_{a}(1))\otimes\mathcal{O}_{{\mathds{P}(V_{c})}}[-k]\longrightarrow\mathcal{D}^{c}_{a}(1)\longrightarrow\mathbb{L}_{\mathcal{O}_{\mathds{P}(V_{c})}}\mathcal{D}^{c}_{a}(1)$$
induces a long exact sequence
\begin{equation}\label{equ:projectmutaion}
\longrightarrow\bigoplus_{k}\mathrm{Ext}_{\mathds{P}(V_{c})}^{k}(\mathcal{O}_{\mathds{P}(V_{c})},\mathcal{D}^{c}_{a}(1))\otimes \mathrm{Ext}_{\mathds{P}(V_{c})}^{i+k}(\mathcal{O}_{\mathds{P}(V_{c})},\mathcal{D}_{a-1+k}^{c}
\otimes\nu^{-1})\longrightarrow\end{equation}
$$\mathrm{Ext}_{\mathds{P}(V_{c})}^{i}(\mathcal{D}^{c}_{a}(1),\mathcal{D}^{c}_{a-1+k}
\otimes\nu^{-1})
\longrightarrow\mathrm{Ext}_{\mathds{P}(V_{c})}^{i}(\mathcal{D}^{c}_{a+1},\mathcal{D}^{c}_{a-1+k}
\otimes\nu^{-1})\longrightarrow$$
noticing Corollary \ref{cor:dualvanishing}, we have $\mathrm{Ext}_{\mathds{P}(V_{c})}^{*}(\mathcal{O}_{\mathds{P}(V_{c})},\mathcal{D}^{c}_{a-1+k}
\otimes\nu^{-1})=0$, and
$$\mathrm{Ext}_{\mathds{P}(V_{c})}^{\leq 0}(\mathcal{D}^{c}_{a+1},\mathcal{D}^{c}_{a-1+k}
\otimes\nu^{-1})=\begin{cases}
     \mathrm{Ext}_{\mathds{P}(V_{c})}^{\leq -1}(\mathrm{\Omega}_{\mathds{P}(V_{c})}^{a+1}(a+1),\mathrm{\Omega}_{\mathds{P}(V_{c})}^{a}(a)\otimes\nu^{-1})=0  & \text{if $k=1$}\\
    \mathrm{Ext}_{\mathds{P}(V_{c})}^{\leq 0}(\mathrm{\Omega}_{\mathds{P}(V_{c})}^{a+1}(a+1),\mathrm{\Omega}_{\mathds{P}(V_{c})}^{a+1}(a+1)\otimes\nu^{-1})=0 & \text{if $k=2$}\\
     0 \text{\quad\quad by Corollary \ref{cor:sodofwps}} & \text{if $k>2$}\\
     \end{cases}$$
then we can see the midterms of the long exact sequence (\ref{equ:projectmutaion}) also vanish at all non-positive degrees.\\

If $a=c-1$,  $\mathcal{D}^{c}_{a}\simeq \mathcal{O}_{\mathds{P}(V_{c})}(-1)[c-1]$ then $$\mathrm{Ext}_{\mathds{P}(V_{c})}^{*}(\mathcal{D}^{c}_{a}(1),\mathcal{D}^{c}_{a-1+k}\otimes\nu^{-1})=\mathrm{Ext}_{\mathds{P}(V_{c})}^{*-c+1}(\mathcal{O}_{\mathds{P}(V_{c})},\mathcal{D}^{c}_{c-2+k}\otimes\nu^{-1})=0$$  for any $k$ by Corollary \ref{cor:sodofwps}.
\end{proof}
So now we can confirm that last term of (\ref{equ:3term}) vanishes for any $a$ and $k$.\\

Let us summarize our observations above, in the classic case \cite[Proposition 11.12]{H}, $\mathcal{H}^{-\lambda} (p^{*}\mathrm{k}_{x})=\iota_{x*}\mathrm{\Omega}^{\lambda}(\lambda)$ is given by the $-\lambda$ term of the Koszul resolution of local blow-up (permitting non-trivial extensions). Since the equivariant pullback functor of blow-up is also fully-faithful and the extension sequence is uniquely non-trivial, we have such morphisms between complexes as inherent are also uniquely non-trivial, and all possible contributions come from
$$\mathrm{Hom}(\iota_{x*}\mathrm{\Omega}^{\lambda}(\lambda)\otimes\nu^{i},\iota_{x*}\mathrm{\Omega}^{\lambda-1}(\lambda-1)\otimes\nu^{i-1})$$
which each of them stems from a unique morphism in (by \cite[Lemma 2]{Bei} and our previous discussion),
\begin{equation}\label{eq:cf}
  \mathrm{Hom}_{[\mathds{P}(V_{c})/\mu_{s}]}(\mathrm{\Omega}^{\lambda}(\lambda)\otimes\nu^{i},\mathrm{\Omega}^{\lambda-1}(\lambda-1)\otimes\nu^{i-1})\simeq [(V_{c-1}\bigoplus\mathrm{k}\otimes\nu)\otimes\nu^{-1}]^{\mu_{s}}\simeq\mathrm{k}
\end{equation}
between sheaves on equivariant projective space.\\

Step 5: We must show for a fixed $i$ each unique non-trivial extension factor in $$\mathrm{Hom}(\iota_{x*}\mathrm{\Omega}^{\lambda}(\lambda)\otimes\nu^{i},\iota_{x*}\mathrm{\Omega}^{\lambda-1}(\lambda-1)\otimes\nu^{i-1})$$
for any $\lambda$ in the range really contributes to the unique morphism of our extension in (\ref{equ:extension}).\\

We take $i=0$ for convenience, there is a heart filtration for $p^{*} l_{*}\mathrm{k}_{x}$:
\begin{equation}\label{eq:hfil}
\begin{adjustbox}{max width=0.98\linewidth}\(\begin{tikzcd}[column sep=0.5em]
 & 0\arrow{rr}&& E^{c-1}\arrow{dl}\arrow{rr}&& E^{c-2}  \arrow{dl} \\
&&\iota_{x*}\mathrm{\Omega}^{c-1}(c-1)[c-1]\arrow[ul,dashed,"\Delta"]&& \iota_{x*}\mathrm{\Omega}^{c-2}(c-2)[c-2]\arrow[ul,dashed,"\Delta"]
\end{tikzcd}...\begin{tikzcd}[column sep=0.5em]
 &  E^{2}\arrow{rr}&& E^{1}\arrow{rr} \arrow{dl}&&  p^{*} l_{*}\mathrm{k}_{x}\arrow{dl} \\
 && \iota_{x*}\mathrm{\Omega}^{1}(1)[1]\arrow[ul,dashed,"\Delta"]  && \iota_{x*}\mathcal{O}\arrow[ul,dashed,"\Delta"]
\end{tikzcd}\)\end{adjustbox}\end{equation}
and a similar filtration for $p^{*} l_{*}\mathrm{k}_{x}\otimes\nu^{-1}$:
$$\begin{tikzcd}[column sep=0.5em]
 & 0\arrow{rr}&& E^{c-1}\otimes\nu^{-1}\arrow{dl} \\
&&\iota_{x*}\mathrm{\Omega}^{c-1}(c-1)\otimes\nu^{-1}[c-1]\arrow[ul,dashed,"\Delta"]
\end{tikzcd}...\begin{tikzcd}[column sep=0.5em]
 &   E^{1}\otimes\nu^{-1}\arrow{rr} &&  p^{*} l_{*}\mathrm{k}_{x}\otimes\nu^{-1}\arrow{dl} \\
 && \iota_{x*}\mathcal{O}\otimes\nu^{-1}\arrow[ul,dashed,"\Delta"]
\end{tikzcd}$$

We consider the extension $\mathrm{Ext}^{1}(\iota_{x*}\mathrm{\Omega}^{i}(i)[i],E^{i+1})$ in each step of (\ref{eq:hfil}), noticing for any $i<j$
$$\mathrm{Ext}^{*}(\iota_{x*}\mathrm{\Omega}^{i}(i)[i],\iota_{x*}\mathrm{\Omega}^{j}(j)[j])= \mathrm{Ext}^{*}(\iota^{*}\iota_{*}\iota'_{x*}\mathrm{\Omega}^{i}(i)[i],\iota'_{x*}\mathrm{\Omega}^{j}(j)[j])$$
by excess distinguished triangle (\ref{equa:exceedofomega1}), it fits in long exact sequence
$$\longrightarrow\mathrm{Ext}^{k}(\iota'_{x*}\mathrm{\Omega}^{i}(i)[i],\iota'_{x*}\mathrm{\Omega}^{j}(j)[j])\longrightarrow\mathrm{Ext}^{k}(\iota^{x*}\iota_{x*}\iota'_{x*}\mathrm{\Omega}^{i}(i)[i],\iota'_{x*}\mathrm{\Omega}^{j}(j)[j])$$
\begin{equation}\label{equ:longexact11}
\longrightarrow\mathrm{Ext}^{k-1}(\iota'_{x*}\mathrm{\Omega}^{i}(i+1)[i],\iota'_{x*}\mathrm{\Omega}^{j}(j)[j])\longrightarrow...
\end{equation}
Similarly to (\ref{equ:speccoh1}) and (\ref{equ:speccoh2}), for the first and third terms of the long exact sequence we have the spectral sequences
$$E^{p,q}_{2}:=\mathrm{Ext}^{p}(\bigwedge^{q}V_{n-c}\otimes\mathrm{\Omega}^{i}(i)[i],\mathrm{\Omega}^{j}(j)[j])\Longrightarrow\mathrm{Ext}^{p+q}(\iota'^{*}_{x}\iota'_{x*}\mathrm{\Omega}^{i}(i)[i],\mathrm{\Omega}^{j}(j)[j])$$
$$E'^{p,q}_{2}:=\mathrm{Ext}^{p}(\bigwedge^{q}V_{n-c}\otimes\mathrm{\Omega}^{i}(i+1)[i],\mathrm{\Omega}^{j}(j)[j])\Longrightarrow\mathrm{Ext}^{p+q}(\iota'^{*}_{x}\iota'_{x*}\mathrm{\Omega}^{i}(i+1)[i],\mathrm{\Omega}^{j}(j)[j])$$
noticing by our assumption $i<j$, $\mathrm{Ext}_{\mathds{P}(V_{c})}^{*}(\mathrm{\Omega}_{\mathds{P}(V_{c})}^{i}(i)[i],\mathrm{\Omega}_{\mathds{P}(V_{c})}^{j}(j)[j])=0$, so the first terms of (\ref{equ:longexact11}) vanish. Following almost the same argument as in Lemma (\ref{lemma:vanshingcoho1}) we have
\begin{equation}\label{eua:omegaomega1}
 \mathrm{Ext}^{*}(\mathrm{\Omega}^{i}(i+1)[i],\mathrm{\Omega}^{j}(j)[j]) =
    \begin{cases}
    \mathrm{k} & \text{if $j=i+1$}\\
      0  & \text{if $j>i+1$}\\
    \end{cases}
\end{equation}
so for the third terms in (\ref{equ:longexact11})
\begin{equation}
  \mathrm{Ext}^{k-1}(\iota'_{x*}\mathrm{\Omega}^{i}(i+1)[i],\iota'_{x*}\mathrm{\Omega}^{j}(j)[j]) =
    \begin{cases}
      \bigwedge^{k-1}V_{n-c} & \text{if $j=i+1$}\\
      0  & \text{if $j>i+1$}\\
    \end{cases}
\end{equation}
and \begin{equation}
  \mathrm{Ext}^{k}(\iota_{x*}\mathrm{\Omega}^{i}(i)[i],\iota_{x*}\mathrm{\Omega}^{j}(j)[j]) =
    \begin{cases}
       \bigwedge^{k-1}V_{n-c} & \text{if $j=i+1$}\\
      0  & \text{if $j>i+1$}\\
    \end{cases}
\end{equation}
furthermore recall $E^{i+1}$ can be realized as iterated extension of $\iota_{x*}\mathrm{\Omega}^{k}(k)[k]$ for $k\geq i+1$, so a easy long exact sequence argument shows
$\mathrm{Ext}^{1}(\iota_{x*}\mathrm{\Omega}^{i}(i)[i],E^{i+1})=\mathrm{k}$ in each step, and if we slightly modify the form of (\ref{eua:omegaomega1}) the extension comes from $\mathrm{Ext}^{1}(\mathrm{\Omega}^{i}(i),\mathrm{\Omega}^{i+1}(i))$ for any $i$ in the range.\\

Then we proceed with a proof similar to that in Lemma \ref{lem:pullbackskycraper}, where we show that every extension is nontrivial. Firstly we can compute $\mathrm{Hom}(\iota_{x*}\mathrm{\Omega}^{i}(i)[i],\iota_{x*}\mathrm{\Omega}^{i}(i)[i])=\mathrm{k}$ similarly by (\ref{equ:longexact11}). Then by spectral sequence for extension
$$\mathrm{Hom}(E^{i},E^{i})= \mathrm{Ker}\, \big\{\mathrm{Hom}(E^{i+1},E^{i+1})\bigoplus\mathrm{Hom}(\iota_{x*}\mathrm{\Omega}^{i}(i)[i],\iota_{x*}\mathrm{\Omega}^{i}(i)[i])\longrightarrow \mathrm{Ext}^{1}(\iota_{x*}\mathrm{\Omega}^{i}(i)[i],E^{i+1})\simeq \mathrm{k}\big\}$$
so $\mathrm{Hom}(E^{i},E^{i})$ equals to $\mathrm{Hom}(E^{i+1},E^{i+1})$ or $\mathrm{Hom}(E^{i+1},E^{i+1})\bigoplus\mathrm{k}$ which depends on whether the extension is non-trivial or not. At last by the fact $p^{*}$ is fully-faithful, we have
$$\mathrm{Hom}(p^{*} l_{*}\mathrm{k}_{x},p^{*} l_{*}\mathrm{k}_{x})=\mathrm{Hom}(l_{*}\mathrm{k}_{x},l_{*}\mathrm{k}_{x})=\mathrm{k}=\mathrm{Hom}(\iota_{x*}\mathrm{\Omega}^{c-1}(c-1)[c-1])=\mathrm{Hom}(E^{c-1},E^{c-1})$$
which tells us each step of extension in the heart filtration (\ref{eq:hfil}) of $p^{*} l_{*}\mathrm{k}_{x}$ should be non-trivial.\\

On the other hand side, let's recall the result in Step 4, $\mathrm{Ext}^{1}(p^{*}(\mathrm{k}_{x}),p^{*}(\mathrm{k}_{x}\otimes\nu^{-1}))$ derives from the contribution of $$\mathrm{Hom}(\iota_{x*}\mathrm{\Omega}^{\lambda}(\lambda),\iota_{x*}\mathrm{\Omega}^{\lambda-1}(\lambda-1)\otimes\nu^{-1})$$ for all $\lambda$,  we can consider a \textit{net} of our extension by combining this fact and our heart filtration in (\ref{eq:hfil}) as following:
$$\begin{tikzcd}[ampersand replacement=\&]
\& E^{\lambda+1}\arrow[d,"e_{\lambda+1}"] \arrow[r]\& E^{\lambda}\arrow[d,"e_{\lambda}"]\arrow[r]\& \iota_{x*}\mathrm{\Omega}^{\lambda}(\lambda)[\lambda]\arrow[d,"\alpha_{\lambda}"] \arrow[r]\& E^{\lambda+1}[1]\arrow[d,"{e_{\lambda+1}[1]}"] \\
 \&E^{\lambda}\otimes\nu^{-1}[1]\arrow[r] \& E^{\lambda-1}\otimes\nu^{-1}[1]\arrow[r] \&\iota_{x*}\mathrm{\Omega}^{\lambda-1}(\lambda-1)[\lambda]\otimes\nu^{-1}\arrow[r] \& E^{\lambda}\otimes\nu^{-1}[2]
\end{tikzcd}$$
with the indices restricted to the nonzero terms of the filtrations.
Also by Step 5, we know the extensions between the horizontal distinguished triangles above are induced by
$$\mathrm{Ext}^{1}(\iota_{x*}\mathrm{\Omega}^{\lambda}(\lambda)[\lambda],E^{\lambda+1})=\mathrm{Ext}^{1}(\iota_{x*}\mathrm{\Omega}^{\lambda}(\lambda)[\lambda],\iota_{x*}\mathrm{\Omega}^{\lambda+1}(\lambda+1)[\lambda+1])$$
so we have a \textit{more sophisticated net} of the extension as following:
$$\begin{tikzcd}[ampersand replacement=\&]
\& \iota_{x*}\mathrm{\Omega}^{\lambda}(\lambda)[\lambda]\arrow[r,"\beta_{\lambda}"]\arrow[d,"\alpha_{\lambda}"] \&\iota_{x*}\mathrm{\Omega}^{\lambda+1}(\lambda+1)[\lambda+1][1]\arrow[d,"\alpha_{\lambda+1}"] \\
 \&\iota_{x*}\mathrm{\Omega}^{\lambda-1}(\lambda-1)[\lambda]\otimes\nu^{-1} \arrow[r,"\beta_{\lambda-1}"] \& \iota_{x*}\mathrm{\Omega}^{\lambda}(\lambda)[\lambda+1]\otimes\nu^{-1}[1]
\end{tikzcd}$$
via calculation on cohomology in Step 4 and Step 5, we know that this \textit{more sophisticated net} is induced by the following commutative diagram on $[\mathds{P}(V_{c})/\mu_{s}]$:
$$\begin{tikzcd}[ampersand replacement=\&]
\& {\mathrm{\Omega}^{\lambda}(\lambda)}\arrow[r,"\beta_{\lambda}"]\arrow[d,"\alpha_{\lambda}"] \&{\mathrm{\Omega}^{\lambda+1}(\lambda)}[1]\arrow[d,"\alpha_{\lambda+1}"] \\
 \&{\mathrm{\Omega}^{\lambda-1}(\lambda-1)\otimes\nu^{-1}} \arrow[r,"\beta_{\lambda-1}"] \& {\mathrm{\Omega}^{\lambda}(\lambda-1)\otimes\nu^{-1}[1]}
\end{tikzcd}$$
The indices are restricted to the nonzero terms. The nonvanishing of
$\alpha_\lambda$ can be checked directly in the Koszul model
\eqref{equ:kouzcenter}. On a blow-up chart $x_1=e$, $x_j=ez_j$, write
$\eta_j=\varepsilon_j-z_j\varepsilon_1$ for its exterior generators.
The root direction $x_c$ gives the boundary map by contraction with
$\varepsilon_c^{\vee}$: it sends $\eta_c$ to $1$ and the other
$\eta_j$ to zero. Thus $\alpha_\lambda\ne0$ for
$1\leq\lambda\leq c-1$; the endpoint $\alpha_0$ is zero.
The twist $\nu^{-1}$ makes this map equivariant. Using the same Koszul
model for the horizontal boundaries retains their compatibility, and
pushforward by $q$ gives the induced maps used in
Corollary~\ref{cor:dualdes}.\\

Step 6: Pushforward $p^{*}\rho^{*}(\mathrm{k}_{x})$ along $q$, using the commutative diagram
$$\begin{tikzcd}[ampersand replacement=\&]
\& {[\mathds{P}(V_{c})/\mu_{s}]}\arrow[r,"\iota_{x}"]\arrow[d,"q"] \&{[wBl_{S}V/\mu_{s}]}\arrow[d,"q"] \\
 \&{\mathcal{P}(1^{c-1},s)} \arrow[r,"\iota_{x}"] \& {\widetilde{wBl_{S}U}}
\end{tikzcd}$$
and Beilinson type result Theorem \ref{thm:beilinson}, there is $$q_{*}\mathrm{\Omega}_{[\mathds{P}(V_{c})/\mu_{s}]}^{\lambda}(\lambda)\otimes\nu^{i}=\mathrm{\Omega}_{\mathcal{P}}^{\lambda}(\mathrm{log}\,D_{\nu^{i^{\sharp}}})(\lambda-i^{\sharp})$$
where $i^{\sharp}:=i\,\mathrm{mod}\,s$.
We keep the factors separate, with the characters and ordering of
Step 1. For $0\leq i<s$, put
\begin{equation}\label{equ:sharp}
\begin{aligned}
\mathcal{F}_{\lambda,i}
&:=\mathcal{H}^{-\lambda}\bigl(q_*p^*(l_*\mathrm{k}_x\otimes\nu^i)\bigr)\\
&=q_*\bigl(\iota_{x*}
\mathrm{\Omega}_{[\mathds{P}(V_c)/\mu_s]}^\lambda(\lambda)
\otimes\nu^i\bigr)\\
&=\iota_{x*}\mathrm{\Omega}_{\mathcal{P}}^\lambda
(\mathrm{log}\,D_{\nu^{i^\sharp}})(\lambda-i^\sharp).
\end{aligned}
\end{equation}
These are factors of the cohomological refinement of the root-filtration
triangles, not summands of a direct-sum decomposition. More explicitly,
applying $q_*$ to \eqref{filtration22} gives the distinguished triangles
\[
q_*p^*E^j\longrightarrow q_*p^*E^{j+1}
\longrightarrow q_*p^*(l_*\mathrm{k}_x\otimes\nu^{j+1})[1]
\longrightarrow q_*p^*E^j[1],
\]
where $E^0=\rho^*\mathrm{k}_x$, $E^s=0$, and the character index is
read modulo $s$. Rotating these triangles back and refining each
pulled-back point sheaf by its cohomology gives the factors
$\mathcal{F}_{\lambda,i}[\lambda]$; equivalently, the displayed
outgoing factors have refinements by
$\mathcal{F}_{\lambda,i}[\lambda+1]$.
The induced filtration of
$\mathcal{H}^{-\lambda}(q_*p^*\rho^*\mathrm{k}_x)$ is computed by
the associated spectral sequence: the sheaves
$\mathcal{F}_{\lambda,i}$ occur on its first page, whereas its limiting
terms are the factors of the total cohomology sheaf. The connecting
maps from Steps 4 and 5 are retained throughout.\\

Step 7: Noticing the extension contributions  from Step 4 and Step 5
  $$\mathrm{Hom}_{[\mathds{P}(V_{c})/\mu_{s}]}(\mathrm{\Omega}^{\lambda}(\lambda)\otimes\nu^{i},\mathrm{\Omega}^{\lambda-1}(\lambda-1)\otimes\nu^{i-1})$$
pushforward along $q$ coincidence with the extension in the dual exceptional collection of $\mathcal{P}(1^{c-1},s)$,
    $$\mathrm{Hom}_{\mathcal{P}(1^{c-1},s)}(\mathrm{\Omega}_{\mathcal{P}}^{\lambda}(\mathrm{log}\,D_{\nu^{i^{\sharp}}})(\lambda-i^{\sharp}),\mathrm{\Omega}_{\mathcal{P}}^{\lambda-1}(\mathrm{log}\,D_{\nu^{(i-1)^{\sharp}}})(\lambda-1-(i-1)^{\sharp}))$$
as we obtained in Corollary \ref{cor:dualdes}, what happen for a general \textit{more sophisticated net} of the extension
$$\begin{tikzcd}[ampersand replacement=\&]
\& \iota_{x*}\mathrm{\Omega}^{\lambda}(\lambda)[\lambda]\otimes\nu^{i}\arrow[r,"\beta_{\lambda}"]\arrow[d,"\alpha_{\lambda}"] \&\iota_{x*}\mathrm{\Omega}^{\lambda+1}(\lambda+1)[\lambda+1][1]\otimes\nu^{i}\arrow[d,"\alpha_{\lambda+1}"] \\
 \&\iota_{x*}\mathrm{\Omega}^{\lambda-1}(\lambda-1)[\lambda]\otimes\nu^{i-1} \arrow[r,"\beta_{\lambda-1}"] \& \iota_{x*}\mathrm{\Omega}^{\lambda}(\lambda)[\lambda+1]\otimes\nu^{i-1}[1]
\end{tikzcd}$$
after pushforward along $q$ is it becomes the following \textit{elementary net}:
$$\begin{tikzcd}[ampersand replacement=\&]
\& \iota_{x*}\mathrm{\Omega}_{\mathcal{P}}^{\lambda}(\mathrm{log}\,D_{\nu^{i^{\sharp}}})(\lambda-i^{\sharp})[\lambda]\arrow[r,"q_{*}\beta_{\lambda}"]\arrow[d,"q_{*}\alpha_{\lambda}"] \&\iota_{x*}\mathrm{\Omega}_{\mathcal{P}}^{\lambda+1}(\mathrm{log}\,D_{\nu^{i^{\sharp}}})(\lambda+1-i^{\sharp})[\lambda+1][1]\arrow[d,"q_{*}\alpha_{\lambda+1}"] \\
 \&\iota_{x*}\mathrm{\Omega}_{\mathcal{P}}^{\lambda-1}(\mathrm{log}\,D_{\nu^{(i-1)^{\sharp}}})(\lambda-1-(i-1)^{\sharp})[\lambda] \arrow[r,"q_{*}\beta_{\lambda-1}"] \& \iota_{x*}\mathrm{\Omega}_{\mathcal{P}}^{\lambda}(\mathrm{log}\,D_{\nu^{(i-1)^{\sharp}}})(\lambda-(i-1)^{\sharp})[\lambda+1][1]
\end{tikzcd}$$
we can see the extended vertical complex of the \textit{elementary net} is just $\mathcal{D}_{\lambda+i}$ and $\mathcal{D}_{\lambda+i+1}$ on the weighted projective stack, hence we can get a filtration of  $q_{*}p^{*}\rho^{*} \mathrm{k}_{x}$ with complex of inherent morphisms as the classic case by continuously applying the triangulated category version of nine lemma (two octahedral axioms) to our \textit{elementary net}:
$$\begin{tikzcd}[column sep=0.5em]
 & 0 \arrow{rr}&& E_{-(c+s-2)}\arrow{dl}\arrow{rr}&& E_{-(c+s-3)}  \arrow{dl} \\
&& \iota_{x*}\mathcal{D}_{(c+s-2)}\arrow[ul,dashed,"\Delta"]&& \iota_{x*}\mathcal{D}_{(c+s-3)}\arrow[ul,dashed,"\Delta"]
\end{tikzcd}...
\begin{tikzcd}[column sep=0.5em]
 &    E_{-1}\arrow{rr}&& q_{*}p^{*}\rho^{*} \mathrm{k}_{x}\arrow{dl} \\
 && \iota_{x*}\mathcal{D}_{0}\arrow[ul,dashed,"\Delta"]
\end{tikzcd}$$
noticing
$$q_{*}p^{*}\rho^{*} (\mathrm{k}_{x})\simeq q_{*}q^{*}\pi^{*} (\mathrm{k}_{x})\simeq\pi^{*}\mathrm{k}_{x}$$
so if changing our notation we have the proposition holds, which finishes our proof.
\end{proof}

Then we will generalize our previous results to varieties admit iterated covering (ramified along a smooth simple normal crossing divisor).\\

We assume $U$ is a $n$ dimensional smooth quasi-projective variety, and $S$ is a smooth codimension $c$ subvariety contained in $U$ which can be realized as complete intersection of smooth divisors $D_{1},...,D_{b},...,D_{c}$ where $1<b\leq c$, we temporarily require that the first $b$ divisors $D_{1},...,D_{b}$ are also permutable, i.e. for any $\sigma$ belonging to the symmetry group $S_{b}$ the divisors sequence $D_{\sigma_{1}},...,D_{\sigma_{b}},...,D_{c}$ is also a locally regular sequence. Specially $U$ admits a $\mu:=\prod^{b}_{i=1}\mu_{s_{i}}$ cyclic covering
$$\rho: V\longrightarrow U$$
branched along divisors $D_{1}$,...,$D_{b}$ with respect to weights $s_{1}$,...,$s_{b}$. We naturally consider its stacky smoothing and denote it also by $\rho$:
$$\rho: [V/\mu]\longrightarrow U$$
We are expected to show how certain skyscraper sheaf changes under pullback along $\rho$ as before, and we still use the character system oppositely defined in \cite[Section 4.1]{KAC} which is compatible with (\ref{equ:actclassical}) for equivariant cyclic covering.\\

\begin{lemma}\label{lem:globalpullbackskycraper}
In the above simple-normal-crossing situation, let $\mathcal M_j$ be the tautological positive root line and let $\mathcal I_j=\mathcal M_j^{-1}$ be the ideal of the reduced ramification divisor corresponding to $D_j$.
\begin{enumerate}
\item\label{enu:globalpullbackskycraper1}
There is an identification
\[
[V/\mu]\simeq
\sqrt[s_1]{D_1/U}\times_U\cdots\times_U\sqrt[s_b]{D_b/U}.
\]
The order of the root constructions can be permuted. For an order $\sigma$, applying Lemma~\ref{lem:pullbackskycraper}(\ref{num:pullbackskycraper1}) at each stage gives the corresponding semi-orthogonal decomposition
\[
\mathrm{D}^{b}([V/\mu])=
\left\langle
\mathrm\Upsilon^{\sigma_b}_{s_{\sigma_b}-1},\ldots,\mathrm\Upsilon^{\sigma_b}_{1},
\ldots,
\mathrm\Upsilon^{\sigma_1}_{s_{\sigma_1}-1},\ldots,\mathrm\Upsilon^{\sigma_1}_{1},
\rho^*\mathrm{D}^{b}(U)
\right\rangle.
\]
Here the component $\mathrm\Upsilon^{\sigma_j}_i$ is obtained from $\mathrm{D}^{b}(\widehat D_{\sigma_j})$ by the root-divisor embedding at stage $j$, twisted by the $i$-th power of its root line and then pulled back through the remaining root constructions. The divisor $\widehat D_{\sigma_j}$ is the inverse image of $D_{\sigma_j}$ before this root is taken.
\item\label{enu:globalpullbackskycraper2}
For a geometric point $x\in S$, let $\mathcal G_x$ denote the trivial-character skyscraper on the residual gerbe over $x$, pushed forward to $[V/\mu]$. The sheaf $\rho^*\mathrm k_x$ has the finite multifiltration
\[
F^{\mathbf a}=\mathcal I_1^{a_1}\cdots\mathcal I_b^{a_b}\rho^*\mathrm k_x,
\qquad \mathbf a=(a_1,\ldots,a_b)\in\mathbb N^b,
\]
with associated factors
\[
\frac{F^{\mathbf a}}{\sum_{j=1}^bF^{\mathbf a+\mathbf e_j}}
\simeq
\mathcal G_x\otimes\nu_1^{-a_1}\otimes\cdots\otimes\nu_b^{-a_b},
\qquad 0\le a_j<s_j.
\]
This multifiltration is compatible with permutations of the root directions. Choosing an order refines it to a finite filtration, and rotation of its short exact sequences gives a left Postnikov system. No uniqueness of all such refinements is asserted. Its factors are associated factors of this filtration, not projections of $\rho^*\mathrm k_x$ to the root components.
\end{enumerate}
\end{lemma}
\begin{proof}
An object of the product of root stacks is a collection of roots of the pairs defining $D_j$. Their independent specification proves that the product is the successive root construction in any order. The divisors remain smooth after each root, because they meet with simple normal crossings. Lemma~\ref{lem:pullbackskycraper}(\ref{num:pullbackskycraper1}) therefore applies at every stage. Pulling the decompositions through the later fully faithful root pullbacks proves \ref{enu:globalpullbackskycraper1}.

For \ref{enu:globalpullbackskycraper2}, work near $x$ with parameters $u_1,\ldots,u_b,z_{b+1},\ldots,z_n$ such that $D_j=(u_j=0)$. Set $u_j=t_j^{s_j}$. The simultaneous root cover is finite flat, and its fibre is the equivariant module
\[
\rho^*\mathrm k_x=
\frac{\mathrm k[t_1,\ldots,t_b]}{(t_1^{s_1},\ldots,t_b^{s_b})}
\simeq
\bigotimes_{j=1}^b\mathrm k[t_j]/(t_j^{s_j}).
\]
The root coordinate $t_j$ has character $\nu_j^{-1}$. The numerator $F^{\mathbf a}$ is the monomial ideal generated by $t_1^{a_1}\cdots t_b^{a_b}$ in this module. Quotienting by the indicated next ideals leaves exactly this monomial, and proves the character formula. Refining by monomial ideals gives the asserted finite filtration. Equivalently, its total-degree filtration has associated quotients equal to the direct sums of these factors for fixed $a_1+\cdots+a_b$. The ideal-sheaf definition shows that the multifiltration glues independently of the chosen parameters. Reordering the factors in the displayed tensor product proves compatibility with permutations.

For clarity about the extensions, multiplication by distinct $t_j$ commutes at the sheaf level. The degree-one boundary classes on the smooth root chart satisfy the Koszul signs. If $\xi_j$ denotes the normal-direction class associated with $t_j$, with its character twists understood, then
\[
\xi_i[1]\circ\xi_j=-\xi_j[1]\circ\xi_i
\qquad(i\ne j).
\]
Indeed, the Koszul resolution of the residual skyscraper for these regular parameters has exterior generators, so its Yoneda algebra is exterior in these normal directions. Tensoring the one-variable filtrations uses the corresponding signs in the total complexes. This specifies the extension data coming from the actual root module; it is not inferred from uniqueness of a multi-step Postnikov system or from unsigned commutative squares.
\end{proof}

\begin{remark}
We considered some cohomology calculations on skyscrapersheaves above, in fact, we only need to do these calculations on local rings (and the quotient stack of local rings). In a Noetherian Cohen--Macaulay local ring, a proper regular sequence is permutable and extends to a system of parameters; it is itself a system of parameters only when its length is the dimension. Thus a separate permutability assumption is unnecessary here. The regularity assumption is retained. In the smooth coordinate situation, the divisors $D_{1},\ldots,D_{b}$ meet with simple normal crossings after shrinking around the centre.
\end{remark}

The following proposition is a natural generalization of Proposition \ref{prop:singlecover} in the iterated cover's cases:

\begin{proposition}\label{prop:pullbackskycaperalongsmoothcenter}
Consider the following commutative diagram, we have:
$$\begin{tikzcd}[ampersand replacement=\&]
\&{[Bl_{S}V/\mu]}\arrow[r,"q"]\arrow[d,"p"] \&\widetilde{wBl_{S}U}\arrow[d,"\pi"] \\
 \&{[V/\mu]} \arrow[r,"{\rho}"] \&  U
\end{tikzcd}$$
for any geometric point $x$ in $S$, $\pi^{*} (\mathrm{k}_{x})$ has a unique Postnikov filtration:
\[\begin{adjustbox}{max width=0.98\linewidth}\(\begin{tikzcd}[column sep=0.5em]
 & 0 \arrow{rr}&& E_{-t}\arrow{dl}\arrow{rr}&& E_{-t+1}  \arrow{dl} \\
&& \iota_{x*}\mathcal{D}_{t}\arrow[ul,dashed,"\Delta"]&& \iota_{x*}\mathcal{D}_{t-1}\arrow[ul,dashed,"\Delta"]
\end{tikzcd}
...\begin{tikzcd}[column sep=0.5em]
 &   E_{-2}\arrow{rr}&&  E_{-1}\arrow{rr}\arrow{dl}&& \pi^{*} (\mathrm{k}_{x})\arrow{dl} \\
 && \iota_{x*}\mathcal{D}_{1}\arrow[ul,dashed,"\Delta"]  &&\iota_{x*}\mathcal{D}_{0}\arrow[ul,dashed,"\Delta"]
\end{tikzcd}\)\end{adjustbox}\]
where $\mathcal{D}_{i}$ is the dual exceptional object for $\mathcal{P}(1^{c-b},s_{b},...,s_{1})$, and $t:=c-b+\sum_{i=1}^{b}s_{i}-1$.\\
\end{proposition}

\begin{proof}
Since the proof is just a boring repetition of the previous proof of Proposition \ref{prop:singlecover}, so instead of writing down all the mundane details we just point out where our previous proof is needed to be adjusted.\\

Step 1: Use the multifiltration of
Lemma~\ref{lem:globalpullbackskycraper}(\ref{enu:globalpullbackskycraper2})
on the root module
\[
\rho^*\mathrm k_x=
\mathrm k[t_1,\ldots,t_b]/(t_1^{s_1},\ldots,t_b^{s_b}).
\]
Its associated sheaves have characters
$\nu_1^{-m_1}\otimes\cdots\otimes\nu_b^{-m_b}$,
where $0\le m_j<s_j$. Choose an order and rotate its quotient sequences
when forming the left Postnikov system. The tensor complexes retain the
Koszul signs; no uniqueness of all refinements is used.\\

Step 2: We have an equivariant Koszul resolution for center $S$:

$$\mathcal{E}^{\vee}:=\bigoplus^{c}_{i=b+1}\mathcal{O}(-\widehat{D}_{i})\bigoplus_{j=1}^{b}\mathcal{M}^{-1}_{\widehat{D}_{j}}\longrightarrow\mathcal{O}_{[V/\mu]}\longrightarrow\mathcal{O}_{[S/\mu]}\longrightarrow0$$
then consider the diagram:

 $$\begin{tikzcd}[column sep=0.5em]
 & {[Bl_{S}V/\mu]} \arrow[dr]\arrow{rr}&& {{\mathds{P}_{[V/\mu]}(\mathcal{E})}}\arrow{dl}\\
&&{[V/\mu]}
\end{tikzcd}$$
for any geometric point $x$ in $[S/\mu]$, we have $$\mathcal{H}^{-\lambda} (p^{*} \mathrm{k}_{x})=\iota_{x*}\mathrm{\Omega}_{[\mathds{P}(V_{c})/\mu]}^{\lambda}\otimes\mathcal{O}_{[\mathds{P}(V_{c})/\mu]}(\lambda)$$

Step 3: We act $p^{*}$ on filtration system of $\rho^{*}(\mathrm{k}_{x})$ in Step 1.\\

Step 4:
The factors that might contribute to the unique extension are
$$\mathrm{Ext}^{1}(p^{*}(\mathrm{k}_{x}\otimes\nu_{j}^{i}),p^{*}(\mathrm{k}_{x}\otimes\nu_{j}^{i-1}))=\mathrm{k}$$
which come from the contribution of
$$\mathrm{Hom}_{[\mathds{P}(V_{c})/\mu]}(\mathrm{\Omega}^{a}(a)\otimes\nu_{j}^{i},\mathrm{\Omega}^{a-1}(a-1)\otimes\nu_{j}^{i-1})=\mathrm{k}$$
for all $a$ and $j$, further, they all contribute to the extension.\\

Step 5: Pushforward $p^{*}\rho^{*}(\mathrm{k}_{x})$ along $q$, using the commutative diagram
$$\begin{tikzcd}[ampersand replacement=\&]
\& {[\mathds{P}(V_{c})/\mu]}\arrow[r,"\iota_{x}"]\arrow[d,"q"] \&{[wBl_{S}V/\mu]}\arrow[d,"q"] \\
 \&{\mathcal{P}(1^{c-b},s_{b},..,s_{1})} \arrow[r,"\iota_{x}"] \& {\widetilde{wBl_{S}U}}
\end{tikzcd}$$
By Theorem \ref{thm:beilinson} and \eqref{equ:sharp}, the individual
factors are now indexed by the deck characters $\chi$:
\[
\mathcal{F}_{\lambda,\chi}:=\iota_{x*}
\mathrm{\Omega}^\lambda(\mathrm{log}\,D_{-\chi})(\lambda-|-\chi|).
\]
As in Step 6 of Proposition \ref{prop:singlecover}, these are the
cohomology sheaves of the successive pulled-back root factors.
Their shifts $[\lambda]$, together with the connecting maps, refine
the filtration triangles of $q_*p^*\rho^*\mathrm{k}_x$.
The filtration on its total cohomology is obtained from the limiting
terms of the associated spectral sequence, not by taking a direct sum
over $\chi$.\\

Step 6:
We can see the \textit{elementary net} of our extension is of form
$$\begin{tikzcd}[ampersand replacement=\&]
\& \iota_{x*}\mathrm{\Omega}_{\mathcal{P}}^{\lambda}(\mathrm{log}\,D_{\nu})(\lambda-|\nu|)[\lambda]\arrow[r,"q_{*}\beta_{\lambda}"]\arrow[d,"q_{*}\alpha_{\lambda}"] \&\iota_{x*}\mathrm{\Omega}_{\mathcal{P}}^{\lambda+1}(\mathrm{log}\,D_{\nu})(\lambda+1-|\nu|)[\lambda+1][1]\arrow[d,"q_{*}\alpha_{\lambda+1}"] \\
 \&\iota_{x*}\mathrm{\Omega}_{\mathcal{P}}^{\lambda-1}(\mathrm{log}\,D_{\nu-\nu_{i}})(\lambda-1-|\nu-\nu_{i}|)[\lambda] \arrow[r,"q_{*}\beta_{\lambda-1}"] \& \iota_{x*}\mathrm{\Omega}_{\mathcal{P}}^{\lambda}(\mathrm{log}\,D_{\nu-\nu_{i}})(\lambda-|\nu-\nu_{i}|)[\lambda+1][1]
\end{tikzcd}$$
In the pulled-back Koszul complex, write $x_i=u^{a_i}\xi_i$ and filter $u^m e_I$ by weight $m+\sum_{i\in I}a_i$, where $e_I$ is the Koszul basis element indexed by $I$.
The weight-$k$ quotient on the exceptional fibre has the signed Koszul differential in the $\xi_i$. Taking derived Hom from $\mathcal O_{\mathcal P}(j)$ gives the degree-$(k-j)$ part of the polynomial Koszul complex for $0\le j,k\le t$, so the quotient is $\mathcal D_k$ by Corollary~\ref{cor:dualvanishing}.
\end{proof}

\begin{theorem}\label{thm:smoohtweightedblowup}
If $X$ is a smooth variety, after $(a_{1},..,a_{c})$ \footnote{\, A sequence of positive integers has no non-trivial common divisor \ref{no:ai}.\label{foot:ai}}-weighted blowing up at a codimension $c$ smooth center $S$, we get another variety (algebraic space) $Y$ with only cyclic quotient singularities and we denote its canonical stack by $\widetilde{Y}$.
$$\begin{tikzcd}[ampersand replacement=\&]
E_{\widetilde{Y}}\arrow[r,"\iota"]\arrow[d,"\pi_{E_{\widetilde{Y}}}"] \&\widetilde{Y}\arrow[d,"\pi"]  \\
S \arrow[r,"\kappa"] \&X
\end{tikzcd}$$
Then we have a semi-orthogonal decomposition of $\widetilde{Y}$ with respect to $X$ as following:
$$\mathrm{D}^{b}(\widetilde{Y})\cong\langle \mathrm{Im}\,\mathrm{\Psi}_{-\sum_{i}{a_{i}}+1},...,\mathrm{Im}\,\mathrm{\Psi}_{-1},\mathrm{Im}\, \mathrm{\Psi}\rangle$$
where $\mathrm{\Psi}$ and $\mathrm{\Psi}_{i}$ are defined by $\mathrm{\Psi}(-):=\pi^{*}(-)$ and $\mathrm{\Psi}_{i}(-):=\iota_{*}\bigl(\pi_{E_{\widetilde{Y}}}^{*}(-)\otimes\mathcal{O}_{E_{\widetilde{Y}}}(i)\bigr)$ for any integer $i$.

\end{theorem}
\begin{proof}

Firstly, we adjust a global to local argument to suit Proposition  \ref{prop:pullbackskycaperalongsmoothcenter} to the local situation, e.g. the same techniques we will use in the proof of Proposition \ref{prop:pullbackcenterlocal}. Then we replace \cite[Proposition 11.12]{H} with our Proposition \ref{prop:pullbackskycaperalongsmoothcenter}, then use the same proof as in \cite[Proposition 11.18]{H} to arrive at this conclusion.
\end{proof}

\begin{remark}
We note that the above theorem can be easily generalized to the case of weighted blow-up of a smooth Deligne-Mumford stack $X$ at a smooth substack center $S$, when $S$ is defined locally as a quotient compatible or independent of the direction we take root in our construction of the weighted blow-up. Since the representations of our iterated roots construction are independent of the representations of our local stabilizer group at points on $S$.  Another method is that we take advantage of the presentation of the local quotient structure and use the equivariant results in \cite{El}.
\end{remark}

\subsection{Derived category of Kawamata blow-up}

For positive-dimensional centres we assume that the quotient gerbe is globally
split, that the character line $\nu$ is fixed globally, and that the chosen
weighted normal data glue. These are the data used in the diagrams below.
At an isolated quotient point over $\mathrm k=\mathbb C$, we fix these data in a quotient chart.
Without them the local character calculations do not by themselves give
the displayed decompositions into untwisted copies of $\mathrm{D}^{b}(S)$.

In this section, we first consider the derived category structure of Kawamata blow-up for a local ring with cyclic quotient singularity on its singular locus. Geometrically, we assume $U$ is a smooth quasi-projective $n$-dimensional variety that admits a cyclic group $\mu_{r}$ action, such that the fixed locus $S$ of the action is small and completely intersected of codimension $c$. So the normal sheaf $\mathcal{N}_{S|U}$ is locally free and inherit the natural action of $\mu_{r}$ from tangent space of $U$ restricting on $S$, we assume there is a characteristic decomposition of $\mathcal{N}_{S|U}$ after shrinking $U$:
 $$\mathcal{N}_{S|U}=\bigoplus^{c}_{i=1}\mathcal{N}_{S|U,i}$$
where $\mathcal{N}_{S|U,i}$ is invertible for any $i$ with a $\mu_{r}$ action (representation) by weight $a_{i}$.\\

After further shrinking of $U$, we also assume each component $\mathcal{N}_{S|U,i}$ correspondences a smooth ray divisor $D_{i}$ and sum of them $\sum_{i} D_{i}$ is SNC and $\mu_{r}$ also act on divisor $D_{i}$  by weight $a_{i}$ respectively.  If we replace $U$ by an iterated cyclic coverings $V$ branched overall $D_{i}$  with degree $a_{i}$ for $i=1,...,c$, then we have two compatible group actions on $V$.
$$\mu:=\mu_{a_{1}}\times,..,\times\mu_{a_{c}}\text{\quad and\quad} \mu_{r}$$
$$G:=\mu_{r}\times\mu$$
under a local coordinate $(x_{1},..,x_{n})$ of $V$ we have compatible action:
\begin{equation}\label{equ:actclassicalrs}
  (\zeta_{a_{1}},..,\zeta_{a_{c}})\times(\zeta_{r})\times(x_{1},..,x_{n})\longmapsto(\zeta_{r}\cdot\zeta_{a_{1}}x_{1},...,\zeta_{r}\cdot\zeta_{a_{c}}x_{c},x_{c+1},...,x_{n})
\end{equation}
specially $\mu_{r}$ act on the ramified ray divisor $M_{D_{i}}$ on $V$ with weight $1$. This choice is because, to satisfy all cyclic coverings, our group action and character system is consistent with (\ref{equ:actclassical}).

\begin{proposition}\label{prop:pullbackcenterglobal}
Consider the following commutative diagram:
$$\begin{tikzcd}[ampersand replacement=\&]
\&{[Bl_{S}V/G]}\arrow[r,"q"]\arrow[d,"p"] \&{[\widetilde {wBl_{S}U}/\mu_{r}]}\arrow[d,"\pi"] \arrow[r,"\varpi"]\&\widetilde{w_{1/r}Bl_{S}U}\\
 \&{[V/G]} \arrow[r,"{\rho}"] \&  {[U/\mu_{r}]}
\end{tikzcd}$$
we have for any generalized geometrical point $\mathrm{k}_{x}\otimes\nu^{i}$ on $[S/\mu_{r}]$, $\overline{\mathcal{D}}_{x,i}:=\varpi_{*}\pi^{*}(\mathrm{k}_{x}\otimes\nu^{i})$ has the following finite Postnikov filtration:
\[\begin{adjustbox}{max width=0.98\linewidth}\(\begin{tikzcd}[column sep=0.5em]
 & 0 \arrow{rr}&& E_{-i-kr}\arrow{dl} \\
&& \iota_{x*}\mathcal{D}_{i+kr}\arrow[ul,dashed,"\Delta"]
\end{tikzcd}......\begin{tikzcd}[column sep=0.5em]
 &   E_{-i-r}\arrow{rr}&&  E_{-i}\arrow{rr}\arrow{dl}&& \overline{\mathcal{D}}_{x,i}\arrow{dl} \\
 && \iota_{x*}\mathcal{D}_{i+r}\arrow[ul,dashed,"\Delta"]  && \iota_{x*}\mathcal{D}_{i}\arrow[ul,dashed,"\Delta"]
 \end{tikzcd}\)\end{adjustbox}\]
For $\sum_{j=1}^c a_j\le r$ and $0\le i<r$, the image is $\iota_{x*}\mathcal D_i$ if $i<\sum_j a_j$, and is zero otherwise.\\

\end{proposition}

\begin{proof}
Its proof is the same as the previous argument in Proposition \ref{prop:pullbackskycaperalongsmoothcenter} and \ref{prop:singlecover}, so we just point out which parts need to be changed.\\

Step 1: Pull back the point sheaf through the simultaneous root cover,
retaining both the deck characters and the additional $\mu_r$-action.
Write $\mathcal J_j$ for the full equivariant ideal of the $j$-th reduced
root divisor. In the local coordinates above its defining sequence is
\[
0\longrightarrow
\mathcal J_j=\mathcal M_{[D_j/\mu_r]}^{-1}\otimes\nu^{-1}
\longrightarrow\mathcal{O}_{[V/G]}
\longrightarrow\mathcal{O}_{[D_j/G]}\longrightarrow0.
\]
The actual ideal-power multifiltration has factors
\[
\mathrm k_x\otimes
\bigotimes_{j=1}^b
\left(\mathcal J_j|_{\mathcal G_x}\right)^{m_j},
\qquad 0\le m_j<s_j,
\]
where $\mathcal G_x$ is the residual gerbe over $x$.
This is Lemma~\ref{lem:globalpullbackskycraper} with the full
equivariant line bundles retained. Reducing a deck-character index
modulo $s_j$ does not reduce the additional $\mu_r$-character or erase
the corresponding base-line twist.\\

Step 2: We have an equivariant Koszul resolution for center $S$:
\begin{equation}\label{equ:equacenter}
\mathcal{E}^{\vee}:=\big (\bigoplus^{c}_{i=b+1}\mathcal{O}(-\widehat{D}_{i})\bigoplus_{j=1}^{b}\mathcal{M}^{-1}_{\widehat{D}_{j}}\big)\otimes\nu^{-1}\longrightarrow\mathcal{O}_{[V/G]}\longrightarrow\mathcal{O}_{[S/G]}\longrightarrow0
\end{equation}
then consider the diagram:
 $$\begin{tikzcd}[column sep=0.5em]
 & {[Bl_{S}V/G]} \arrow[dr]\arrow{rr}&& {{\mathds{P}_{[V/G]}(\mathcal{E})}}\arrow{dl}\\
&&{[V/G]}
\end{tikzcd}$$

For any geometric point $x$ in $[S/G]$, we have a heart filtration
\[
\mathcal{H}^{-\lambda} (p^{*} \mathrm{k}_{x}) = \iota_{x*}\Omega_{[\mathds{P}(V_{c})/G]}^{\lambda} \otimes \mathcal{O}_{[\mathds{P}(V_{c})/G]}(\lambda)
\]
noticing $[\mathds{P}(V_{c})/G]$ can be realized as a $\mu_{r}$-gerbe over $[\mathds{P}(V_{c})/\mu]$. In fact the equivariant presentation of $\Omega_{[\mathds{P}(V_{c})/G]}(1)$ comes from the following equivariant filtration of $p^{*} \mathcal{O}_{S}$:
\[
\begin{tikzcd}[column sep=normal]
\cdots \arrow[r] &
\bigwedge^{2}\mathcal{E}^{\vee}(E)\otimes\nu^{-1} \arrow[r] \arrow[d, two heads] &
\mathcal{E}^{\vee}\otimes\nu^{-1} \arrow[r,"v"] \arrow[d, two heads] &
\mathcal{O} \arrow[r] \arrow[d, equal] &
\mathcal{F}^1(p^{*} \mathcal{O}_{S}) \arrow[d, two heads] \\
\cdots \arrow[r] &
\bigwedge^{2}\mathcal{E}^{\vee}(2E) \arrow[r] \arrow[d, hook] &
\mathcal{E}^{\vee}(E) \arrow[r] \arrow[d, hook] &
\mathcal{O} \arrow[r] \arrow[d] &
0 \arrow[d, hook] \\
\cdots \arrow[r] &
\bigwedge^{2}\mathcal{E}^{\vee}(2E)|_E \arrow[r] &
\mathcal{E}^{\vee}(E)|_E \arrow[r] &
0 \arrow[r] &
\operatorname{gr}^1(p^{*} \mathcal{O}_{S})
\end{tikzcd}
\]
where $v$ is determined by the composite
$\mathcal{E}^{\vee}\otimes\nu^{-1} \to \mathcal{O}(-E)\otimes\nu^{-1} \to \mathcal{O}$, then
\[
\begin{tikzcd}[column sep=normal]
\cdots \arrow[r] &
\bigwedge^{2}\mathcal{E}^{\vee}(2E)\otimes\nu^{-2} \arrow[r] \arrow[d, hook] &
\mathcal{E}^{\vee}\otimes\nu^{-1} \arrow[r] \arrow[d, hook] &
\mathcal{O} \arrow[r] \arrow[d, equal] &
\mathcal{F}^2(p^{*} \mathcal{O}_{S}) \arrow[d, hook] \\
\cdots \arrow[r] &
\bigwedge^{2}\mathcal{E}^{\vee}(E)\otimes\nu^{-1} \arrow[r] \arrow[d, two heads] &
\mathcal{E}^{\vee}\otimes\nu^{-1} \arrow[r,"v"] \arrow[d, two heads] &
\mathcal{O} \arrow[r] \arrow[d] &
\mathcal{F}^1(p^{*} \mathcal{O}_{S}) \arrow[d, two heads] \\
\cdots \arrow[r] &
\bigwedge^{2}\mathcal{E}^{\vee}(E)|_E\otimes\nu^{-1} \arrow[r] &
0 \arrow[r] &
0 \arrow[r] &
\operatorname{gr}^2(p^{*} \mathcal{O}_{S})
\end{tikzcd}
\]
Also note that the exact sequence:
\[
\begin{tikzcd}[column sep=normal]
\cdots \arrow[r] &
\bigwedge^{2}\mathcal{E}^{\vee}(2E)|_E \arrow[r] &
\mathcal{E}^{\vee}(E)|_E \arrow[r] &
\mathcal{O}|_E \arrow[r] &
0
\end{tikzcd}
\]
we have equivariantly $\Omega_{[\mathds{P}(V_{c})/G]}(1)$ can be realized as a submodule of $\mathcal{E}^{\vee}\otimes\nu^{-1}|_E$ hence it carries a nontrivial representation $\nu^{-1}$. So we get an equivariant version of filtration
$$\mathcal{H}^{-\lambda} (p^{*} \mathrm{k}_{x})=\iota_{*}\mathrm{\Omega}_{[\mathds{P}(V_{c})/\mu]}^{\lambda}\otimes\mathcal{O}_{[\mathds{P}(V_{c})/\mu]}(\lambda)\otimes\nu^{-\lambda}$$

Step 3: We act $p^{*}$ on filtration system we obtained in Step 1.\\

Step 4: Pushforward $p^{*}\rho^{*}(\mathrm{k}_{x})$ along $q$, using the commutative diagram
$$\begin{tikzcd}[ampersand replacement=\&]
\& {[\mathds{P}(V_{c})/G]}\arrow[r,"\iota_{x}"]\arrow[d,"q"] \&{[wBl_{S}V/G]}\arrow[d,"q"] \\
 \&{[\mathcal{P}(1^{c-b},s_{b},..,s_{1})/\mu_{r}]} \arrow[r,"\iota_{x}"] \& {[\widetilde{wBl_{S}U}/\mu_{r}]}
\end{tikzcd}$$
By Theorem \ref{thm:beilinson} and \eqref{equ:sharp}, the individual
factors, with their additional $\mu_r$-characters, are
\[
\mathcal{F}_{\lambda,\chi}:=\iota_{x*}
\mathrm{\Omega}^\lambda(\mathrm{log}\,D_{-\chi})(\lambda-|-\chi|)
\otimes\nu^{-\lambda-|\chi|}.
\]
Here $\chi$ ranges over the deck characters. Their shifts $[\lambda]$
are the factors in the refinement of the root-filtration triangles
for $q_*p^*\rho^*\mathrm{k}_x\simeq\pi^*\mathrm{k}_x$.
We retain the connecting maps and the additional characters; as above,
the factors of total cohomology are the limiting terms of the
associated spectral sequence.\\

The weighted Koszul filtration in Proposition~\ref{prop:pullbackskycaperalongsmoothcenter} now gives the factors $\iota_{x*}\mathcal D_k\otimes\nu^{-k}$, since $u^m e_I$ has character $\nu^{-(m+\sum_{i\in I}a_i)}$. Thus we have the filtration
\[\begin{adjustbox}{max width=0.98\linewidth}\(\begin{tikzcd}[column sep=0.5em]
 & 0 \arrow{rr}&& E_{-t+1}\arrow{dl}\arrow{rr}&& E_{-t+2}  \arrow{dl} \\
&& \iota_{x*}\mathcal{D}_{(t-1)}\otimes\nu^{-(t-1)}\arrow[ul,dashed,"\Delta"]&& \iota_{x*}\mathcal{D}_{t-2}\otimes\nu^{-(t-2)}\arrow[ul,dashed,"\Delta"]
\end{tikzcd}
...\begin{tikzcd}[column sep=0.5em]
 &   E_{-2}\arrow{rr}&&  E_{-1}\arrow{rr}\arrow{dl}&& \pi^{*} \mathrm{k}_{x}\arrow{dl} \\
 && \iota_{x*}\mathcal{D}_{1}\otimes\nu^{-1}\arrow[ul,dashed,"\Delta"]  && \iota_{x*}\mathcal{D}_{0}\arrow[ul,dashed,"\Delta"]
\end{tikzcd}\)\end{adjustbox}\]
where $t:=\sum_{i=1}^{b}s_{i}+c-b$.\\

Step 5: Pushforward along morphism $\varpi$,
$$\begin{tikzcd}[ampersand replacement=\&]
\& {[\mathcal{P}/\mu_{r}]}\arrow[r,"\iota_{x}"]\arrow[d,"\varpi"] \&{[wBl_{S}U/\mu_{r}]}\arrow[d,"\varpi"] \\
 \&{\mathcal{P}} \arrow[r,"\iota_{x}"] \& {\widetilde{wBl_{S}U/\mu_{r}}}
\end{tikzcd}$$
which is a $\mu_{r}$ gerbe when restricting on $[\mathcal{P}/\mu_{r}]$, so pushforward cancels all the terms with nontrivial representation. That is:
 \[\begin{adjustbox}{max width=0.98\linewidth}\(\begin{tikzcd}[column sep=0.5em]
 & 0 \arrow{rr}&& E_{-i-kr}\arrow{dl}\arrow{rr}&& E_{-i-kr+r}  \arrow{dl} \\
&& \iota_{x*}\mathcal{D}_{i+kr}\arrow[ul,dashed,"\Delta"]&& \iota_{x*}\mathcal{D}_{i+kr-r}\arrow[ul,dashed,"\Delta"]
\end{tikzcd}...\begin{tikzcd}[column sep=0.5em]
 &   E_{-i-r}\arrow{rr}&&  E_{-i}\arrow{rr}\arrow{dl}&& \varpi_{*}\pi^{*} (\mathrm{k}_{x}\otimes\nu^{i})\arrow{dl} \\
 && \iota_{x*}\mathcal{D}_{i+r}\arrow[ul,dashed,"\Delta"]  && \iota_{x*}\mathcal{D}_{i}\arrow[ul,dashed,"\Delta"]
\end{tikzcd}\)\end{adjustbox}\]
we finish the proof.
\end{proof}

Then we give a global to local argument in the derived category, to suit the previous result to general projective varieties' cases.\\

If $X$ is a $n$ dimensional (projective) variety admits a $\mu_{r}$ cyclic quotient small center $S$, which is also smooth and locally complete intersected of codimension $c$. Furthermore, we assume the cyclic group action is of weight $(a_{1},..,a_{c})$ at center $S$, which means the action at analytical neighborhood is a cyclic quotient of affine space $\mathds{A}^{n}$ with weight $(a_{1},..,a_{c})$ and these neighborhoods are compatible with each other.\\

We roughly explain  Vistoli's construction here. By the definition of quotient singularity, we have an isomorphism between formal stalks $\widehat{\mathcal{O}}_{X,x}$ and $\widehat{\mathcal{O}}_{\mathds{A}^{n}/\mu_{r},x'}$ for any point $x$ in the cyclic quotient center $S$, so by Artin approximation we have a scheme (variety) $W$ and a common $\acute{e}$tale diagram
 $$\begin{tikzcd}[column sep=0.5em]
 & && W\arrow[dl]\arrow[dr] &&  \\
&& X&& \mathds{A}^{n}/\mu_{r}
\end{tikzcd}$$
we define $U:=\mathds{A}^{n}\times_{\mathds{A}^{n}/\mu_{r}}W$ which is a smooth variety with a natural $\mu_{r}$ action, and it's not difficult to see $U/\mu_{r}=W$, especially the normal sheaf of fixed locus inherent the weight of group action on $\mathds{A}^{n}/\mu_{r}$. The canonical stack lift gives an $\acute{e}$tale morphism
$$\beta: [U/\mu_{r}]\longrightarrow\widetilde{X}$$
Let's consider Vistoli's $\acute{e}$tale neighborhood near a quotient singular point $x$, saying $\alpha:U/\mu_{r}\rightarrow W'\subseteq X$, since $\acute{e}$tale is locally quasi-finite, by replacing $U$ with its Zariski neighborhood of a point in the fiber, they are just $m$ points $t_{1},...,t_{m}$ in $T:=\alpha^{-1}(x)$.\\

The following result is a global to local generalization of the Proposition \ref{prop:pullbackcenterglobal}.
\begin{proposition}\label{prop:pullbackcenterlocal}
If we Kawamata blow up $X$ along the cyclic quotient center $S$ and consider its stacky lift picture:
$$\begin{tikzcd}[ampersand replacement=\&]
 \&\sqrt[r]{E_{\widetilde{Y}}/\widetilde{Y}}\arrow[d,"\pi"] \arrow[r,"\varpi"]\&\widetilde{Y}:=\widetilde{w_{1/r}Bl_{S}X}\\
 \&  {\widetilde{X}}
\end{tikzcd}$$
where $E_{\widetilde{Y}}$ is the stacky lift of the exceptional divisor of Kawamata blow-up. Then we have for any generalized geometrical point $\mathrm{k}_{x}\otimes\nu^{i}$ on $[S/\mu_{r}]$, $\overline{\mathcal{D}}_{x,i}:=\varpi_{*}\pi^{*}(\mathrm{k}_{x}\otimes\nu^{i})$ has the following finite Postnikov filtration:
\[\begin{adjustbox}{max width=0.98\linewidth}\(\begin{tikzcd}[column sep=0.5em]
 & 0 \arrow{rr}&& E_{-i-kr}\arrow{dl} \\
&& \iota_{x*}\mathcal{D}_{i+kr}\arrow[ul,dashed,"\Delta"]
\end{tikzcd}......\begin{tikzcd}[column sep=0.5em]
 &   E_{-i-r}\arrow{rr}&&  E_{-i}\arrow{rr}\arrow{dl}&& \overline{\mathcal{D}}_{x,i}\arrow{dl} \\
 && \iota_{x*}\mathcal{D}_{i+r}\arrow[ul,dashed,"\Delta"]  && \iota_{x*}\mathcal{D}_{i}\arrow[ul,dashed,"\Delta"]
 \end{tikzcd}\)\end{adjustbox}\]
and $\mathcal{D}_{x,i}$ is the $i$-th dual exceptional object on fiber over $x$. For $\sum_{j=1}^c a_j\le r$ and $0\le i<r$, the image is $\iota_{x*}\mathcal D_i$ if $i<\sum_j a_j$, and is zero otherwise.\\
\end{proposition}
\begin{proof}
All complexes and morphisms in the following argument are derived.
Fix a closed point $x\in S$ and put
$\mathcal P_x=f^{-1}(x)$. We first record why the calculation may be
carried out after completion at $x$. If $F,G\in\mathrm{D}^{b}(\widetilde Y)$
have cohomology supported on $\mathcal P_x$, then
\[
\operatorname{Ext}_{\widetilde Y}^q(F,G)
\longrightarrow
\operatorname{Ext}_{\widetilde Y_{\widehat x}}^q
 (F_{\widehat x},G_{\widehat x})
\]
is an isomorphism for every $q$. Here the subscript $\widehat x$ denotes
base change to $\operatorname{Spec}\widehat{\mathcal{O}}_{X,x}$.
Indeed $\widetilde Y$ is smooth, so $F$ is perfect. Properness of $f$
therefore makes the cohomology of
$f_*\mathcal RHom(F,G)$ coherent. It is supported at $x$, and hence is
of finite length over $\mathcal{O}_{X,x}$. Flat base change and completion
do not change these finite-length modules. The same argument applies
to any finite tower of such complexes. In particular, completion is
fully faithful on this subcategory and detects isomorphisms there.

Choose an $\acute e$tale quotient chart above $x$ and a point of that chart with
the same residue field. Its completed local ring agrees with
$\widehat{\mathcal{O}}_{X,x}$. The compatible weighted normal data identify
the completed root-stack and blow-up diagrams with the corresponding
completed diagram of Proposition~\ref{prop:pullbackcenterglobal}.
Flat base change gives the natural identification
\[
\bigl(\varpi_*\pi^*(\mathrm k_x\otimes\nu^i)\bigr)_{\widehat x}
\simeq
\varpi_{\widehat x*}\pi_{\widehat x}^*
 (\mathrm k_x\otimes\nu^i).
\]
Apply the finite Postnikov construction of that proposition on this
diagram. Its factors are the completions of the actual fibre objects
$\iota_{x*}\mathcal D_{i+kr}$ appearing in the statement; indices outside
$0\le i+kr<\sum_j a_j$ contribute zero.

This finite tower descends without choosing a free algebra on the
reduced fibre. Starting with its first factor, lift each attaching map
by the Ext isomorphism above and take its cone. Induction constructs
a tower on $\widetilde Y$ with exactly the stated factors. Every partial
object remains supported on $\mathcal P_x$. Its final object has the same
completion as $\varpi_*\pi^*(\mathrm k_x\otimes\nu^i)$, and that
isomorphism lifts by full faithfulness. The cone of the lifted morphism
has zero completion and is supported on $\mathcal P_x$, so faithful
flatness over $\mathcal{O}_{X,x}$ makes it zero. Thus the final object is
$\varpi_*\pi^*(\mathrm k_x\otimes\nu^i)$ itself.

All comparison maps of this prescribed tower descend by the same
Ext isomorphisms. This proves descent of the stated construction; it does
not assert uniqueness of arbitrary filtrations with the same factors.
For $\sum_j a_j\le r$, the single permitted factor gives
$\iota_{x*}\mathcal D_i$ when $0\le i<\sum_j a_j$, and the image is
zero when $\sum_j a_j\le i<r$.
\end{proof}

For the definition of a spanning class, we refer to \cite[Definition 1.47]{H}. Once spanning classes and admissible functors appear together, they have a powerful interaction. First, we need the following simple lemmas about the application of spanning class.

\begin{lemma}
If $\widetilde{X}$ is a smooth Delinge-Mumford stack, the set $$\{\mathrm{k}_{x}\otimes\chi\}_{x\in X(\mathrm{k}),\,\chi\in \mathrm{Irr.Rep}(G_{x})}$$ forms a spanning class for $\mathrm{D}^{b}(\widetilde{X})$.

\end{lemma}
\begin{proof}
By spectral sequence, for any element $F^{*}$ in $\mathrm{D}^{b}(\widetilde{X})$,
$$E^{p,q}_{2}:=\mathrm{Ext}^{p}(\mathcal{H}^{-q}(F^{*}),\mathrm{k}_{x}\otimes\chi) \Longrightarrow \mathrm{Ext}^{p+q}(F^{*}, \mathrm{k}_{x}\otimes\chi)$$
and Serre duality for embedding morphism of point $x$ which is proper, e.g. \cite[Corollary 2.10]{Nir}, it's only enough to verify the cases $F$ is a sheaf in one direction. So if $\mathrm{Hom}(F,\mathrm{k}_{x}\otimes\chi)=0$ for any $x,\chi$, then $F|_{x}=0$ for any $x$, hence $F\cong 0$.
\end{proof}

\begin{lemma}[Fully-faithfulness {\cite[Proposition 1.49]{H}}]\label{lemma:ff}
If $F: \mathcal{D}\longrightarrow \mathcal{D}'$ be an exact admissible functor between triangulated categories.
Suppose $\Omega$ is a spanning class of $\mathcal{D}$ such that for all objects $a, b$ in $\Omega$ and all the natural homomorphisms
$$\mathrm{Ext}_{\mathcal{D}}^{*}(a,b)=\mathrm{Ext}_{\mathcal{D}'}^{*}(F(a),F(b))$$
are bijective. Then $F$ is fully faithful.
\end{lemma}

\begin{lemma}[Orthogonality]\label{lemma:o}
If $F_{i}: \mathcal{D}_{i}\longrightarrow \mathcal{D}$ $i=1,2$ are two admissible functors between triangulated categories.
Suppose $\Omega_{i}$ is a spanning class of $\mathcal{D}_{i}$  for $i=1,2$, if
\begin{enumerate}

  \item  $F_{1}^{*}F_{2}(b)=0$ for any $b$ in $\Omega_{2}$,
  \item or $F_{2}^{!}F_{1}(a)=0$ for any $a$ in $\Omega_{1}$,
\end{enumerate}
then $\langle \mathrm{Im}\,F_{1}, \mathrm{Im}\,F_{2}\rangle$ is semi-orthogonal, furthermore all these conditions are equivalent.
\end{lemma}
\begin{proof}
  The orthogonality of $\mathrm{Im}\,F_{1}, \mathrm{Im}\,F_{2}$ is expected to show for any element $x_{i}$ in $\mathcal{D}_{i}$ $i=1,2$, we have
  $$\mathrm{Ext}^{*}(F_{2}(x_{2}),F_{1}(x_{1}))=0$$
  consider $x_{2}$ traversing all elements in $\Omega_{2}$, according to the definition of spanning class which is equivalent to $F_{2}^{!}F_{1}(x_{1})=0$ for any $x_{1}$.\\

  By the definition of spanning class, $F_{2}^{!}F_{1}(x_{1})=0$ if and only if for any $b$ in $\Omega_{2}$ we have
  $$\mathrm{Ext}^{*}(b,F_{2}^{!}F_{1}(x_{1}))=\mathrm{Ext}^{*}(F_{1}^{*}F_{2}(b),x_{1})=0$$
  consider $x_{1}$ traversing all elements in $\Omega_{1}$, according to the definition of spanning class this is equivalent to $F_{1}^{*}F_{2}(b)=0$  for any $b$ in $\Omega_{2}$.

\end{proof}

\begin{proposition}\label{prop:singularff} Under the assumption  $r\leq\sum_{i=1}^{c}a_{i}$, the functor
$$\mathrm{\Theta}: \mathrm{D}^{b}(\widetilde{X})\longrightarrow \mathrm{D}^{b}(\widetilde{Y})$$
$$(-)\longmapsto\varpi_{*}\pi^{*}(-)  $$
is admissible and fully-faithful.
\end{proposition}

\begin{proof}

The functor is admissible since $\pi^{*}$ and $\varpi_{*}$ are all admissible by \cite[Theorem 2.22]{Nir}.\\

For the fully-faithfulness, we are expected by Lemma \ref{lemma:ff} to show
$$\mathrm{Ext}^{*}(\mathrm{k}_{x_{i}}\otimes\nu^{a},\mathrm{k}_{x_{j}}\otimes\nu^{b})=\mathrm{Ext}^{*}(\mathrm{\Theta}(\mathrm{k}_{x_{i}}\otimes\nu^{a}),\mathrm{\Theta}(\mathrm{k}_{x_{j}}\otimes\nu^{b}))$$
for any $x_{i}$, $x_{j}$ and  $\nu^{a}$ ($\nu^{b}$) belong in the representation of stabilizer of  $x_{i}$ ($x_{j}$). Since $\mathrm{\Theta}$ is admissible, it's enough for us to show
$$\mathrm{\Theta}^{!}\mathrm{\Theta}(\mathrm{k}_{x}\otimes\nu^{a})=\mathrm{k}_{x}\otimes\nu^{a}$$
for any $x$ contained in $\widetilde{X}$ and any $\nu^{a}$ belong in the representation of stabilizer of  $x$.\\

\begin{enumerate}
\item If $x$ is contained in $S$, we already know that $\mathrm{\Theta}(\mathrm{k}_{x}\otimes\nu^{a})=\overline{\mathcal{D}}_{x,a}$ and $\mathrm{\Theta}^{!}(-)=\pi_{*}\varpi^{!}(-)$, so
$$\mathrm{\Theta}^{!}\mathrm{\Theta}(\mathrm{k}_{x}\otimes\nu^{a})=\pi_{*}\varpi^{!}\overline{\mathcal{D}}_{x,a}$$

Firstly, we compute $\varpi^{*}\overline{\mathcal{D}}_{x,a}$ using a similar strategy in Lemma \ref{lem:pullbackskycraper} (\ref{num:pullbackskycraper2}).\\

We have a semi-orthogonal decomposition for its root structure,
\begin{equation}\label{euq:sodofexp}
\mathrm{D}^{b}(\sqrt[r]{E_{\widetilde{Y}}/\widetilde{Y}})\simeq \langle \mathrm{\Upsilon}_{r-1},...,\mathrm{\Upsilon}_{1},\mathrm{\Upsilon}:=\mathrm{Im}\,\rho^{*}\rangle
\end{equation}
where embedding functor of $\mathrm{\Upsilon}_{i}$ is defined by the following commutative diagram:
$$\begin{tikzcd}[ampersand replacement=\&]
\&{[E_{\widetilde{Y}}/\mu_{r}]}\arrow[r,"\iota"]\arrow[d,"\pi"] \& {\sqrt[r]{E_{\widetilde{Y}}/\widetilde{Y}}}\\
 \&{E_{\widetilde{Y}}}
\end{tikzcd}$$
$i_{\mathrm{\Upsilon}_{i}}(-):=\iota_{*}\pi^{*}(-)\otimes \mathcal{M}_{E_{\widetilde{Y}}}^{\otimes i}$, its left adjoint follows $i_{\mathrm{\Upsilon}_{i}}^{*}(-)=\pi_{*}\iota^{*}(\mathcal{M}_{E_{\widetilde{Y}}}^{\otimes -i}\otimes (-))$. \\

Hence $\varpi^{*}\mathcal{D}_{x,a}=\varpi^{*}\varpi_{*}\mathcal{D}_{x,a}$ has a Postnikov filtration as following:
$$
\begin{tikzcd}[column sep=0.5em]
 & \varpi^{*}\mathcal{D}_{x,a} \arrow{rr}&& E^{1}\arrow{dl}\arrow{rr}&& E^{2}  \arrow{dl} \\
&&\mathcal{D}_{x,a}(r-1)\otimes\nu[1]\arrow[ul,dashed,"\Delta"]&& \mathcal{D}_{x,a}(r-2)\otimes\nu^{2}[1]\arrow[ul,dashed,"\Delta"]
\end{tikzcd}...
$$
$$
...\begin{tikzcd}[column sep=0.5em]
 &  E^{r-2}\arrow{rr}&& \mathcal{D}_{x,a}\arrow{dl} \arrow{rr}&& 0\arrow{dl} \\
 && \mathcal{D}_{x,a}(1)\otimes\nu^{r-1}[1] \arrow[ul,dashed,"\Delta"]  &&  \mathcal{D}_{x,a}[1]\arrow[ul,dashed,"\Delta"]
\end{tikzcd}
$$
while we also know the expression of $\overline{\mathcal{D}}_{x,a}$ as in Proposition \ref{prop:pullbackcenterlocal}:
$$\begin{tikzcd}[column sep=0.5em]
 & 0 \arrow{rr}&& E_{-a-k'r}\arrow{dl} \\
&& \iota_{x*}\mathcal{D}_{a+k'r}\arrow[ul,dashed,"\Delta"]
\end{tikzcd}......\begin{tikzcd}[column sep=0.5em]
 &   E_{-a-r}\arrow{rr}&&  E_{-a}\arrow{rr}\arrow{dl}&& \overline{\mathcal{D}}_{x,a}\arrow{dl} \\
 && \iota_{x*}\mathcal{D}_{a+r}\arrow[ul,dashed,"\Delta"]  && \iota_{x*}\mathcal{D}_{a}\arrow[ul,dashed,"\Delta"]
 \end{tikzcd}$$
so by noticing $\mathcal{M}_{E_{\widetilde{Y}}}|_{E_{\widetilde{Y}}}\simeq\mathcal{O}_{E}(-1)\otimes\nu$ and $\varpi^{!}(-)\simeq \mathcal{M}_{E_{\widetilde{Y}}}^{\otimes r-1}\otimes\varpi^{*}(-)$ e.g. \cite[Proposition 3.4]{Nir}, so $\pi_{*}\varpi^{!}\overline{\mathcal{D}}_{x,a}$ has a filtration with factors as following for variables $k$ and $\lambda$:
Here $\operatorname{gr}^{\mathrm{root}}_j$ denotes the outgoing
factor of the chosen root-ideal quotient filtration, including the
displayed rotation shift. It is not an SOD projection.
\begin{align*}
\pi_{*}\big(\mathcal{M}_{E_{\widetilde{Y}}}^{\otimes r-1}
\otimes\operatorname{gr}^{\mathrm{root}}_{r-k}(\varpi^{*}\mathcal{D}_{x,a+\lambda r})\big)
&= \pi_{*}\big(\mathcal{M}_{E_{\widetilde{Y}}}^{\otimes r-1}
\otimes \mathcal{D}_{x,a+\lambda r}(k)\otimes\nu^{r-k}[1]\big) \\
&= \pi_{*}\big(\mathcal{D}_{x,a+\lambda r}(k-r+1)\otimes\nu^{-k-1}[1]\big)
\stepcounter{equation}\tag{\theequation}\label{equ:kxextension}
\end{align*}
where $0\leq k\leq r-1$ and $0\leq \lambda\leq k'$,
then by Corollary \ref{cor:dualvanishing} we have
$$\pi_{*}(\iota_{x*}\mathcal{D}_{a+\lambda r}(k-r+1)\otimes\nu^{-k-1}[1])=\mathrm{H}^{*}_{\mathcal{P}}(\mathcal{P},\mathcal{D}_{a+\lambda r}(k-r+1))\otimes\nu^{-k-1}[1]\\
= \mathrm{k}_{x}\otimes\nu^{a}[1]\cdot\delta_{a+k,(1-\lambda)r-1}$$
the range of  $k$, $\lambda$ and  non-vanishing condition $a+k=(1-\lambda)r-1\geq 0$ implies $\lambda=0$ and $k=r-1-a$. So (\ref{equ:kxextension}) has only one non-trivial factor and we have an expression as follows
$$\pi_{*}\varpi^{!}\overline{\mathcal{D}}_{x,a}={ \pi_{*}(\mathcal{D}_{x,a}(-a)\otimes\nu^{a-r}[1])}[-1]=\mathrm{k}_{x}\otimes\nu^{a}[1][-1]$$
by our left  Postnikov filtration under the action of $\mathrm{\Theta}^{!}$, that is
$$\mathrm{\Theta}^{!}\mathrm{\Theta}(\mathrm{k}_{x}\otimes\nu^{a})=\pi_{*}\varpi^{!}\overline{\mathcal{D}}_{x,a}=\mathrm{k}_{x}\otimes\nu^{a}$$

\item If $x$ is not contained in $S$, consider following diagram:
$$\begin{tikzcd}[ampersand replacement=\&]
\widehat{x} \arrow[r]\arrow[d,equal] \&\sqrt[r]{E_{\widetilde{Y}}/\widetilde{Y}}\arrow[d,"\pi"]\&{\varpi^{-1}(E_{\widetilde{Y}})}\arrow[l,"\iota"]\arrow[d,"\pi_{\varpi^{-1}(E_{\widetilde{Y}})}"] \\
x \arrow[r] \& \widetilde{X}\& S\arrow[l]
\end{tikzcd}$$
we can see easily that $\pi^{*}\mathrm{k}_{x}=\mathrm{k}_{\widehat{x}}$  using the torsion independent base. Hence we have
$$\mathrm{\Theta}(\nu^{a}_{x})=\varpi_{*}\pi^{*}\nu^{a}_{x}=\nu^{a}_{\widehat{x}}$$
by similar argument,
$$\mathrm{\Theta}^{!}\mathrm{\Theta}(\mathrm{k}_{x}\otimes\nu^{a})=\mathrm{k}_{x}\otimes\nu^{a}$$
which finishes our proof.
\end{enumerate}

The computations above identify the actual adjunction unit as well.
Write $P=\mathrm k_x\otimes\nu^a$ and let
$\eta_P:P\to\Theta^!\Theta(P)$ be that unit. Its target is isomorphic
to $P$. Since $\Theta(P)\ne0$, the triangle identity
\[
\varepsilon_{\Theta(P)}\circ\Theta(\eta_P)
=\operatorname{id}_{\Theta(P)}
\]
shows that $\eta_P$ is nonzero. As
$\operatorname{End}(P)=\mathrm k$, it is an isomorphism. Under adjunction,
the natural map on Ext between any two generalized points is composition
with this unit, and hence is an isomorphism. Lemma~\ref{lemma:ff} applies.
\end{proof}

\begin{proposition}  Under the assumption  $r<\sum_{i=1}^{c}a_{i}$, the functors
$$\mathrm{\Theta}_{i}: \mathrm{D}^{b}(S)\longrightarrow \mathrm{D}^{b}(\widetilde{Y})$$
$$(-)\longmapsto \iota_{*}\pi^{*}(-)\otimes\mathcal{O}_{E_{\widetilde{Y}}}(i)$$
fit in the following commutative diagram:
$$\begin{tikzcd}[ampersand replacement=\&]
\&{E_{\widetilde{Y}}}\arrow[r,"\iota"]\arrow[d,"\pi"] \& \widetilde{Y} \\
 \&{S}
\end{tikzcd}$$
are admissible and fully-faithful, for all integer $i$.
\end{proposition}
\begin{proof}
The functors are admissible since $\pi^{*}$, $\iota{*}$ and $\mathcal{O}_{E_{\widetilde{Y}}}(i)\otimes(-)$ are all admissible by \cite[Theorem 2.22]{Nir}.\\

By Lemma \ref{lemma:ff} it is enough for us to check,
$$\mathrm{Ext}^{*}(\mathrm{k}_{x_{i}},\mathrm{k}_{x_{j}})=\mathrm{Ext}^{*}(\mathrm{\Theta}_{i}(\mathrm{k}_{x_{i}}),\mathrm{\Theta}_{i}(\mathrm{k}_{x_{j}}))$$
for any $x_{i}$, $x_{j}$ on $S$, since $\mathcal{O}_{E_{\widetilde{Y}}}(-E_{\widetilde{Y}})\simeq\mathcal{O}_{E_{\widetilde{Y}}}(r)$ by excess distinguished triangle,
$$\mathcal{O}_{E_{\widetilde{Y}}}(r)\otimes(-)[1]\longrightarrow \iota^{*}\iota_{*}(-)\longrightarrow (-)\longrightarrow\mathcal{O}_{E_{\widetilde{Y}}}(r)\otimes(-)[2]$$
we have,
\begin{align*}
  \mathrm{Ext}^{*}(\mathrm{\Theta}_{i}(\mathrm{k}_{x_{i}}),\mathrm{\Theta}_{i}(\mathrm{k}_{x_{j}})) &= \mathrm{Ext}^{*}(\iota_{*}\big(\pi^{*}(\mathrm{k}_{x_{i}})\otimes \mathcal{O}_{E_{\widetilde{Y}}}(i)\big),\iota_{*}\big(\pi^{*}(\mathrm{k}_{x_{j}})\otimes \mathcal{O}_{E_{\widetilde{Y}}}(i)\big)) \\
   &= \mathrm{Ext}^{*}(\iota^{*}\iota_{*}\big(\pi^{*}(\mathrm{k}_{x_{i}})\otimes \mathcal{O}_{E_{\widetilde{Y}}}(i)\big),\big(\pi^{*}(\mathrm{k}_{x_{j}})\otimes \mathcal{O}_{E_{\widetilde{Y}}}(i)\big))
\end{align*}
fits in long exact sequence,
$$\mathrm{Ext}^{*}(\pi^{*}(\mathrm{k}_{x_{i}})\otimes \mathcal{O}_{E_{\widetilde{Y}}}(i+r)[1],\pi^{*}(\mathrm{k}_{x_{j}})\otimes \mathcal{O}_{E_{\widetilde{Y}}}(i))\longrightarrow\mathrm{Ext}^{*}(\mathrm{\Theta}_{i}(\mathrm{k}_{x_{i}}),\mathrm{\Theta}_{i}(\mathrm{k}_{x_{j}}))$$
$$\longrightarrow\mathrm{Ext}^{*}(\pi^{*}(\mathrm{k}_{x_{i}})\otimes \mathcal{O}_{E_{\widetilde{Y}}}(i),\pi^{*}(\mathrm{k}_{x_{j}})\otimes \mathcal{O}_{E_{\widetilde{Y}}}(i))\longrightarrow...$$
since we assume $r<\sum_{i=1}^{c}a_{i}$, we have $\mathrm{H}^{*}(\mathcal{P},\mathcal{O}_{\mathcal{P}}(-r))=0$ by the relative version Beilinson type theorem for weighted projective bundle over $S$, e.g. Lemma \ref{lemma:wpscoho}, the first term of above long exact sequence vanishes. While the third term is
$$\mathrm{Ext}^{*}(\pi^{*}(\mathrm{k}_{x_{i}})\otimes \mathcal{O}_{E_{\widetilde{Y}}}(i),\pi^{*}(\mathrm{k}_{x_{j}})\otimes \mathcal{O}_{E_{\widetilde{Y}}}(i))=\mathrm{Ext}^{*}(\mathrm{k}_{x_{i}},\mathrm{k}_{x_{j}})$$ we can see our proposition holds.

\end{proof}

\begin{proposition} \label{prop:a-r}Under the assumption  $r<\sum_{i=1}^{c}a_{i}$,  the functors
$$\mathrm{\Theta}_{-(\sum_{i=1}^{c}a_{i}-r)},...,\mathrm{\Theta}_{-1},\mathrm{\Theta} $$  form a well-defined semi-orthogonal decomposition:
$$\mathcal{D}:=\langle\mathrm{Im}\,\Theta_{-(\sum_{i=1}^{c}a_{i}-r)},...,\mathrm{Im}\,\Theta_{-1},\mathrm{Im}\, \mathrm{\Theta}\rangle$$
and they generate the whole derived category, that is $\mathcal{D}\simeq\mathrm{D}^{b}(\widetilde{Y})$.
\end{proposition}

\begin{proof}
Firstly, we prove  $\mathrm{Im}\,\mathrm{\Theta}_{-i}$ and $\mathrm{Im}\,\mathrm{\Theta} $  are semi-orthogonal, and it can be reduced by Lemma \ref{lemma:o} for us to prove for any point $x$ in $S$,
$$\mathrm{\Theta}^{!}\mathrm{\Theta}_{-i}(\mathrm{k}_{x})=0$$
since $\mathrm{\Theta}$ is admissible. Recall
$$\mathrm{\Theta}^{!}\mathrm{\Theta}_{-i}(\mathrm{k}_{x})=\pi_{*}\varpi^{!}(\iota_{*}\pi^{*}\mathrm{k}_{x}\otimes\mathcal{O}_{E_{\widetilde{Y}}}(-i))$$
by a similar argument, $\varpi^{*}\iota_{*}\pi^{*}\mathrm{k}_{x}\otimes\mathcal{O}_{E_{\widetilde{Y}}}(-i)=\varpi^{*}\varpi_{*}\iota_{*}\pi^{*}\mathrm{k}_{x}\otimes\mathcal{O}_{E_{\widetilde{Y}}}(-i)$ has a left  Postnikov filtration as following:
$$\operatorname{gr}^{\mathrm{root}}_{k}(\varpi^{*}\iota_{*}\pi^{*}\mathrm{k}_{x}\otimes\mathcal{O}_{E_{\widetilde{Y}}}(-i))=\iota_{*}\pi^{*}\mathrm{k}_{x}\otimes\mathcal{O}_{E_{\widetilde{Y}}}(-i+r-k)\otimes\nu^{k}[1]$$
for $1\leq k\leq r$, noticing $\omega_{\varpi}=\mathcal{M}_{E_{\widetilde{Y}}}^{r-1}$, $\mathcal{M}_{E_{\widetilde{Y}}}^{r-1}|_{\mathcal{P}}=\mathcal{O}_{\mathcal{P}}(-r+1)\otimes\nu^{r-1}$ and $\varpi^{!}(-)\simeq \mathcal{M}_{E_{\widetilde{Y}}}^{\otimes r-1}\otimes\varpi^{*}(-)$ e.g. \cite[Proposition 3.4]{Nir}
\begin{align*}
  \omega_{\varpi}\otimes\operatorname{gr}^{\mathrm{root}}_{k}(\varpi^{*}\iota_{*}\pi^{*}\mathrm{k}_{x}\otimes\mathcal{O}_{E_{\widetilde{Y}}}(-i)) &= \iota_{*}\pi^{*}\mathrm{k}_{x}\otimes\mathcal{O}_{E_{\widetilde{Y}}}(-i-k+1)\otimes\nu^{k+r-1}[1] \\
   &=\mathrm{H}_{\mathcal{P}}^{*}(\mathcal{P},\mathcal{O}_{E_{\widetilde{Y}}}(-i-k+1))\otimes\nu^{k+r-1}[1]
\end{align*}
considering the ranges of all variables $1\leq k\leq r$ and $0<i\leq\sum_{i=1}^{c}a_{i}-r$, we have
$$-\sum_{i=1}^{c}a_{i}<-i-k+1<0$$
so all the cohomology vanish by Lemma \ref{lemma:wpscoho}, we get $\mathrm{\Theta}^{!}\mathrm{\Theta}_{-i}(\mathrm{k}_{x})=0$.\\

Secondly, we prove  $\mathrm{Im}\,\mathrm{\Theta}_{-i}$ and $\mathrm{Im}\,\mathrm{\Theta}_{-j} $  are semi-orthogonal if $j>i$ in this range. By Lemma \ref{lemma:o} it can be reduced to for us to prove for any point $x$ in $S$,
$$\mathrm{\Theta}_{-i}^{!}\mathrm{\Theta}_{-j}(\mathrm{k}_{x})=0$$
recall,
\begin{align*}
{\mathrm{\Theta}_{-i}^{!}\mathrm{\Theta}_{-j}(\mathrm{k}_{x})}&={ \pi_{*}\iota^{!}(\iota_{*}\pi^{*}\mathrm{k}_{x}\otimes\mathcal{O}_{E_{\widetilde{Y}}}(i-j))} \\
&={ \pi_{*}\iota^{*}(\iota_{*}\pi^{*}\mathrm{k}_{x}\otimes\mathcal{O}_{E_{\widetilde{Y}}}(i-j-r))}[-1]
\end{align*}
and notice the excess distinguished triangle,
$$\pi^{*}\mathrm{k}_{x}\otimes\mathcal{O}_{E_{\widetilde{Y}}}(i-j)[1]\longrightarrow \iota^{*}(\iota_{*}\pi^{*}\mathrm{k}_{x}\otimes\mathcal{O}_{E_{\widetilde{Y}}}(i-j-r))\longrightarrow\pi^{*}\mathrm{k}_{x}\otimes\mathcal{O}_{E_{\widetilde{Y}}}(i-j-r)$$
pushforward along $\pi$,
$$\mathrm{H}^{*}_{\mathcal{P}}(\mathcal{P},\mathcal{O}_{\mathcal{P}}(i-j))[1]\longrightarrow \pi_{*}\iota^{*}(\iota_{*}\pi^{*}\mathrm{k}_{x}\otimes\mathcal{O}_{E_{\widetilde{Y}}}(i-j-r))\longrightarrow\mathrm{H}^{*}_{\mathcal{P}}(\mathcal{P},\mathcal{O}_{\mathcal{P}}(i-j-r))$$
considering the ranges of all variables, we have $-\sum_{i=1}^{c}a_{i}<i-j<0$ and $-\sum_{i=1}^{c}a_{i}<i-j-r<0$, so all the cohomology vanish, we get $\mathrm{\Theta}_{-i}^{!}\mathrm{\Theta}_{-j}(\mathrm{k}_{x})=0$.\\

Finally, we prove they generate the whole category. By combinatorial property, it's not
difficult to see  $\mathcal{D}$ pullback to $\mathcal{P}$ is essential surjective, that means: \\

Firstly, by Corollary \ref{cor:okfmtest} $\mathcal{O}_{\mathcal{P}}(-1)$ under Fourier-Mukai transformation respect to complex in Theorem \ref{thm:beilinson} i.e. $p_{2*}(p_{1}^{*}\mathcal{O}_{\mathcal{P}}(-1)\otimes \mathbf{e})$ is $\mathcal{D}_{\sum_{i} a_{i}-1}$ up to certain shift, so $\iota_{x*}\mathcal{O}_{\mathcal{P}}(-1)$ is also isomorphic to $\iota_{x*}\mathcal{D}_{\sum_{i} a_{i}-1}$ up to certain shift. Then we test $\mathcal{O}_{\mathcal{P}}(-2)$ similarly, which is a complex
$$V\otimes\mathcal{D}_{\sum_{i} a_{i}-1}\rightarrow\mathcal{D}_{\sum_{i} a_{i}-2}[1]$$
up to certain shift, and $\iota_{x*}\mathcal{O}_{\mathcal{P}}(-2)$ is isomorphic to
$$\iota_{x*}V\otimes\mathcal{D}_{\sum_{i} a_{i}-1}\rightarrow\iota_{x*}\mathcal{D}_{\sum_{i} a_{i}-2}[1]$$
continue to repeat this test for $\mathcal{O}_{\mathcal{P}}(-k)$ by Corollary \ref{cor:okfmtest}, where $1\leq k\leq\sum_{i=1}^{c}a_{i}-r$, we can see
$$\iota_{x*}\mathcal{D}_{r},...,\iota_{x*}\mathcal{D}_{\sum_{i} a_{i}-1}$$
are all contained in $\langle\mathrm{Im}\,\Theta_{-(\sum_{i=1}^{c}a_{i}-r)},...,\mathrm{Im}\,\Theta_{-1}\rangle$ for any $x$, while combine the result in Proposition \ref{prop:pullbackcenterlocal}, for any $0\leq i<r$, $\varpi_{*}\pi^{*}(\mathrm{k}_{x}\otimes\nu^{i})$ has the following finite Postnikov filtration:
\[\begin{adjustbox}{max width=0.98\linewidth}\(\begin{tikzcd}[column sep=0.5em]
 & 0 \arrow{rr}&& E_{-i-kr}\arrow{dl} \\
&& \iota_{x*}\mathcal{D}_{i+kr}\arrow[ul,dashed,"\Delta"]
\end{tikzcd}......\begin{tikzcd}[column sep=0.5em]
 &   E_{-i-r}\arrow{rr}&&  E_{-i}\arrow{rr}\arrow{dl}&& \varpi_{*}\pi^{*}(\mathrm{k}_{x}\otimes\nu^{i})\arrow{dl} \\
 && \iota_{x*}\mathcal{D}_{i+r}\arrow[ul,dashed,"\Delta"]  && \iota_{x*}\mathcal{D}_{i}\arrow[ul,dashed,"\Delta"]
 \end{tikzcd}\)\end{adjustbox}\]
so
$$\iota_{x*}\mathcal{D}_{0},...,\iota_{x*}\mathcal{D}_{r-1}$$
are also contained in $$\mathcal{D}:=\langle\mathrm{Im}\,\Theta_{-(\sum_{i=1}^{c}a_{i}-r)},...,\mathrm{Im}\,\Theta_{-1},\mathrm{Im}\, \mathrm{\Theta}\rangle$$
 for any $x$.\\

Summary any generalized skyscrapersheaf concentrates at exceptional divisor of Kawamata blow-up is contained in $\mathcal{D}$.
Similarly, any generalized skyscrapersheaf concentrates outside exceptional divisor is contained in $\mathrm{Im}\,\mathrm{\Theta}$. So a spanning class of $\mathrm{D}^{b}(\widetilde{Y})$ is contained in $\mathcal{D}$.\\

On the other side, note that $\mathcal{D}$ is a admissible subcategory of $\mathrm{D}^{b}(\widetilde{Y})$ since $\Theta$ and $\Theta_{-i}$ are all admissible. If $G$ is contained in $\mathcal{D}^{\perp}$, we have
$$\mathrm{Ext}_{\widetilde{Y}}^{*}(a,G)=0$$ for a spanning class of $\mathrm{D}^{b}(\widetilde{Y})$, hence $G\simeq0$ and we finish the proof.

\end{proof}

\begin{proposition}\label{prop:ffofemeb}  Under the assumption  $\sum_{i=1}^{c}a_{i}\leq r$, the functor
$$\mathrm{\Phi}: \mathrm{D}^{b}(\widetilde{Y})\longrightarrow \mathrm{D}^{b}(\widetilde{X})$$
$$(-)\longmapsto\pi_{*}\varpi^{*}(-)  $$
is admissible and fully-faithful.
\end{proposition}
\begin{proof}
Put $e=\sum_i a_i$. For each test object
$W=\iota_{x*}\mathcal D_a(1-r)$ we have
$\operatorname{End}_{\widetilde Y}(W)=\mathrm k$.
This can be checked on the completed local model, using the Ext
comparison in Proposition~\ref{prop:pullbackcenterlocal}. That model is
the completion of the zero section in
$\operatorname{Tot}_{\mathcal P}\mathcal{O}_{\mathcal P}(-r)$, with the
smooth centre parameters adjoined. On the uncompleted total space let
\[
T=\bigoplus_{j=0}^{e-1}\mathcal{O}(j).
\]
For $r\ge e$ and $n\ge0$,
$-e<j-i+nr$ for all $0\le i,j<e$. Thus
$\mathrm H^q(\mathcal P,\mathcal{O}(j-i+nr))=0$ for $q>0$.
The bundle $T$ consequently has no higher self-Ext. It generates because
the projection to $\mathcal P$ is affine and
$\mathcal{O},\ldots,\mathcal{O}(e-1)$ generate on $\mathcal P$.
Under the resulting tilting equivalence,
$\iota_{x*}\mathcal D_a$ corresponds to the one-dimensional simple
module at the $a$-th vertex: this is exactly
Corollary~\ref{cor:dualvanishing}. Its endomorphism space is therefore
$\mathrm k$. A common line-bundle twist gives the same result for $W$.

The functor is admissible since $\pi_{*}$ and $\varpi^{*}$ are all admissible by \cite[Theorem 2.22]{Nir}.\\

The idea of proof is the same as above, by Lemma \ref{lemma:ff} we only need to prove that
$$\mathrm{\Phi}^{*}\mathrm{\Phi}(\iota_{x*}\mathcal{D}_{a}(-r+1))=\iota_{x*}\mathcal{D}_{a}(-r+1)$$
holds for any $a$ such that $0\leq a<\sum_{i}a_{i} $ and $\iota_{x*}\mathcal{D}_{a}(-r+1)$ supporting on fiber of any point $x$ on $S$. Recall
$$\mathrm{\Phi}^{*}\mathrm{\Phi}(\iota_{x*}\mathcal{D}_{a}(-r+1))=\varpi_{!}\pi^{*}\pi_{*}\varpi^{*}(\iota_{x*}\mathcal{D}_{a}(-r+1))$$
and by semi-orthogonal decomposition (\ref{euq:sodofexp}), $\varpi^{*}(\iota_{x*}\mathcal{D}_{a}(-r+1))$ has a left  Postnikov filtration as following:
$$\operatorname{gr}^{\mathrm{root}}_{k}(\varpi^{*}(\iota_{x*}\mathcal{D}_{a}(-r+1)))=\iota_{x*}\mathcal{D}_{a}(-k+1)\otimes\nu^{k}[1]$$
for $1\leq k\leq r$.\\

Then pushforward it along $\pi_{*}$, and separate it into two parts
$$\pi_{*}\operatorname{gr}^{\mathrm{root}}_{k}\varpi^{*}(\iota_{x*}\mathcal{D}_{a}(-r+1))=\pi_{*}\iota_{x*}\mathcal{D}_{a}(-k+1)\otimes\nu^{k}[1]$$
\begin{enumerate}
  \item\label{enuks2} $0\leq k-1\leq \sum_{i}a_{i}-1$:
\begin{align*}
  \pi_{*}\iota_{x*}\mathcal{D}_{a}(-k+1)\otimes\nu^{k}[1]&=\mathrm{H}^{*}(\mathcal{P},\mathcal{D}_{a}(-k+1))\otimes\nu^{k}[1]\\
 &=\mathrm{Ext}_{\mathcal{P}}^{*}(\mathcal{O}_{\mathcal{P}}(k-1),\mathcal{D}_{a})\otimes\nu^{k}[1]\cdot\delta_{k-1,a}\\
 &=\mathrm{Ext}_{\mathcal{P}}^{*}(\mathcal{O}_{\mathcal{P}}(a),\mathcal{D}_{a})\otimes\nu^{a+1}[1]\cdot\delta_{k-1,a}=\mathrm{k}_{x}\otimes \nu^{a+1}[1]\cdot\delta_{k-1,a}
\end{align*}
by Corollary \ref{cor:dualvanishing}.\\
  \item $\sum_{i}a_{i}\leq k-1\leq r-1$:
$$\pi_{*}\iota_{x*}\mathcal{D}_{a}(-k+1)\otimes\nu^{k}[1]=F^{*}\otimes\nu^{k} $$
where $F^{*}$ is certain complex with trivial representation supporting on $x$.
\end{enumerate}

We pullback this two components separately along $\pi$ and then Verdier pushforward along $\varpi$, the second part vanishes since we have a filtration for  $\pi^{*}(\mathrm{k}_{x}\otimes\nu^{k})$ as Proposition \ref{prop:pullbackcenterlocal},
\[\begin{adjustbox}{max width=0.98\linewidth}\(\begin{tikzcd}[column sep=0.5em]
 & 0 \arrow{rr}&& E_{-t+1}\arrow{dl} \\
&& \iota_{x*}\mathcal{D}_{(t-1)}\otimes\nu^{k-(t-1)}\arrow[ul,dashed,"\Delta"]&&
\end{tikzcd}
...\begin{tikzcd}[column sep=0.5em]
 &   E_{-2}\arrow{rr}&&  E_{-1}\arrow{rr}\arrow{dl}&& \pi^{*} \mathrm{k}_{x}\otimes\nu^{k}\arrow{dl} \\
 && \iota_{x*}\mathcal{D}_{1}\otimes\nu^{k-1}\arrow[ul,dashed,"\Delta"]  && \iota_{x*}\mathcal{D}_{0}\otimes\nu^{k}\arrow[ul,dashed,"\Delta"]
\end{tikzcd}\)\end{adjustbox}\]
then tensor it with $\mathcal{M}_{E_{\widetilde{Y}}}^{\otimes(r-1)}$ and pushforward along $\varpi$, noticing $\mathcal{M}_{E_{\widetilde{Y}}}^{\otimes(r-1)}|_{E_{\widetilde{Y}}}\simeq\mathcal{O}_{E_{\widetilde{Y}}}(-r+1)\otimes\nu^{r-1}$ and the range of variables $k$ in the second part, we have all presentations above are among
$$\nu,...,\nu^{r-1}$$
pushforward along $\varpi$ cancels all nontrivial presentations, hence
$$\varpi_{!}\pi^{*}(\mathrm{k}_{x}\otimes\nu^{k})=0$$
if $\sum_{i}a_{i}\leq k-1\leq r-1$. So only $k=a+1$ term in part (\ref{enuks2}) is left and
$$\varpi_{!}\pi^{*}\pi_{*}\varpi^{*}(\iota_{x*}\mathcal{D}_{a}(-r+1))=\varpi_{!}\pi^{*}\mathrm{k}_{x}\otimes\nu^{a+1}$$
similarly noticing $\pi^{*}\mathrm{k}_{x}\otimes\nu^{a+1}$ has a filtration by the argument in Proposition \ref{prop:pullbackcenterlocal}
\[\begin{adjustbox}{max width=0.98\linewidth}\(\begin{tikzcd}[column sep=0.5em]
 & 0 \arrow{rr}&& E_{-t}\arrow{dl}\arrow{rr}&& E_{-t+1}  \arrow{dl} \\
&& \iota_{x*}\mathcal{D}_{t}\otimes \nu^{a+1-t}\arrow[ul,dashed,"\Delta"]&& \iota_{x*}\mathcal{D}_{t-1}\otimes \nu^{a+1-(t-1)}\arrow[ul,dashed,"\Delta"]
\end{tikzcd}...\)\end{adjustbox}\]

\begin{equation}\label{equ:syzoffir}\begin{tikzcd}[column sep=0.5em]
 &   E_{-2}\arrow{rr}&&  E_{-1}\arrow{rr}\arrow{dl}&& \pi^{*} (\mathrm{k}_{x}\otimes \nu^{a+1})\arrow{dl} \\
 && \iota_{x*}\mathcal{D}_{1}\otimes \nu^{a}\arrow[ul,dashed,"\Delta"]  && \iota_{x*}\mathcal{D}_{0}\otimes \nu^{a+1}\arrow[ul,dashed,"\Delta"]
\end{tikzcd}\end{equation}
where $t:=\sum a_{i}-1$. Finally we tensor (\ref{equ:syzoffir}) with $\mathcal{M}_{E_{\widetilde{Y}}}^{\otimes{r-1}}$ and pushforward along $\varpi$, which cancels all terms with non-trivial representation, hence we have
$$\mathrm{\Phi}^{*}\mathrm{\Phi}(\iota_{x*}\mathcal{D}_{a}(-r+1))=\iota_{x*}\mathcal{D}_{a}(-r+1)$$
and we get the result.

Let $\epsilon_W:\Phi^*\Phi(W)\to W$ be the counit for
$\Phi^*\dashv\Phi$. The calculation identifies its source with $W$.
It also implies $\Phi(W)\ne0$. The triangle identity
\[
\Phi(\epsilon_W)\circ\eta_{\Phi(W)}
=\operatorname{id}_{\Phi(W)}
\]
then shows that $\epsilon_W$ is nonzero. Since
$\operatorname{End}(W)=\mathrm k$, it is an isomorphism. On generalized
points away from the exceptional divisor the unit and counit are the
identity under the geometric isomorphism. For any two members $W,W'$
of the stated spanning class, adjunction identifies the natural map on
Ext with composition by $\epsilon_W$, so it is an isomorphism.
Lemma~\ref{lemma:ff} proves full faithfulness.
\end{proof}

\begin{proposition}  Under the assumption  $\sum_{i=1}^{c}a_{i}< r$, the functors
$$\mathrm{\Phi}_{i}: \mathrm{D}^{b}(S)\longrightarrow \mathrm{D}^{b}(\widetilde{X})$$
$$(-)\longmapsto \varsigma_{*}(\tau^{*}(-)\otimes \nu^{i})$$
fit in diagram
$$\begin{tikzcd}[ampersand replacement=\&]
\&{[S/\mu_{r}]}\arrow[r,"\varsigma"]\arrow[d,"\tau"] \& \widetilde{X} \\
 \&{S}
\end{tikzcd}$$
are admissible and fully-faithful.
\end{proposition}

\begin{proof}
The functor are admissible since $\tau^{*}$,  $\varsigma_{*}$ and $\nu^{i}\otimes(-)$ are all admissible by \cite[Theorem 2.22]{Nir}.\\

By Lemma \ref{lemma:ff} we only need to prove
$$\mathrm{\Phi}_{i}^{*}\mathrm{\Phi}_{i}(\mathrm{k}_{x})=\mathrm{k}_{x}$$
for any geometrical point $x$ on $S$, since $\mathrm{\Phi}_{i}$ is admissible. Recall
$$\mathrm{\Phi}_{i}^{*}\mathrm{\Phi}_{i}(\mathrm{k}_{x})= \tau_{*}(\varsigma^{*}\varsigma_{*}(\tau^{*}(\mathrm{k}_{x})\otimes \nu^{i})\otimes\nu^{-i})$$
noticing $\varsigma^{*}\varsigma_{*}(\tau^{*}(\mathrm{k}_{x})\otimes \nu^{i})$ has a heart filtration as following:
$$\mathcal{H}^{-k}(\varsigma^{*}\varsigma_{*}(\tau^{*}(\mathrm{k}_{x})\otimes \nu^{i}))=\bigwedge^{k}\mathcal{N}_{[S/\mu_{r}]|\widetilde{X}}^{\vee}\otimes(\tau^{*}(\mathrm{k}_{x})\otimes \nu^{i})$$
for any $k$. Hence:
$$\nu^{-i}\otimes\mathcal{H}^{-k}(\varsigma^{*}\varsigma_{*}(\tau^{*}(\mathrm{k}_{x})\otimes \nu^{i}))=\bigwedge^{k}\mathcal{N}_{[S/\mu_{r}]|\widetilde{X}}^{\vee}\big |_{x}$$
by assumption $S$ is contained in $X$ as a $\mu_{r}$ cyclic quotient center with weight $a_{1},...,a_{c}$, which $i$-th coordinate contributes to $\mathcal{N}_{[S/\mu_{r}]|\widetilde{X}}^{\vee}\big |_{x}$ a presentation $\nu^{-a_{i}}$.\\

If $k>0$ and the representation term in $\bigwedge^{k}\mathcal{N}_{[S/\mu_{r}]|\widetilde{X}}^{\vee}\big |_{x}$ is $\nu^{-j_{k}}$ for $1\leq k\leq c$, then $0<j_{k}\leq\sum_{1}^{c}a_{i}$ and by our assumption $\sum_{1}^{c}a_{i}<r$ they are all non-trivial, so
$$\tau_{*}\varsigma^{*}\varsigma_{*}(\tau^{*}(\mathrm{k}_{x})\otimes \nu^{i})\otimes\nu^{-i}=\nu^{-i}\otimes\mathcal{H}^{0}(\varsigma^{*}\varsigma_{*}(\tau^{*}(\mathrm{k}_{x})\otimes \nu^{i}))=\mathrm{k}_{x}$$
since $\tau$ is a gerbe, we have
$$\mathrm{\Phi}_{i}^{*}\mathrm{\Phi}_{i}(\mathrm{k}_{x})=\mathrm{k}_{x}$$

\end{proof}

\begin{proposition}\label{prop:r-a}
Assume $e:=\sum_{i=1}^c a_i\le r$, with the compatible global character data fixed above. Then
\[
\mathrm{D}^{b}(\widetilde X)
=\langle\operatorname{Im}\Phi_{r-e},\ldots,
\operatorname{Im}\Phi_1,\operatorname{Im}\Phi\rangle,
\qquad \Phi=\pi_*\varpi^*.
\]
For $e=r$ the exceptional range is empty.
\end{proposition}
\begin{proof}
The full faithfulness of $\Phi$ and, when $e<r$, of the $\Phi_i$
was established above. For semiorthogonality,
$\Phi^!=\varpi_*\pi^!$ and
\[
\omega_\pi=\mathcal M_E^{e-1},\qquad
\omega_\pi|_{[\mathcal P/\mu_r]}
=\mathcal{O}_{\mathcal P}(1-e)\otimes\nu^{e-1}.
\]
The filtration of $\pi^*(\mathrm k_x\otimes\nu^i)$ has factors
$\iota_{x*}\mathcal D_k\otimes\nu^{i-k}$ for $0\le k<e$.
For $1\le i\le r-e$, all characters $i-k+e-1$ lie in $[1,r-1]$,
so $\Phi^!\Phi_i(\mathrm k_x)=0$.
For $1\le j<i\le r-e$, use the left adjoint $\Phi_i^*$:
its factors have characters $j-i-\sum_{h\in A}a_h$, for subsets
$A\subset\{1,\ldots,c\}$. Every such exponent lies in $[-r+1,-1]$.
Thus $\Phi_i^*\Phi_j(\mathrm k_x)=0$, giving the required direction of
semiorthogonality.

For generation, write $\mathcal D$ for the subcategory displayed above.
For $0\le a<e$, apply $\Phi$ to
$\iota_{x*}\mathcal D_a(1-e)$. The root-stack filtration has pushforward
factors
\[
\operatorname{Ext}^*_{\mathcal P}
  (\mathcal{O}_{\mathcal P}(\ell_k),\mathcal D_a)
\otimes\nu^k[1],\qquad
\ell_k=e+k-r-1,\quad1\le k\le r.
\]
If $k\le r-e$, then $\ell_k<0$, and every such factor is a complex of
copies of $\mathrm k_x\otimes\nu^k$, already in $\mathcal D$.
If $r-e+1\le k\le r$, then $0\le\ell_k<e$. Duality with the collection
$\mathcal{O},\ldots,\mathcal{O}(e-1)$ leaves exactly the factor
\[
\mathrm k_x\otimes\nu^{r+1+a-e}[1].
\]
The total filtered object and all other factors lie in $\mathcal D$,
so this last factor does as well. Varying $a$ supplies the characters
$r-e+1,\ldots,r$, while the exceptional components supply $1,\ldots,r-e$.
Hence all $r$ generalized point sheaves over $x$ lie in $\mathcal D$.
Away from the centre the functor is an isomorphism. Admissibility and the
spanning class of generalized points therefore imply
$\mathcal D=\mathrm{D}^{b}(\widetilde X)$.

\end{proof}

\begin{theorem}\label{thm:singularweightedblowup}
If $X$ is a variety with some codimension $c$ and $(a_{1},..,a_{c})/r$  type \footnote{\, Same as above \ref{foot:ai}, additionally condition $a_{i}<r$ is unnecessary.} cyclic quotient center $S$, we denotes $\widetilde{X}$ as its canonical stack, with the global split character and weighted normal data specified in Section~4.3, after weighted blowing up with the compatible weight at $S$, we get another variety (algebraic space) $Y$ with some cyclic quotient center, and denotes its canonical stack by $\widetilde{Y}$.\\
\begin{enumerate}
\item If $r<\sum_{i=1}^{c} a_{i}$:\\

There is a semi-orthogonal decomposition of $\widetilde{Y}$ with respect to $\widetilde{X}$ as following:
$$\mathrm{D}^{b}(\widetilde{Y})\cong\langle \mathrm{Im}\,\mathrm{\Theta}_{-\sum_{i}{a_{i}}+r},..., \mathrm{Im}\,\mathrm{\Theta}_{-1},\mathrm{Im}\, \mathrm{\Theta}\rangle$$
where $\mathrm{\Theta}$ and $\mathrm{\Theta}_{i}$ are defined by considering the following diagrams:
$$\begin{tikzcd}[ampersand replacement=\&]
E_{\widetilde{Y}} \arrow[r,"\iota"]\arrow[d,"f_{E_{\widetilde{Y}}}"] \&\widetilde{Y}\arrow[d,"f"]  \&\sqrt[ r]{E_{\widetilde{Y}}/\widetilde{Y}} \arrow[r,"\varpi"]\arrow[d,"\pi"] \&\widetilde{Y}\arrow[d,"f"] \\
S \arrow[r,"\kappa"] \&X  \&\widetilde{X} \arrow[r,"\theta"] \& X
\end{tikzcd}$$
we have $\mathrm{\Theta}(-):=\varpi_{*}\pi^{*}(-)$ and $\mathrm{\Theta}_{i}(-):=\iota_{*}\bigl(f_{E_{\widetilde{Y}}}^{*}(-)\otimes\mathcal{O}_{E_{\widetilde{Y}}}(i)\bigr)$
for any integer $i$.\\

\item If $r=\sum_{i=1}^{c} a_{i}$:\\

There is an equivalence induced by either $\mathrm{\Theta}$ or $\mathrm{\Phi}$:

$$\mathrm{D}^{b}(\widetilde{Y})\simeq \mathrm{D}^{b}(\widetilde{X})$$

\item If $r>\sum_{i=1}^{c} a_{i}$:\\

There is a semi-orthogonal decomposition of $\widetilde{X}$ with respect to $\widetilde{Y}$ as following:
$$\mathrm{D}^{b}(\widetilde{X})=\langle\mathrm{Im}\,\Phi_{r-\sum_{i=1}^{c}a_i},\ldots,\mathrm{Im}\,\Phi_{1},\mathrm{Im}\, \mathrm{\Phi}\rangle$$
where $\mathrm{\Phi}$ and $\mathrm{\Phi}_{i}$ are defined by considering the following diagrams:
$$\begin{tikzcd}[ampersand replacement=\&]
[S/\mu_{r}] \arrow[r,"\varsigma"]\arrow[d,"\tau"] \&\widetilde{X}  \&\sqrt[r]{E_{\widetilde{Y}}/\widetilde{Y}} \arrow[r,"\varpi"]\arrow[d,"\pi"] \&\widetilde{Y}\arrow[d,"f"] \\
S  \& \&\widetilde{X} \arrow[r,"\theta"] \& X
\end{tikzcd}$$
we have $\mathrm{\Phi}(-):=\pi_{*}\varpi^{*}(-)$ and $\mathrm{\Phi}_{i}(-):=\varsigma_{*}(\tau^{*}(-)\otimes \nu^{i})$ for any integer $i$.
\end{enumerate}
\end{theorem}
\begin{proof}

By combining the results of Proposition \ref{prop:singularff}, Proposition \ref{prop:a-r}, Proposition \ref{prop:ffofemeb}, and Proposition \ref{prop:r-a} we get the theorem.

\end{proof}

\subsection{Derived category of Francia's flip}

In this section we will discuss the behaviors of $3$-dimensional locally non-commutative Francia's flip, the result can be easily adjusted to the higher dimension cases with arbitrary weight as a weighted version of standard flip (without discussion in here). In the global case, we refer to \cite[Section 4]{KaFf}. Instead, we consider an analytically local flip to have this property.\footnote{\, $3$-folds we are dealing with in this section could be Moishezon (algebraic space).}\\

If $X$ is $3$-fold with only terminal quotient singularities which admits $\acute{e}$tale locally $3$-dimensional Francia's flip of type $(1, a+b; a, b)$, where $a$ and $b$ are coprime positive integers, by  a similar argument as the previous case, we can firstly reduce our argument to an $\acute{e}$tale cover of certain classical Francia's flip, that is

 $$\begin{tikzcd}[column sep=0.5em]
 & &X\arrow[dr]& &X^{+}\arrow[dl]&  \\
&& &X^{0}&
\end{tikzcd}$$
where $X^{0}$ could be realize as $U//\mathds{G}_{m}$ for certain smooth dimension $4$ variety $U$ with a coordinate defined by regular sequence
$x^{0}_{1},x^{0}_{2},x^{0}_{3},x^{0}_{4}$ (which corresponds to the intrinsic coordinate of the affine space through an $e\acute{}$tale mapping),
while the action of $U//\mathds{G}_{m}$ is given by
$$\mathds{G}_{m}\times U\longrightarrow U,\quad \lambda\times (x^{0}_{1},x^{0}_{2},x^{0}_{3},x^{0}_{4})\longmapsto (\lambda x^{0}_{1},\lambda^{(a+b)}x^{0}_{2},\lambda^{-a}x^{0}_{3},\lambda^{-b}x^{0}_{4}) $$
while two different directions of VGIT give the variety $X$ and $X^{+}$.\\

If we consider a $G:=\mu_{a+b}\times\mu_{a}\times\mu_{b}$ Galois covering $W$ of $U$ which is defined by steps of cyclic covering along divisor $\mathrm{div}(x^{0}_{i})$, for $i=1,..,4$ to weight $(1,a+b,a,b)$, then we have a natural Galois covering of the flip diagram,

 $$\begin{tikzcd}[column sep=0.5em]
 & &Y:=W/\mathds{G}_{m}\arrow[dr]& &Y^{+}:=W/_{+}\mathds{G}_{m}\arrow[dl]&  \\
&& &Y^{0}:=W//\mathds{G}_{m}&
\end{tikzcd}$$
where the action of $\mathds{G}_{m}$ on $W$ has weight $(1,1,-1,-1)$, especially since it is a standard flip, and noticing the compatible of actions of group $\mathds{G}_{m}$ and $G$, it yields a finite group $G$ action on the covering diagram.\\

We consider the quotient of $G$, then we pass this picture to equivariant situations and associate it with their canonical stacks,

 $$\begin{tikzcd}[column sep=0.5em]
  & && {[Bl_{L}Y/G]}\arrow[dl,"\phi"]\arrow[dr,"\varphi"] &&  \\
 & &{[Y/G]}\arrow[dr]& &{[Y^{+}/G]}\arrow[dl]&  \\
&& &{[Y^{0}/G]}&
\end{tikzcd}$$

\begin{proposition}[{\cite[Section 4]{KaFf}}]\label{prop:rootpropflip}

\begin{enumerate}
\item The quotient morphism $Y\longrightarrow X$ is ramified at $S_{a+b}$  of degree $a+b$, $S_{a}$  of degree $a$ and $S_{b}$  of degree $b$ in codimension 1 locus, so
$$[Y/G]=\sqrt[a+b]{{S_{a+b}}/\widetilde{X}}\times_{\widetilde{X}}\sqrt[a]{{S_{a}}/\widetilde{X}}\times_{\widetilde{X}}\sqrt[ b]{{S_{b}}/\widetilde{X}}$$
\item  The quotient morphism $Y^{+}\longrightarrow X^{+}$ is ramified at $S^{+}_{a+b}$  of degree $a+b$ , $S^{+}_{a}$  of degree $a$ and  $S^{+}_{b}$  of degree $b$ in codimension 1 locus, so
$$[Y^{+}/G]=\sqrt[a+b]{{S^{+}_{a+b}}/\widetilde{X}^{+}}\times_{\widetilde{X}^{+}}\sqrt[ a]{{S^{+}_{a}}/\widetilde{X}^{+}}\times_{\widetilde{X}^{+}}\sqrt[b]{{S^{+}_{b}}/\widetilde{X}^{+}}$$
\item  The quotient morphism $Y\longrightarrow X$ is ramified at point $t_{a+b}$ of degree $a+b$ in small locus. \\
\item  The quotient morphism $Y^{+}\longrightarrow X^{+}$ is ramified at point $t^{+}_{b}$ of degree $b$ and  $t^{+}_{a}$ of degree $a$ in small locus.\end{enumerate}
\end{proposition}

The ramified divisors' behaviors under local flip:\\
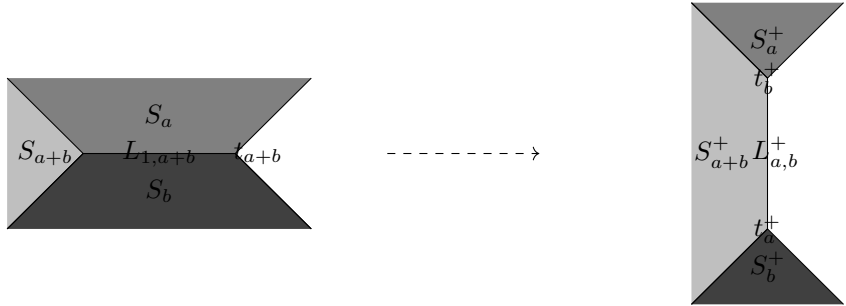
$$\begin{tikzpicture}
\coordinate(A)at(0,0);
\coordinate(B)at(1,1);
\coordinate(C)at(0,2);
\coordinate(D)at(3,1);
\coordinate(E)at(4,0);
\coordinate(F)at(4,2);

\fill[lightgray][label=left:A](A)--(B)--(C)--cycle;
\fill[white](D)--(E)--(F)--cycle;
\fill[darkgray](A)--(B)--(D)--(E)--cycle;
\fill[gray](F)--(D)--(B)--(C)--cycle;
\draw (A)--(B)--(C);
\draw (F)--(D)--(E);
\draw(A)--(B)--(D)--(E);
\node at (3.3,1) {$t_{a+b}$};
\node at (2,1) {$L_{1,a+b}$};
\node at (0.5,1) {$S_{a+b}$};
\node at (2,1.5) {$S_{a}$};
\node at (2,0.5) {$S_{b}$};

\coordinate(A)at(9,-1);
\coordinate(B)at(10,0);
\coordinate(C)at(11,-1);
\coordinate(D)at(10,2);
\coordinate(E)at(9,3);
\coordinate(F)at(11,3);

\fill[darkgray][label=left:A](A)--(B)--(C)--cycle;
\fill[gray](E)--(D)--(F)--cycle;
\fill[lightgray](A)--(B)--(D)--(E)--cycle;
\fill[white](C)--(B)--(D)--(F)--cycle;
\draw (A)--(B)--(C);
\draw (F)--(D)--(E);
\draw(A)--(B)--(D)--(E);
\draw [->,dashed] (5,1) --(7,1);

\node at (D) {$t^{+}_{b}$};
\node at (B) {$t^{+}_{a}$};
\node at (9.4,1) {$S^{+}_{a+b}$};
\node at (10,2.5) {$S^{+}_{a}$};
\node at (10,-0.5) {$S^{+}_{b}$};
\node at (10.1,1) {$L^{+}_{a,b}$};

\end{tikzpicture}$$

If we pass this equivariant standard flop to the picture of its canonical stack of quotient space, by comparing the exceptional locus we have:

 $$\begin{tikzcd}[column sep=0.5em]
  & && {Z:={wBl^{a,b}_{\widetilde{L}_{1,a+b}}\widetilde{X}}\simeq{wBl^{1,a+b}_{\widetilde{L}^{+}_{a,b}}\widetilde{X}^{+}}}\arrow[dl,"\phi"]\arrow[dr,"\varphi"] &&  \\
 & &{\widetilde{X}}& &{\widetilde{X}^{+}}&
\end{tikzcd}$$
In addition, due to a lack of reference, we need the following calculations for the normal sheaves of these flipping (flipped) loci like the standard flip cases.
\begin{proposition}
We have $$\mathcal{N}_{\widetilde{L}_{1,a+b}|\widetilde{X}}\simeq \mathcal{O}_{\widetilde{L}_{1,a+b}}(-a)\bigoplus\mathcal{O}_{\widetilde{L}_{1,a+b}}(-b)$$ and $$\mathcal{N}_{\widetilde{L}^{+}_{a,b}|\widetilde{X}^{+}}\simeq \mathcal{O}_{\widetilde{L}^{+}_{a,b}}(-1)\bigoplus\mathcal{O}_{\widetilde{L}^{+}_{a,b}}(-a-b)$$
where we view $\widetilde{L}_{1,a+b}$ (or $\widetilde{L}^{+}_{a,b}$) as weighted projective stack $\mathcal{P}(1,a+b)$ (or $\mathcal{P}(a,b)$).
\end{proposition}

\begin{proof}
Noticing the above construction, we can see $\widetilde{L}_{1,a+b}$ is a complete intersection of $S_{a}$ and $S_{b}$, we assume
$S_{a}|_{\widetilde{L}_{1,a+b}}=\mathcal{O}_{\widetilde{L}_{1,a+b}}(m)$ for certain integer $m$. \\

First, we consider the first step of root construction in Proposition \ref{prop:rootpropflip},
$$\pi_{1}:\sqrt[a]{{S_{a}}/\widetilde{X}}\longrightarrow\widetilde{X}$$
we have
$$\big(\mathcal{M}_{S_{a}}|_{\sqrt[a]{\widetilde{L}_{1,a+b}}}\big)^{\otimes a}=\pi_{1}^{*}\big(S_{a}|_{\widetilde{L}_{1,a+b}}\big)=\mathcal{O}_{\sqrt[a]{\widetilde{L}_{1,a+b}}}(m)$$
so $\mathcal{M}_{S_{a}}|_{\sqrt[a]{\widetilde{L}_{1,a+b}}}=\mathcal{O}_{\sqrt[a]{\widetilde{L}_{1,a+b}}}(m/a)$, where $\sqrt[a]{\widetilde{L}_{1,a+b}}\simeq[\widetilde{L}_{1,a+b}/\mu_{a}]$.\\

Then follows the second root structure:
$$\pi_{2}:\sqrt[a]{{S_{a}}/\widetilde{X}}\times_{\widetilde{X}}\sqrt[ b]{{S_{b}}/\widetilde{X}}\simeq \sqrt[ b]{{\widehat{S_{b}}}/\sqrt[a]{{S_{a}}/\widetilde{X}}}\longrightarrow\sqrt[a]{{S_{a}}/\widetilde{X}}$$
but this time the pullback and strict transformation of $\mathcal{M}_{S_{a}}$ are consistent, we have
$$\widehat{\mathcal{M}_{S_{a}}}|_{\sqrt[b]{\sqrt[a]{\widetilde{L}_{1,a+b}}}}=\pi_{2}^{*}\big(\mathcal{M}_{S_{a}}|_{\sqrt[a]{\widetilde{L}_{1,a+b}}}\big)=\mathcal{O}_{\sqrt[b]{\sqrt[a]{\widetilde{L}_{1,a+b}}}}(m/a)$$
where ${\sqrt[b]{\sqrt[a]{\widetilde{L}_{1,a+b}}}}\simeq[\widetilde{L}_{1,a+b}/\mu_{a}\times\mu_{b}]$.\\

Finally, we consider the third root structure:
$$\pi_{3}:\sqrt[a+b]{{S_{a+b}}/\widetilde{X}}\times_{\widetilde{X}}\sqrt[a]{{S_{a}}/\widetilde{X}}\times_{\widetilde{X}}\sqrt[ b]{{S_{b}}/\widetilde{X}}\longrightarrow\sqrt[a]{{S_{a}}/\widetilde{X}}\times_{\widetilde{X}}\sqrt[ b]{{S_{b}}/\widetilde{X}}$$
restricted to the flipping locus, it is
$$\big[[\mathds{P}^{1}/\mu_{a+b}]/\mu_{a}\times\mu_{b}\big]\longrightarrow\big[\mathcal{P}(1,a+b)/\mu_{a}\times\mu_{b}\big]$$
it is a $\mu_{a}\times\mu_{b}$ gerbe over
\begin{equation}\label{equ:coarstop1a}[\mathds{P}^{1}/\mu_{a+b}]\longrightarrow\mathcal{P}(1,a+b)\end{equation}
we defined in Lemma \ref{lemma:Kawambi} and Theorem \ref{thm:beilinson}. And similar to our second root, the pullback and strict transformation of $\widehat{\mathcal{M}_{S_{a}}}$ along $\pi_{3}$ are consistent, we have
$$\widehat{\mathcal{M}_{S_{a}}}|_{\big[[\mathds{P}^{1}/\mu_{a+b}]/\mu_{a}\times\mu_{b}\big]}=\pi_{3}^{*}\big(\widehat{\mathcal{M}_{S_{a}}}|_{\big[\mathcal{P}(1,a+b)/\mu_{a}\times\mu_{b}\big]}\big)=\pi_{3}^{*}\mathcal{O}_{\big[\mathcal{P}(1,a+b)/\mu_{a}\times\mu_{b}\big]}(m/a)$$
and by our previous argument, we know after this 3 steps of root, it's an equivariant standard flop, so $\widehat{\mathcal{M}_{S_{a}}}|_{\big[[\mathds{P}^{1}/\mu_{a+b}]/\mu_{a}\times\mu_{b}\big]}=\mathcal{O}_{\big[[\mathds{P}^{1}/\mu_{a+b}]/\mu_{a}\times\mu_{b}\big]}(-1)=\pi_{3}^{*}\mathcal{O}_{\big[\mathcal{P}(1,a+b)/\mu_{a}\times\mu_{b}\big]}(m/a)$, so by morphism (\ref{equ:coarstop1a}), Lemma \ref{lemma:Kawambi} and Theorem \ref{thm:beilinson}, we have $m/a=-1$ hence $m=-a$.\\

We have $S_{a}|_{\widetilde{L}_{1,a+b}}=\mathcal{O}_{\widetilde{L}_{1,a+b}}(-a)$, similarly $S_{b}|_{\widetilde{L}_{1,a+b}}=\mathcal{O}_{\widetilde{L}_{1,a+b}}(-b)$ and we finish the proof.
\end{proof}
So we can adjust the classical non-commutative standard flip's arguments and apply them from the ordinary Orlov's non-commutative blow-up formula to our weighted version, we have an evident generalization of a classical result:

\begin{theorem}\label{thm:fraciaflip}
If $X$  is  a $3$-fold with only terminal quotient singularities, and it admits a Francia's type flip of weight $(1,a+b; a,b)$ to $X^{+}$, then we have a semi-orthogonal decomposition of $\mathrm{D}^{b}(\widetilde{X})$ respect to $\mathrm{D}^{b}(\widetilde{X}^{+})$:
$$\mathrm{D}^{b}(\widetilde{X})\simeq\langle \mathcal{O}_{L_{1,a+b}}(-1),\phi_{*}\varphi^{*}\mathrm{D}^{b}(\widetilde{X}^{+})\rangle$$
where we view $\widetilde{L}_{1,a+b}$ as weighted projective stack $\mathcal{P}(1,a+b)$.
\end{theorem}
\begin{proof}
  See \cite[Theorem A]{BFR}, while noticing that their proof depends only on the geometry structure of blow-up blow-down and their derived category formulas, the proof is nearly the same.
\end{proof}

\begin{example}[Shokurov link for $(1,2:1,1)$-type Francia's flip]\label{ex:fanciaflip}
If $X^{+}$ is a smooth $3$-dimensional projective variety with a $(-1,-2)$ smooth rational curve $C^{+}$, we blow up $X^{+}$ along the curve $C^{+}$ then run the minimal model program, which consists of a standard flop (flop the zero section of the exceptional divisor) and a blow-down (contract the strict transformation of the exceptional divisor (a plane) to a $(1,1,1)/2$-type singularity).  Specially the strict transformation of  flopped curve to $X$ is a stacky curve $C=\mathcal{P}^{1}_{2}$ (or $\mathcal{P}(1,2))$, with normal sheaf $(-1/2,-1/2)$ (or $(-1,-1)$).\\
$$\begin{tikzcd}[ampersand replacement=\&]
Y=w_{1/2}Bl^{(1,1,1)}_{p}X \arrow[r,dashrightarrow] \arrow[d,"f"] \& Y^{+}\arrow[d,"f^{+}"] \\
X\& X^{+}
\end{tikzcd}$$

We can calculate the behaviors of normal sheaves by exact sequences.

$$\begin{tikzcd}[ampersand replacement=\&]
Z:=\sqrt[2]{{E/Y}}\arrow[r,"\varpi"] \arrow[d,"\pi"] \& Y:=w_{1/2}Bl^{(1,1,1)}_{p} X \\
\widetilde{X}\&
\end{tikzcd}\hspace{2cm}
\begin{tikzcd}[column sep=0.5em]
 & && W\arrow[dl,"\overline{\phi}"]\arrow[dr,"\overline{\varphi}"] &&  \\
&& Y&& Y^{+}
\end{tikzcd}$$
Noticing the flopped curve $\widehat{C_{Y}}$ is a smooth rational curve with $\mathcal{N}_{\widehat{C_{Y}}| Y}$ isomorphic to $ \mathcal{O}_{\widehat{C_{Y}}}(-1)\oplus\mathcal{O}_{\widehat{C_{ Y}}}(-1)$, so
$$\mathrm{\Omega}_{ Y}|_{\widehat{C_{Y}}}\simeq \mathcal{O}_{\widehat{C_{Y}}}(1)^{2}\oplus\mathcal{O}_{\widehat{C_{Y}}}(-2)$$
and $\varpi$ is stacky smoothing of a cyclic cover branched at exceptional divisor $E$, then we restrict the exact sequence
$$0\longrightarrow \varpi^{*}\mathrm{\Omega}_{ Y}\longrightarrow\mathrm{\Omega}_{Z}\longrightarrow
\mathcal{O}_{M_{E}}(-M_{E})\longrightarrow0$$
on $\widehat{C_{Z}}$ (the strict transformation of flopping curve at $Z$), then
$$\mathrm{\Omega}_{ Z}|_{\widehat{C_{Z}}}\simeq \mathcal{O}_{\widehat{C_{Z}}}(1)^{\oplus2}\oplus\mathcal{O}_{\widehat{C_{Z}}}(-2+\frac{1}{2})$$
similarly, we restrict the exact sequence of differential
$$0\longrightarrow \pi^{*}\mathrm{\Omega}_{ \widetilde{X}}\longrightarrow\mathrm{\Omega}_{Z}\longrightarrow
\mathrm{\Omega}_{M_{E}/[p/\mu_{2}]}\longrightarrow0$$
on $\widehat{C_{Z}}$,  and since $M_{E}$ is a $\mu_{2}$ gerbe over plane, we have
$$\mathrm{\Omega}_{ \widetilde{X}}|_{C}\simeq \mathcal{O}_{C}(\frac{1}{2})^{\oplus 2}\oplus\mathcal{O}_{C}(-2+\frac{1}{2})$$
and $C$ (the strict transformation of flopped curve at $\widetilde{X}$) is a stacky curve $\mathcal{P}^{1}_{2}$ with $\mathrm{K}_{C}\simeq \mathcal{O}_{C}(-3/2)$, so
$$\mathcal{N}_{\widetilde{X}|C}\simeq \mathcal{O}_{C}(-\frac{1}{2})\oplus\mathcal{O}_{C}(-\frac{1}{2})$$
We can also calculate in the opposite direction to get the relationship between their normal sheaves.\\

Conversely, if $X$ is a $3$-fold with only one $(1,1,1)/2$ singular point, and a stacky $(-1,-1)$ curve $C\simeq\mathcal{P}(1,2)$ passing it, then we can do a Kawamata blow-up along this $(1,1,1)/2$ point, then flop the strict transformation of this stacky curve $C$, and blow down the strict transformation of the exceptional divisor of Kawamata blow-up (which is a $\mathbb{F}_{1}$ surface), we finally get a smooth $3$-fold with a $(-1,-2)$ smooth rational curve $C^{+}$ (strict transformation of flopped curve). Also noting that for the general $(1,n,-1,-1)$ type Francia's flip we have similar birational descriptions. In these cases, all (coarse moduli of) objects are projective varieties if $C$ or $C^{+}$ is contractible.\qed
\end{example}

\begin{example}[Francia's non-commutative projection]\label{ex:fanciaprojection}
Let us continue to consider the situation of the previous example,

$$\begin{tikzcd}[ampersand replacement=\&]
Y=w_{1/2}Bl^{(1,1,1)}_{p}X\arrow[dr,"g"] \arrow[d,"f"] \& \& Y^{+}\arrow[d,"f^{+}"] \arrow[dl,"g^{+}"]\\
X\& Z\& X^{+}
\end{tikzcd}$$

It's not difficult to see $\mathrm{\Theta}(\mathcal{O}_{C}(-1/2))$ fits in a distinguished triangle as following:
\begin{equation}\label{equ:distingcc}\mathrm{\Theta}(\mathcal{O}_{C}(-1/2))\longrightarrow\mathcal{O}_{\widehat{C_{Y}}}(-1)\longrightarrow\mathcal{O}_{E}(-1)[2]\end{equation}
where $\mathrm{\Theta}(-):=\varpi_{*}\pi^{*}(-)$, $\widehat{C_{Y}}$ is the flopping curve, since
\begin{equation}\label{equ:wpit}\mathrm{\Theta}^{!}(\mathcal{O}_{\widehat{C_{Y}}}(-1))=\pi_{*}(\varpi^{*}\mathcal{O}_{\widehat{C_{Y}}}(-1)\otimes \mathcal{M}_{E})=\pi_{*}\mathcal{O}_{\widehat{C_{Y}}}(-1/2)=\mathcal{O}_{C}(-1/2)\end{equation}
 and $$\mathrm{\Theta}_{-1}^{*}(\mathcal{O}_{\widehat{C_{Y}}}(-1))=f_{E*}\big(\iota^{*}(\mathcal{O}_{\widehat{C_{Y}}}(-1))(-2)\big)[2]=\mathrm{k}[2]$$
 recall Theorem \ref{thm:singularweightedblowup} we get the result.\\

 Noticing $g$ is an extremal contraction which contracts $\widehat{C_{Y}}$ to an ordinary double point, we are expected to pushforward the semi-orthogonal decomposition of
 $\mathrm{D}^{b}(\mathrm{Y})$ obtained in Theorem \ref{thm:singularweightedblowup} and Theorem \ref{thm:fraciaflip}
 $$\langle \mathcal{O}_{E}(-1), \mathrm{\Theta}(\mathrm{D}^{b}(\widetilde{X}))\rangle\simeq\langle \mathcal{O}_{E}(-1), \mathrm{\Theta}(\mathcal{O}_{C}(-1/2)),\mathrm{\Theta}(\mathcal{A})\rangle$$
 along $g$, where $\mathcal A$ denotes the orthogonal component in $\mathrm{D}^{b}(\widetilde X)\simeq\langle\mathcal{O}_C(-1/2),\mathcal A\rangle$. We do not identify this embedding with the unmodified block for the normalization having exceptional object $\mathcal{O}_C(-1)$.\\

Put $D:=\widehat{E_Z}$ and let $p\in D$ be the ordinary double point of $Z$.
The local pair at $p$ is $(xy-zw=0,\ x=z=0)$, so $D$ is regular there.
Away from $p$, the map $g$ is an isomorphism near $D$.
The restriction $g|_E:E\to D$ is proper, quasi-finite and birational, hence
is an isomorphism because $D$ is normal. In particular $D\simeq\mathds P^2$.
Since $g$ is crepant and $\mathcal{O}_E(E)=\mathcal{O}_E(-2)$,
adjunction on $Y$ gives
\[
\omega_Z|_D\simeq\omega_Y|_E
\simeq\omega_E\otimes\mathcal{O}_E(-E)
\simeq\mathcal{O}_D(-1).
\]
Consequently $\mathcal E xt_Z^1(\mathcal{O}_D,\mathcal{O}_Z)
\simeq\omega_D\otimes\omega_Z^{-1}|_D\simeq\mathcal{O}_D(-2)$.

The self-extension sheaves can be computed directly. At $p$, choose the
completed local presentation
\[
R=\mathrm k[[x,y,z,w]]/(xy-zw),\qquad
M=R/(x,z)=\mathrm k[[y,w]].
\]
The free resolution of $M$ starts with $d_1=(x\ \ z)$ and is then
$2$-periodic, with alternating matrices
\[
d_2=\begin{pmatrix}-z&y\\x&-w\end{pmatrix},\qquad
d_3=\begin{pmatrix}-w&-y\\-x&-z\end{pmatrix}.
\]
Here $d_2d_3=d_3d_2=-(xy-zw)I_2$ over the regular power-series ring.
This is the stabilization of the Koszul presentation of $(x,z)$, so it is
exact over $R$. Applying $\operatorname{Hom}_R(-,M)$ gives the differentials
\[
0,\qquad
\begin{pmatrix}0&0\\y&-w\end{pmatrix},\qquad
\begin{pmatrix}-w&0\\-y&0\end{pmatrix},\quad\ldots.
\]
The first extension module is free of rank one. More precisely,
\[
\begin{aligned}
\operatorname{Ext}_R^1(M,M)&=M(w,y),\\
\operatorname{Ext}_R^{2j}(M,M)&=\mathrm k,\qquad
\operatorname{Ext}_R^{2j+1}(M,M)=0\quad(j\ge1).
\end{aligned}
\]
Away from $p$, the first self-extension sheaf is the normal bundle
$\mathcal{O}_D(-2)$. Since the local computation shows that it is locally free
also at $p$, these calculations identify
\[
\mathcal E xt_Z^q(\mathcal{O}_D,\mathcal{O}_D)=
\begin{cases}
\mathcal{O}_D,&q=0,\\
\mathcal{O}_D(-2),&q=1,\\
\mathrm k_p,&q\ge2\text{ even},\\
0,&q\ge3\text{ odd}.
\end{cases}
\]
The same self-extension sheaves are obtained after twisting both arguments
by $\mathcal{O}_D(-1)$. The vanishing
\[
\mathrm H^*(D,\mathcal{O}_D(-2))=0
\]
and the local-to-global spectral sequence give
\[
\operatorname{Ext}_Z^q(\mathcal{O}_D(-1),\mathcal{O}_D(-1))=
\begin{cases}
\mathrm k,&q=0\text{ or }q\ge2\text{ even},\\
0,&\text{otherwise}.
\end{cases}
\]
In particular, the pushforward $g_*\mathcal{O}_E(-1)=\mathcal{O}_D(-1)$ is
not exceptional, despite the exceptionality of $\mathcal{O}_E(-1)$ on $Y$.
\end{example}

The following proposition is another proof of the derived category of Francia's flip, by changing the path of the birational correction. There could be more semi-orthogonal decompositions, but we only consider one of them here.

\begin{proposition}\label{prop:secondshoukufr}
For $(1,2:1,1)$ type Francia's flip, we have another categorical description under Shokurov's link:
\begin{equation}\mathrm{D}^{b}(\widetilde{X})\simeq\langle\mathcal{O}_{C}(-1), \mathrm{Im}\,\Gamma\big(D^{b}(X^{+})\big) \rangle\end{equation}
where $\Gamma(-):=\pi_{*}\varpi^{!}\overline{\phi}_{*}\overline{\varphi}^{!}f^{+*}(-)$  is a fully-faithful functor, and we view $C$ as $\mathcal{P}_{2}^{1}$.
\end{proposition}

\begin{proof}
\begin{enumerate}
\item Consider $\mathrm{D}^{b}(Y^{+})$ by blowing up formula of derived category e.g. Theorem \ref{thm:singularweightedblowup}, we have
$$\mathrm{D}^{b}(Y^{+})\simeq\langle \iota^{+}_{*}\big( f^{+*}_{E}\mathrm{D}^{b}(C^{+})\otimes\mathcal{O}_{E^{+}}(-1)\big), f^{+*}\mathrm{D}^{b}(X^{+})\rangle$$
where $f^{+}_{E}$ is $f^{+}$ restriction on exception locus and $\iota^{+}$ is the embedding of exception locus of $f^{+}$.
\item By the construction of flip, we know $\mathrm{K}_{X}C=-1/2$ and $\mathrm{K}_{X^{+}}C^{+}=1$, so we have a description for
$$\mathrm{D}^{b}(C^{+})\simeq\langle\mathcal{O}_{C^{+}}(a),\mathcal{O}_{C^{+}}(a+1)\rangle\simeq\langle\mathcal{O}_{C^{+}}(a\mathrm{K}_{X^{+}}),\mathcal{O}_{C^{+}}((a+1)\mathrm{K}_{X^{+}})\rangle $$
for any integer $a$.
\begingroup\small
\begin{align*}
  \mathrm{D}^{b}(Y^{+})&\simeq\langle \iota^{+}_{*}\big(f^{+*}_{E}\mathcal{O}_{C^{+}}(a\mathrm{K}_{X^{+}})\otimes\mathcal{O}_{E^{+}}(-1)\big),\iota^{+}_{*} \big(f^{+*}_{E}\mathcal{O}_{C^{+}}((a+1)\mathrm{K}_{X^{+}})\otimes\mathcal{O}_{E^{+}}(-1)\big), f^{+*}\mathrm{D}^{b}(X^{+})\rangle \\
   &\simeq \langle \mathcal{O}_{E^{+}}(a\mathrm{K}_{X^{+}}+E^{+}), \mathcal{O}_{E^{+}}((a+1)\mathrm{K}_{X^{+}}+E^{+}), f^{+*}\mathrm{D}^{b}(X^{+})\rangle\stepcounter{equation}\tag{\theequation}\label{sod:fr1}
\end{align*}
\endgroup
\item Mutate the first two terms in (\ref{sod:fr1}) to the far right, noticing $\mathrm{K}_{Y^{+}}\sim \mathrm{K}_{X^{+}}+E^{+}$
$$ \mathrm{D}^{b}(Y^{+})\simeq\langle f^{+*}\mathrm{D}^{b}(X^{+}), \mathcal{O}_{E^{+}}((a-1)\mathrm{K}_{X^{+}}), \mathcal{O}_{E^{+}}(a\mathrm{K}_{X^{+}})\rangle$$
and we pick $a=1$ for our convenient.
\item Act standard flop derived functor $\overline{\phi}_{*}\overline{\varphi}^{!}(-)=\overline{\phi}_{*}(\mathcal{O}_{W}(E_{W})\otimes\overline{\varphi}^{*}(-))$ by e.g.\cite{H}, where $E_{W}$ is the exceptional divisor of the standard flop. From the position of the flopping (flopped) curve we know $E^{+}\sim E+E_{W}$ and $\mathrm{K}_{Y^{+}}\sim \mathrm{K}_{Y}$, so\\
\begin{enumerate}

      \item\label{equ:acompu0} By exact sequence:
    $$0\longrightarrow\mathcal{O}_{W}\longrightarrow\mathcal{O}_{W}(E_{W})\longrightarrow\mathcal{O}_{E_{W}}(E_{W})\longrightarrow0$$
    noticing the last term pushforward along $\overline{\phi}$ is zero, since $\overline{\phi}$ is a $\mathds{P}^{1}$ fibration at exceptional locus, so we have
$\overline{\phi}_{*}(\mathcal{O}_{W}(E_{W})\otimes\overline{\varphi}^{*}\mathcal{O}_{Y^{+}})\simeq \overline{\phi}_{*}\mathcal{O}_{W}(E_{W})\simeq  \overline{\phi}_{*}\mathcal{O}_{W}\simeq\mathcal{O}_{Y}$.\\

      \item and $\overline{\phi}_{*}(\mathcal{O}_{W}(E_{W})\otimes\overline{\varphi}^{*}\mathcal{O}_{Y^{+}}(-E^{+}))\simeq \overline{\phi}_{*}(\mathcal{O}_{W}(E_{W}-E-E_{W}))\simeq  \overline{\phi}_{*}\mathcal{O}_{W}(-E)\simeq\mathcal{O}_{Y}(-E)$\\

      \item combining the above observations and the exact sequence
    $$0\longrightarrow\mathcal{O}_{Y^{+}}(-E^{+})\longrightarrow\mathcal{O}_{Y^{+}}\longrightarrow\mathcal{O}_{E^{+}}\longrightarrow0$$
     $$0\longrightarrow\mathcal{O}_{Y}(-E)\longrightarrow\mathcal{O}_{Y}\longrightarrow\mathcal{O}_{E}\longrightarrow0$$
    we have $\overline{\phi}_{*}(\mathcal{O}_{W}(E_{W})\otimes\overline{\varphi}^{*}\mathcal{O}_{E^{+}})\simeq \mathcal{O}_{E}$.\\

      \item Similarly, for
\begin{equation}\label{equ:shoreqideal1}
0\longrightarrow\mathcal{O}_{Y^{+}}(-E^{+}-\mathrm{K}_{X^{+}})\longrightarrow\mathcal{O}_{Y^{+}}(-\mathrm{K}_{X^{+}})\longrightarrow\mathcal{O}_{E^{+}}(-\mathrm{K}_{X^{+}})\longrightarrow0
\end{equation}
we have
\begin{align*}
  \overline{\phi}_{*}(\mathcal{O}_{W}(E_{W})\otimes\overline{\varphi}^{*}(\mathcal{O}_{Y^{+}}(-\mathrm{K}_{X^{+}}))) &\simeq\overline{\phi}_{*}\mathcal{O}_{W}(E_{W}-\mathrm{K}_{Y^{+}}+E^{+}) \\
   &\simeq \overline{\phi}_{*}\mathcal{O}_{W}(E_{W}-\mathrm{K}_{Y}+E+E_{W}) \\
    &\simeq \overline{\phi}_{*}\mathcal{O}_{W}(2E_{W})\otimes\mathcal{O}_{X}(-\mathrm{K}_{Y}+E)\stepcounter{equation}\tag{\theequation}\label{equ:diviosre}
\end{align*}
by short exact sequence
$$0\longrightarrow\mathcal{O}_{W}(E_{W})\longrightarrow\mathcal{O}_{W}(2E_{W})\longrightarrow\mathcal{O}_{E_{W}}(2E_{W})\longrightarrow0$$
$\overline{\phi}_{*}\mathcal{O}_{E_{W}}(2E_{W})=\overline{\phi}_{*}\mathcal{O}_{\mathds{P}^{1}\times\mathds{P}^{1}}(-2,-2)=\mathcal{O}_{\widehat{C_{Y}}}(-2)[-1]$ and $(-\mathrm{K}_{Y}+E)\widehat{C_{Y}}=E\widehat{C_{Y}}=1$, we can see (\ref{equ:diviosre}) fits in a distinguished triangle
$$\mathcal{O}_{Y}(-\mathrm{K}_{Y}+E)\longrightarrow\overline{\phi}_{*}\mathcal{O}_{W}(2E_{W})\otimes\mathcal{O}_{Y}(-\mathrm{K}_{Y}+E)\longrightarrow\mathcal{O}_{\widehat{C_{Y}}}(-1)[-1]\longrightarrow$$
where $\widehat{C_{Y}}$ is the strict transformation flipping curve at $Y$ which is also the flopping curve.
    \begin{align*}
  \overline{\phi}_{*}(\mathcal{O}_{W}(E_{W})\otimes\overline{\varphi}^{*}(\mathcal{O}_{Y^{+}}(-E^{+}-\mathrm{K}_{X^{+}}))) &\simeq\overline{\phi}_{*}\mathcal{O}_{W}(E_{W}-\mathrm{K}_{Y^{+}}) \\
    &\simeq \overline{\phi}_{*}\mathcal{O}_{W}(E_{W})\otimes\mathcal{O}_{X}(-\mathrm{K}_{Y})\simeq\mathcal{O}_{X}(-\mathrm{K}_{Y})\stepcounter{equation}\tag{\theequation}\label{equ:diviosre1}
\end{align*}
combining  (\ref{equ:shoreqideal1}), (\ref{equ:diviosre}) and (\ref{equ:diviosre1}), we have $t:=\overline{\phi}_{*}(\mathcal{O}_{W}(E_{W})\otimes\overline{\varphi}^{*}(\mathcal{O}_{E^{+}}(-\mathrm{K}_{X^{+}})))$ has a filtration:
\begin{equation}\label{equ:diextension1}\begin{tikzcd}[column sep=0.5em]
 &  {\mathcal{O}_{Y}(-\mathrm{K}_{Y}+E)}\arrow{rr}&& \overline{\phi}_{*}(\mathcal{O}_{W}(E_{W})\otimes\overline{\varphi}^{*}(\mathcal{O}_{Y^{+}}(-\mathrm{K}_{X^{+}}))) \arrow{rr}\arrow{dl}&& t\arrow{dl} \\
 && {\mathcal{O}_{\widehat{C_{Y}}}(-1)[-1]}\arrow[ul,dashed,"\Delta"]  &&{\mathcal{O}_{Y}(-\mathrm{K}_{Y})[1]}\arrow[ul,dashed,"\Delta"]
 \end{tikzcd}\end{equation}
   \item Let's compute the extension
 $$\mathrm{Ext}^{*}_{Y}(\mathcal{O}_{Y}(-\mathrm{K}_{Y})[1],{\mathcal{O}_{\widehat{C_{Y}}}(-1)[-1]})=0$$
 since $\mathrm{K}_{Y}{\widehat{C_{Y}}}=0$, and
\begin{equation}\label{equ:acompu1}\mathrm{Ext}^{1}_{Y}(\mathcal{O}_{Y}(-\mathrm{K}_{Y})[1],{\mathcal{O}_{Y}(-\mathrm{K}_{Y}+E)})=\mathrm{H}^{0}(Y,{\mathcal{O}_{Y}(E)})=\mathrm{k}\end{equation}
which is due to we have a short exact sequence
 $$0\longrightarrow\mathcal{O}_{Y}\longrightarrow\mathcal{O}_{Y}(E)\longrightarrow\mathcal{O}_{E}(E)\longrightarrow0$$
 and $\mathcal{O}_{E}(E)\simeq\mathcal{O}_{\mathds{P}^{2}}(-2)$ which cohomology is trivial, we can see that the contribution of this non-trivial extension in (\ref{equ:diextension1}) comes from the unique morphism between $\mathcal{O}_{Y}(-\mathrm{K}_{Y})$ and $\mathcal{O}_{Y}(-\mathrm{K}_{Y}+E)$, and it yields a short exact sequence
 $$0\longrightarrow\mathcal{O}_{Y}(-\mathrm{K}_{Y})\longrightarrow\mathcal{O}_{Y}(-\mathrm{K}_{Y}+E)\longrightarrow\mathcal{O}_{E}(-\mathrm{K}_{Y}+E)\longrightarrow0$$
 then the previous filtration (\ref{equ:diextension1}) can be simplified to
 $$\mathcal{O}_{E}(-\mathrm{K}_{Y}+E)\longrightarrow t\longrightarrow {\mathcal{O}_{\widehat{C_{Y}}}(-1)[-1]}$$
 notice $(-\mathrm{K}_{Y}+E)|_{E}=(-\mathrm{K}_{E}+2E)|_{E}$, $\mathcal{O}_{E}(E)\simeq\mathcal{O}_{E}(-2)$ and $\mathcal{O}_{E}(\mathrm{K}_{E})\simeq\mathcal{O}_{E}(-3)$, so we can see $t:=\overline{\phi}_{*}(\mathcal{O}_{W}(E_{W})\otimes\overline{\varphi}^{*}(\mathcal{O}_{E^{+}}(-\mathrm{K}_{X^{+}})))$ has a clear description as following:
\begin{equation}\label{equ:keyexte}\mathcal{O}_{E}(-1)\longrightarrow t\longrightarrow {\mathcal{O}_{\widehat{C_{Y}}}(-1)[-1]}\end{equation}
 \\
Combining all the above calculations, we get
\begin{equation}\label{sod:2}
\mathrm{D}^{b}(Y)\simeq\langle \overline{\phi}_{*}\overline{\varphi}^{!}f^{+*}\mathrm{D}^{b}(X^{+}), t, \mathcal{O}_{E}\rangle
\end{equation}
where there is a filtration of $t$ in (\ref{equ:keyexte}).

\end{enumerate}
\item Right mutate $\mathcal{O}_{E}$ in (\ref{sod:2}) to the far left, we have
\begin{equation}\label{sod:fr3}\mathrm{D}^{b}(Y)\simeq\langle \mathcal{O}_{E}(-1),  \overline{\phi}_{*}\overline{\varphi}^{!}f^{+*}\mathrm{D}^{b}(X^{+}), t\rangle\end{equation}
since $\mathrm{K}_{Y}|_{E}=\mathrm{K}_{X}+E/2|_{E}=E|_{E}/2$, we have $\mathcal{O}_{E}(\mathrm{K}_{Y})\simeq\mathcal{O}_{\mathds{P}^{2}}(-1)$.
\item Notice $Y$ can also be obtained by Kawamata blow-up $X$ along its singular point, by our weighted blow-up formula in Theorem \ref{thm:singularweightedblowup} we have
\begin{equation}\label{sod:fr4}\mathrm{D}^{b}(Y)\simeq\langle \mathcal{O}_{E}(-1), \varpi_{*}\pi^{*}\mathrm{D}^{b}(\widetilde{X})\rangle\end{equation}
\item Let's compare the semi-orthogonal decomposition in (\ref{sod:fr3}) and (\ref{sod:fr4}), we have an equivalence induced by projection along $\pi_{*}\varpi^{!}$, that is
\begin{equation}\mathrm{D}^{b}(\widetilde{X})\simeq\langle \pi_{*}\varpi^{!} \overline{\phi}_{*}\overline{\varphi}^{!}f^{+*}\mathrm{D}^{b}(X^{+}), \pi_{*}\varpi^{!}t\rangle\end{equation}
\item We compute $\pi_{*}\varpi^{!}t$ by filtration in (\ref{equ:keyexte}), by the orthogonality of semi-orthogonal decomposition (\ref{sod:fr4}), we have
$\pi_{*}\varpi^{!}t\simeq \pi_{*}\varpi^{!}{\mathcal{O}_{\widehat{C_{Y}}}(-1)[-1]}$, where the latter term we have already calculated in (\ref{equ:wpit}) is just $\mathcal{O}_{C}(-1/2)[-1]$, so we have
\begin{equation}\label{sod:fr10}\mathrm{D}^{b}(\widetilde{X})\simeq\langle \pi_{*}\varpi^{!} \overline{\phi}_{*}\overline{\varphi}^{!}f^{+*}\mathrm{D}^{b}(X^{+}), \mathcal{O}_{C}(-1/2)\rangle\end{equation}
\item Finally, we mutate the second item in (\ref{sod:fr10}) to the far left
\begin{equation}\mathrm{D}^{b}(\widetilde{X})\simeq\langle \mathcal{O}_{C}(-1), \pi_{*}\varpi^{!} \overline{\phi}_{*}\overline{\varphi}^{!}f^{+*}\mathrm{D}^{b}(X^{+})\rangle\end{equation}
since $\mathrm{K}_{X}C=-1/2$, we get our result.
\end{enumerate}
\end{proof}
\begin{remark}
  We give a less precise guess here: we can use Shokurov's method to simplify some sufficiently complicated Francia's flips into simpler Francia's flips through weighted blow-up and weighted blow-down, and using a generalized argument in Proposition \ref{prop:secondshoukufr}, we can prove the Theorem  \ref{thm:fraciaflip}.
\end{remark}

\section{Takagi's hourglasses}\label{sec:Tak}

In this chapter, we refer to a series of articles \cite{Tak} \cite{Tak1} \cite{Tak2} by H.Takagi, we will briefly review these works first and then adapt them to our application, the interested readers could also read the above references further. In this chapter, we only consider the case variety lies in Takagi's classification Table 1 \cite[p.121]{Tak} or as shown in Table (\ref{table1}) with only orbifold singularities. The category classification of Table 2-5 in Takagi's classification can also be carried out similarly (if we consider their (non-commutative) crepant resolution for varieties with $cA_{n}$ singularity), we can use \cite[Theorem 0.11-0.21]{Tak1} to similarly construct their non-commutative structures. But because the reverse construction for varieties lies in Table 2-5 is not complete, and more appropriately we need to consider more exotic derived category structures of themselves (rather than their resolutions in terms of category), we ignore this part of the discussion.

\subsection{Classification by Takagi}
\begin{definition}\begin{enumerate}
\item Takagi's varieties are \lq\lq non-trivial" \footnote{\, Which means $h^{0}(-\mathrm{K}_{\mathrm{X}})\geq 4$, this number is related to the codimension ($h^{0}(-\mathrm{K}_{\mathrm{X}})+N(\mathrm{X})-4$ in these cases) of $\mathds{Q}$-Fano 3-fold, if they are small enough, we have structural theorems for their equations, see e.g. \cite{BKR}.\label{foot:sturthm}} prime $\mathds{Q}$-Fano 3-fold $\mathrm{X}$ over $\mathrm{k}$, with Fano index $1/2$,
also there exists a point $\mathrm{p}$ analytically isomorphic to
$$(xy+z^{2}+u^{a}=0)\Big/\mu_{2}^{(1,1,1,0)},a\in\mathds{N}$$
\item Specially, we call an hourglass is a \lq\lq non-trivial"\footnote{\, See above \ref{foot:sturthm}.} prime $\mathds{Q}$-Fano 3-fold $\mathrm{X}$ over $\mathrm{k}$, with Fano index $\mathrm{\frac{1}{2}}$ and only
$\mathrm{\frac{(1^{3})}{2}}$ singularities.
\end{enumerate}
\end{definition}

Takagi's program is to play a two-ray game after a weighted blow-up of the singular point $\mathrm{p}$, the picture is like this:

\begin{equation}\label{sqm}\begin{tikzcd}[ampersand replacement=\&]
\mathrm{Y}=w_{\frac{1}{2}}Bl^{(1,1,1,2)}_{\mathrm{p}}\mathrm{X} \arrow[r,dashrightarrow] \arrow[d,"\mathrm{f}"]\& \mathrm{Y}^{1}\arrow[r,dashrightarrow] \& \mathrm{Y}^{2}=\mathrm{Y'}\arrow[d,"\mathrm{f'}"] \\
\mathrm{X}\&\& \mathrm{X'}
\end{tikzcd}\end{equation}
in which the $\mathrm{SQM}$s consist of a step of flop first and then a flip. He extended Takeuchi's method and converted the classification problem of algebraic varieties into a series of Diophantine equations through the classification of Mori-type local extremal contractions and the calculation of the intersection number of divisors. By solving these Diophantine equations, he obtained tables of all possible classifications of these varieties e.g. Table (\ref{table1}) .\\

For the convenience of readers, we need to explain how to read the information in Takagi's table.\\

\begin{enumerate}
  \item $\mathrm{X}_{1,i}$ is a prime (with Picard number 1) normal projective variety with terminal singularity and $-\mathrm{K}_{\mathrm{X}_{1,i}}$ ample.
  \item Gorenstein index $I(\mathrm{X}_{1,i}):=\mathrm{min}\,\{I| \,I\mathrm{K}_{\mathrm{X}_{1,i}} \text{is Cartier}\}=2$.
  \item Cartier Fano index $F_{c}(\mathrm{X}_{1,i}):=r(\mathrm{X}_{1,i})/I(\mathrm{X}_{1,i})=1/2$, where $r(\mathrm{X}_{1,i})=1$ is the unique positive integer such that $$-I_{\mathrm{X}_{1,i}}\mathrm{K}_{\mathrm{X}_{1,i}}\sim r(\mathrm{X}_{1,i})H_{c} $$
      where $H_{c}$ is a primitive Cartier divisor.
  \item Weil Fano index $F_w(\mathrm X_{1,i})=1$, so that
$-\mathrm{K}_{\mathrm X_{1,i}}\sim H_w$, where $H_w$ is a primitive ample Weil divisor.
  \item Genus $g(\mathrm{X}_{1,i}):=h^{0}(-\mathrm{K}_{\mathrm{X}_{1,i}})-2$.
  \item Degree $\mathrm{Deg}(\mathrm{X}_{1,i}):=(-\mathrm{K}_{\mathrm{X}_{1,i}})^{3}$.
  \item\label{axialweight} Axial weight $\mathrm{N}(\mathrm{X}_{1,i})$ is defined to be the sum of all \textit{axial weight} of all terminal singularity points of singular index $>1$, where \textit{axial weight} of a singularity point $p$ means the sum of cyclic singularities up to a general deformation of $p$. If $\mathrm{X}_{1,i}$ is hourglass $\mathrm{N}(\mathrm{X}_{1,i})$ is the sum of orbifold points.
\item \label{flope}$e(\mathrm{X}_{1,i}):=\mathrm{Exc}(\mathrm{f})^{3}-\widehat{\mathrm{Exc}(\mathrm{f})_{\mathrm{Y}^{1}_{1,i}}}^{3}$, where $\widehat{\mathrm{Exc}(\mathrm{f})_{\mathrm{Y}^{1}_{1,i}}}$ is the strict transformation of exceptional divisor of $\mathrm{f}$ $\mathrm{Exc}(\mathrm{f})$ in $\mathrm{Y}^{1}_{1,i}$. If $\mathrm{X}_{1,i}$ is hourglass $e(\mathrm{X}_{1,i})$ is the number of flopping curves.
  \item \label{equ:fanoindex}$z(\mathrm{X}_{1,i})$ and $u(\mathrm{X}_{1,i})$ are integers defined by considering the unique decomposition in $\mathrm{Y}'_{\mathrm{X}_{1,i}}$
  $$\mathrm{Exc}(\mathrm{f}')\sim z(\mathrm{X}_{1,i})(-\mathrm{K}_{\mathrm{Y}'_{1,i}})-u(\mathrm{X}_{1,i})\widehat{\mathrm{Exc}(\mathrm{f})}$$
where $\mathrm{Exc}(\mathrm{f}')$ is the exceptional divisor of $\mathrm{f}'$, and $\widehat{\mathrm{Exc}(\mathrm{f})}$ is the strict transformation of  $\mathrm{f}$. If $\mathrm{X}_{1,i}$ is hourglass (in Table 1), we have $u(\mathrm{X}_{1,i})=z(\mathrm{X}_{1,i})+1$, specially it tells us $$F_{w}(\mathrm{X}'_{1,i})=m(z+1)$$
for some integer $m$.
\item \label{flipn}$n(\mathrm{X}_{1,i}):=2\big[(-\mathrm{K}_{\mathrm{Y}_{1,i}})^{3}-(-\mathrm{K}_{\mathrm{Y}'_{1,i}})^{3}\big]$. If $\mathrm{X}_{1,i}$ is hourglass $n(\mathrm{X}_{1,i})$ is the number of flipping curves.
\item $g(\mathrm{C_{1,i}} )$ is the genus of curve $\mathrm{C_{1,i}} $ which is the center of $\mathrm{f}'$.  If $\mathrm{X}_{1,i}$ is hourglass then $\mathrm{C_{1,i}} $ is smooth.

\item $\mathrm{deg}(\mathrm{C_{1,i}} ):=H_{c}(\mathrm{X}'_{1,i})\mathrm{C_{1,i}} $ is the degree of curve $\mathrm{C_{1,i}} $, where $H_{c}(\mathrm{X}'_{1,i})$ is the primitive ample Cartier divisor on $\mathrm{X}'_{1,i}$.
\item \label{clssX}$\mathrm{X}'_{1,i}$  is one of the following algebraic varieties, all of which are quasi-smooth if we consider $\mathrm{X}_{1,i}$ is an hourglass.
\begin{enumerate}\label{targetX}
  \item $\mathds{P}^{3}$, $F_{w}(\mathds{P}^{3})=4$
  \item $Q_{3}:=(2)\subset\mathds{P}^{4}$, $F_{w}(Q_{3})=3$
  \item $B_{3}:=(3)\subset\mathds{P}^{4}$, $F_{w}(B_{3})=2$
  \item $B_{4}:=(2)^{2}\subset\mathds{P}^{5}$, $F_{w}(B_{4})=2$
  \item $B_{5}:=(1)^{3}\subset\mathbf{Gr}(2,5)$, $F_{w}(B_{5})=2$
  \item $[2]:=(4)\subset\mathds{P}(1^{3},2^{2})$, $F_{w}([2])=3$
  \item $[3]:=(3)\subset\mathds{P}(1^{4},2)$, $F_{w}([3])=3$
  \item $[5]:=\mathds{P}(1^{3},2)$, $F_{w}([5])=5$
\end{enumerate}
and we can see  $m$ in (\ref{equ:fanoindex}) are all $1$ and $F_{w}(\mathrm{X}'_{1,i})=z(\mathrm{X}_{1,i})+1$.

\item\label{targetZ} Gorenstein terminal non-$\mathds{Q}$-factorial varieties $A_{\mathrm{Deg}(\mathrm{X}_{1,i})-N(\mathrm{X}_{1,i})/2}$ which could serve as the global midpoint of Sarkisov links are as following:
\begin{enumerate}
  \item\label{a6} $A_{6}:=(2)\cap(3)\subset\mathds{P}^{5}$, $F_{c}(A_{6})=1$
  \item\label{a8} $A_{8}:=(2)^{3}\subset\mathds{P}^{6}$, $F_{c}(A_{8})=1$
  \item\label{a10} $A_{10}:=(1)^{2}\cap(2)\subset\mathbf{Gr}(2,5)$, $F_{c}(A_{10})=1$
  \item\label{a12} $A_{12}:=(1)^{7}\subset\mathbf{OGr}_{+}(5,10)$, $F_{c}(A_{12})=1$
  \item\label{a14} $A_{14}:=(1)^{5}\subset\mathbf{Gr}(2,6)$, $F_{c}(A_{14})=1$
\end{enumerate}
and two special Picard number $2$ cases that may serve as examples of our extended global midpoint
\begin{enumerate}
  \item\label{c6} $C_{6}:$ a double cover of $\mathds{P}^{1}\times\mathds{P}^{2}$ branched along a $(2,4)$ divisor,
  \item\label{c8} $C_{8}:=Bl_{L}A_{2}$  where $L$ is a special elliptic curve.
\end{enumerate}
\end{enumerate}
\begin{table}  \label{table1}\centering
        \begin{tabular}{c c c c c c c c c c}
            \toprule
            \textbf{Hourglass} & \textbf{Genus} & \textbf{Degree} & \textbf{N}  & \textbf{e}  & \textbf{z} & \textbf{n} & \textbf{d(C)} & \textbf{g(C)}& \textbf{Target/X'}\\
            \midrule
            $\mathrm{X_{1.1}}$      & 4   & 7         & 2     & 7    & 4  & 0   & 7   & 8  & $\mathrm{[5]}$   \\
            $\mathrm{X_{1.2}}$      & 4   & 7.5       & 3     & 7    & 2   & 0  & 3   & 0   & $\mathrm{[2] } $  \\
            $\mathrm{X_{1.4}}$      & 5   & 8.5        & 1     & 6    & 3  & 0   & 9   & 9  & ${\mathds{P}}^{3}$   \\
            $\mathrm{X_{1.5}}$      & 5   & 9          & 2     & 6    & 2  & 0   & 6   & 3  & $\mathrm{[3]}$   \\
            $\mathrm{X_{1.9}}$      & 6   & 10.5      & 1     & 6    & 1  & 0   & 3   & 0  & $B_{3}$   \\
            $\mathrm{X_{1.10}}$      & 6   & 10.5       & 1     & 5    & 2  & 0   & 9   & 6 & $Q_{3}$    \\
            $\mathrm{X_{1.12}}$      & 7   & 12.5       & 1     & 5    & 1  & 0   & 5   & 1  & $B_{4}$   \\
            $\mathrm{X_{1.13}}$      & 8   & 14.5     & 1     & 4    & 1   & 0  & 7   & 2 & $B_{5}$    \\
            \bottomrule
        \end{tabular}
        \caption{part of Takagi's tables (Table 1 and $\mathrm{f'}$ is of $(2,1)$ type and $n=0$)}
\quad\\
        \begin{tabular}{c c c c c c c c c c}
            \toprule
            \textbf{Hourglass} & \textbf{Genus} & \textbf{Degree} & \textbf{N}  & \textbf{e}  & \textbf{z} & \textbf{n} & \textbf{d(C)} & \textbf{g(C)}& \textbf{Target/X'}\\
            \midrule
            $\mathrm{X_{1.3}}$       & 4   & 7.5     & 3     & 6    & 4   & 1  & 6   & 3  & $\mathrm{[5]}$   \\
            $\mathrm{X_{1.6}}$      & 5   & 9        & 2     & 5    & 3  & 1   & 8   & 5   & ${\mathds{P}}^{3}$  \\
            $\mathrm{X_{1.7}}$      & 5   & 9.5      & 3     & 5    & 2  & 1   & 5   & 0  & $\mathrm{[3]}$   \\
            $\mathrm{X_{1.8}}$      & 5   & 9.5     & 3     & 4    & 3   & 2  & 7   & 1   & ${\mathds{P}}^{3}$  \\
            $\mathrm{X_{1.11}}$      & 6   & 11       & 2     & 4    & 2  & 1   & 8   & 3  & $Q_{3}$   \\
            $\mathrm{X_{1.14}}$      & 8   & 15      & 2     & 3    & 1   & 1  & 6   & 0  & $B_{5}$   \\
            \bottomrule
        \end{tabular}
        \caption{part of Takagi's tables (Table 1 and $\mathrm{f'}$ is of $(2,1)$ type and $n\neq0$)}
    \end{table}
\subsubsection{Geometry of hourglasses}\label{Geometryofhourglasses}
A general hourglass is characterized by the following features:
\begin{enumerate}
  \item a $\mathrm{\frac{(1^{3})}{2}}$ type singular point $\mathrm{p}$,
  \item a general smooth element in $|-2\mathrm{K}_{\mathrm{X}}|$, which serves as a base of the hourglass, we denote it (and its strict transformation in each point on the birational links) as $\mathrm{B}$,
  \item $e$ lines passing through the marked singular point $\mathrm{p}$, we denote them by $L_{1},...,L_{e}$,\\
and after blowing up the singular point, we get an exceptional plane,
  \item we denote this exceptional plane (and its strict transformation in each point of the birational links) as $\mathrm{A}$.
\end{enumerate}

We have a clear description of these SQMs in (\ref{sqm}) if we only consider just hourglasses.

\begin{remark}
Unless otherwise stated, we work under the generality conditions for which
all flops under consideration are standard, with normal bundle
$\mathcal{O}_{\mathds{P}^{1}}(-1)^{\oplus2}$. In the hourglass construction,
we use the transverse-intersection hypothesis in
Proposition~\ref{prop:flopstand}. We expect the $(-2,0)$ case to admit
an analogous treatment, but do not give that argument in this article.
\end{remark}

\begin{proposition}\label{prop:flopstand}
Suppose that the strict transforms $\mathrm{A}$ and $\mathrm{B}$
meet transversely and cut out the reduced flopped locus $F$.
Then all flops in (\ref{sqm}) are standard flops.
\end{proposition}

\begin{proof}
By \cite[Proposition 5.1]{Tak1}, the normal bundle of an irreducible
flopping curve is either $\mathcal{O}(-1)^{\oplus2}$ or
$\mathcal{O}\oplus\mathcal{O}(-2)$. We use the normal-sheaf
sequence on the flopped side to distinguish these possibilities.
By \cite[Claim 4.1]{Tak1}, $\mathrm{A}F_k<0$ for each component $F_k$.
Since $\mathrm{B}\sim-2\mathrm{K}_{\mathrm{Y}^{1}}+\mathrm{A}$ and
$\mathrm{K}_{\mathrm{Y}^{1}}F_k=0$, we also have
$\mathrm{B}F_k=\mathrm{A}F_k<0$.
The transverse intersection in the hypothesis identifies
$\mathcal{N}_{F_k/\mathrm{A}}$ with
$\mathcal{O}_{F_k}(\mathrm{B})$ and gives
\begin{equation}\label{exactnormal}
0\longrightarrow\mathcal{O}_{F_k}(\mathrm{B})
\longrightarrow\mathcal{N}_{F_k/\mathrm{Y}^{1}}
\longrightarrow\mathcal{O}_{F_k}(\mathrm{A})\longrightarrow0.
\end{equation}
Both outer terms have the same negative integral degree
$a_k:=\mathrm{A}F_k=\mathrm{B}F_k$.
By adjunction,
\[
\det\mathcal{N}_{F_k/\mathrm{Y}^{1}}
\simeq\mathcal{O}_{\mathds{P}^{1}}(-2),\qquad 2a_k=-2.
\]
Thus $a_k=-1$. Since $\mathrm{H}^{1}(\mathds{P}^{1},\mathcal{O})=0$,
the extension \eqref{exactnormal} splits and
$\mathcal{N}_{F_k/\mathrm{Y}^{1}}\simeq\mathcal{O}(-1)^{\oplus2}$.
The transverse intersection is the geometric condition required here;
it does not follow merely from the two negative intersection numbers.
\end{proof}

\begin{corollary}\label{cor:syzygyofflocus}
Under the transverse-intersection hypothesis of Proposition~\ref{prop:flopstand},
the flopped locus $F$ of each flop in $\mathrm{SQM}$s of (\ref{sqm}) admits a syzygy:
        $$0\longrightarrow\mathcal{O}_{\mathrm{Y}^{1}}(-\mathrm{A}-\mathrm{B})\longrightarrow\mathcal{O}_{\mathrm{Y}^{1}}(-\mathrm{A})\oplus\mathcal{O}_{\mathrm{Y}^{1}}(-\mathrm{B})\longrightarrow \mathcal{O}_{\mathrm{Y}^{1}}\longrightarrow  \mathcal{O}_{F}\longrightarrow 0$$
        where $\mathrm{A}$ is the strict transformation of a plane, $\mathrm{B}$ is the strict transformation of a base of the hourglass.

\end{corollary}
\begin{proof}
Since we already know $\mathrm{AB}=\sum_{i}F_{i}=F$ and the Koszul resolution given by the intersection equation of divisors $\mathrm{A}$ and $\mathrm{B}$ gives this syzygy.
\end{proof}

\begin{proposition}\label{prop:flipstandard}
Each flip in $\mathrm{SQM}$s of (\ref{sqm}) is a $(1,2:-1,-1)$ type Francia's flip.
\end{proposition}
\begin{proof}
See \cite[Proposition 0.8]{Tak1}.
\end{proof}

\begin{theorem}[Hourglass construction I]\label{thm:hourglassconstrucion}
Fix a row of Takagi's Table 1, as listed above, and let $\mathrm{X}'$,
$n,z,\deg C$ and $g(C)$ be the corresponding data. Put
$\mathrm{H}:=(-\mathrm{K}_{\mathrm{X}'})/(z+1)$. Assume that:
\begin{enumerate}
\item $C\subset\mathrm{X}'$ is smooth with the prescribed degree and
genus. The $n$ lines $l_i$ are mutually disjoint, lie in
$\operatorname{Reg}\mathrm{X}'$, and are $(z+2)$-secant to $C$ with
respect to $\mathrm{H}$.
\item The base locus is
\[
\operatorname{Bs}|-\mathrm{K}_{\mathrm{X}'}-C|
=C\cup\bigcup_{i=1}^{n}l_i\cup
\operatorname{Sing}_{1/2}(\mathrm{X}'),
\]
where $\operatorname{Sing}_{1/2}(\mathrm{X}')$ denotes the
$\frac12(1,1,1)$ points. This is an equality of base loci as sets.
\item There is a surface $\mathrm{S}\subset\mathrm{X}'$ containing $C$
with $\mathrm{S}\equiv z\mathrm{H}$.
\end{enumerate}
For the orbifold hourglasses considered here, restrict to the data for
which the final contraction produces a $\frac12(1,1,1)$ point, rather
than a QODP, and impose the standard-flop generality convention above.
Then the following operations reconstruct an hourglass in the same row:
\begin{enumerate}
\item Blow up $\mathrm{X}'$ along $C$, and then along the strict
transforms of the lines $l_i$.
\item Perform the standard flops along the $e$ rational curves in the
strict transform of $\mathrm{S}$ and the $n$ relevant sections in the
exceptional divisors over the lines.
\item Contract the strict transforms of $\mathrm{S}$ and those exceptional
divisors to points.
\end{enumerate}
Conversely, the hourglasses under consideration arise from these data.
\end{theorem}
\begin{proof}
Apply \cite[Proposition 0.8 and Theorem 0.10(A), (B)]{Tak1}.
The blow-ups along the transformed lines resolve the inverse Francia
flips and give the three-step presentation above. Its hypotheses give the
irreducibility of $\mathrm{S}$ and $C\not\subset\operatorname{Sing}\mathrm{S}$;
smoothness of $\mathrm{S}$ is not an additional assumption here.
The restriction on the last contraction excludes the QODP possibility
in \cite[Remark 0.9]{Tak1}. Standardness is understood under the
generality convention preceding Proposition~\ref{prop:flopstand}.
\end{proof}

\begin{remark}
   If we know the first midpoint $\mathrm{Z}$ of the Sarkisov link, we directly construct the (Weil but not Cartier) divisor corresponding to $\mathrm{A}$ and $\mathrm{B}$ in the projective bundle $\mathds{P}_{\mathrm{Z}}(\mathcal{O}(\mathrm{A})\oplus\mathcal{O}(\mathrm{B}))$ noticing even in $\mathrm{A}$ and $\mathrm{B}$  are not Cartier but their difference is, and then shrink the divisor $\mathrm{A}$. Hence $\mathrm{X}$ has $\mathrm{Cone(v^{(-2\mathrm{K}_{\mathrm{Z}})}\mathrm{Z})}$ as an ambient space.\\

   The idea of the proof is just using the result in  Corollary $\ref{cor:syzygyofflocus}$, consider the surjective map
   $$0\longrightarrow\mathcal{O}_{\mathrm{Y}^{1}}(-\mathrm{A}-\mathrm{B})\longrightarrow\mathcal{O}_{\mathrm{Y}^{1}}(-\mathrm{A})\oplus\mathcal{O}_{\mathrm{Y}^{1}}(-\mathrm{B})\longrightarrow \mathcal{I}_{F/\mathrm{Y}^{1}}\longrightarrow0$$
   which induces a closed immersion of $Bl_{F}\mathrm{Y}^{1}$ to $\mathds{P}_{\mathrm{Y}^{1}}(\mathcal{O}_{\mathrm{Y}^{1}}(\mathrm{A})\oplus\mathcal{O}_{\mathrm{Y}^{1}}(\mathrm{B}))$ as the section corresponding to the first map of above exact sequence, then we blow down $Bl_{F}\mathrm{Y}^{1}$ along another ray, correspondingly $\mathrm{Y}^{1}$ shrinks to $\mathrm{Z}$, we can see there is also an embedding of $\mathrm{Y}$ to $\mathds{P}_{\mathrm{Z}}(\mathcal{O}_{\mathrm{Z}}(\mathrm{A})\oplus\mathcal{O}_{\mathrm{Z}}(\mathrm{B}))$ as a tautological section. And we know $\mathrm{B}-\mathrm{A}\sim -2\mathrm{K}_{\mathrm{Z}}$, its a Veronese cone after contraction a section of $\mathds{P}_{\mathrm{Z}}(\mathcal{O}_{\mathrm{Z}}(\mathrm{A})\oplus\mathcal{O}_{\mathrm{Z}}(\mathrm{B}))$.
\end{remark}

\subsection{Noncommutative hourglasses}

In this section, we first review the derived category classification for smooth Fano $3$-folds, we begin with a famous theorem by D.Orlov.

\begin{theorem}[{\cite[Thm. 2.5]{Orlh}}]\label{thm:orlov1}
Let $A$ be a connected graded Noetherian algebra that is Gorenstein with Gorenstein parameter $a>0$, then we have the following semi-orthogonal decomposition:
$$\mathrm{D}^{b}(\mathrm{qgr}\,A)\simeq\langle A(-a+1),..,A(0), \mathrm{D}^{\mathrm{gr}}_{\mathrm{sg}}(A)\rangle$$
where $\mathrm{qgr}\,A$ is the $\mathrm{gr}\,A$ quotient a torsion subcategory and $\mathrm{D}^{\mathrm{gr}}_{\mathrm{sg}}(A)$ is the singularity category of $\mathrm{D}^{b}(\mathrm{gr}\,A)$.
\end{theorem}

In later work, D.Orlov and others proved that if $A$ is also a complete intersected graded algebra, then $\mathrm{D}^{\mathrm{gr}}_{\mathrm{sg}}(A)$ is equivalent to a (higher) matrix factorization category. So according to the convention of this article, we call the left orthogonal component of
$$\langle A(-a+1),..,A(0)\rangle$$
in $\mathrm{D}^{b}(\mathrm{qgr}\,A)$ as \textit{matrix factorization category} of $A$ and denote it by $\mathcal{MF}(A)$, if $A$ is a complete intersected Gorenstein connected graded Noetherian algebra with Gorenstein parameter $a>0$.\\

A simple corollary of the above theorem is the derived category decomposition of a completely intersected substack in a weighted projective stack. If
$\mathcal{X}\subset\mathcal{P}:=\mathcal{P}(a_{0},...,a_{n})$ is a complete intersection of $m$ hypersurface $H_{1},...,H_{m}$ of grade $d_{1}$,...,$d_{m}$, then $$A_{\mathcal{X}}:=\bigoplus_{i\geq 0}\mathrm{H}^{0}(\mathcal{X},\mathcal{O}_{\mathcal{X}}(i))$$
is a graded complete intersected Gorenstein algebra with Gorenstein parameter $a:=\sum_{i=0}^{i=n}a_{i}-\sum_{i=1}^{i=m}d_{i}$, we assume it's a positive number which is equivalent to $-\mathrm{K}_{\mathcal{X}}>0$, i.e. $\mathcal{X}$ is Fano.\\

On the other hand, we have an equivalence
$$\mathrm{qgr}\,A_{\mathcal{X}}\simeq \mathrm{Coh} (\mathcal{X})$$
by (\ref{equ:serreequa}), and easy to see this equivalence sends $A_{\mathcal{X}}(i)$ to $\mathcal{O}_{\mathcal{X}}(i)$ for any $i$. So we have the following corollary:
\begin{corollary}
If $\mathcal{X}$ is the Gorenstein connected Fano stack constructed above, we have
$$\mathrm{D}^{b}(\mathrm{Coh} (\mathcal{X}))\simeq\langle \mathcal{O}_{\mathcal{X}}(-a+1),..,\mathcal{O}_{\mathcal{X}}, \mathcal{MF}(\mathcal{X})\rangle$$
where $a:=\sum_{i=0}^{i=n}a_{i}-\sum_{i=1}^{i=m}d_{i}$ and $\mathcal{MF}(\mathcal{X}):= \mathcal{MF}(A_{\mathcal{X}})$.
\end{corollary}

According to the classification of $\mathrm{X}'_{1,i}$ in (\ref{targetX}), we find that except for the case of $B_{5}$ they all meet the conditions of the above corollary, so we have

\begin{corollary}\label{cor:mfx'}
If $\mathrm{\widetilde{X}}'_{1,i}$ is the canonical stack of $\mathrm{X}'_{1,i}$ where $i\neq 13,14$, we have
$$\mathrm{D}^{b}(\mathrm{\widetilde{X}}'_{1,i})\simeq\langle\mathcal{MF}(\mathrm{\widetilde{X}}'_{1,i}), \mathrm{Exc}_{i},\mathcal{O}_{\mathrm{\widetilde{X}}'_{1,i}}(-1),\mathcal{O}_{\mathrm{\widetilde{X}}'_{1,i}}\rangle$$
where $\mathrm{Exc}_{i}:=\langle \mathcal{O}_{\mathrm{\widetilde{X}}'_{1,i}}(-F_{w}(\mathrm{X}'_{1,i})+1),..,\mathcal{O}_{\mathrm{\widetilde{X}}'_{1,i}}(-2)\rangle$ is an exceptional collection\footnote{\,Allowed to be empty.} and we mutate the matrix factorization category $\mathcal{MF}(\mathrm{\widetilde{X}}'_{1,i})$ to the far left.
\end{corollary}

For the remaining examples, we need to deal with $B_{5}$ which is a smooth quintic del Pezzo $3$-fold, although the conclusion is already very mature, we give a few more inspiring ideas to illustrate that these remaining examples are not special:
\begin{proposition}\label{prop:dcb5}
If $\mathrm{\widetilde{X}}'_{1,i}$ is the canonical stack of $\mathrm{X}'_{1,i}$ where $i= 13,14$, we have
$$\mathrm{D}^{b}(\mathrm{\widetilde{X}}'_{1,i})\simeq\langle \mathcal{MF}(\mathrm{\widetilde{X}}'_{1,i}), \mathrm{Exc}_{i},\mathcal{O}_{\mathrm{\widetilde{X}}'_{1,i}}(-1),\mathcal{O}_{\mathrm{\widetilde{X}}'_{1,i}}\rangle$$
where $\mathrm{Exc}_{i}$ is empty and $\mathcal{MF}(\mathrm{\widetilde{X}}'_{1,i})$ denotes the residual block $\langle\mathcal{Q}_{B_{5}}(-2),\mathcal{U}_{B_{5}}(-1)\rangle$ in this case.
\end{proposition}
\begin{proof}
\begin{enumerate}
  \item\label{B5mua} For the exceptionality of pair $\mathcal{O}_{B_{5}}(-1)$, $\mathcal{O}_{B_{5}}$, we can easily see from the famous Sarkisov link by V.lskovskikh:
    $$\begin{tikzcd}[column sep=0.5em]
  & && Bl_L{B_{5}}= Bl_{C}Q_{3}\arrow[dl]\arrow[dr] &&  \\
 & &B_{5}&  &Q_{3}&
\end{tikzcd}$$
where $L$ is a line on $B_{5}$ and $C$ is a cubic rational curve on $Q_{3}$, and further we can see this Sarkisov link can give all the exception objects we want by mere mutation (and using some simple divisor relationships).
  \item Similarly, we are expected to consider an extended Sarkisov link like \cite[Conjecture 4.12]{KaFf}.
  \item\label{B5mua1} More directly, we can calculate this exceptional collection through the sheaves on the key variety of $B_{5}$ as \cite{OrB5}, we have
  \begin{equation}\label{equ:b5sod1}\mathrm{D}^{b}(B_{5})\simeq\langle \mathcal{Q}_{B_{5}}(-2),\mathcal{U}_{B_{5}}(-1),\mathcal{O}_{B_{5}}(-1),\mathcal{O}_{\mathrm{B_{5}}}\rangle\end{equation}
since $B_{5}\subset\mathbf{Gr}(2,5)$.
\item Slightly changing (\ref{equ:b5sod1}), we have
  \begin{equation}\label{equ:b5sod2}\mathrm{D}^{b}(B_{5})\simeq\langle \mathcal{U}_{B_{5}}(-1),\mathcal{O}_{B_{5}}(-1),\mathcal{U}_{B_{5}},\mathcal{O}_{\mathrm{B_{5}}}\rangle\end{equation}
\item By comparing (\ref{equ:b5sod1}) with the mutation construction in (\ref{B5mua}), we can easily check that this exceptional collection comes from the birational correlation of
$$\mathrm{D}^{b}(Q_{3})=\langle\mathcal{MF}(Q_{3}), \mathcal{O}_{Q_{3}}(-2),\mathcal{O}_{Q_{3}}(-1),\mathcal{O}_{\mathrm{Q_{3}}}\rangle$$
by modifying one element $\mathcal{O}_{Q_{3}}(-2)$ with $\mathcal{U}_{B_{5}}(-1)$  and retaining the other parts in an equivalent sense.
\end{enumerate}
\end{proof}

Now that we have all the necessary preparations, our goal is as follows.

\begin{theorem}[Semi-orthogonal decomposition for hourglasses]\label{thm:SODh}
For any $i$, the canonical stack $\widetilde{\mathrm{X}}_{1,i}$ of hourglass $\mathrm{X}_{1,i}$ admits semi-orthogonal decomposition as following:
        $$\mathrm{D}^{b}(\widetilde{\mathrm{X}}_{1,i})\cong \langle\,  {\mathcal{R}}es_{1,i},\mathrm{Exc}_{1,i} \,\rangle$$
        where $\mathrm{Exc}_{1,i}$ is an exceptional collection,
        and the residue part ${\mathcal{R}}es_{1,i} $ is equivalent to a span of $\mathrm{D}^{b}(\mathrm{C_{1,i}} )$ and the matrix factorization category $\mathcal{MF}(\mathrm{\widetilde{X}}'_{1,i})$ as in Corollary \ref{cor:mfx'}, that is
        $${\mathcal{R}}es_{1,i}\simeq\langle\mathrm{D}^{b}(\mathrm{C_{1,i}} ),\mathcal{MF}(\mathrm{\widetilde{X}}'_{1,i})\rangle $$
Specifically, if $g(\mathrm{C_{1,i}} )=0$ and the matrix factorization category is degenerated, the noncommutative hourglass admits a full exceptional collection.
\end{theorem}
\begin{proof}
Fix $i$ and omit the subscripts $1,i$. We retain the notation
$\widetilde{\mathrm{X}}$, $\widetilde{\mathrm{Y}}$,
$\widetilde{\mathrm{Y}}^{1}$ and $\widetilde{\mathrm{Y}}'$ for the canonical
stacks in the hourglass diagram. Recall that
$\mathrm{\Theta}=\varpi_*\pi^*$ is the fully-faithful functor of the Kawamata
blow-up, with right adjoint $\mathrm{\Theta}^{!}=\pi_*\varpi^{!}$.
Write $\mathrm{A}$ for its exceptional divisor and its strict transform on $\widetilde{\mathrm{Y}}^{1}$,
and $\mathrm{A}^+$ for the strict transform on $\widetilde{\mathrm{Y}}'$. All functors are derived.
We retain cohomological shifts in object calculations and suppress them
only when writing the generated components of a semi-orthogonal decomposition.

\begin{enumerate}
\item The Kawamata formula gives
\begin{equation}\label{thsod1}
\mathrm{D}^{b}(\widetilde{\mathrm{Y}})=\langle\mathcal{O}_{\mathrm{A}}(-1),\mathrm{\Theta}\,\mathrm{D}^{b}(\widetilde{\mathrm{X}})\rangle.
\end{equation}
Since $\mathcal{O}_{\mathrm{A}}(-\mathrm{K}_{\widetilde{\mathrm{Y}}})=\mathcal{O}_{\mathrm{A}}(1)$, rotating the first component
to the end gives
\begin{equation}\label{thsod2}
\mathrm{D}^{b}(\widetilde{\mathrm{Y}})=\langle\mathrm{\Theta}\,\mathrm{D}^{b}(\widetilde{\mathrm{X}}),\mathcal{O}_{\mathrm{A}}\rangle.
\end{equation}
Here the rotation is the inverse Serre functor; its shift has no effect
on the generated component.

\item\label{mg:step3} Apply the standard-flop equivalence
$\varphi_*\phi^*$ associated with the same correspondence as above:
\[
\begin{tikzcd}[column sep=small]
&\widetilde{\mathrm{W}}\arrow[dl,"\phi"']\arrow[dr,"\varphi"]&\\
\widetilde{\mathrm{Y}}&&\widetilde{\mathrm{Y}}^{1}.
\end{tikzcd}
\]
The divisor sequence gives
$\varphi_*\phi^*\mathcal{O}_{\mathrm{A}}\simeq\mathcal{O}_{\mathrm{A}}$; also
$\varphi_*\phi^*\mathrm{\Theta}\mathcal{O}_{\widetilde{\mathrm{X}}}\simeq\mathcal{O}_{\widetilde{\mathrm{Y}}^{1}}$.
Thus \eqref{thsod2} becomes
\begin{equation}\label{thsod3}
\mathrm{D}^{b}(\widetilde{\mathrm{Y}}^{1})=
\left\langle\varphi_*\phi^*\mathrm{\Theta}\,\mathrm{D}^{b}(\widetilde{\mathrm{X}}),\mathcal{O}_{\mathrm{A}}\right\rangle.
\end{equation}
We compare this decomposition with the one obtained from $\widetilde{\mathrm{X}}'$.

\item The blow-up formula for $f'$ and the chosen decomposition on
$\widetilde{\mathrm{X}}'$ give
\begin{equation}\label{XSOD4}
\mathrm{D}^{b}(\widetilde{\mathrm Y}')=\big\langle \iota'_*(f_{\mathrm E}'^*\mathrm{D}^{b}(\mathrm C)(-1)), f'^*\mathcal{MF}(\widetilde{\mathrm X}'), f'^*\mathrm{Exc}_i, \mathcal{O}(-\mathrm H'), \mathcal{O}\big\rangle.
\end{equation}
Here $(-1)$ inside $\iota'_*$ means tensoring on $\mathrm E'$ by
$\mathcal{O}_{\mathrm E'}(-1)$. In the formulas below a common functor is
applied to the explicitly displayed ordered span of the first three
components, including all Hom complexes between them.
Right-mutate the first three components through $\mathcal{O}(-\mathrm{H}')$,
retaining their order:
\begin{equation}\label{XSOD5}
\mathrm{D}^{b}(\widetilde{\mathrm Y}')=\big\langle \mathcal{O}(-\mathrm H'), \mathbb{R}_{\mathcal{O}(-\mathrm H')}\langle \iota'_*(f_{\mathrm E}'^*\mathrm{D}^{b}(\mathrm C)(-1)), f'^*\mathcal{MF}(\widetilde{\mathrm X}'), f'^*\mathrm{Exc}_i\rangle, \mathcal{O}\big\rangle.
\end{equation}
Rotate the first component to the end. Since
$-\mathrm{K}_{\widetilde{\mathrm{Y}}'}-\mathrm{H}'=\mathrm{A}^+$, its line-bundle class becomes
$\mathcal{O}(\mathrm{A}^+)$. Thus
\begin{equation}\label{XSOD6}
\mathrm{D}^{b}(\widetilde{\mathrm Y}')=\big\langle \mathbb{R}_{\mathcal{O}(-\mathrm H')}\langle \iota'_*(f_{\mathrm E}'^*\mathrm{D}^{b}(\mathrm C)(-1)), f'^*\mathcal{MF}(\widetilde{\mathrm X}'), f'^*\mathrm{Exc}_i\rangle, \mathcal{O}, \mathcal{O}(\mathrm A^+)\big\rangle.
\end{equation}
\item Tensor \eqref{XSOD6} by $\mathcal{O}(-\mathrm{A}^+)$. This gives
\begin{equation}\label{XSOD7}
\mathrm{D}^{b}(\widetilde{\mathrm Y}')=\big\langle \bigl(\mathbb{R}_{\mathcal{O}(-\mathrm H')}\langle \iota'_*(f_{\mathrm E}'^*\mathrm{D}^{b}(\mathrm C)(-1)), f'^*\mathcal{MF}(\widetilde{\mathrm X}'), f'^*\mathrm{Exc}_i\rangle\bigr)(-\mathrm A^+), \mathcal{O}(-\mathrm A^+), \mathcal{O}\big\rangle.
\end{equation}
Right mutation and tensoring restrict to equivalences on these three
components and preserve all Hom complexes between them.

\item Let $\Gamma_1,\ldots,\Gamma_n\subset \widetilde{\mathrm{Y}}^{1}$ be the disjoint Francia
flipping curves and put $\Gamma=\coprod_j\Gamma_j$. Thus
$\mathcal{O}_\Gamma(-\mathrm A)=\bigoplus_j\mathcal{O}_{\Gamma_j}(-\mathrm A)$
and
\[
\langle\mathcal{O}_\Gamma(-\mathrm A)\rangle
=\langle\mathcal{O}_{\Gamma_1}(-\mathrm A),\ldots,
\mathcal{O}_{\Gamma_n}(-\mathrm A)\rangle.
\]
An entry $\mathcal{O}_\Gamma(-\mathrm A)$ in a decomposition denotes this
span of its $n$ direct summands. The same convention applies to its
images under the displayed additive functors. Each
$\Gamma_j\simeq\mathcal P(1,2)$ has normal bundle
$\mathcal{O}_{\Gamma_j}(-1)^{\oplus2}$ in weighted notation. The divisor
$\mathrm{A}$ meets it at one ordinary point, so
$\mathcal{O}_{\Gamma_j}(-\mathrm{A})=\mathcal{O}_{\mathcal P(1,2)}(-2)$.
This is coarse degree $-1$, not weighted degree $-1$.
On the common weighted resolution write
$\phi':\widetilde{\mathrm{V}}\to \widetilde{\mathrm{Y}}^{1}$ and $\varphi':\widetilde{\mathrm{V}}\to \widetilde{\mathrm{Y}}'$, and let $\mathrm{E}_{\mathrm{V}}$ be the sum
of its exceptional divisors. Thus $\phi'$ is the
ordinary blow-up of the smooth stacky curves $\Gamma_j$, and $\varphi'$ is the
weighted blow-up of weights $(1,2)$ on the other side. Use the functor
\[
\Phi_1(-)=\phi'_*\bigl(\varphi'^*(-)\otimes\mathcal{O}_{\widetilde{\mathrm{V}}}(\mathrm{E}_{\mathrm{V}})\bigr).
\]
The weighted flip formula gives the decomposition
\[
\mathrm{D}^{b}(\widetilde{\mathrm{Y}}^{1})=
\left\langle\mathcal{O}_{\Gamma_1}(-\mathrm{A}),\ldots,
\mathcal{O}_{\Gamma_n}(-\mathrm{A}),
\Phi_1\mathrm{D}^{b}(\widetilde{\mathrm{Y}}')\right\rangle.
\]
Indeed, the untwisted functor $\phi'_*\varphi'^*$ is the one in
Theorem~\ref{thm:fraciaflip}. The canonical divisor identities are
\[
\mathrm{K}_{\widetilde{\mathrm{V}}}
=\phi'^*\mathrm{K}_{\widetilde{\mathrm{Y}}^{1}}+\mathrm{E}_{\mathrm{V}}
=\varphi'^*\mathrm{K}_{\widetilde{\mathrm{Y}}'}+2\mathrm{E}_{\mathrm{V}}.
\]
Consequently
\[
\Phi_1=(-\otimes\omega_{\widetilde{\mathrm{Y}}^{1}})
\circ\phi'_*\varphi'^*\circ
(-\otimes\omega_{\widetilde{\mathrm{Y}}'}^{-1}).
\]
Full faithfulness is preserved by these tensor autoequivalences. Since
$\omega_{\widetilde{\mathrm{Y}}^{1}}|_{\Gamma_j}
=\mathcal{O}_{\mathcal P(1,2)}(-1)$, the exceptional objects become
$\mathcal{O}_{\mathcal P(1,2)}(-2)=\mathcal{O}_{\Gamma_j}(-\mathrm{A})$.

In local Cox coordinates the common correspondence is
\[
\phi'(\tau,s,t,u,v)=(s,t;\tau u,\tau v),\qquad
\varphi'(\tau,s,t,u,v)=(u,v;\tau s,\tau^2t).
\]
At its ordinary intersection with the flipping curve, $\mathrm{A}$ has
equation $t+\alpha s^2=0$. The target equation pulls back to
$\tau^2(t+\alpha s^2)$, giving
\[
\varphi'^*\mathrm{A}^+=\phi'^*\mathrm{A}+2\mathrm{E}_{\mathrm{V}},\qquad
\phi'_*\mathcal{O}(\mathrm{E}_{\mathrm{V}})=\mathcal{O},\qquad
\phi'_*\mathcal{O}(-\mathrm{E}_{\mathrm{V}})=\mathcal I_\Gamma.
\]
It follows that
\[
\Phi_1\mathcal{O}=\mathcal{O},\qquad
\Phi_1\mathcal{O}(-\mathrm{A}^+)=\mathcal I_\Gamma(-\mathrm{A}).
\]
Applying $\Phi_1$ to \eqref{XSOD7} therefore gives
\begingroup\small
\begin{equation}\label{XSOD8}
\mathrm{D}^{b}(\widetilde{\mathrm Y}^{1})=\big\langle \mathcal{O}_\Gamma(-\mathrm A), \Phi_1\bigl(\bigl(\mathbb{R}_{\mathcal{O}(-\mathrm H')}\langle \iota'_*(f_{\mathrm E}'^*\mathrm{D}^{b}(\mathrm C)(-1)), f'^*\mathcal{MF}(\widetilde{\mathrm X}'), f'^*\mathrm{Exc}_i\rangle\bigr)(-\mathrm A^+)\bigr), \mathcal I_\Gamma(-\mathrm A), \mathcal{O}\big\rangle.
\end{equation}
\endgroup
If $n=0$, the curve entries are absent, $\Phi_1$ is the identity,
and $\mathcal I_\Gamma(-\mathrm{A})=\mathcal{O}(-\mathrm{A})$.

\item Move each of the three displayed images through
$\langle\mathcal{O}_\Gamma(-\mathrm A)\rangle$ by left mutation.
The resulting semi-orthogonal decomposition of
$\mathrm{D}^{b}(\widetilde{\mathrm Y}^{1})$ is
\begingroup\small
\begin{equation}\label{XSOD9}
\big\langle \mathbb{L}_{\langle\mathcal{O}_\Gamma(-\mathrm A)\rangle}\Phi_1\bigl(\bigl(\mathbb{R}_{\mathcal{O}(-\mathrm H')}\langle \iota'_*(f_{\mathrm E}'^*\mathrm{D}^{b}(\mathrm C)(-1)), f'^*\mathcal{MF}(\widetilde{\mathrm X}'), f'^*\mathrm{Exc}_i\rangle\bigr)(-\mathrm A^+)\bigr), \mathcal{O}_\Gamma(-\mathrm A), \mathcal I_\Gamma(-\mathrm A), \mathcal{O}\big\rangle.
\end{equation}
\endgroup
The left mutation is through the span generated by $\mathcal{O}_\Gamma(-\mathrm A)$.
Since the curves are disjoint, these sheaves are completely orthogonal;
the mutation is the composite of the individual left mutations in any
fixed order. For $n=0$ it is the identity. Its restriction to the three
components above is an equivalence and preserves their internal
semi-orthogonal decomposition.

The next mutations follow from the ideal sequence
\[
0\longrightarrow\mathcal I_\Gamma(-\mathrm{A})
\longrightarrow\mathcal{O}(-\mathrm{A})
\longrightarrow\bigoplus_{j=1}^n\mathcal{O}_{\Gamma_j}(-\mathrm{A})
\longrightarrow0.
\]
The weighted normal determinant and the required vanishing are
\[
\det\mathcal{N}_{\Gamma_j/\widetilde{\mathrm{Y}}^{1}}
=\mathcal{O}_{\Gamma_j}(-2),\qquad
\mathrm{H}^*(\mathcal P(1,2),\mathcal{O}(-2))=0.
\]
Codimension-two duality gives
\[
\begin{aligned}
\operatorname{RHom}(\mathcal{O}_{\Gamma_j}(-\mathrm{A}),
                         \mathcal{O}(-\mathrm{A}))=0,\quad
\operatorname{RHom}(\mathcal{O}_{\Gamma_j}(-\mathrm{A}),
                         \mathcal I_\Gamma(-\mathrm{A}))=\mathrm k[-1].
\end{aligned}
\]
The connecting map of the ideal sequence is the mutation evaluation map.
Hence
\[
\mathbb{L}_{\langle\mathcal{O}_{\Gamma_j}(-\mathrm{A})\rangle_{j=1}^n}
\mathcal I_\Gamma(-\mathrm{A})=\mathcal{O}(-\mathrm{A}).
\]
On the last components, the two successive mutations read
\[
\langle\mathcal{O}_\Gamma(-\mathrm A),\mathcal I_\Gamma(-\mathrm A),\mathcal{O}\rangle
=\langle\mathcal{O}(-\mathrm A),\mathcal{O}_\Gamma(-\mathrm A),\mathcal{O}\rangle
=\langle\mathcal{O}(-\mathrm A),\mathcal{O},
\mathbb{R}_{\mathcal{O}}\mathcal{O}_\Gamma(-\mathrm A)\rangle.
\]
For each $j$ the right-mutated object is specified by the triangle
\[
\mathbb{R}_{\mathcal{O}}\mathcal{O}_{\Gamma_j}(-\mathrm A)
\longrightarrow\mathcal{O}_{\Gamma_j}(-\mathrm A)
\longrightarrow\mathcal{O}[2]
\longrightarrow\bigl(\mathbb{R}_{\mathcal{O}}\mathcal{O}_{\Gamma_j}(-\mathrm A)\bigr)[1],
\]
because
$\operatorname{RHom}(\mathcal{O}_{\Gamma_j}(-\mathrm{A}),\mathcal{O})
=\mathrm k[-2]$.

\item\label{mutacohcomp} The remaining line-bundle mutation is computed
on $\widetilde{\mathrm{Y}}^{1}$, not on $\widetilde{\mathrm{Y}}'$. Here $\mathrm{A}$ is a rational surface obtained by blowing
up the plane in the first flops, with
$\mathcal{O}_{\mathrm{A}}(\mathrm{A})=\mathcal{O}_{\mathrm{A}}(-2\mathrm{h}+\sum \mathrm{e}_j)$, where the $\mathrm{e}_j$ are its
disjoint exceptional curves. The sequences adding each $\mathrm{e}_j$ have quotient
$\mathcal{O}_{\mathrm{e}_j}(-1)$, so
\[
\mathrm{H}^*(\mathrm{A},\mathcal{O}_{\mathrm{A}}(\mathrm{A}))
=\mathrm{H}^*(\mathds P^2,\mathcal{O}(-2))=0.
\]
Also $\mathrm{H}^*(\widetilde{\mathrm{Y}}^{1},\mathcal{O})=\mathrm k$: it follows from Fano
vanishing on $\widetilde{\mathrm{X}}$, the Kawamata extraction and the standard flops.
The divisor sequence therefore gives
$\operatorname{RHom}(\mathcal{O}(-\mathrm{A}),\mathcal{O})=\mathrm k$.
Its nonzero map is the divisor equation, whence
\[
\mathbb{R}_{\mathcal{O}}\mathcal{O}(-\mathrm{A})=\mathcal{O}_{\mathrm{A}}[-1].
\]
Combining the three mutations gives the following decomposition of $\mathrm{D}^{b}(\widetilde{\mathrm Y}^{1})$:
\begingroup\small
\begin{equation}\label{XSOD9'}
\big\langle \mathbb{L}_{\langle\mathcal{O}_\Gamma(-\mathrm A)\rangle}\Phi_1\bigl(\bigl(\mathbb{R}_{\mathcal{O}(-\mathrm H')}\langle \iota'_*(f_{\mathrm E}'^*\mathrm{D}^{b}(\mathrm C)(-1)), f'^*\mathcal{MF}(\widetilde{\mathrm X}'), f'^*\mathrm{Exc}_i\rangle\bigr)(-\mathrm A^+)\bigr), \mathcal{O}, \mathcal{O}_{\mathrm A}, \mathbb{R}_{\mathcal{O}}\mathcal{O}_\Gamma(-\mathrm A)\big\rangle.
\end{equation}
\endgroup
The shift $[-1]$ of $\mathcal{O}_{\mathrm{A}}$ is suppressed only in the
generated component. For $n=0$, $\widetilde{\mathrm Y}'=
\widetilde{\mathrm Y}^{1}$, and this is the cohomological calculation
used in the cubic example below.

\item Rotate the final span of $n$ direct summands to the front. On the threefold
$\widetilde{\mathrm Y}^{1}$, Serre rotation is tensoring by
$\omega_{\widetilde{\mathrm Y}^{1}}$ followed by $[3]$. Thus
\begingroup\small
\begin{multline}\label{XSOD10}
\mathrm{D}^{b}(\widetilde{\mathrm Y}^{1})=\big\langle (\mathbb{R}_{\mathcal{O}}\mathcal{O}_\Gamma(-\mathrm A))\otimes\omega_{\widetilde{\mathrm Y}^{1}}[3],\\[2pt]
\mathbb{L}_{\langle\mathcal{O}_\Gamma(-\mathrm A)\rangle}\Phi_1\bigl(\bigl(\mathbb{R}_{\mathcal{O}(-\mathrm H')}\langle \iota'_*(f_{\mathrm E}'^*\mathrm{D}^{b}(\mathrm C)(-1)), f'^*\mathcal{MF}(\widetilde{\mathrm X}'), f'^*\mathrm{Exc}_i\rangle\bigr)(-\mathrm A^+)\bigr), \mathcal{O}, \mathcal{O}_{\mathrm A}\big\rangle.
\end{multline}
\endgroup
The final component is the same embedded
$\mathcal{O}_{\mathrm{A}}$ as in \eqref{thsod3}. Their complements are
both $\mathcal{O}_{\mathrm{A}}^{\perp}$. Removing that common component
identifies the preceding entries with
$\varphi_*\phi^*\mathrm{\Theta}\,\mathrm{D}^{b}(\widetilde{\mathrm X})$.
On this image, apply its right adjoint
$\mathrm{\Theta}^{!}\phi_*\varphi^{!}$, which is its inverse. We obtain
\begingroup\small
\begin{multline}\label{XSOD10'}
\mathrm{D}^{b}(\widetilde{\mathrm X})=\big\langle \mathrm\Theta^!\phi_*\varphi^!((\mathbb{R}_{\mathcal{O}}\mathcal{O}_\Gamma(-\mathrm A))\otimes\omega_{\widetilde{\mathrm Y}^{1}})[3],\\[2pt]
\mathrm\Theta^!\phi_*\varphi^!\mathbb{L}_{\langle\mathcal{O}_\Gamma(-\mathrm A)\rangle}\Phi_1\bigl(\bigl(\mathbb{R}_{\mathcal{O}(-\mathrm H')}\langle \iota'_*(f_{\mathrm E}'^*\mathrm{D}^{b}(\mathrm C)(-1)), f'^*\mathcal{MF}(\widetilde{\mathrm X}'), f'^*\mathrm{Exc}_i\rangle\bigr)(-\mathrm A^+)\bigr), \mathcal{O}_{\widetilde{\mathrm X}}\big\rangle.
\end{multline}
\endgroup
All three transported components retain their Hom complexes and the
maps between them: each mutation restricts to an equivalence on its
component, and the last right adjoint is inverse on the entire common
complement.

\item\label{mg:step11} Rotate the first $n$ components of
\eqref{XSOD10'} to the end. On $\widetilde{\mathrm X}$ the inverse Serre
functor is tensoring by
$\omega_{\widetilde{\mathrm X}}^{-1}
=\mathcal{O}_{\widetilde{\mathrm X}}(-\mathrm{K}_{\widetilde{\mathrm X}})$
followed by $[-3]$. Its shift cancels the preceding $[3]$, since the
intervening functors are exact. Explicitly, for every $j$ the last object is
\begingroup\small
\begin{equation}\label{XSOD11'}
\begin{gathered}
\mathrm{\Theta}^{!}\phi_*\varphi^{!}
\bigl((\mathbb{R}_{\mathcal{O}}\mathcal{O}_{\Gamma_j}(-\mathrm{A}))
\otimes\omega_{\widetilde{\mathrm Y}^{1}}\bigr)
\otimes\mathcal{O}_{\widetilde{\mathrm X}}(-\mathrm{K}_{\widetilde{\mathrm X}}).
\end{gathered}
\end{equation}
\endgroup
The canonical bundles are tensored on their respective varieties in this
order. In particular, the resulting decomposition is
\begingroup\small
\begin{multline}\label{XSOD11}
\mathrm{D}^{b}(\widetilde{\mathrm X})=\big\langle \mathrm\Theta^!\phi_*\varphi^!\mathbb{L}_{\langle\mathcal{O}_\Gamma(-\mathrm A)\rangle}\Phi_1\bigl(\bigl(\mathbb{R}_{\mathcal{O}(-\mathrm H')}\langle \iota'_*(f_{\mathrm E}'^*\mathrm{D}^{b}(\mathrm C)(-1)), f'^*\mathcal{MF}(\widetilde{\mathrm X}'), f'^*\mathrm{Exc}_i\rangle\bigr)(-\mathrm A^+)\bigr),\\[2pt]
\mathcal{O}_{\widetilde{\mathrm X}}, \mathrm\Theta^!\phi_*\varphi^!((\mathbb{R}_{\mathcal{O}}\mathcal{O}_\Gamma(-\mathrm A))\otimes\omega_{\widetilde{\mathrm Y}^{1}})\otimes\mathcal{O}_{\widetilde{\mathrm X}}(-\mathrm \mathrm{K}_{\widetilde{\mathrm X}})\big\rangle.
\end{multline}
\endgroup
Let $\mathcal{R}es_{1,i}$ be generated by the transported curve and
matrix-factorization components in the explicitly displayed inner span.
Let $\mathrm{Exc}_{1,i}$ consist, in order, of the transported
$\mathrm{Exc}_i$, the structure sheaf, and the $n$ direct summands of
the final object, equivalently the individual objects in \eqref{XSOD11'}. Every one is
exceptional, and their semi-orthogonality is preserved by the displayed
mutations and rotations. Therefore
\begin{equation}\label{XSOD12}
\begin{gathered}
\mathrm{D}^{b}(\widetilde{\mathrm X})=
\langle\mathcal{R}es_{1,i},\mathrm{Exc}_{1,i}\rangle,\qquad
\mathcal{R}es_{1,i}\simeq
\langle\mathrm{D}^{b}(\mathrm C),\mathcal{MF}(\widetilde{\mathrm X}')\rangle.
\end{gathered}
\end{equation}
If $\mathrm C$ is rational and the indicated matrix-factorization
component has a full exceptional collection, these components give a
full exceptional collection on $\widetilde{\mathrm X}$.
\end{enumerate}
\end{proof}

For the later curve notation, if $L\subset\widetilde X$ has strict
transform $\widehat L\subset Y$, put
$\mathcal{O}_L(-a)^\circ:=\pi_*\varpi^!\mathcal{O}_{\widehat L}(-a)$.
This defines the indicated derived modification.

We record the degree changes for rational curves along the link of Theorem~\ref{thm:hourglassconstrucion}.

\begin{proposition}[Rational curves on hourglasses]\label{hourglass:curves}
Fix $i$ and omit its subscripts. Let $L\subset\mathrm X$ be an integral rational curve whose strict transform $\Gamma\subset\mathrm Y$ is not a flopping or flipping curve. Write $\Gamma'\subset\mathrm Y'$ for its strict transform and put
\[
k:=-2\mathrm{K}_{\mathrm X}\cdot L,\qquad a:=\mathrm{A}\cdot\Gamma,\qquad d:=\mathrm f'^*\mathrm{H}'\cdot\Gamma'.
\]
Let $s,t$ be the total intersection multiplicities of the strict transform of $\Gamma$ with the exceptional divisors over the flop and Francia flip centres on a common weighted resolution. Then
\begin{equation}\label{eq:hourglass-curve-degree}
d=\frac{k-3a}{2}-s-t,\qquad a\in\mathbb Z_{\geq0},\quad s,t\in\mathbb Q_{\geq0}.
\end{equation}
In particular, $d\geq0$, with equality precisely when $\mathrm f'$ contracts $\Gamma'$. For the low degrees we have:
\begin{enumerate}
\item\label{equ1curve} If $k=1$, then $a=0$ and $d=\frac12-s-t$.
\item\label{equ2curve} If $k=3$, then either $a=0$ and $d=\frac32-s-t$, or $a=1$, $s=t=0$ and $d=0$.
\item\label{equ3curve} If $k=5$, then either $a=0$ and $d=\frac52-s-t$, or $a=1$ and $d=1-s-t$.
\end{enumerate}
If $\Gamma'$ is not contained in $\mathrm{E}_C:=\operatorname{Exc}(\mathrm f')$, put $m:=\mathrm{E}_C\cdot\Gamma'$. Then
\begin{equation}\label{eq:hourglass-target-degree}
k=(3z+2)d-3m-s-4t.
\end{equation}
Here $\mathrm{H}'$ is the primitive ample Weil class on $\mathrm X'$, and $z$ is determined by $-\mathrm{K}_{\mathrm X'}\sim(z+1)\mathrm{H}'$.
\end{proposition}
\begin{proof}
The blow-up and the final contraction give
\[
\mathrm f^*(-2\mathrm{K}_{\mathrm X})\sim-2\mathrm{K}_{\mathrm Y}+\mathrm{A},\qquad
\mathrm f'^*\mathrm{H}'\sim_{\mathbb Q}-\mathrm{K}_{\mathrm Y'}-\mathrm{A}'.
\]
Thus $-\mathrm{K}_{\mathrm Y}\cdot\Gamma=(k-a)/2$.

For a standard flop, use the common resolution $\widetilde{\mathrm W}$
with maps $\phi,\varphi$ as above, and let $\mathrm{E}_{\mathrm W}$ be its exceptional divisor.
The exceptional divisor $\mathrm{A}$ satisfies $\mathrm{A}\cdot C=1$ on each
flopping curve. The local model gives
\[
\varphi^*(-\mathrm{K}_{\widetilde{\mathrm Y}^{1}})
=\phi^*(-\mathrm{K}_{\widetilde{\mathrm Y}}),\qquad
\varphi^*\mathrm{A}=\phi^*\mathrm{A}+\mathrm{E}_{\mathrm W}.
\]
For a Francia flip, the source curve has type $\mathcal P^1_2$ and meets
$\mathrm{A}$ transversely at an ordinary point. Use the common weighted
resolution $\widetilde{\mathrm V}$ with maps $\phi',\varphi'$ as above,
and write $\mathrm{E}_{\mathrm V}$ for its exceptional divisor. The two discrepancies are
$1,2$, and the local weights $(1,2:-1,-1)$ give
\[
\begin{aligned}
\mathrm{K}_{\widetilde{\mathrm V}}
&=\phi'^*\mathrm{K}_{\widetilde{\mathrm Y}^{1}}+\mathrm{E}_{\mathrm V}
=\varphi'^*\mathrm{K}_{\widetilde{\mathrm Y}'}+2\mathrm{E}_{\mathrm V},\\
\varphi'^*\mathrm{A}^+&=\phi'^*\mathrm{A}+2\mathrm{E}_{\mathrm V}.
\end{aligned}
\]
Adding these contributions along the link yields
\[
-\mathrm{K}_{\mathrm Y'}\cdot\Gamma'=\frac{k-a}{2}+t,\qquad
\mathrm{A}'\cdot\Gamma'=a+s+2t.
\]
Subtracting proves \eqref{eq:hourglass-curve-degree}. The divisor $\mathrm{A}$ is Cartier and effective on $\mathrm Y$, so $a$ is a nonnegative integer. Since $\mathrm{H}'$ is ample, $d$ is nonnegative and vanishes exactly on contracted curves. The three cases follow.

Finally, $\mathrm{E}_C\sim_{\mathbb Q}z\,\mathrm f'^*\mathrm{H}'-\mathrm{A}'$, so $m=zd-a-s-2t$. Substitution gives \eqref{eq:hourglass-target-degree}.
\end{proof}

The multiplicities $s,t$ are not generally ordinary point counts; a curve through an index-two point can contribute a fractional intersection. Nor does $a$ count distinct orbifold points: it records contact with the chosen exceptional divisor. In the general construction, the flopping curves give degree-one curves on $\mathrm X$, whereas a Francia flipping curve gives a degree-two curve, since $-\mathrm{K}_{\mathrm Y^1}\cdot C=1/2$ and $\mathrm A_{\mathrm Y^1}\cdot C=1$.

For a target curve not contained in $\mathrm C$ whose strict transform avoids the small modification centres, \eqref{eq:hourglass-target-degree} reduces to $k=(3z+2)d-3m$. Thus a target line gives a degree-five curve only when it is $(z-1)$-secant to $\mathrm C$. These incidence conditions must be retained when comparing curve spaces.

\begin{proposition}
For any $i$, $\mathrm{B}_{1,i}$ is a minimal surface of general type and $\mathrm{D}^{b}(\mathrm{B}_{1,i})$ is indecomposable.
\end{proposition}
\begin{proof}
Adjunction gives
$\mathrm{K}_{\mathrm B_{1,i}}=(-\mathrm{K}_{\mathrm X_{1,i}})|_{\mathrm B_{1,i}}$,
which is ample. Thus the smooth surface $\mathrm B_{1,i}$ is minimal of
general type. By \cite[Corollary 3.2]{Tak1}, the anticanonical sections
are the weight-one coordinates in the weighted projective embedding.
Their base locus is contained in the index-two points.
A general member of $|-2\mathrm{K}_{\mathrm X_{1,i}}|$ avoids these points, so
their restrictions generate $\mathrm{K}_{\mathrm B_{1,i}}$.
The indecomposability now follows from
\cite[Corollary 1.5]{KawataniOkawa}. Very ampleness of the canonical
bundle is not needed.
\end{proof}

\subsubsection{Cubic type $\mathrm{X}_{1.9}$}\label{sub:cubictype}

Of all the examples, we are particularly interested in $\mathrm{X}_{1.9}$ because it may explain more relevant cubic geometries.
We start from a general smooth cubic $3$-fold $\mathrm{X'}_{1.9}$, it contains a smooth cubic surface $\mathrm{S}$ and a smooth twisted cubic $\mathrm{C}$ on $\mathrm{S}$.\\
$$\begin{tikzcd}[ampersand replacement=\&]
\mathrm{Y} \arrow[dr,"\mathrm{g}"] \arrow[d,"\mathrm{f}"]\& \& \mathrm{Y'}\arrow[dl,"\mathrm{g'}"]\arrow[d,"\mathrm{f'}"] \\
\mathrm{X}_{1.9}\&\mathrm{Z}\& \mathrm{X'}_{1.9}
\end{tikzcd}$$
Firstly, we blow up $\mathrm{X'}$ along the twisted cubic $\mathrm{C}$ (we denote its exceptional divisor by $\mathrm{E}$), and then easy to see there are six $(-1,-1)$ lines bi-secant to the twisted cubic which is $\mathrm{K}_{\mathrm{Y'}}$ trivial.  Then we flop these lines, noticing the strict transformation (we denote this divisor by $\mathrm{A}$) of $\mathrm{S}$ will also contract six disjoint lines that a just a plane, whose normals sheaf $\mathcal{O}_{\mathrm{A}}(\mathrm{A})\simeq \mathcal{O}_{\mathrm{A}}(-2)$. Finally, we contract this plane to a $(1^{3})/2$ point to get our target variety $\mathrm{X}_{1.9}$.\\

We can translate it into projective geometry. First, we project $\mathrm{X}_{1.9}$ along the singular point, and we get a Gushel $3$-fold with a $\sigma_{2,2}$ plane and six ODPs. Then we project the Gushel $3$-fold along this $\sigma_{2,2}$ plane, we get our target cubic $3$-fold $\mathrm{X'}_{1.9}$.
\[\begin{adjustbox}{max width=0.98\linewidth}\(\begin{tikzcd}[ampersand replacement=\&]
\sqrt{\mathrm{Y}}:=\sqrt[2]{{\mathrm{A}/\mathrm{Y}}}\arrow[r,"\varpi"] \arrow[d,"\pi"] \& \mathrm{Y}\arrow[d,"\mathrm{g}"]\arrow[r,dashed]\& \mathrm{Y}'\arrow[d,"\mathrm{f}'"] \\
\widetilde{\mathrm{X}}_{1.9}\subset \mathrm{Cone}(\mathrm{v}^{(2)}\mathbf{Gr}(2,5)\cap(1)^{2})\subset\mathds{P}(1^{8},2)\arrow[r,dashed,"\mathrm{pr}"]\& \mathrm{Z}\subset \mathbf{Gr}(2,5)\cap(1)^{2}\subset\mathds{P}(1^{8})\arrow[r,dashed,"\mathrm{pr}'"]\& \mathrm{X}'\subset\mathds{P}(1^{5})\subset\mathds{P}(1^{6})
\end{tikzcd}\)\end{adjustbox}\]

For the construction of other examples, we refer to \cite{Tak2}. The different positions of the unprojection plane on the midpoint of the Sarkisov link determine the different hourglasses with the same degree.

\begin{lemma}\label{lemma:lf}
If $\mathrm{U}_{\mathrm{Z}}$ denotes the tautological sheaf of $\mathbf{Gr}(2,5)$ restricts at mid point of the Sarkisov link $\mathrm{Z}$, we have
$$\pi_{*}\varpi^{*}g^{*}\mathrm{U}_{\mathrm{Z}}\quad\text{and}\quad\pi_{*}\varpi^{*}g^{*}\mathrm{U}^{\vee}_{\mathrm{Z}}$$
are  rank $2$ sheaves on $\widetilde{\mathrm{X}}_{1.9}$, without ambiguity we denote them by $\mathrm{U}_{\widetilde{\mathrm{X}}_{1.9}}$ and $\mathrm{U}^{\vee}_{\widetilde{\mathrm{X}}_{1.9}}$.
\end{lemma}

\begin{proof}
Firstly we prove that $\pi_{*}\varpi^{*}g^{*}\mathrm{U}_{\mathrm{Z}}$ is concentrated in cohomological degree zero, to prove this, we only need to prove that as any Zariski point $x$ on $\widetilde{\mathrm{X}}$, we have ${\mathbf{R}^{i}\pi_{*}\varpi^{*}g^{*}\mathrm{U}_{\mathrm{Z}}}|_{x}$ is zero up to formal completion for any $i\neq 0$, which trivial for any smooth point. So we only need to consider the situation on the singular point $\mathrm{p}$, the idea is to take a presentation of the local $\acute{e}$tale atlas of $\pi$, simplify the problem to the scheme and use the theorem on formal functions. First notice we unproject the $\sigma_{2,2}$ plane:
$$\{M|M\subset V_{3},\, \mathrm{rank} (M)=2\}=\mathbf{Gr}(2,3)\subset\{M|M\subset V_{5},\, \mathrm{rank} (M)=2\}=\mathbf{Gr}(2,5)$$
so $\mathrm{U}_{\mathrm{Z}}|_{\mathbf{Gr}(2,3)}$ is the tautological bundle on $\mathbf{Gr}(2,3)$, and under its isomorphism to a plane  $\mathbf{Gr}(2,3)\simeq\mathrm{A}$, we have $\mathrm{U}_{\mathrm{Z}}|_{\mathrm{A}}=\Omega_{\mathrm{A}}(1)$, further we have $\varpi^{*}g^{*}\mathrm{U}_{\mathrm{Z}}|_{M_{\mathrm{A}}}=\Omega_{M_{\mathrm{A}}}(1)$, where $M_{\mathrm{A}}$ is the tautological divisor of our root of $\mathrm{A}$, i.e. exceptional divisor of local equivariant blow-up, which is also a plane.\\

Then we are expected to compute \begin{equation}\label{equ:h1a}\mathrm{H}^{i}(\sqrt[2]{{\mathrm{A}/\mathrm{Y}}},(\varpi^{*}g^{*}\mathrm{U}_{\mathrm{Z}})_{(n)})\end{equation}
where subscript $(n)$ means their restriction on $n$-th thicken neighbourhood of exceptional divisor of $\pi$. Specially it's easy to see $\mathrm{H}^{i}(M_{\mathrm{A}},\Omega_{M_{\mathrm{A}}}(j))=0$ for any $i\neq0$ and $j\geq1$, so by exact sequences
$$0\longrightarrow\mathcal{O}_{M_{\mathrm{A}}}(n)\longrightarrow\mathcal{O}_{M_{\mathrm{A}},(n+1)}\longrightarrow\mathcal{O}_{M_{\mathrm{A}},(n)}\longrightarrow0$$
for any $n\geq 1$, tensor with locally free sheaf $\varpi^{*}g^{*}\mathrm{U}_{\mathrm{Z}}$,
$$0\longrightarrow\Omega_{M_{\mathrm{A}}}(n+1)\longrightarrow{\varpi^{*}g^{*}\mathrm{U}_{\mathrm{Z}}}_{(n+1)}\longrightarrow{\varpi^{*}g^{*}\mathrm{U}_{\mathrm{Z}}}_{(n)}\longrightarrow0$$
by repeating the cohomology calculation on the above short exact sequences,  we can see all terms in (\ref{equ:h1a}) are zero if $i\neq 0$. So ${\mathbf{R}^{i}\pi_{*}\varpi^{*}g^{*}\mathrm{U}_{\mathrm{Z}}}|_{\mathrm{p}}$ is zero for any $i\neq 0$ by the theorem on formal functions and Nakayama's lemma.
If we restrict outside the singular point, $\pi_{*}\varpi^{*}g^{*}\mathrm{U}_{\mathrm{Z}}|_{\widetilde{X}_{1.9}-\mathrm{p}}$ is a
locally free sheaf defined by pullback of $\mathrm{U}_{\mathrm{Z}}$ along the projection, it's generic rank $2$.\\

Similar result holds for $\pi_{*}\varpi^{*}g^{*}\mathrm{U}^{\vee}_{\mathrm{Z}}$, noticing $\mathrm{U}^{\vee}_{\mathrm{Z}}|_{\mathrm{A}}=\mathrm{T}_{\mathrm{A}}(-1)$ and $\mathrm{H}^{i}(M_{\mathrm{A}},\mathrm{T}_{\mathrm{A}}(j))=0$ for any $i\neq0$ and $j\geq-1$.
\end{proof}

\begin{remark}
The sheaf $\mathrm{U}^{\vee}_{\widetilde{\mathrm{X}}_{1.9}}$ is locally free away from $\mathrm p$. Proposition~\ref{prop:ECof1.9} shows that it is not locally free at $\mathrm p$, by its Ext groups with the stacky lines $L_i$.
\end{remark}

Applying the Theorem \ref{thm:SODh}, we know that the \textit{matrix factorization category} of cubic $3$-fold $\mathrm{B}_{3}$ is retained in $\mathrm{X}_{1.9}$, and we are more curious about the meaning of the elements in the exceptional collection of $\mathrm{X}_{1.9}$. We have the following description,

\begin{proposition}\label{prop:X1.9}
We have a semi-orthogonal decomposition of $\widetilde{\mathrm{X}}_{1.9}$:
$$\mathrm{D}^{b}(\widetilde{\mathrm{X}}_{1.9})\simeq \langle \mathcal{O}_{\widetilde{\mathrm{X}}_{1.9}}(-2\mathrm{p})\rightarrow \bigoplus_{i=1}^{6}\mathcal{O}_{L_{i}}(-1),\mathcal{O}_{\widetilde{\mathrm{X}}_{1.9}},\mathrm{U}^{\vee}_{\widetilde{\mathrm{X}}_{1.9}},\mathcal{MF}(\widetilde{\mathrm{X}}_{1.9})\rangle$$
and there is an equivalence between the residue component $\mathcal{MF}(\widetilde{\mathrm{X}}_{1.9})$ and $\mathcal{MF}(B_{3})$ via a functor
$$\Gamma_{5}:\mathcal{MF}(B_{3})\longrightarrow\mathcal{MF}(\widetilde{\mathrm{X}}_{1.9})$$
$$\Gamma_{5}(-):=\pi_{*}\varpi^{*}\big[\phi_{*}\varphi^{*}\mathbb{R}_{\langle\mathcal{O}_{\mathrm{S}}(\mathrm{S}),\mathcal{O}_{\mathrm{Y}'},g'^{*}\mathrm{U}_{\mathrm{Z}}^{\vee}\rangle}\big(\mathrm{f}'^{*}(-)(\mathrm{H})\big)\big]$$
\end{proposition}

\begin{proof}

\begin{enumerate}

            \item Consider  $\mathrm{D}^{b}(\mathrm{Y}')$ under blowing up formula of derived category e.g. Theorem \ref{thm:singularweightedblowup}:
          $$\begin{tikzcd}[ampersand replacement=\&]
\mathrm{E}\arrow[r,"\iota'"]\arrow[d,"\mathrm{f}'_{\mathrm{E}}"] \&\mathrm{Y}'\arrow[d,"\mathrm{f}'"] \\
\mathrm{C} \arrow[r,"\kappa'"] \&\mathrm{X}'
\end{tikzcd}$$
          we have
           \begin{equation} \mathrm{D}^{b}(\mathrm{\widetilde{Y}'})\simeq \langle \iota'_{*}\big(\mathrm{f}'^{*}_{\mathrm{E}}\mathrm{D}^{b}(\mathrm{C})\otimes\mathcal{O}_{\mathrm{E}'}(-1)\big), \mathrm{f}'^{*}\mathrm{D}^{b}(\mathrm{X'})\rangle\end{equation}
           where $\mathrm{X}'$ is a cubic $3$-fold
          $$\mathrm{D}^{b}(\mathrm{X}')\simeq \langle \mathcal{MF}(\mathrm{X}'),\mathcal{O}_{\mathrm{X}'}(-\mathrm{H}),\mathcal{O}_{\mathrm{X}'}\rangle$$
          by Corollary \ref{cor:mfx'}, and $\mathrm{C}$ is a twisted cubic, its derived category is
          $$\mathrm{D}^{b}(\mathrm{C})\simeq \langle\mathcal{O}_{\mathrm{C}}(ah),\mathcal{O}_{\mathrm{C}}((a+1)h)\rangle$$
          for any integer $a$ and $h$ is the point class on rational curve $\mathrm{C}$, we have a relationship $\mathrm{H}|_{\mathrm{C}}=3h$.\\

           So
           \begin{equation}\label{sodc1}\mathrm{D}^{b}(\mathrm{Y'})\simeq \langle \mathcal{O}_{\mathrm{E}}(ah+\mathrm{E}), \mathcal{O}_{\mathrm{E}}((a+1)h+\mathrm{E}),\mathrm{f}'^{*}\mathcal{MF}(\mathrm{X}'),\mathcal{O}_{\mathrm{Y}'}(-\mathrm{H}),\mathcal{O}_{\mathrm{Y}'}\rangle\end{equation}
for any integer $a$.
           \item  Mutate terms  $\mathcal{O}_{\mathrm{E}}(ah+\mathrm{E}), \mathcal{O}_{\mathrm{E}}((a+1)h+\mathrm{E})$ in (\ref{sodc1}) to the far right
         \begin{equation}\label{sodc2}\mathrm{D}^{b}(\mathrm{Y'})\simeq \langle \mathrm{f}'^{*}\mathcal{MF}(\mathrm{X}'),\mathcal{O}_{\mathrm{Y}'}(-\mathrm{H}),\mathcal{O}_{\mathrm{Y}'},\mathcal{O}_{\mathrm{E}}((a+6)h), \mathcal{O}_{\mathrm{E}}((a+7)h)\rangle\end{equation}
  noticing $-\mathrm{K}_{\mathrm{Y}'}=-\mathrm{K}_{\mathrm{X}'}-\mathrm{E}=2\mathrm{H}-\mathrm{E}$, and $\mathrm{H}|_{\mathrm{C}}=3h$.

            \item For the convenience we pick $a=-6$,  and left mutate term $\mathcal{O}_{\mathrm{E}}$ respect to $\mathcal{O}_{\mathrm{Y}'}$ in (\ref{sodc2}), since by mutation triangle:
            $$\mathbb{L}_{\mathcal{O}_{\mathrm{Y}'}}\mathcal{O}_{E}[-1]\longrightarrow \mathrm{Ext}^{*}(\mathcal{O}_{\mathrm{Y}'},\mathcal{O}_{\mathrm{E}})\otimes\mathcal{O}_{\mathrm{Y}'}\longrightarrow\mathcal{O}_{\mathrm{E}}\longrightarrow\mathbb{L}_{\mathcal{O}_{\mathrm{Y}'}}\mathcal{O}_{\mathrm{E}}$$
            we have $\mathbb{L}_{\mathcal{O}_{\mathrm{Y}'}}\mathcal{O}_{\mathrm{E}}\simeq\mathcal{O}_{\mathrm{Y'}}(-\mathrm{E})[1]$, so
             \begin{equation}\label{sodc3}\mathrm{D}^{b}(\mathrm{Y'})\simeq \langle \mathrm{f}'^{*}\mathcal{MF}(\mathrm{X}'),\mathcal{O}_{\mathrm{Y}'}(-\mathrm{H}),\mathcal{O}_{\mathrm{Y'}}(-\mathrm{E}),\mathcal{O}_{\mathrm{Y}'}, \mathcal{O}_{\mathrm{E}}(h)\rangle\end{equation}

            \item Left mutate object $\mathcal{O}_{\mathrm{E}}(h)$ respect to $\mathcal{O}_{\mathrm{Y'}}$ in (\ref{sodc3}), noticing the mutation triangle
            $$\mathbb{L}_{\mathcal{O}_{\mathrm{Y}'}}\mathcal{O}_{\mathrm{E}}(h)[-1]\longrightarrow \mathrm{Ext}^{*}(\mathcal{O}_{\mathrm{Y}'},\mathcal{O}_{\mathrm{E}}(h))\otimes\mathcal{O}_{\mathrm{Y}'}\longrightarrow\mathcal{O}_{\mathrm{E}}(h)\longrightarrow\mathbb{L}_{\mathcal{O}_{\mathrm{Y}'}}\mathcal{O}_{\mathrm{E}}(h)$$
             it's not difficult to see $\mathrm{Ext}^{*}(\mathcal{O}_{\mathrm{Y}'},\mathcal{O}_{\mathrm{E}}(h))=\mathrm{H}^{*}(\mathrm{C},\mathcal{O}_{\mathrm{C}}(1))=\mathrm{k}^{2}$, then take advantage of short exact sequence from \cite[Proposition 5.5]{KuGM} and the following is an affine base change of it, we have
             $$0\longrightarrow g'^{*}\mathrm{U}_{\mathrm{Z}}^{\vee}(-\mathrm{H})\longrightarrow\mathcal{O}_{\mathrm{Y}'}^{\oplus 2}\longrightarrow\mathcal{O}_{\mathrm{E}}(h)\longrightarrow0$$
so we can see $\mathbb{L}_{\mathcal{O}_{\mathrm{Y}'}}\mathcal{O}_{\mathrm{E}}(h)[-1]=g'^{*}\mathrm{U}_{\mathrm{Z}}^{\vee}(-\mathrm{H})$, and
   \begin{equation}\label{sodc4}\mathrm{D}^{b}(\mathrm{Y'})\simeq \langle \mathrm{f}'^{*}\mathcal{MF}(\mathrm{X}'),\mathcal{O}_{\mathrm{Y}'}(-\mathrm{H}),\mathcal{O}_{\mathrm{Y'}}(-\mathrm{E}),g'^{*}\mathrm{U}_{\mathrm{Z}}^{\vee}(-\mathrm{H}),\mathcal{O}_{\mathrm{Y}'}\rangle\end{equation}

            \item Tensor $\mathcal{O}_{\mathrm{Y}'}(\mathrm{H})$ on the both hand side of (\ref{sodc4}),
                   \begin{equation}\label{sodc5}\mathrm{D}^{b}(\mathrm{Y'})\simeq \langle \mathrm{f}'^{*}\mathcal{MF}(\mathrm{X}')(\mathrm{H}),\mathcal{O}_{\mathrm{Y}'},\mathcal{O}_{\mathrm{Y'}}(\mathrm{H}-\mathrm{E}),g'^{*}\mathrm{U}_{\mathrm{Z}}^{\vee},\mathcal{O}_{\mathrm{Y}'}(\mathrm{H})\rangle\end{equation}
         we know the cubic surface $\mathrm{S}$ is linear equivalent to $\mathrm{H}-\mathrm{E}$ in $\mathrm{Y}'$.
           \item Left mutate term  $\mathcal{O}_{\mathrm{Y'}}(\mathrm{H-E})$ respect to $\mathcal{O}_{\mathrm{Y}'}$ in (\ref{sodc5}), noticing mutation triangle
             $$\mathbb{L}_{\mathcal{O}_{\mathrm{Y}'}}\mathcal{O}_{\mathrm{Y'}}(\mathrm{H-E})[-1]\longrightarrow \mathrm{Ext}^{*}(\mathcal{O}_{\mathrm{Y}'},\mathcal{O}_{\mathrm{Y'}}(\mathrm{H-E}))\otimes\mathcal{O}_{\mathrm{Y}'}\longrightarrow\mathcal{O}_{\mathrm{Y'}}(\mathrm{H-E})$$
           and  $\mathrm{Ext}^{*}(\mathcal{O}_{\mathrm{Y}'},\mathcal{O}_{\mathrm{Y'}}(\mathrm{H-E}))=\mathrm{H}^{*}(\mathrm{Y},\mathcal{O}_{\mathrm{Y}}(\mathrm{A}))=\mathrm{k}$ by a totally same argument as in (\ref{mutacohcomp}), so $\mathbb{L}_{\mathcal{O}_{\mathrm{Y}'}}\mathcal{O}_{\mathrm{Y'}}(\mathrm{H-E})=\mathcal{O}_{\mathrm{S}}(\mathrm{S})$ and

          \begin{equation}\label{sodc6}\mathrm{D}^{b}(\mathrm{Y'})\simeq \langle \mathrm{f}'^{*}\mathcal{MF}(\mathrm{X}')(\mathrm{H}),\mathcal{O}_{\mathrm{S}}(\mathrm{S}),\mathcal{O}_{\mathrm{Y}'},g'^{*}\mathrm{U}_{\mathrm{Z}}^{\vee},\mathcal{O}_{\mathrm{Y}'}(\mathrm{H})\rangle\end{equation}

            \item Mutate term $\mathcal{O}_{\mathrm{Y'}}(\mathrm{H})$ in (\ref{sodc6}) to the far left, we have
                \begin{equation}\label{sodc61}\mathrm{D}^{b}(\mathrm{Y'})\simeq \langle \mathcal{O}_{\mathrm{Y}'}(-\mathrm{S}),\mathrm{f}'^{*}\mathcal{MF}(\mathrm{X}')(\mathrm{H}),\mathcal{O}_{\mathrm{S}}(\mathrm{S}),\mathcal{O}_{\mathrm{Y}'},g'^{*}\mathrm{U}_{\mathrm{Z}}^{\vee}\rangle\end{equation}
     since $-\mathrm{K}_{\mathrm{Y}'}\sim2\mathrm{H}-\mathrm{E}$ and $\mathrm{S}\sim\mathrm{H}-\mathrm{E}$.
            \item Mutate term $\mathrm{f}'^{*}\mathcal{MF}(\mathrm{X}')(\mathrm{H})$ in (\ref{sodc61}) to the far right, we have
                \begin{equation}\label{sodc7}\mathrm{D}^{b}(\mathrm{Y'})\simeq \langle \mathcal{O}_{\mathrm{Y}'}(-\mathrm{S}),\mathcal{O}_{\mathrm{S}}(\mathrm{S}),\mathcal{O}_{\mathrm{Y}'},g'^{*}\mathrm{U}_{\mathrm{Z}}^{\vee},\Gamma_{1}\mathcal{MF}(\mathrm{X}')\rangle\end{equation}
   where $\Gamma_{1}(-):=\mathbb{R}_{\langle\mathcal{O}_{\mathrm{S}}(\mathrm{S}),\mathcal{O}_{\mathrm{Y}'},g'^{*}\mathrm{U}_{\mathrm{Z}}^{\vee}\rangle}\mathrm{f}'^{*}(-)(\mathrm{H})$

            \item Act the standard flop functor $\phi_{*}\varphi^{*}(-)$ on both hand sides of (\ref{sodc7})
             $$\begin{tikzcd}[column sep=0.5em]
  & && \mathrm{W}\arrow[dl,"\phi"]\arrow[dr,"\varphi"] &&  \\
 & &\mathrm{Y}&  &\mathrm{Y}'&
\end{tikzcd}$$
noticing $\mathrm{S}-\mathrm{E_{W}}\sim \mathrm{A}$, where $\mathrm{E_{W}}$ is the exceptional locus on $\mathrm{W}$, since $\mathrm{S}$ contains all the flopped curves. It's not difficult to see
\begin{enumerate}
  \item\label{equ:comp1} $\phi_{*}\varphi^{*}\mathcal{O}_{\mathrm{Y}'}=\mathcal{O}_{\mathrm{Y}}$
  \item\label{equ:comp2} $\phi_{*}\varphi^{*}\mathcal{O}_{\mathrm{Y}'}(\mathrm{S})=\mathcal{O}_{\mathrm{Y}}(\mathrm{A})$
  \item\label{equ:comp3} $\phi_{*}\varphi^{*}\mathcal{O}_{\mathrm{S}}(\mathrm{S})=\mathcal{O}_{\mathrm{A}}(\mathrm{A})$
  \item\label{equ:comp4} $\phi_{*}\varphi^{*}\mathcal{O}_{\mathrm{Y}'}(-\mathrm{S})=\mathcal{O}_{\mathrm{Y}}(-\mathrm{A})\otimes\mathcal{I}_{\widehat{L_{\mathrm{Y}}}}$
  \item\label{equ:comp5} $\phi_{*}\varphi^{*}g'^{*}\mathrm{U}_{\mathrm{Z}}^{\vee}=g^{*}\mathrm{U}_{\mathrm{Z}}^{\vee}$

\end{enumerate}
for (\ref{equ:comp1})-(\ref{equ:comp3}) the argument is the same as in (\ref{mutacohcomp}), (\ref{equ:comp4}) is due to $\phi_{*}\mathcal{O}_{\mathrm{W}}(-\mathrm{E_{W}})=\mathcal{I}_{\widehat{L_{\mathrm{Y}}}}$, where $\widehat{L_{\mathrm{Y}}}$ is the union of all flopping curves, (\ref{equ:comp5}) is due to the commutative of flop diagram. So
 \begin{equation}\label{sodc8}\mathrm{D}^{b}(\mathrm{Y})\simeq\phi_{*}\varphi^{*}\mathrm{D}^{b}(\mathrm{Y'})\simeq \langle \mathcal{O}_{\mathrm{Y}}(-\mathrm{A})\otimes\mathcal{I}_{\widehat{L_{\mathrm{Y}}}},\mathcal{O}_{\mathrm{A}}(\mathrm{A}),\mathcal{O}_{\mathrm{Y}},g^{*}\mathrm{U}_{\mathrm{Z}}^{\vee},\Gamma_{2}\mathcal{MF}(\mathrm{X}')\rangle\end{equation}
   where $\Gamma_{2}(-):=\phi_{*}\varphi^{*}\mathbb{R}_{\langle\mathcal{O}_{\mathrm{S}}(\mathrm{S}),\mathcal{O}_{\mathrm{Y}'},g'^{*}\mathrm{U}_{\mathrm{Z}}^{\vee}\rangle}\mathrm{f}'^{*}(-)(\mathrm{H})$

\item Tensor with $\mathcal{O}_{\mathrm{Y}}(-\mathrm{K}_{\mathrm{Y}})$ on each term of (\ref{sodc8}),
 \begin{equation}\label{sodc9}\mathrm{D}^{b}(\mathrm{Y})\simeq\langle \mathcal{O}_{\mathrm{Y}}(-\mathrm{A}-\mathrm{K}_{\mathrm{Y}})\otimes\mathcal{I}_{\widehat{L_{\mathrm{Y}}}},\mathcal{O}_{\mathrm{A}}(-1),\mathcal{O}_{\mathrm{Y}}(-\mathrm{K}_{\mathrm{Y}}),g^{*}\mathrm{U}_{\mathrm{Z}}^{\vee}(-\mathrm{K}_{\mathrm{Y}}),\Gamma_{3}\mathcal{MF}(\mathrm{X}')\rangle\end{equation}
since $-\mathrm{K}_{\mathrm{Y}}|_{\mathrm{A}}\sim-\mathrm{A}/2|_{\mathrm{A}}$, where $\Gamma_{3}(-):=\big[\phi_{*}\varphi^{*}\mathbb{R}_{\langle\mathcal{O}_{\mathrm{S}}(\mathrm{S}),\mathcal{O}_{\mathrm{Y}'},g'^{*}\mathrm{U}_{\mathrm{Z}}^{\vee}\rangle}[\mathrm{f}'^{*}(-)(\mathrm{H})]\big]\otimes\mathcal{O}_{\mathrm{Y}}(-\mathrm{K}_{\mathrm{Y}})$.

\item Right mutate term $\mathcal{O}_{\mathrm{Y}}(-\mathrm{A}-\mathrm{K}_{\mathrm{Y}})\otimes\mathcal{I}_{\widehat{L_{\mathrm{Y}}}}$ respect to $\mathcal{O}_{\mathrm{A}}(-1)$ in (\ref{sodc9}), it's just
 \begin{equation}\label{sodc10}\mathrm{D}^{b}(\mathrm{Y})\simeq \langle \mathcal{O}_{\mathrm{A}}(-1), \mathbb{R}_{\mathcal{O}_{\mathrm{A}}(-1)}\bigl(\mathcal{O}_{\mathrm{Y}}(-\mathrm{A}-\mathrm{K}_{\mathrm{Y}})\otimes\mathcal{I}_{\widehat{L_{\mathrm{Y}}}}\bigr),\mathcal{O}_{\mathrm{Y}}(-\mathrm{K}_{\mathrm{Y}}),g^{*}\mathrm{U}_{\mathrm{Z}}^{\vee}(-\mathrm{K}_{\mathrm{Y}}),\Gamma_{3}\mathcal{MF}(\mathrm{X}')\rangle\end{equation}

\item Recall our weighted blow-up formula for $\mathrm{D}^{b}(\mathrm{Y})$ e.g. Theorem \ref{thm:singularweightedblowup}, we have
 \begin{equation}\label{sodc11}\mathrm{D}^{b}(\mathrm{Y})\simeq \langle\mathcal{O}_{\mathrm{A}}(-1),\varpi_{*}\pi^{*}\mathrm{D}^{b}(\widetilde{\mathrm{X}})\rangle\end{equation}
comparing (\ref{sodc11}) with (\ref{sodc10}), we have a semi-orthogonal decomposition of $\mathrm{D}^{b}(\widetilde{\mathrm{X}})$ by projection along $\mathcal{O}_{\mathrm{A}}(-1)$, that is
 \begin{equation}\label{sodc12}\mathrm{D}^{b}(\widetilde{\mathrm{X}})\simeq \langle \pi_{*}\varpi^{!}\mathbb{R}_{\mathcal{O}_{\mathrm{A}}(-1)}\bigl(\mathcal{O}_{\mathrm{Y}}(-\mathrm{A}-\mathrm{K}_{\mathrm{Y}})\otimes\mathcal{I}_{\widehat{L_{\mathrm{Y}}}}\bigr),\pi_{*}\varpi^{!}\mathcal{O}_{\mathrm{Y}}(-\mathrm{K}_{\mathrm{Y}}),\pi_{*}\varpi^{!}g^{*}\mathrm{U}_{\mathrm{Z}}^{\vee}(-\mathrm{K}_{\mathrm{Y}}),\Gamma_{4}\mathcal{MF}(\mathrm{X}')\rangle\end{equation}
where $\Gamma_{4}(-):=\pi_{*}\varpi^{!}\big\{\big[\phi_{*}\varphi^{*}\mathbb{R}_{\langle\mathcal{O}_{\mathrm{S}}(\mathrm{S}),\mathcal{O}_{\mathrm{Y}'},g'^{*}\mathrm{U}_{\mathrm{Z}}^{\vee}\rangle}[\mathrm{f}'^{*}(-)(\mathrm{H})]\big]\otimes\mathcal{O}_{\mathrm{Y}}(-\mathrm{K}_{\mathrm{Y}})\big\}$.\\

We have the following simple calculation:
\begin{enumerate}
  \item \label{enu:exp1} $\pi_{*}\varpi^{!}\big(g^{*}\mathrm{U}_{\mathrm{Z}}^{\vee}(-\mathrm{K}_{\mathrm{Y}})\big)=\pi_{*}\varpi^{*}g^{*}\mathrm{U}_{\mathrm{Z}}^{\vee}(-\mathrm{K}_{\widetilde{\mathrm{X}}})=\mathrm{U}_{\widetilde{\mathrm{X}}}^{\vee}(-\mathrm{K}_{\widetilde{\mathrm{X}}})$
  \item $\pi_{*}\varpi^{!}\mathcal{O}_{\mathrm{Y}}(-\mathrm{K}_{\mathrm{Y}})=\mathcal{O}_{\widetilde{\mathrm{X}}}(-\mathrm{K}_{\widetilde{\mathrm{X}}})$
  \item \label{easycomp}$\varpi^{*}\mathcal{O}_{\widehat{L_{\mathrm{Y}}}}=\mathcal{O}_{\widehat{L_{\mathrm{\sqrt{Y}}}}}$, it is due to a computation using the semi-orthogonal decomposition in (\ref{euq:sodofexp}), we have $\mathcal{O}_{\widehat{L_{\mathrm{\sqrt{Y}}}}}\otimes\mathcal{M}^{-1}_{\mathrm{A}}|_{\mathrm{A}}=\mathrm{k}_{M_{\mathrm{A}}\cap \widehat{L_{\mathrm{\sqrt{Y}}}}}\otimes\nu$ and $\varpi_{*}\mathcal{O}_{\widehat{L_{\mathrm{\sqrt{Y}}}}}=\mathcal{O}_{\widehat{L_{\mathrm{Y}}}}$.
  \item Noticing the mutation triangle
$$\mathbb{R}_{\mathcal{O}_{\mathrm{A}}(-1)}\bigl(\mathcal{O}_{\mathrm{Y}}(-\mathrm{A}-\mathrm{K}_{\mathrm{Y}})\otimes\mathcal{I}_{\widehat{L_{\mathrm{Y}}}}\bigr)\longrightarrow\mathcal{O}_{\mathrm{Y}}(-\mathrm{A}-\mathrm{K}_{\mathrm{Y}})\otimes\mathcal{I}_{\widehat{L_{\mathrm{Y}}}}\longrightarrow V^{*}\otimes\mathcal{O}_{\mathrm{A}}(-1)$$
and $\pi_{*}\varpi^{!}\mathcal{O}_{\mathrm{A}}(-1)=0$ by orthogonality in (\ref{sodc11}), we have
$$\pi_{*}\varpi^{!}\big(\mathbb{R}_{\mathcal{O}_{\mathrm{A}}(-1)}\bigl(\mathcal{O}_{\mathrm{Y}}(-\mathrm{A}-\mathrm{K}_{\mathrm{Y}})\otimes\mathcal{I}_{\widehat{L_{\mathrm{Y}}}}\bigr)\big)=\pi_{*}\varpi^{!}\big(\mathcal{O}_{\mathrm{Y}}(-\mathrm{A}-\mathrm{K}_{\mathrm{Y}})\otimes\mathcal{I}_{\widehat{L_{\mathrm{Y}}}}\big)$$
it's easy to see after (\ref{easycomp})
$$\pi_{*}\varpi^{!}\big(\mathcal{O}_{\mathrm{Y}}(-\mathrm{A}-\mathrm{K}_{\mathrm{Y}})\otimes\mathcal{I}_{\widehat{L_{\mathrm{Y}}}}\big)=\pi_{*}\bigl(\mathcal{O}_{\mathrm{\sqrt{Y}}}(-\mathrm{A}-\mathrm{K}_{\widetilde{\mathrm{X}}})\otimes\mathcal{I}_{\widehat{L_{\mathrm{\sqrt{Y}}}}}\bigr)$$
since $\varpi^{!}(-)=\varpi^{*}(-)\otimes\mathcal{M}_{\mathrm{A}}$, and $\mathrm{K}_{\mathrm{Y}}\sim \mathrm{K}_{\mathrm{X}}+M_{\mathrm{A}}$ in $\sqrt{\mathrm{Y}}$.
\end{enumerate}
combine all,
\begin{equation}\label{sodc13}\mathrm{D}^{b}(\widetilde{\mathrm{X}})\simeq \langle\pi_{*}\bigl(\mathcal{O}_{\mathrm{\sqrt{Y}}}(-\mathrm{A}-\mathrm{K}_{\widetilde{\mathrm{X}}})\otimes\mathcal{I}_{\widehat{L_{\mathrm{\sqrt{Y}}}}}\bigr), \mathcal{O}_{\widetilde{\mathrm{X}}}(-\mathrm{K}_{\widetilde{\mathrm{X}}}),\mathrm{U}_{\widetilde{\mathrm{X}}}^{\vee}(-\mathrm{K}_{\widetilde{\mathrm{X}}}),\Gamma_{4}\mathcal{MF}(\mathrm{X}')\rangle\end{equation}
where $\Gamma_{4}(-):=\pi_{*}\varpi^{!}\big\{\big[\phi_{*}\varphi^{*}\mathbb{R}_{\langle\mathcal{O}_{\mathrm{S}}(\mathrm{S}),\mathcal{O}_{\mathrm{Y}'},g'^{*}\mathrm{U}_{\mathrm{Z}}^{\vee}\rangle}[\mathrm{f}'^{*}(-)(\mathrm{H})]\big]\otimes\mathcal{O}_{\mathrm{Y}}(-\mathrm{K}_{\mathrm{Y}})\big\}$.\\

\item Tensor with $\mathcal{O}_{\widetilde{\mathrm{X}}}(\mathrm{K}_{\widetilde{\mathrm{X}}})$ on both hand sides of (\ref{sodc13}),
and notice that $\pi_{*}\mathcal{I}_{\widehat{L_{\mathrm{\sqrt{Y}}}}}(-A)$ is quasi-isomorphic to a two terms complex
$$\mathcal{O}_{\widetilde{\mathrm{X}}}(-2\mathrm{p})\rightarrow \mathcal{O}_{L}(-1):=\bigoplus_{i}\mathcal{O}_{L_{i}}(-1)$$
The first term is placed in degree $0$ and the second in degree $1$;
the map is induced by the ideal sequence on the root space.
Here $\mathcal{O}_{\widetilde{\mathrm{X}}}(-2\mathrm{p})$ denotes the
squared ideal $\mathcal{I}_{\mathrm p}^{\,2}$.
Thus
\begin{equation}\label{sodc14}\mathrm{D}^{b}(\widetilde{\mathrm{X}})\simeq \langle\mathcal{O}_{\widetilde{\mathrm{X}}}(-2\mathrm{p})\rightarrow \bigoplus_{i}\mathcal{O}_{L_{i}}(-1), \mathcal{O}_{\widetilde{\mathrm{X}}},\mathrm{U}_{\widetilde{\mathrm{X}}}^{\vee},\Gamma_{5}\mathcal{MF}(\mathrm{X}')\rangle\end{equation}
where $\Gamma_{5}(-):=\pi_{*}\varpi^{*}\big[\phi_{*}\varphi^{*}\mathbb{R}_{\langle\mathcal{O}_{\mathrm{S}}(\mathrm{S}),\mathcal{O}_{\mathrm{Y}'},g'^{*}\mathrm{U}_{\mathrm{Z}}^{\vee}\rangle}\big(\mathrm{f}'^{*}(-)(\mathrm{H})\big)\big]$,
we finish the proof.
\end{enumerate}
\end{proof}

In the following proposition, $J(\mathrm X)$ is defined from the ordinary Betti cohomology of the smooth canonical stack:
\[
J(\mathrm X):=
\frac{\mathrm H^3(\widetilde{\mathrm X},\mathbb C)}
{F^2\mathrm H^3(\widetilde{\mathrm X},\mathbb C)
+\mathrm H^3(\widetilde{\mathrm X},\mathbb Z)_{\mathrm{free}}}.
\]
Here $F^\bullet$ is the Hodge filtration. The Abel--Jacobi map of the degree-five family is computed on its strict transforms in $\mathrm Y$, with a fixed general member as base cycle.

\begin{proposition}\label{prop:cubic-curve-families}
For $\mathrm X=\mathrm X_{1.9}$, let $\mathcal F_3,\mathcal F_5$ be the reduced closures in the Chow variety of the families obtained from general fibres of $\operatorname{Exc}(\mathrm f')\to\mathrm C$ and from general lines on $\mathrm X'$, respectively. Then $\mathrm X$ has six integral degree-one rational curves. The normalization $\mathcal F_3^{\nu}$ is isomorphic to $\mathrm C\simeq\mathbb P^1$, and strict transform gives a birational correspondence $\mathcal F_5\dashrightarrow F_1(\mathrm X')$. Over $\mathbb C$, if $S_5$ is a smooth projective model of $\mathcal F_5$, the quintic family induces an Abel--Jacobi isomorphism
\[
\mathrm{AJ}:\operatorname{Alb}(S_5)\xrightarrow{\sim}J(\mathrm X_{1.9}).
\]
\end{proposition}
\begin{proof}
Here $z=1$, there are no flips, and $d,s$ are integers. For $k=1$, \eqref{eq:hourglass-curve-degree} forces $a$ to be odd, hence positive, and then $d<0$. Thus the six flop images are the only integral degree-one rational curves.

A general fibre of $\operatorname{Exc}(\mathrm f')\to\mathrm C$ avoids the flop centres and has $a=1$, $d=s=0$, so its image has degree three. Strict transform identifies the general fibres with the general members of $\mathcal F_3$. Hence $\mathcal F_3^\nu\simeq\mathrm C\simeq\mathbb P^1$.

A general line on $\mathrm X'$ avoids $\mathrm C$ and the six bisecant lines. Its transform has $d=1$ and $m=s=0$, so \eqref{eq:hourglass-target-degree} gives degree five. The inverse strict transform recovers the line, proving $\mathcal F_5\dashrightarrow F_1(\mathrm X')$ birational.

Put $U=\mathrm X\setminus\{\mathrm p\}=\mathrm Y\setminus\mathrm A$. The Gysin sequences give
\[
\mathrm H^3(\widetilde{\mathrm X},\mathbb Z)
\xrightarrow{\sim}\mathrm H^3(U,\mathbb Z)
\xleftarrow{\sim}\mathrm H^3(\mathrm Y,\mathbb Z).
\]
For the second isomorphism, the Gysin map $\mathrm H^2(\mathrm A)(-1)\to\mathrm H^4(\mathrm Y)$ is injective: if $\iota:\mathrm A\hookrightarrow\mathrm Y$ and $h$ is the hyperplane class, then $\iota^*\iota_*h=-2h^2$. The identifications preserve the Hodge structures and pairings, computed using compactly supported representatives on $U$. Thus $J(\mathrm X)\simeq J(\mathrm Y)$.

Blowing up $\mathrm C\simeq\mathbb P^1$ and the common blow-ups of the six standard flops add no degree-three cohomology, so $J(\mathrm Y)\simeq J(\mathrm X')$. Proper pushforward from the common resolution identifies the Abel--Jacobi homomorphisms of the two families. Birational invariance of the Albanese variety and \cite[Theorem~11.19]{CG1972} now give the asserted isomorphism for $S_5$.
\end{proof}

\begin{corollary}
The base $\mathrm B$ of $\mathrm X_{1.9}$ and the cubic Fano surface $F_1(\mathrm X'_{1.9})$ are minimal surfaces of general type and are not birational.
\end{corollary}
\begin{proof}
Put $H:=-\mathrm{K}_{\mathrm X_{1.9}}$. Since $\mathrm B\in|2H|$ is smooth and avoids the quotient point, adjunction gives $\mathrm{K}_{\mathrm B}=H|_{\mathrm B}$, which is ample. The divisor sequence and vanishing on the Fano threefold give
\[
\mathrm H^1(\mathrm B,\mathcal{O}_{\mathrm B})\simeq\mathrm H^2(\mathrm X_{1.9},\mathcal{O}_{\mathrm X_{1.9}}(-2H))\simeq\mathrm H^1(\mathrm X_{1.9},\mathcal{O}_{\mathrm X_{1.9}}(H))^{\vee}=0.
\]
The cubic Fano surface is minimal of general type and has irregularity $5$ \cite{CG1972}. Hence the two surfaces are not birational.
\end{proof}

On the other hand, we have the following exceptional collections.
\begin{proposition}\label{prop:ECof1.9}
For any stacky lines $L_{i}$, $1\leq i\leq 6$ on $\widetilde{\mathrm{X}}_{1.9}$,
\begin{enumerate}
  \item\label{enu:ECof1.91} the order sequence
$$\langle\mathcal{O}_{L_{i}}(-1),\mathcal{O}_{\widetilde{\mathrm{X}}_{1.9}},\mathrm{U}_{\widetilde{\mathrm{X}}_{1.9}}^{\vee}\rangle$$
forms an exceptional collection.
  \item\label{enu:ECof1.92} while
  $$\langle\mathcal{O}_{L_{i}}(-1/2),\mathcal{O}_{\widetilde{\mathrm{X}}_{1.9}},\mathrm{U}_{\widetilde{\mathrm{X}}_{1.9}}^{\vee}\rangle$$
   is not semi-orthogonal.
   \item\label{enu:ECof1.93} and $\mathrm{U}_{\widetilde{\mathrm{X}}_{1.9}}^{\vee}$ is not a locally free sheaf.
\end{enumerate}

\end{proposition}
\begin{proof}
For (\ref{enu:ECof1.91}), by  Proposition \ref{prop:X1.9} the only non-trivial part left for us is to show the orthogonality between $\mathcal{O}_{L_{i}}(-1)$ and $\mathrm{U}_{\widetilde{\mathrm{X}}_{1.9}}^{\vee}$. Noticing $-\mathrm{K}_{\mathrm{X}_{1.9}}L_{i}=1/2$, after tensoring $\mathcal{O}_{\widetilde{\mathrm{X}}}(-\mathrm{K}_{\widetilde{\mathrm{X}}_{1.9}})$, it's enough for us to show the orthogonality of pair,
$$\langle\mathcal{O}_{L_{i}}(-1/2),\mathrm{U}_{\widetilde{\mathrm{X}}_{1.9}}^{\vee}(-\mathrm{K}_{\widetilde{\mathrm{X}}_{1.9}})\rangle$$
since $\mathrm{\Theta}:=\varpi_{*}\pi^{*}$ is fully-faithful by Theorem \ref{thm:singularweightedblowup}, it can be reduced to show the orthogonality of pair,

$$\langle\mathrm{\Theta}(\mathcal{O}_{L_{i}}(-1/2)),\mathrm{\Theta}(\mathrm{U}_{\widetilde{\mathrm{X}}_{1.9}}^{\vee}(-\mathrm{K}_{\widetilde{\mathrm{X}}_{1.9}}))\rangle$$
recall (\ref{enu:exp1}) in Proposition \ref{prop:X1.9}, we have $\mathrm{U}_{\widetilde{\mathrm{X}}_{1.9}}^{\vee}(-\mathrm{K}_{\widetilde{\mathrm{X}}_{1.9}})=\mathrm{\Theta}^{!}\big(g^{*}\mathrm{U}_{\mathrm{Z}}^{\vee}(-\mathrm{K}_{\mathrm{Y}})\big)$, and $g^{*}\mathrm{U}_{\mathrm{Z}}^{\vee}(-\mathrm{K}_{\mathrm{Y}})$ is also orthogonal to $\mathcal{O}_{\mathrm{A}}(-1)$ by  (\ref{sodc9}), so we can see

$$\mathrm{\Theta}(\mathrm{U}_{\widetilde{\mathrm{X}}_{1.9}}^{\vee}(-\mathrm{K}_{\widetilde{\mathrm{X}}_{1.9}}))=\mathrm{\Theta}\,\mathrm{\Theta}^{!}\big(g^{*}\mathrm{U}_{\mathrm{Z}}^{\vee}(-\mathrm{K}_{\mathrm{Y}})\big)=g^{*}\mathrm{U}_{\mathrm{Z}}^{\vee}(-\mathrm{K}_{\mathrm{Y}})$$
and recall Example \ref{ex:fanciaprojection}, we have an exact sequence: $$\mathrm{\Theta}(\mathcal{O}_{L_{i}}(-1/2))\longrightarrow\mathcal{O}_{\widehat{L_{i,\mathrm{Y}}}}(-1)\longrightarrow\mathcal{O}_{\mathrm{A}}(-1)[2]$$
so $$\mathrm{Ext}^{*}(\mathrm{\Theta}(\mathrm{U}_{\widetilde{\mathrm{X}}_{1.9}}^{\vee}(-\mathrm{K}_{\widetilde{\mathrm{X}}_{1.9}})),\mathrm{\Theta}(\mathcal{O}_{L_{i}}(-1/2)))=\mathrm{Ext}^{*}(g^{*}\mathrm{U}_{\mathrm{Z}}^{\vee}(-\mathrm{K}_{\mathrm{Y}}),\mathcal{O}_{\widehat{L_{i,\mathrm{Y}}}}(-1))=0
$$
since $\mathrm{\Theta}^{!}\mathcal{O}_{\mathrm{A}}(-1)=0$ and $g_{*}\mathcal{O}_{\widehat{L_{i,\mathrm{Y}}}}(-1)=0$, we finish (\ref{enu:ECof1.91}).\\

For (\ref{enu:ECof1.92}), let's consider the ideal sheaf of stack point at
$L_{i}$:
\begin{equation}\label{equ:computaionexcseq1}
0\longrightarrow\mathcal{O}_{L_{i}}(-1/2)\longrightarrow\mathcal{O}_{L_{i}}\longrightarrow\mathcal{O}_{\mathrm{p}}\longrightarrow0
\end{equation}
act functor $\mathrm{\Theta}$, we have distinguished triangle:
\begin{equation}\label{equ:computaionexcseq2}\mathrm{\Theta}(\mathcal{O}_{L_{i}}(-1/2))\longrightarrow\mathrm{\Theta}(\mathcal{O}_{L_{i}})\longrightarrow\mathrm{\Theta}(\mathcal{O}_{\mathrm{p}})\end{equation}
about the first term by our computation in Example \ref{ex:fanciaprojection}, it fits in distinguished triangle
\begin{equation}\label{equ:computaionexcseq3}\mathrm{\Theta}(\mathcal{O}_{L_{i}}(-1/2))\longrightarrow\mathcal{O}_{\widehat{L_{i,\mathrm{Y}}}}(-1)\longrightarrow\mathcal{O}_{\mathrm{A}}(-1)[2]\end{equation}
and by Proposition \ref{prop:pullbackcenterlocal}, we have
\begin{equation}\label{equ:computaionexcseq4}\mathcal{O}_{\mathrm{A}}(-1)[2]\longrightarrow\mathrm{\Theta}(\mathcal{O}_{\mathrm{p}})\longrightarrow\mathcal{O}_{\mathrm{A}}\end{equation}
We adopt the same strategy as above (\ref{enu:ECof1.91}), the semi-orthogonality of pair
$$\langle\mathcal{O}_{L_{i}}(-1/2),\mathrm{U}_{\widetilde{\mathrm{X}}_{1.9}}^{\vee}\rangle$$
is  equivalent to seeing the semi-orthogonality of the pair
$$\langle\mathrm{\Theta}(\mathcal{O}_{L_{i}}),\mathrm{\Theta}(\mathrm{U}_{\widetilde{\mathrm{X}}_{1.9}}^{\vee}(-\mathrm{K}_{\widetilde{\mathrm{X}}_{1.9}}))\rangle$$
the second term is just $g^{*}\mathrm{U}_{\mathrm{Z}}^{\vee}(-\mathrm{K}_{\mathrm{Y}})$, so the only contribution term to the extension is
\begin{equation}\label{equ:computaionexcseq5}\mathrm{Ext}^{*}(\mathrm{\Theta}(\mathrm{U}_{\widetilde{\mathrm{X}}_{1.9}}^{\vee}(-\mathrm{K}_{\widetilde{\mathrm{X}}_{1.9}})),\mathrm{\Theta}(\mathcal{O}_{L_{i}}))=\mathrm{Ext}^{*}(g^{*}\mathrm{U}_{\mathrm{Z}}^{\vee}(-\mathrm{K}_{\mathrm{Y}}),\mathcal{O}_{\mathrm{A}})=\mathrm{Ext}^{*}(\mathrm{T}_{\mathds{P}^{2}},\mathcal{O}_{\mathds{P}^{2}})=\mathrm{k}[-1]\neq 0
\end{equation}
by considering our distinguished triangles in (\ref{equ:computaionexcseq2}) (\ref{equ:computaionexcseq3}) and (\ref{equ:computaionexcseq4}), we finish (\ref{enu:ECof1.92}).\\

For (\ref{enu:ECof1.93}), we assume $\mathrm{U}_{\widetilde{\mathrm{X}}_{1.9}}^{\vee}$ is locally free, hence by \cite{MT} its restriction on any stacky line $L_{i}$ is a direct sum of line bundles, i.e. $$\mathrm{U}_{\widetilde{\mathrm{X}}_{1.9}}^{\vee}|_{L_{i}}=\mathcal{O}_{L_{i}}(-a/2)\bigoplus\mathcal{O}_{L_{i}}(-b/2)$$
for two integers $a$ and $b$. By (\ref{enu:ECof1.91}) and (\ref{equ:computaionexcseq5}) we have
$$\mathrm{Ext}_{L_{i}}^{*}(\mathcal{O}_{L_{i}}(-a/2)\bigoplus\mathcal{O}_{L_{i}}(-b/2),\mathcal{O}_{L_{i}}(-1))=0$$
while $$\mathrm{Ext}_{L_{i}}^{*}(\mathcal{O}_{L_{i}}(-a/2)\bigoplus\mathcal{O}_{L_{i}}(-b/2),\mathcal{O}_{L_{i}}(-1/2))=\mathrm{k}[-1]$$
an easy cohomology computation using Lemma \ref{lemma:wpscoho} yields a contradiction.
\end{proof}

By Lemma~\ref{lemma:lf}, $\mathrm{U}_{\widetilde{\mathrm{X}}_{1.9}}^{\vee}$ is a rank-two sheaf, locally free away from $\mathrm p$. It is torsion-free: it is the pushforward along the dominant morphism $\pi$ of the locally free sheaf $\varpi^*g^*\mathrm U_{\mathrm Z}^{\vee}$, and multiplication by a nonzero local function remains injective after pushforward. Its natural map to its reflexive hull is therefore injective. Put
$${\mathrm{U}_{\widetilde{\mathrm{X}}_{1.9}}^{\vee}}_{,\circ}:=\bigl(\mathrm{U}_{\widetilde{\mathrm{X}}_{1.9}}^{\vee}\bigr)^{**}.$$
We have an exact sequence
$$0\longrightarrow\mathrm{U}_{\widetilde{\mathrm{X}}_{1.9}}^{\vee}\longrightarrow{\mathrm{U}_{\widetilde{\mathrm{X}}_{1.9}}^{\vee}}_{,\circ}\longrightarrow\tau_2\longrightarrow0,$$
where $\tau_2$ is supported at $\mathrm p$. The hull is reflexive and locally free away from $\mathrm p$; reflexivity alone does not imply local freeness at $\mathrm p$. In particular, $\mathrm{U}_{\widetilde{\mathrm{X}}_{1.9}}^{\vee}$ is reflexive precisely when $\tau_2=0$.
\begin{definition}\label{def:Ai}
We denote the left orthogonal component of exceptional collection $$\langle\mathcal{O}_{L_{i}}(-1),\mathcal{O}_{\widetilde{\mathrm{X}}_{1.9}},\mathrm{U}_{\widetilde{\mathrm{X}}_{1.9}}^{\vee}\rangle$$ by $\mathcal{MF}_{i}(\mathrm{X}_{1.9})$ for any $i$, so we have another semi-orthogonal decomposition of  $\widetilde{\mathrm{X}}_{1.9}$,
$$\mathrm{D}^{b}(\widetilde{\mathrm{X}}_{1.9})\simeq \langle\mathcal{O}_{L_{i}}(-1), \mathcal{O}_{\widetilde{\mathrm{X}}_{1.9}},\mathrm{U}_{\widetilde{\mathrm{X}}_{1.9}}^{\vee},\mathcal{MF}_{i}(\mathrm{X}_{1.9})\rangle$$
\end{definition}

\begin{proposition}
For every $i$ there is a noncanonical isomorphism of abelian groups
\[
\mathrm K_0(\mathcal{MF}_i(\mathrm X_{1.9}))
\simeq
\mathrm K_0(\mathcal{MF}(\mathrm X'_{1.9})).
\]
\end{proposition}
\begin{proof}
The decompositions in Proposition~\ref{prop:X1.9} and
Definition~\ref{def:Ai} have three exceptional objects in addition to
their respective residual categories. Additivity of the Grothendieck
group for a semi-orthogonal decomposition therefore gives
\[
\begin{aligned}
\mathrm K_0(\mathcal{MF}_i(\mathrm X_{1.9}))\oplus\mathds Z^3
&\simeq\mathrm K_0(\mathrm{D}^{b}(\widetilde{\mathrm X}_{1.9}))\\
&\simeq\mathrm K_0(\mathcal{MF}(\mathrm X'_{1.9}))\oplus\mathds Z^3.
\end{aligned}
\]
A finitely generated abelian direct summand can be cancelled by the
Cohn--Walker cancellation theorem \cite{Cohn1956}. This proves the
abstract isomorphism; it does not identify an equivalence of the residual
categories.
\end{proof}

\begin{remark}
The six curve objects $\mathcal{O}_{L_i}(-1/2)$ have Euler pairing
$\chi(\mathcal{O}_{L_i}(-1/2),\mathcal{O}_{L_j}(-1/2))=1$ for all $i,j$,
by the Ext calculation below. Thus their displayed Euler matrix has rank
one. This is a numerical statement. It does not compute the intrinsic
Grothendieck group of their generated category, nor the torsion in its
image in $\mathrm K_0(\widetilde{\mathrm X}_{1.9})$.
\end{remark}

So it is natural for us to speculate that
\begin{conjecture}\label{conj:cubic}
For any $i$, $\mathcal{MF}_{i}(\mathrm{X}_{1.9})$ is equivalent to $\mathcal{MF}(\mathrm{X}'_{1.9})$.
\end{conjecture}

\subsection{Non-commutative projection}

This section is motivated by work on geometric unprojection, for example \cite{BKR}, we consider the relationship between the derived category of a general hourglass and the derived category of a midpoint of its Sarkisov link. The $e(\mathrm{X}_{1,i})$ flopping curves give stacky lines on $\mathrm{X}_{1,i}$ through the cyclic quotient centre $\mathrm p$ of the weighted blow-up. We denote them by $L_{i}$, for $1\leq i\leq e$, and their union as $L$.  We consider the ideal sheaf of the singularity point $\mathrm{p}$ on $L_{i}$, and denote them by $\mathcal{O}_{L_{i}}(-1/2)$. We omit the subscript $(1,i)$ of $\mathrm{X}_{1,i}$ in the following.
\begin{lemma}\begin{enumerate}
\item\label{equ:ecpOl} $\mathcal{O}_{L_{i}}(-1/2)$ is exceptional for any $i$.
\item\label{equ:ecpO2} For any $i$, $j$ such that $i\neq j$,

\begin{equation}
 \mathrm{Ext}_{\widetilde{\mathrm{X}}}^{l}(\mathcal{O}_{L_{i}}(-1/2),\mathcal{O}_{L_{j}}(-1/2)) =
    \begin{cases}
     \mathrm{k}  & \text{if\, $l=2$}\\
      0 & \text{else}\\
    \end{cases}
\end{equation}
so $\mathcal{O}_{L_{i}}(-1/2)$ and $\mathcal{O}_{L_{j}}(-1/2)$ are not mutually orthogonal if $i\neq j$.
\end{enumerate}
\end{lemma}
\begin{proof}
  For (\ref{equ:ecpOl}), noticing after the weighted blow-up it follows a standard flop which contracts the strict transformation of $L_{i}$, its geometry is the same as in Example \ref{ex:fanciaflip},  so we follow a same argument,
  $$\mathcal{N}_{L_{i}/\widetilde{\mathrm{X}}}=\mathcal{O}_{L_{i}}(-1/2)\oplus\mathcal{O}_{L_{i}}(-1/2)$$
  combine with spectral sequence
  $$E^{p,q}_{2}:=\mathrm{H}^{p}(L_{i},\bigwedge^{q}\mathcal{N}_{L_{i}/\widetilde{\mathrm{X}}})\Longrightarrow\mathrm{Ext}^{*}(\mathcal{O}_{L_{i}}(-1/2),\mathcal{O}_{L_{i}}(-1/2))$$
and $L_{i}=\mathcal{P}^{1}_{2}$ with $\mathrm{K}_{L_{i}}=\mathcal{O}_{L_{i}}(-3/2)$, we get the result.\\

 For (\ref{equ:ecpO2}), since we have a distinguished triangle,
\begin{equation}\label{eq:lstirct}
\mathrm{\Theta}(\mathcal{O}_{L_{i}}(-1/2))\longrightarrow\mathcal{O}_{\widehat{L_{i}}}(-1)\longrightarrow\mathcal{O}_{\mathrm{A}}(-1)[2]
\end{equation}
and $\mathrm{\Theta}$ is fully-faithful, we have

$$\mathrm{Ext}^{*}(\mathcal{O}_{L_{i}}(-1/2),\mathcal{O}_{L_{j}}(-1/2))=\mathrm{Ext}^{*}(\mathrm{\Theta}(\mathcal{O}_{L_{i}}(-1/2)),\mathrm{\Theta}(\mathcal{O}_{L_{j}}(-1/2)))$$
then a simple spectral sequence and
\begin{enumerate}
  \item $\mathrm{Ext}_{\mathrm{Y}}^{*}(\mathcal{O}_{\widehat{L_{i}}}(-1),\mathcal{O}_{\widehat{L_{j}}}(-1))=0$, since flopping curves are disjoint in $\mathrm{Y}$.
  \item $\mathrm{Ext}_{\mathrm{Y}}^{*}(\mathcal{O}_{\mathrm{A}}(-1)[2],\mathcal{O}_{\mathrm{A}}(-1)[2])=\mathrm{k}$, since $\mathcal{O}_{\mathrm{A}}(-1)$ is exceptional.
  \item  $\mathrm{Ext}_{\mathrm{Y}}^{*}(\mathcal{O}_{\mathrm{A}}(-1)[2],\mathcal{O}_{\widehat{L_{j}}}(-1))=\mathrm{k}[-3]$,
   \item  $\mathrm{Ext}_{\mathrm{Y}}^{*}(\mathcal{O}_{\widehat{L_{i}}}(-1),\mathcal{O}_{\mathrm{A}}(-1)[2])=\mathrm{k}$,
\end{enumerate}
give us the result.
\end{proof}

\begin{definition} For a general noncommutative hourglass $\mathrm{D}^{b}(\widetilde{\mathrm{X}})$,\\
(1)  we denote the subcategory generated by $\mathcal{O}_{L_{k}}(-1/2)$ as $\mathcal{L}_{k}$, its left orthogonal component by $\mathcal{B}_{k}$,\\
(2) the subcategory generated by $\mathcal{O}_{L_{i}}(-1/2)$ for $1\leq i\leq e$ by $\mathcal{L}$, its left orthogonal component by $\mathcal{B}$.
\end{definition}

\begin{proposition}[Non-commutative projection]\label{prop:projection}
We consider the first midpoint $\mathrm{Z}$ of our Sarkisov link, it admits an equivalence up to quotient via projection:
$$g_{*}\mathrm{\Theta}:\mathrm{D}^{b}(\widetilde{\mathrm{X}})/\mathcal{L}\simeq\mathrm{D}^{b}(\mathrm{\widetilde{Z}})/\langle\mathcal{O}_{\mathrm{A}}(-1)\rangle $$
where $\mathrm{\widetilde{Z}}$ is stacky smoothing of $\mathrm{Z}$ at orbifold points.
\end{proposition}

\begin{proof}
The Kawamata blow-up formula gives
\[
\mathrm{D}^{b}(\widetilde{\mathrm Y})=
\langle\mathcal{O}_{\mathrm A}(-1),
\mathrm{\Theta}\,\mathrm{D}^{b}(\widetilde{\mathrm X})\rangle.
\]
Thus $\mathrm{\Theta}^{!}$ identifies the quotient by
$\langle\mathcal{O}_{\mathrm A}(-1)\rangle$ with
$\mathrm{D}^{b}(\widetilde{\mathrm X})$. Here brackets denote the thick
triangulated subcategory generated by the displayed objects.
By the triangles \eqref{eq:lstirct}, the images of
$\mathcal{O}_{\widehat{L_{i,\mathrm Y}}}(-1)$ in this quotient are
$\mathcal{O}_{L_i}(-1/2)$, which generate $\mathcal L$.

On the other hand, the contraction $g$ has only the finitely many ordinary
exceptional curves $\widehat{L_{i,\mathrm Y}}$, and
$g_*\mathcal{O}_{\widetilde{\mathrm Y}}\simeq
\mathcal{O}_{\widetilde{\mathrm Z}}$. Hence
\cite[Proposition 2.10]{Xie} gives
\[
\mathrm{D}^{b}(\widetilde{\mathrm Y})/
\langle\mathcal{O}_{\widehat{L_{i,\mathrm Y}}}(-1)\rangle_i
\simeq\mathrm{D}^{b}(\widetilde{\mathrm Z}).
\]
The cited localization applies on the ordinary double point charts.
Away from these charts $g$ is an isomorphism; the unchanged quotient
charts are treated equivariantly. This gives the same localization for
the canonical stacks. Moreover,
$g_*\mathcal{O}_{\mathrm A}(-1)=
\mathcal{O}_{\mathrm A_{\widetilde{\mathrm Z}}}(-1)$.
Transitivity of Verdier quotients now gives
\[
\begin{aligned}
\mathrm{D}^{b}(\widetilde{\mathrm X})/\mathcal L
&\simeq
\frac{\mathrm{D}^{b}(\widetilde{\mathrm Y})}
{\langle\mathcal{O}_{\mathrm A}(-1),
\mathcal{O}_{\widehat{L_{i,\mathrm Y}}}(-1)\rangle_i}\simeq\mathrm{D}^{b}(\widetilde{\mathrm Z})/
\langle\mathcal{O}_{\mathrm A_{\widetilde{\mathrm Z}}}(-1)\rangle.
\end{aligned}
\]
The induced functor is $g_*\mathrm{\Theta}$. This uses transitivity of
Verdier quotients without assuming that $\mathcal L$ or
$\langle\mathcal{O}_{\mathrm A_{\widetilde{\mathrm Z}}}(-1)\rangle$
is admissible.
\end{proof}

Since $\mathcal{O}_{L_{k}}(-1/2)$ is exceptional for any $k$, so $\mathcal{B}_{k}$ is admissible, $\mathrm{D}^{b}(\widetilde{\mathrm{X}})\simeq\langle\mathcal{L}_{k},\mathcal{B}_{k}\rangle$. But we don't know whether $\mathcal{L}$ or $\mathcal{B}$ is admissible or not, since $i$ is generally larger than $1$, we have $\langle\mathcal{L},\mathcal{B}\rangle$ is merely a subcategory of $\mathrm{D}^{b}(\widetilde{\mathrm{X}})$.\\

The restriction to the orthogonal component must retain the perfectness condition on the midpoint. We use the following form of perfect localization, rather than infer local freeness of an arbitrary complex from its restriction to the reduced exceptional fibre.
\begin{lemma}[{\cite[Prop. 5.5 and Prop. 6.1]{KuS}}]
The functor $g^*$ induces an equivalence
\[
\mathrm{D}^{\mathrm{perf}}(\widetilde{\mathrm Z})
\xrightarrow{\sim}{}^\perp\ker g_*
\subset\mathrm{D}^{b}(\widetilde{\mathrm Y}),
\]
with inverse $g_*$.
\end{lemma}
\begin{proof}
The contraction is a categorical localization by
\cite[Proposition 2.10]{Xie}. Since $\widetilde{\mathrm Y}$ is smooth,
\cite[Proposition 5.5(iii), Proposition 6.1 and Example 6.3]{KuS}
identify the essential image of $g^*$ with the indicated left orthogonal.
The calculation is on schemes near the ODPs; at the unchanged quotient
charts the morphism is an isomorphism.
\end{proof}

The preceding lemma identifies the relevant perfect orthogonal on the midpoint. The resulting unprojection statement is as follows:

\begin{proposition}\label{prop:sgequi}
If we restrict the equivalent functor
$$g_{*}\mathrm{\Theta}:\mathrm{D}^{b}(\widetilde{\mathrm{X}})/\mathcal{L}\isomto \mathrm{D}^{b}(\mathrm{\widetilde{Z}})/\langle\mathcal{O}_{\mathrm{A}}(-1)\rangle$$
in Proposition \ref{prop:projection} on subcategory $\mathcal{B}$ it yields another equivalence
$$g_{*}\mathrm{\Theta}:\mathcal{B}\isomto \mathcal{B}_{\mathrm{\widetilde{Z}}}^{perf}$$
where $\mathcal{B}_{\mathrm{\widetilde{Z}}}^{perf}$  is equivalent to the subcategory of $^{\perp}\langle\mathcal{O}_{\mathrm{A}}(-1)\rangle$ consists of perfect complexes on $\widetilde{\mathrm{Z}}$.
\end{proposition}
\begin{proof}
The perfect localization theorem
\cite[Proposition 5.5(iii), Proposition 6.1 and Example 6.3]{KuS} gives
\[
g^*:\mathrm{D}^{\mathrm{perf}}(\widetilde{\mathrm Z})
\xrightarrow{\sim}{}^\perp\langle\mathcal{O}_{\widehat{L_{i,\mathrm Y}}}(-1)\rangle_i
\subset\mathrm{D}^{b}(\widetilde{\mathrm Y}),
\]
with inverse $g_*$. Here $\widetilde{\mathrm Y}$ is smooth, so every bounded
coherent complex on it is perfect. The result applies near each ordinary
double point; at the remaining quotient charts $g$ is an isomorphism.

If $b\in\mathcal B$, then
$\mathrm\Theta b\in{}^\perp\langle\mathcal{O}_{\mathrm A}(-1)\rangle$. Applying
$\mathrm{Ext}^*(\mathrm\Theta b,-)$ to \eqref{eq:lstirct} and using full
faithfulness of $\mathrm\Theta$ gives
$\mathrm\Theta b\in{}^\perp\langle\mathcal{O}_{\widehat{L_{i,\mathrm Y}}}(-1)\rangle_i$ as well. Hence
$\mathrm\Theta b\simeq g^*F$ for a unique
$F\in\mathrm{D}^{\mathrm{perf}}(\widetilde{\mathrm Z})$. By adjunction,
\[
\mathrm{Ext}^*(F,\mathcal{O}_{\mathrm A_{\widetilde{\mathrm Z}}}(-1))
\simeq
\mathrm{Ext}^*(g^*F,\mathcal{O}_{\mathrm A}(-1))=0.
\]
Thus $g_*\mathrm\Theta$ takes $\mathcal B$ into
$\mathcal B^{\mathrm{perf}}_{\widetilde{\mathrm Z}}$.

Conversely, take a perfect $F$ satisfying this vanishing.
The perfect localization theorem gives
$g^*F\in{}^\perp\langle\mathcal{O}_{\widehat{L_{i,\mathrm Y}}}(-1)\rangle_i$, and adjunction gives
$g^*F\in{}^\perp\langle\mathcal{O}_{\mathrm A}(-1)\rangle$.
The Kawamata decomposition therefore places $g^*F$ in the image of
$\mathrm\Theta$. Writing $b:=\mathrm\Theta^!g^*F$, the triangles
\eqref{eq:lstirct} show that
$\mathrm{Ext}^*(b,\mathcal{O}_{L_i}(-1/2))=0$ for every $i$, so
$b\in\mathcal B$. The two constructions are inverse by the adjunction
identities for $g^*$ and $\mathrm\Theta$. Hence the claimed equivalence
has inverse $\mathrm\Theta^!g^*$.
\end{proof}

For the multi-centre construction below, we impose the following
geometric hypotheses; the diagrams do not assert the existence of a
common model from the numerical table alone.
\begin{enumerate}
\item The chosen centres are distinct $\frac12(1,1,1)$ points of
$\mathrm{X}_{1,i}$. Their simultaneous projective Kawamata extraction
exists and its canonical stack is smooth. The exceptional planes
corresponding to the centres extracted on that model are disjoint;
on the full extraction they are $\mathrm{A}_1,\ldots,\mathrm{A}_N$. The functor
$\mathrm{\Theta}$ denotes the composite of the extraction functors
already defined.
\item On each extracted model used below, the anticanonical divisor is
semiample. A sufficiently divisible multiple defines a projective small
contraction $g$ to a normal midpoint. Its nontrivial fibres are finitely
many pairwise disjoint smooth rational curves, each with normal bundle
$\mathcal{O}(-1)^{\oplus2}$ and contracted to a distinct ordinary double
point. Away from them, $g$ is an isomorphism, also on the unchanged
quotient charts.
\item Any further centres used in \eqref{exsark3} lie in this common
isomorphism locus. The blow-ups on the two sides are taken at the
corresponding centres in these isomorphic neighbourhoods.
\end{enumerate}
The ordinary double point models give
$g_*\mathcal{O}\simeq\mathcal{O}$, with derived pushforward understood,
and the localization and perfect-orthogonal comparison used in
Propositions~\ref{prop:projection} and~\ref{prop:sgequi} apply.
The plane sheaves on the midpoint mean the images
$g_*\mathcal{O}_{\mathrm{A}_j}(-1)$, as in the one-centre argument.
We assume these hypotheses whenever the following global midpoint
diagrams or the multiple-unprojection corollary are used.

\begin{equation}\label{exsark}
\begin{tikzcd}[column sep=0.5em]
  & && \mathrm{Y}^{G}_{1,i}:=wBl_{\mathrm{p}_{1},...,\mathrm{p}_{N}}\mathrm{\widetilde{X}}_{1,i}\arrow[dl]\arrow[dr,"g"] &&  \\
 & &\mathrm{\widetilde{X}}_{1,i}&  &\mathrm{Z}^{G}_{1,i}&
\end{tikzcd}\end{equation}
and we take some examples from \cite{Tak2}:

\begin{enumerate}
  \item\label{equ:unpo2plan} For an hourglass $\mathrm{X}_{1,1}$ it has two $(\frac{1^{3}}{2})$ type orbifold points, and we weighted blow-up all of them, contract it along $-\mathrm{K}_{\mathrm{Y}}$ to a midpoint $\mathrm{Z}$, it can be realized as a Gorenstein non-factorial $A_{6}$, e.g.(\ref{a6}), and we have two disjoint unprojection planes,
      $$\mathrm{A}_{1}:=\{M|M\subset V_{3},\, \mathrm{rank} (M)=2\}\simeq\mathds P^2\subset\mathbf{Gr}(2,4)$$
      and
      $$\mathrm{A}_{2}:=\{M|V_{1}\subset M\subset V_{4}, \, \mathrm{rank} (M)=2\}\simeq\mathds P^2\subset\mathbf{Gr}(2,4)$$
      where $V_{1}\bigoplus V_{3}=V_{4}$.
  \item\label{equ:unpo2plan1} For $\mathrm{X}_{1,4}$, we weighted blow up the orbifold point then contract it to a Gorenstein non-factorial $A_{10}$ e.g.(\ref{a8}), with a trivial plane in $\mathds{P}^{6}$,
\item\label{equ:unpo2plan2} For $\mathrm{X}_{1,9}$, we weighted blow up the orbifold point then contract it to a Gorenstein non-factorial $A_{10}$ e.g.(\ref{a10}), with a $\sigma_{2,2}$ plane,
       $$\mathrm{A}_{1}:=\{M|M\subset V_{3}, \, \mathrm{rank} (M)=2\}\simeq\mathds P^2\subset\mathbf{Gr}(2,5)$$
  \item\label{equ:unpo2plan3} For $\mathrm{X}_{1,10}$, we weighted blow up the orbifold point then contract it to a Gorenstein non-factorial $A_{10}$ e.g.(\ref{a10}), with a $\sigma_{3,1}$ plane,
       $$\mathrm{A}_{2}:=\{M|V_{1}\subset M\subset V_{4}, \, \mathrm{rank} (M)=2\}\simeq\mathds P^2\subset\mathbf{Gr}(2,5)$$
  \item\label{equ:unpo2plan4} For $\mathrm{X}_{1,13}$, we weighted blow up the orbifold point then contract it to a Gorenstein non-factorial $A_{14}$ e.g.(\ref{a14}), with one unprojection plane,
       $$\mathrm{A}_{2}:=\{M|V_{1}\subset M\subset V_{4}, \, \mathrm{rank} (M)=2\}\simeq\mathds P^2\subset\mathbf{Gr}(2,6)$$
  \item For the remaining hourglasses, the expected midpoint has type
$A_{\mathrm{Deg}(\mathrm{X}_{1,i})-N(\mathrm{X}_{1,i})/2}$, as in
\eqref{targetZ}. We use this model only when the extraction and small
contraction satisfy the hypotheses above; see also
\cite[Proposition 2.6]{Tak2}.
\end{enumerate}
According to the work of S.Mukai and a recent work by A.Bayer, A.Kuznetsov, E.Macrì, we know if $A_{2g-2}$ are smooth we have some semi-orthogonal decompositions of $\mathrm{D}^{b}(A_{2g-2})$:

\begin{proposition}\cite[Lemma 3.6]{KuFano}\label{kufanolemma}
If  $A_{2g-2}$ be a  smooth Fano $3$ fold of genus $g=2t$, there is a locally free sheaf $\mathrm{U}_{2g-2}$ of rank 2 on $A_{2g-2}$  constructed from its key variety, such that
$$\langle\mathrm{U}_{2g-2},\mathcal{O}_{A_{2g-2}}\rangle$$
forms an exceptional pair.
\end{proposition}
In addition, we also have mature similar results for smooth $A_{8}$ and $A_{12}$, e.g.\cite{KuV12}. However, the author is not clear whether the above properties hold for the Gorenstein non-factorial cases or not. One idea is that we expect these midpoints $\mathrm{Z}^{G}$ to be the complete intersection of some smooth high-dimensional varieties as in \cite[Proposition 2.7]{Tak2}, and on these varieties we have very symmetric decompositions of exceptional collections e.g. (\ref{equ:b5sod2}) that induce $\langle\mathrm{U}_{\mathrm{Z}^{G}},\mathcal{O}_{\mathrm{Z}^{G}}\rangle$ is still an exceptional pair on $\mathrm{Z}_{i}^{G}$ for $i=9,10$, which also inherits the spirit of Theorem \ref{thm:orlov1}.\\

 Suppose that the restrictions of the tautological bundles give an
exceptional pair $\langle\mathrm U_{\mathrm Z},\mathcal{O}_{\mathrm Z}\rangle$
on the singular midpoint. The restrictions to the two types of plane give
\[
\begin{aligned}
\mathrm{Ext}^*(\mathrm U_{\mathrm Z},\mathcal{O}_{\mathrm A_1}(-1))
&=\mathrm H^*(\mathcal{T}_{\mathds P^2}(-2))=0,\\
\mathrm{Ext}^*(\mathrm U_{\mathrm Z}^{\vee},\mathcal{O}_{\mathrm A_2}(-1))
&=\mathrm H^*(\mathcal{O}_{\mathds P^2}(-2)\oplus
\mathcal{O}_{\mathds P^2}(-1))=0.
\end{aligned}
\]
Also $\mathrm{Ext}^*(\mathcal{O}_{\mathrm Z},\mathcal{O}_{\mathrm A_j}(-1))=0$.
These are orthogonality statements, not exceptional collections containing
the plane sheaf: $\mathcal{O}_{\mathrm A_j}(-1)$ is nonperfect at an ordinary
double point on the plane and is not exceptional there.

The appropriate perfect pair consequently belongs to the perfect
left orthogonal of the plane sheaf. Proposition~\ref{prop:sgequi} sends it
to an exceptional pair in $\mathcal B$. For one of the selected exceptional
curve objects, this gives respectively
\[
\langle\mathcal{O}_{L_k}(-1/2),
\mathrm\Theta^!g^*\mathrm U_{\mathrm Z},
\mathrm\Theta^!g^*\mathcal{O}_{\mathrm Z}\rangle,
\]
or
\[
\langle\mathcal{O}_{L_k}(-1/2),
\mathrm\Theta^!g^*\mathcal{O}_{\mathrm Z},
\mathrm\Theta^!g^*\mathrm U_{\mathrm Z}^{\vee}\rangle.
\]
The exceptionality of the tautological pair on a singular midpoint is a
separate global calculation. The plane restrictions alone do not establish
it. For several source curve objects, their mutual semi-orthogonality must
also be checked.
To distinguish the number of centres from the number of curve
objects, choose $0\leq N_0\leq N$ of the original orbifold points.
Subject to the geometric hypotheses above, the partial extraction has
the diagram
\begin{equation}\label{exsark2}
\begin{tikzcd}[column sep=0.5em]
  & && \mathrm{Y}^{G}_{1,i}:=wBl_{\mathrm{p}_{1},\ldots,\mathrm{p}_{N_0}}\widetilde{\mathrm{X}}_{1,i}\arrow[dl]\arrow[dr,"g"] && \\
 & &\widetilde{\mathrm{X}}_{1,i}& &\mathrm{Z}^{G}_{1,i}&
\end{tikzcd}
\end{equation}
If the remaining centres lie in the isomorphism locus of this
contraction, the corresponding blow-ups on both sides give
\begin{equation}\label{exsark3}
\begin{tikzcd}[column sep=0.5em]
  & && \overline{\mathrm{Y}}:=wBl_{\mathrm{p}_{N_0+1},\ldots,\mathrm{p}_{N}}\mathrm{Y}^{G}_{1,i}\arrow[dl]\arrow[dr,"g"] && \\
 & &\widetilde{\mathrm{X}}_{1,i}& &\overline{\mathrm{Z}}_{1,i}:=wBl_{\mathrm{p}_{N_0+1},\ldots,\mathrm{p}_{N}}\mathrm{Z}^{G}_{1,i}&
\end{tikzcd}
\end{equation}
Here the lower-left map is the composite extraction. The commuting
modification is justified by the isomorphic neighbourhoods of the
remaining centres; the upper source is smooth on its canonical stack.
We call $\overline{\mathrm{Z}}_{1,i}$ an \textit{extended global midpoint}
when these properties hold. The types in \eqref{targetZ} are candidates
for this construction, not an existence assertion for every row.\\

Under these geometric hypotheses, the one-centre argument gives
the following conditional multiple-unprojection statement.
\begin{corollary}\label{cor:multipleunprojection}
Under the preceding geometric hypotheses, choose mutually semi-orthogonal
exceptional curve objects
$\mathcal{O}_{L_{i_j},\mathrm{p}_{a(j)}}(-1/2)$ whose strict transforms
are exceptional fibres of $g$ with the one-centre local form of
\eqref{eq:lstirct}. Let $\langle e_1,\ldots,e_l\rangle$ be a perfect
exceptional collection on $\overline{\mathrm{Z}}_{1,i}$, left orthogonal
to all $\mathcal{O}_{\mathrm{A}_a}(-1)$. Then
\[
\langle\mathcal{O}_{L_{i_1},\mathrm{p}_{a(1)}}(-1/2),\ldots,
\mathcal{O}_{L_{i_m},\mathrm{p}_{a(m)}}(-1/2),
\mathrm{\Theta}^{!}g^*e_1,\ldots,\mathrm{\Theta}^{!}g^*e_l\rangle
\]
is an exceptional collection on $\widetilde{\mathrm{X}}_{1,i}$.
\end{corollary}
\begin{proof}
The disjoint-centre extraction formula identifies the image of
$\mathrm{\Theta}$ with the left orthogonal to all
$\mathcal{O}_{\mathrm{A}_a}(-1)$. For each perfect $e_b$, adjunction and
the assumed orthogonality put $g^*e_b$ in this image. The localization
argument of Proposition~\ref{prop:sgequi} makes
$\mathrm{\Theta}^{!}g^*$ fully faithful on these objects. Since each selected strict transform is a fibre
of $g$, its sheaf $\mathcal{O}_{\widehat{L_{i_j}}}(-1)$ has zero
derived pushforward. Applying $\mathrm{Ext}^*(g^*e_b,-)$ to the comparison triangle of each selected
curve gives orthogonality to that source curve object, exactly as in
the one-centre proof. Together with the assumed semi-orthogonality of
the curve objects, this proves the assertion.
\end{proof}

Finally, we end with a conjecture for a uniform description.
\begin{conjecture}\label{conj:noncommuativehourglasses}
For a general hourglass $\mathrm{X}_{1,i}$, we have its non-commutative resolution
$$\mathrm{D}^{b}(\mathrm{\widetilde{X}}_{1,i})\simeq\langle\mathcal{O}_{L_{i_{1}},\mathrm{p}_{1}}(-a_{1}),..., \mathcal{O}_{L_{i_{m}},\mathrm{p}_{m}}(-a_{m}), \mathcal{O}_{\mathrm{\widetilde{X}}_{1,i}},\mathrm{U}^{\vee}_{\mathrm{\widetilde{X}}_{1,i}}, \mathcal{MF}_{\underline{k}}(\mathrm{\widetilde{X}_{1,i}}) \rangle$$
for any stacky line $L_{i_{j},\mathrm{p}_{j}}$ passing through the $j$-th orbifold point $\mathrm{p}_{j}$ on $\widetilde{\mathrm{X}}_{1,i}$. \\

\begin{enumerate}

\item $a_i\in\{1/2,1\}$. For intersecting $L_{i_j,\mathrm{p}_j}$, replace $\mathcal{O}_{L_{i_j}}(-a_j)$ by its stacky modification $\mathcal{O}_{L_{i_j}}(-a_j)^{\circ}$.
  \item where $\mathrm{U}^{\vee}_{\mathrm{\widetilde{X}_{1,i}}}$ comes from the key Grassmannian of the \textit{extended global midpoint}  $\mathrm{\overline{Z}}_{1,i}$ by appropriate adjustment,
   \item $\mathcal{MF}_{\underline{k}}(\mathrm{\widetilde{X}_{1,i}})$ is equivalent to each other for any choice of $L_{i_{j},\mathrm{p}_{j}}$, and we denote this equivalent category as $\mathcal{MF}(\mathrm{\widetilde{X}_{1,i}})$,
   \item if $g(\mathrm{C}_{1,i})=0$, there is an equivalence
   $$\mathcal{MF}(\mathrm{\widetilde{X}_{1,i}})\simeq\mathcal{MF}(\mathrm{\widetilde{X}_{1,i}'})$$
   where $\mathcal{MF}(\mathrm{\widetilde{X}_{1,i}'})$ is defined in Corollary \ref{cor:mfx'},
   \item in any case $\mathcal{MF}(\mathrm{\widetilde{X}_{1,i}})$ can always be realized as a (higher) \textit{matrix factorization category} (it is non-commutative completely intersected).
\end{enumerate}

\end{conjecture}

\section{Acknowledgements}
I am grateful to my supervisor Professor Nagai Yasunari for reminding me of Takagi's series works, and also for his continued weekly discussions with me. Especially, I would like to thank Professor Hiromichi Takagi for his interest in my work and for being able to publish his related work in advance. In a private email with him, he showed his unpublished observation of the derived category of varieties on his tables, which inspired me to many ideas in Section \ref{sec:Tak}. Finally, I would also like to thank the members of the Northwest Tokyo Algebraic Geometry Group for their helpful discussions and concern.

\quad\\
\author{HAO XINGBANG}

\newcommand{\Addresses}{{
  \bigskip
  \footnotesize

  HAO XINGBANG, \textsc{Waseda University, Ookubo, Shinjuku-ku, Tokyo, 169-8555, Japan.}\par\nopagebreak
  \textit{E-mail address}: \texttt{hao@fuji.waseda.jp}
}}
\Addresses

\end{document}